# RESERVATION FOR OTHER BACKWARD CLASSES IN INDIAN CENTRAL GOVERNMENT INSTITUTIONS LIKE IITs, IIMs AND AIIMS – A STUDY OF THE ROLE OF MEDIA USING FUZZY SUPER FRM MODELS

**W. B. Vasantha Kandasamy**
**Florentin Smarandache**
**K. Kandasamy**

**2009**

# RESERVATION FOR OTHER BACKWARD CLASSES IN INDIAN CENTRAL GOVERNMENT INSTITUTIONS LIKE IITs, IIMs AND AIIMS – A STUDY OF THE ROLE OF MEDIA USING FUZZY SUPER FRM MODELS

**W. B. Vasantha Kandasamy**
e-mail: vasanthakandasamy@gmail.com
web: http://mat.iitm.ac.in/~wbv
www.vasantha.in

**Florentin Smarandache**
e-mail: smarand@unm.edu

**K. Kandasamy**
e-mail: dr.k.kandasamy@gamil.com

**2009**

# CONTENTS

# DEDICATION

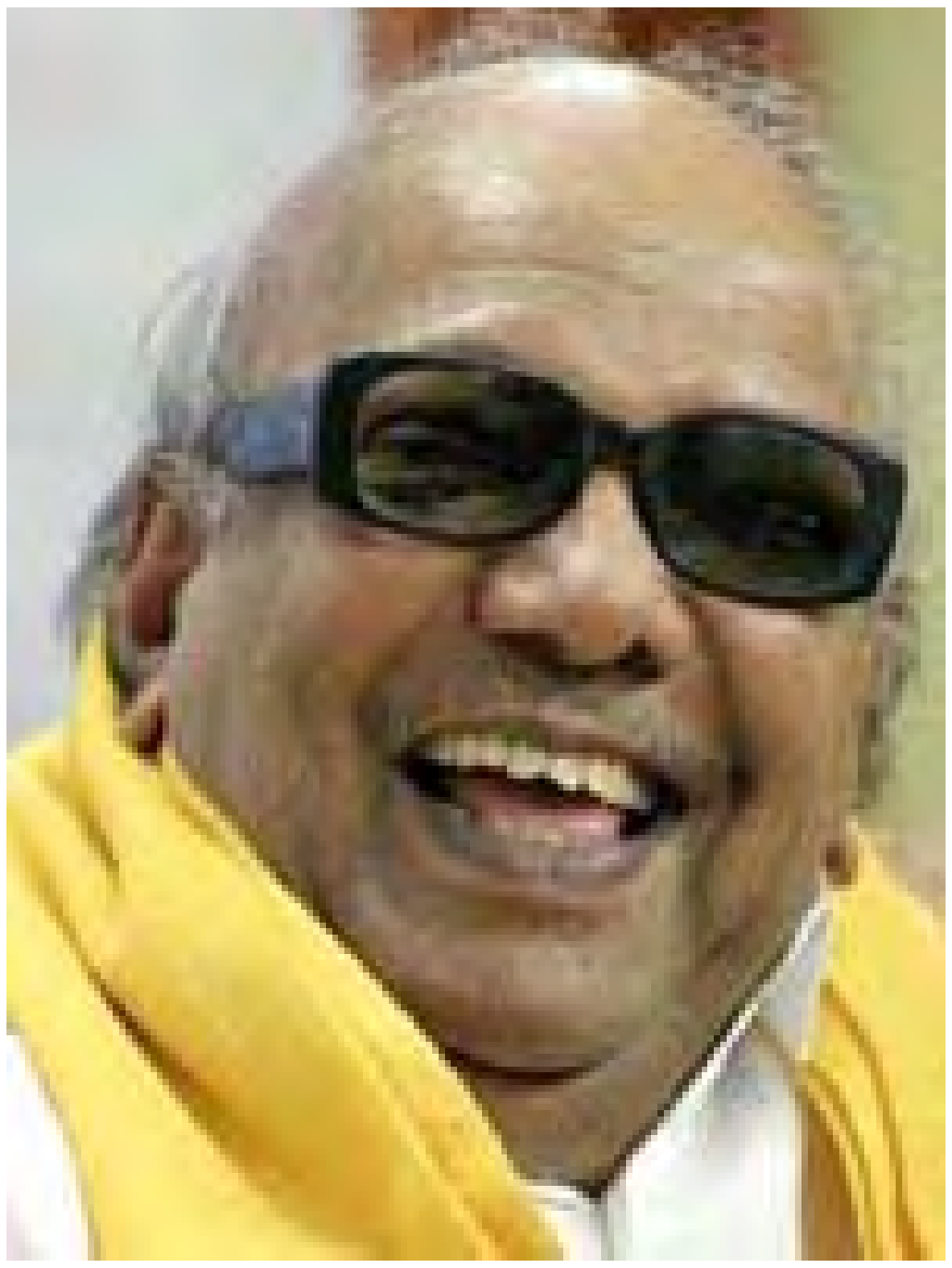

*We dedicate this book to Thanthai Periyar's foremost follower and five-time Tamil Nadu Chief Minister Hon'ble Dr. Kalaignar for his uncompromising struggle to ensure social justice through 27% reservation for the Other Backward Classes (OBC) in Central Government-run higher educational institutions like the IITs, IIMs and AIIMS.*

# PREFACE

The new notions of super column FRM model, super row FRM model and mixed super FRM model are introduced in this book. These three models are introduced specially to analyze the biased role of the print media on 27 percent reservation for the Other Backward Classes (OBCs) in educational institutions run by the Indian Central Government. This book has four chapters. In chapter one the authors introduce the three types of super FRM models. Chapter two uses these three new super fuzzy models to study the role of media which feverishly argued against 27 percent reservation for OBCs in Central Government-run institutions in India. The experts we consulted were divided into 19 groups depending on their profession. These groups of experts gave their opinion and comments on the news-items that appeared about reservations in dailies and weekly magazines, and the gist of these lengthy discussions form the third chapter of this book. The fourth chapter gives the conclusions based on our study. Our study was conducted from April 2006 to March 2007, at which point of time the Supreme Court of India stayed the 27 percent reservation for OBCs in the IITs, IIMs and AIIMS. After the aforesaid injunction from the Supreme Court, the experts did not wish to give their opinion since the matter was *sub-judice*. The authors deeply acknowledge the service of each and every expert who contributed their opinion and thus made this book a possibility. We have analyzed the data using the opinion of the experts who formed a heterogeneous group consisting of administrators, lawyers, OBC/SC/ST students, upper caste students and Brahmin students, educationalists, university vice-chancellors, directors, professors, teachers, retired Judges, principals of colleges, parents, journalists, members of the public, politicians, doctors, engineers, NGOs and government staff.

The authors deeply acknowledge the unflinching support of Kama and Meena.

W.B.VASANTHA KANDASAMY
FLORENTIN SMARANDACHE
K.KANDASAMY

**Chapter One**

# INTRODUCTION TO NEW SUPER FUZZY MODELS

In this chapter the authors introduce some new fuzzy models using supermatrices. These new fuzzy models are used in chapter two to analyze the role of media in the context of OBC (Other Backward Castes/Classes) reservation in the institutions run by the central government. In the first section we recall the definition of supermatrices and then we define the notion of fuzzy supermatrices. Section two defines the new notion of fuzzy super models and shows how they function.

## 1.1 Supermatrices and Fuzzy Supermatrices

Here we just recall the notion of supermatrices and define the new notion of fuzzy supermatrices.

**DEFINITION 1.1.1:** *Let V = ($V_1$ $V_2$ …$V_n$) where each $V_i$ is a row vector* $\left(v_1^i ... v_{t_1}^i\right)$*,* $1 \le i \le n$*,* $v_j^i \in Q$ *or R, then V is denoted by V = ($V_1$ … $V_n$) =* $\left\{\left(v_1^1 v_2^1 \cdots v_{t_1}^1 | v_1^1 v_2^1 \cdots v_{t_2}^1 | \cdots | v_1^n v_2^n \cdots v_{t_n}^n\right)\right\}$ *where* $v_j^i \in Q$ *or R, i=1, 2, …, n and* $i \le j \le t_i$*. If $t_1$, $t_2$, …, $t_n$ are distinct then we call V to be a mixed super row vector. If $t_1$ = $t_2$ = … = $t_n$ we then call V to be a super row vector.*

We just illustrate this by the following example.

***Example 1.1.1*:** Let V= (1 2 3 0 | 5 7 –3 1 | 6 1 2 7 | 0 0 –1 4); V is called the super row vector and V = ($V_1$ $V_2$ $V_3$ $V_4$) where each $V_i$ is a $1 \times 4$ row vector.

**DEFINITION 1.1.2:** *Let $V = (V_1\ V_2\ \ldots V_n)$ be a super row vector, if the entries of each $V_i$ is from [0, 1], the unit interval; i=1, 2,…, n then we call V to be a simple fuzzy super row vector. Let $V= (V_1\ V_2\ \ldots V_n)$ be the mixed super row vector, where each $V_i$ is a $1 \times t_i$ row vector, i =1, 2, …, n and $t_i \neq t_j$; if $i \neq j$ for atleast one i and j: $1 \leq i, j \leq n$ ; and if the entries of each of the row vector $V_i$ is from the unit interval [0, 1]; $1 \leq i \leq n$ ; then we call V to be the mixed simple fuzzy super row vector.*

We illustrate the above definition by the following examples.

***Example 1.1.2*:** Let V = (10.7 0.8 | 0.9 0.2 0.3 0.4 | 0.6 1 0.5 0.2 1 0.8 0). V is a mixed simple fuzzy super row vector, the entry of each $V_i$ is from [0, 1] where V = ($V_1$ $V_2$ $V_3$); i = 1, 2, 3.

***Example 1.1.3*:** Let A = (1 0 1 0 | 0 1 1 1 | 1 1 1 1 | 1 0 0 1 | 1 1 1 0) = ($A_1$ $A_2$ $A_3$ $A_4$ $A_5$). A is a simple fuzzy super row vector for entries of each $A_i$ is from the set $\{0, 1\} \subseteq [0, 1]$.

Now we proceed on to define simple super column vector and mixed simple super column vector.

**DEFINITION 1.1.3:** *Let*

$$M = \begin{bmatrix} M_1 \\ M_2 \\ \vdots \\ M_s \end{bmatrix}$$

*where $M_i$ is a $t_i \times 1$ column vector with entries from Q or R, $1 \leq i \leq s$; then we call M to be a simple mixed super column vector if $t_i \neq t_j$ for at least one $i \neq j$. If $t_1 = t_2 = \ldots = t_s$ then we call M to be a simple super column vector.*

We illustrate this by following examples.

***Example 1.1.4**:* Let

$$M = \begin{bmatrix} 3 \\ 7 \\ \frac{1}{5} \\ 2 \end{bmatrix} = \begin{bmatrix} M_1 \\ M_2 \end{bmatrix}.$$

M is a simple mixed super column vector.

***Example 1.1.5**:* Consider

$$M = \begin{bmatrix} M_1 \\ M_2 \\ M_3 \\ M_4 \end{bmatrix} = \begin{bmatrix} 3 \\ 1 \\ \frac{1}{2} \\ 0 \\ \frac{5}{7} \\ 2 \\ \frac{1}{-3} \\ 1 \\ 2 \end{bmatrix}.$$

M is a simple super column vector.

Now we define simple fuzzy super column vector and mixed simple fuzzy super column vector.

Let

$$M = \begin{bmatrix} M_1 \\ M_2 \\ \vdots \\ M_s \end{bmatrix}$$

be a mixed simple super column vector, if each of the entries are from the unit interval [0, 1] then we call M to be mixed simple fuzzy super column vector.

$$M = \begin{bmatrix} M_1 \\ M_2 \\ \vdots \\ M_t \end{bmatrix}$$

be a simple super column vector where each $M_i$ is a $r \times 1$ column vector i =1, 2, …, t and each entry of $M_i$ is from the unit interval [0, 1] then we call M to be a fuzzy super column vector.

***Example 1.1.6:*** Let

$$M = \left[\begin{array}{c} 1 \\ 0 \\ 1 \\ 1 \\ \hline 0.2 \\ 0.7 \\ 0 \\ 1 \\ 1 \\ 0.5 \\ 0.2 \\ \hline 0.7 \\ 0.5 \\ 1 \\ 0.9 \end{array}\right].$$

M is a fuzzy super column vector.

Next we proceed on to define the notion of super row matrix and mixed super row matrix.

**DEFINITION 1.1.4:** *Let $M = [M_1\ M_2\ \ldots\ M_j]$, where $M_t$ is a $s \times p_t$ matrix with s rows and $p_t$ columns, t = 1, 2, …, j ; where* *elements of each $M_t$ is from Q or R; $1 \leq t \leq j$. Then we call M to be mixed super row matrix. Here at least one $p_t \neq p_i$, if $t \neq i$. If we have on the other hand $p_t = p_i$ for $i \neq t$ and $1 \leq t, i \leq j$, then we call M to be a super row matrix.*

We illustrate this by the following examples.

***Example 1.1.7:*** Consider the super row matrix

$$M = [M_1\ M_2\ M_3\ M_4] =$$

$$\left[\begin{array}{cccc|cccc|cccc|cccc} 3 & 1 & 0 & 1 & 0 & 2 & 1 & 3 & 1 & 2 & 3 & 4 & 3 & 6 & 1 & 0 \\ 0 & 4 & 5 & 7 & 3 & 1 & 0 & 2 & 5 & 6 & 7 & 8 & 1 & 7 & 3 & 7 \\ 2 & 1 & 1 & 5 & 1 & 5 & 7 & 9 & 9 & 0 & 1 & 5 & 1 & 1 & 0 & 2 \end{array}\right].$$

Clearly M is a super row matrix.

***Example 1.1.8:*** Let

$$M = \left[\begin{array}{ccc|cc|c} 3 & 0 & 2 & 3 & 5 & 2 \\ 1 & 0 & 1 & 1 & 0 & 3 \\ 0 & 0 & 1 & 0 & 1 & 4 \\ 2 & 2 & 1 & 1 & 1 & 0 \end{array}\right] = [M_1\ M_2\ M_3],$$

M is a mixed super row matrix.

**DEFINITION 1.1.5:** *Let $M = [M_1\ M_2\ \ldots\ M_s]$ be a super row matrix i.e., each $M_i$ is a $t \times p$ matrix. If each $M_i$ takes its entries from the unit internal [0, 1], then we call M to be fuzzy super row matrix, $I= 1,2,\ldots,s$. Suppose $M = [M_1\ M_2\ \ldots\ M_t]$ be a mixed row matrix and if each $M_i$ is a $n \times p_i$ fuzzy matrix with entries from [0, 1] having n rows and $p_i$ columns with at least one $p_i \neq p_j$, if $i \neq j$, $1 \leq i, j \leq t$, then we call M to be a mixed fuzzy super row matrix.*

We illustrate this situation by the following examples.

***Example 1.1.9:*** Let

$$M = [M_1\ M_2\ M_3\ M_4] =$$

$$\left[\begin{array}{cccc|cccc} 0.3 & 1 & 0.2 & 0.7 & 1 & 0 & 1 & 0 \\ 0.6 & 0.5 & 1 & 0.2 & 0 & 1 & 0 & 1 \\ 0.9 & 0.3 & 0 & 1 & 1 & 1 & 1 & 1 \\ 0.8 & 0.9 & 0.3 & 0.2 & 0 & 0 & 1 & 0 \end{array}\right.$$

$$\left|\begin{array}{cccc|cccc} 0.2 & 0.3 & 0 & 1 & 1 & 0.7 & 0.9 & 1 \\ 0 & 1 & 0.6 & 0.8 & 0 & 0.8 & 1 & 0.6 \\ 1 & 0.4 & 1 & 0.9 & 1 & 0.2 & 0.5 & 1 \\ 0.2 & 1 & 1 & 1 & 0.7 & 0 & 0.6 & 0.2 \end{array}\right].$$

M is a fuzzy super row matrix. Each $M_i$ is a $4 \times 4$ fuzzy matrix for i = 1, 2, 3, 4.

***Example 1.1.10:*** Let

$$M = [M_1|M_2\ |M_3] =$$

$$\left[\begin{array}{cccc|cc|ccccc} 0.9 & 1 & 0.8 & 0.2 & 0.2 & 1 & 1 & 0 & 0.4 & 0.7 & 0.8 \\ 0.2 & 0.6 & 0 & 0.1 & 0.7 & 0.4 & 0.9 & 1 & 0.3 & 0.2 & 0.9 \\ 0.7 & 0.5 & 1 & 0.4 & 0.3 & 0 & 0.2 & 0.5 & 0.6 & 0.1 & 0.3 \end{array}\right],$$

is the mixed fuzzy super matrix; $M_1$ is a $3 \times 4$ fuzzy matrix, $M_2$ is a $3 \times 2$ fuzzy matrix and $M_3$ is a $3 \times 5$ fuzzy matrix.

**DEFINITION 1.1.6:** *Let*

$$S = \begin{bmatrix} S_1 \\ S_2 \\ \vdots \\ S_n \end{bmatrix}$$

*where each $S_i$ is a $p_i \times m$ matrix with entries from Q or R, i = 1, 2, …, n. We call S to be a mixed super column matrix if $p_i \neq p_j$ for at least one $i \neq j$, $1 \leq i, j \leq n$.*

**DEFINITION 1.1.7:** *Let*

$$P = \begin{bmatrix} P_1 \\ P_2 \\ \vdots \\ P_m \end{bmatrix}$$

*where $P_i$ is a $t \times s$ matrix with entries from Q or R; $1 \leq i \leq m$. Then we call P to be a super column matrix.*

We proceed on to illustrate these by the following examples.

***Example 1.1.11:*** Let

$$S = \begin{bmatrix} S_1 \\ S_2 \\ S_3 \\ S_4 \end{bmatrix} = \left[\begin{array}{cccc} 3 & 1 & 5 & 9 \\ 0 & 2 & 1 & -7 \\ \hline 1 & 0 & 2 & 3 \\ 4 & 5 & 6 & 0 \\ -7 & 8 & 9 & 2 \\ \hline 1 & -5 & 0 & 3 \\ 2 & 6 & 1 & 5 \\ 3 & 7 & 7 & 9 \\ \hline 4 & 8 & 2 & 6 \\ 1 & 9 & 0 & 2 \end{array}\right].$$

S is a mixed super column vector for $S_i$ is a $2 \times 4$ matrix. $S_2$ is a $3 \times 4$ matrix. $S_3$ is a $3 \times 4$ matrix and $S_4$ is a $2 \times 4$ matrix.

***Example 1.1.12:*** Let

$$S = \begin{bmatrix} S_1 \\ S_2 \\ S_3 \\ S_4 \end{bmatrix}$$

$$= \left[\begin{array}{ccc} 3 & 0 & 2 \\ 1 & 1 & 0 \\ 5 & 2 & 7 \\ 1 & 0 & -3 \\ \hline 1 & 1 & 1 \\ 2 & 3 & -1 \\ 4 & 5 & 6 \\ 0 & 2 & 7 \\ \hline 1 & 2 & 3 \\ 4 & 5 & 6 \\ 7 & 8 & 9 \\ 0 & -1 & 2 \\ \hline 0 & 2 & 1 \\ 0 & 1 & 0 \\ 0 & 3 & -3 \\ 7 & 5 & 7 \end{array}\right].$$

S is a super column vector and $S_i$ is a $4 \times 3$ matrix.

Now we proceed onto define the notion of super fuzzy column matrix and mixed super fuzzy column matrix.

**DEFINITION 1.1.8:** *Let*

$$S = \begin{bmatrix} S_1 \\ S_2 \\ \vdots \\ S_m \end{bmatrix}$$

*be a super column matrix where each $S_i$ is a $t \times s$ fuzzy matrix; $1 \leq i \leq m$. We call S to be a fuzzy super column matrix if each entry of $S_i$ is from the unit interval [0, 1], $1 \leq i \leq m$.*

**DEFINITION 1.1.9:** *Let*

$$M = \begin{bmatrix} M_1 \\ M_2 \\ \vdots \\ M_m \end{bmatrix}$$

*be a super column matrix where each $M_i$ is a $n_i \times s$ fuzzy matrix, i.e., each of the matrix $M_i$ takes its entries from the unit interval [0, 1]; with atleast one $M_i \neq M_j$ if $i \neq j$; $1 \leq i, j \leq m$; then we call M to be a mixed fuzzy super column matrix.*

We now illustrate the definitions 1.1.8 and 1.1.9 by the following examples.

***Example 1.1.13***: Let

$$S = \begin{bmatrix} S_1 \\ S_2 \\ S_3 \\ S_4 \end{bmatrix} = \left[\begin{array}{cccc} 0.3 & 1 & 0.7 & 0.2 \\ 0.6 & 0.2 & 1 & 0.5 \\ 1 & 0.9 & 0.4 & 1 \\ 0.7 & 0.2 & 0.6 & 0.4 \\ 0.6 & 0 & 1 & 0.3 \\ \hline 0.1 & 0 & 1 & 0.2 \\ 0.7 & 0.6 & 0.9 & 0.8 \\ \hline 0.7 & 1 & 0.5 & 0.2 \\ 0.6 & 0 & 1 & 0.1 \\ 0.9 & 0.8 & 0.4 & 0 \\ \hline 0.9 & 1 & 0 & 0.9 \\ 0.3 & 0 & 1 & 0.2 \end{array}\right],$$

S is a mixed fuzzy super column matrix**.**

***Example 1.1.14:*** Let

$$T = \begin{bmatrix} T_1 \\ T_2 \\ T_3 \\ T_4 \end{bmatrix} = \left[\begin{array}{ccc} 0.3 & 0.2 & 0.1 \\ 1 & 0 & 0.7 \\ 0.4 & 0.5 & 0.9 \\ 0.6 & 1 & 1 \\ \hline 0.9 & 1 & 0 \\ 0.2 & 0 & 0.7 \\ 1 & 0.8 & 0.4 \\ 0 & 0.1 & 1 \\ \hline 0.1 & 0.6 & 0.3 \\ 1 & 0 & 0.9 \\ 0.7 & 0.9 & 1 \\ 0 & 1 & 0.2 \\ \hline 1 & 0.2 & 0.4 \\ 0.5 & 1 & 0.7 \\ 1 & 0 & 1 \\ 0.9 & 1 & 0.2 \end{array}\right].$$

T is a fuzzy super column matrix and each $T_i$ is a 4 × 3 fuzzy matrix, i = 1, 2, 3, 4.

Now we proceed on to define supermatrices.

**DEFINITION 1.1.10:** *Let*

$$M = \begin{bmatrix} M_{11} & M_{12} & \cdots & M_{1n} \\ M_{21} & M_{22} & \cdots & M_{2n} \\ \vdots & \vdots & & \vdots \\ M_{m1} & M_{m2} & \cdots & M_{mn} \end{bmatrix}$$

*where*

*(1) $M_{ij}$ are $t_i \times s_j$ matrices with entries form Q or R, $1 \le i \le m$ and $1 \le j \le n$.*

2) $$\begin{bmatrix} M_{11} & M_{12} & \ldots & M_{1n} \end{bmatrix}$$
$$\begin{bmatrix} M_{21} & M_{22} & \ldots & M_{2n} \end{bmatrix}$$
$$\vdots$$
$$\begin{bmatrix} M_{m1} & M_{n2} & \ldots & M_{mn} \end{bmatrix}$$

*are mixed super row matrices.*

3) $$\begin{bmatrix} M_{11} \\ M_{21} \\ \vdots \\ M_{m1} \end{bmatrix}, \begin{bmatrix} M_{12} \\ M_{22} \\ \vdots \\ M_{m2} \end{bmatrix}, \cdots, \begin{bmatrix} M_{1n} \\ M_{2n} \\ \vdots \\ M_{mn} \end{bmatrix}$$

*are mixed super column matrices.*

*Then we call M to be a supermatrix or we can define M as follows:*

*$M_{11}$, $M_{12}$, ..., $M_{1n}$ are matrices each of $M_{1j}$ ($1 \leq j \leq n$) having same number of rows but different number of columns. $M_{21}$, $M_{22}$, ..., $M_{2n}$ are matrices each of $M_{2k}$ ($1 \leq k \leq n$) having the same number of rows but different number of columns and so on. Thus $M_{m1}$, $M_{m2}$ ,..., $M_{mn}$ are matrices each of $M_{mt}$ ($1 \leq t \leq n$) have the same number of rows but different number of columns.*

*Like wise $M_{11}$, $M_{21}$,..., $M_{m1}$ are matrices each of $M_{i1}$ , ($1 \leq i \leq m$) having same number of columns but different number of rows and so on, i.e., $M_{1n}$, $M_{2n}$ ,..., $M_{mn}$ are matrices each of $M_{jn}$ ($1 \leq j \leq m$) having same number of columns but different number of rows.*

We illustrate this by the following example.

***Example 1.1.15:*** Let

$$M = \left[\begin{array}{cc|cccc} 1 & 2 & 3 & 4 & 5 & 6 \\ 7 & 8 & 9 & 1 & 0 & 1 \\ -2 & 3 & 6 & 0 & 9 & -1 \\ \hline 1 & 2 & 3 & 4 & 5 & 0 \\ 7 & 6 & 5 & 4 & 3 & 2 \\ 1 & 0 & 2 & 1 & 3 & 0 \\ 4 & 5 & 0 & 6 & 7 & 0 \\ \hline 8 & 9 & 1 & 2 & 3 & 1 \\ 9 & 0 & 1 & 0 & 3 & 0 \end{array}\right] = \begin{bmatrix} M_{11} & M_{12} \\ M_{21} & M_{22} \\ M_{31} & M_{32} \end{bmatrix}$$

where $M_{11}$ $M_{12}$ have the same number of rows but different number of columns i.e., each of $M_{11}$ and $M_{12}$ have only 3 rows but $M_{11}$ has 2 columns and $M_{12}$ has 3 rows and 4 columns. $M_{21}$ and $M_{22}$ has same number of rows viz. four but different number of columns i.e.; $M_{21}$ has 4 rows and 2 columns where as $M_{22}$ has 4 rows and 4 columns. Now $M_{31}$ and $M_{32}$ have same number of rows viz., 2 but different number of columns; i.e. $M_{31}$ has 2 rows and 2 columns but $M_{32}$ has two rows and four columns,

Now the set $\{M_{11}\ M_{21}\ M_{31}\}$ have different number of rows but same number of columns i.e., $M_{11}$ has 3 rows and two columns. $M_{21}$ has four rows and two columns and $M_{31}$ has two rows and two columns i.e., all of $M_{11}$, $M_{21}$ and $M_{31}$ have same number of columns viz; two. Similarly $M_{12}$, $M_{22}$ and $M_{32}$ have same number of columns namely four but $M_{12}$ has three rows $M_{22}$ has four rows and $M_{32}$ two rows. Thus we can say ($M_{11}$ $M_{12}$), ($M_{21}$ $M_{22}$) and ($M_{31}$ $M_{32}$) are mixed super column matrices and

$$\begin{bmatrix} M_{11} \\ M_{21} \\ M_{31} \end{bmatrix} \text{ and } \begin{bmatrix} M_{21} \\ M_{22} \\ M_{32} \end{bmatrix}$$

are mixed super row matrices.

We give yet another example so that the reader becomes familiar with it.

***Example 1.1.16:*** Consider the supermatrix

$$M = \left[\begin{array}{c|cc|ccccc|ccc} 1 & 2 & 9 & 8 & 7 & 6 & 4 & 5 & 0 & 2 & 4 \\ 0 & 9 & 3 & -1 & 2 & 0 & -4 & 0 & 3 & -1 & -5 \\ 9 & 7 & 0 & 1 & 5 & 9 & 8 & 7 & 6 & 4 & 5 \\ \hline 6 & 0 & -3 & 2 & 6 & 0 & 1 & 3 & 5 & 7 & 9 \\ 3 & -2 & 0 & 3 & 7 & 1 & 2 & 4 & -6 & 8 & -1 \\ \hline 1 & 4 & 1 & 4 & 8 & 2 & -3 & +3 & 7 & 0 & 2 \end{array}\right] =$$

$$\begin{bmatrix} M_{11} & M_{12} & M_{13} & M_{14} \\ M_{21} & M_{22} & M_{23} & M_{24} \\ M_{31} & M_{32} & M_{33} & M_{34} \end{bmatrix}$$

Clearly ($M_{11}$ $M_{12}$ $M_{13}$ $M_{14}$) is a mixed super column matrix for $M_{11}$ is a $3 \times 1$ matrix, $M_{12}$ is a $3 \times 2$ matrix, $M_{13}$ is a $3 \times 5$ matrix and $M_{14}$ is a $3 \times 3$ matrix. Thus ($M_{11}$ $M_{12}$ $M_{13}$ $M_{14}$) is a mixed super column matrix, hence all of them have only 3 rows. Now ($M_{21}$ $M_{22}$ $M_{23}$ $M_{24}$) is also a mixed super column matrix and each of $M_{2j}$ has only 2 rows for $1 \leq j \leq 4$. Similarly ($M_{31}$ $M_{32}$ $M_{33}$ $M_{34}$) is a mixed super column matrix where each of $M_{3k}$ has only one row i.e., ($M_{31}$ $M_{32}$ $M_{33}$ $M_{34}$) is a mixed super column matrix.

Now consider

$$\begin{bmatrix} M_{11} \\ M_{21} \\ M_{31} \end{bmatrix};$$

each $M_{j1}$ is a matrix with only one column i.e.,

$$\begin{bmatrix} M_{11} \\ M_{21} \\ M_{31} \end{bmatrix}$$

is a mixed super row vector.

$$\begin{bmatrix} M_{21} \\ M_{22} \\ M_{32} \end{bmatrix}$$

is a mixed super row matrix and each of $M_{k2}$ has only 2 columns but different number of rows ; $1 \leq k \leq 3$.

$$\begin{bmatrix} M_{13} \\ M_{23} \\ M_{33} \end{bmatrix}$$

is a mixed super row matrix and each of $M_{r3}$, $1 \leq r \leq 3$ have only 5 columns. Like wise

$$\begin{bmatrix} M_{14} \\ M_{24} \\ M_{34} \end{bmatrix}$$

is a mixed super row matrix and each of $M_{p3}$, $1 \leq p \leq 3$ have only 3 columns. Thus M is a supermatrix; this supermatrix has column vectors, row vectors, square matrices and rectangular matrices as its components.

***Example 1.1.17:*** Suppose

$$M = \left[\begin{array}{ccc|ccc|ccc} -2 & 1 & 0 & 3 & 4 & 5 & 7 & 8 & 9 \\ 1 & 2 & 3 & 9 & 8 & 7 & 6 & 7 & 8 \\ 0 & 1 & 2 & 0 & -3 & 0 & 4 & 0 & 5 \\ \hline 1 & -7 & 1 & -9 & 1 & 5 & 1 & 3 & 2 \\ 2 & 4 & 2 & 6 & 2 & -8 & 2 & 0 & 1 \\ -3 & 1 & 0 & 2 & 3 & 1 & 3 & -1 & 0 \\ \hline 7 & -5 & 3 & 2 & -9 & 4 & 9 & 8 & 1 \\ 3 & 1 & 2 & 0 & 5 & 0 & 1 & 1 & 0 \\ 8 & 0 & 1 & -1 & 0 & 8 & 0 & -1 & 2 \end{array}\right]$$

$$= \begin{bmatrix} M_{11} & M_{12} & M_{13} \\ M_{21} & M_{22} & M_{23} \\ M_{31} & M_{32} & M_{33} \end{bmatrix}.$$

We see M is a supermatrix. Now ($M_{11}$ $M_{12}$ $M_{13}$) is a super row matrix ($M_{21}$ $M_{22}$ $M_{23}$) and ($M_{31}$ $M_{32}$ $M_{33}$) are also super row matrices. None of them is mixed.

$$\begin{bmatrix} M_{11} \\ M_{21} \\ M_{31} \end{bmatrix}, \begin{bmatrix} M_{12} \\ M_{22} \\ M_{32} \end{bmatrix} \text{ and } \begin{bmatrix} M_{13} \\ M_{23} \\ M_{33} \end{bmatrix}$$

are super column matrices. None of them is mixed.

Further we see each of the matrices $M_{ij}$ in M is a $3 \times 3$ square matrix; $1 \leq i, j \leq 3$. We call such supermatrices as perfect square supermatrices.

**DEFINITION 1.1.11:** *Let M be a supermatrix.*

$$M = \begin{bmatrix} M_{11} & M_{12} & \cdots & M_{1n} \\ M_{21} & M_{22} & \cdots & M_{2n} \\ \vdots & \vdots & & \vdots \\ M_{n1} & M_{n2} & \cdots & M_{nn} \end{bmatrix}$$

*If each $M_{ij}$ is a $t \times t$ square matrix $1 \leq i, j \leq n$ then we define M to be a perfect square supermatrix. All supermatrices in general need not be perfect square supermatrices.*

In a perfect square supermatrix M, number of column will be equal to the number of rows. If the number of columns in a supermatrix M is not equal to the number of rows in the super matrix then we cannot have M to be a perfect square matrix. Now we give a method for obtaining a supermatrix from the given matrix.

Suppose we have a matrix

$$M = \begin{bmatrix} m_{11} & m_{12} & \cdots & m_{1n} \\ m_{21} & m_{22} & \cdots & m_{2n} \\ \vdots & \vdots & & \vdots \\ m_{m1} & m_{m2} & \cdots & m_{mn} \end{bmatrix},$$

where $m_{ij}$ are in Q or R. If we draw lines in between two columns say in between the $j^{th}$ and the $(j + 1)^{th}$ column. Suppose we draw lines between $(j + r)^{th}$ column and $(j + r + 1)^{th}$ column and so on say between $(j + t)^{th}$ and $(j + t + 1)^{th}$ column $j + t + 1 < n$ and $j < j + r < \ldots < j + t$ then M is said to become a super column matrix which has m rows. Similarly if the rows are partitioned i.e., lines are drawn between two rows, then M is said to have become a super row matrix that has only n columns. If both lines are drawn in between rows as well as in between columns to the given matrix $M = (m_{ij})$. M becomes a super matrix. Thus any $n \times n$ square matrix n not a prime can be made into a perfect square matrix.

We proceed to give one illustration. The reader is requested to refer [125] for more properties regarding super matrices.

Fuzzy supermatrices are defined only in the book [316].

***Example 1.1.18:*** Let

$$M = \begin{bmatrix} 3 & 4 & 5 & 7 & 8 & 9 \\ 1 & 2 & 3 & 0 & 1 & 2 \\ 0 & 6 & 1 & 4 & 0 & -3 \\ -1 & +1 & 2 & 5 & 1 & 4 \end{bmatrix}$$

be a $4 \times 6$ matrix. Suppose we draw lines between the rows 2 and 3, 3 and 4 only we get

$$M = \left[\begin{array}{cccccc} 3 & 4 & 5 & 7 & 8 & 9 \\ 1 & 2 & 3 & 0 & 1 & 2 \\ \hline 0 & 6 & 1 & 4 & 0 & -3 \\ \hline -1 & 1 & 2 & 5 & 1 & 4 \end{array}\right] = \begin{bmatrix} M_1 \\ M_2 \\ M_3 \end{bmatrix}$$

to be only a mixed super row matrix.

Suppose for the same matrix M we draw lines between the 3rd and 4th column and 5th and 6th column we see

$$M = \left[\begin{array}{ccc|cc|c} 3 & 4 & 5 & 7 & 8 & 9 \\ 1 & 2 & 3 & 0 & 1 & 2 \\ 0 & 6 & 1 & 4 & 0 & -3 \\ -1 & 1 & 2 & 5 & 1 & 4 \end{array}\right] = \left[ M_1' \quad M_2' \quad M_3' \right]$$

where M is mixed super column matrix.

Now for the matrix M we draw lines between the rows 2 and 3 and lines between the columns 2 and 3 and 4 and 5 we get

$$M = \left[\begin{array}{cc|cc|cc} 3 & 4 & 5 & 7 & 8 & 9 \\ 1 & 2 & 3 & 0 & 1 & 2 \\ \hline 0 & 6 & 1 & 4 & 0 & -3 \\ -1 & 1 & 2 & 5 & 1 & 4 \end{array}\right]$$

$$= \begin{bmatrix} M_{11} & M_{12} & M_{13} \\ M_{21} & M_{22} & M_{23} \end{bmatrix}.$$

Thus M is a super matrix in fact we call M to be a square supermatrix. However M is not a perfect square supermatrix.

Now we just define the fuzzy analogue of these concepts.

**DEFINITION 1.1.12:** *Let M be a supermatrix i.e.,*

$$M = \begin{bmatrix} M_{11} & M_{12} & \cdots & M_{1n} \\ M_{21} & M_{22} & \cdots & M_{2n} \\ \vdots & \vdots & & \vdots \\ M_{m1} & M_{m2} & \cdots & M_{mn} \end{bmatrix}.$$

*If each of the matrices $M_{ij}$ are fuzzy matrices for $1 \le i \le m$ and $1 \le j \le n$; i.e., the entries of each of the matrices $M_{ij}$ are taken from the unit interval [0, 1] then we call M to be fuzzy super matrix.*

***Example 1.1.19**:* Let

$$M = \begin{bmatrix} M_{11} & M_{12} \\ M_{21} & M_{22} \\ M_{31} & M_{32} \end{bmatrix} = \left[\begin{array}{ccc|cc} 0.1 & 1 & 0.8 & 0.6 & 0.5 \\ 1 & 0.7 & 1 & 0 & 0.7 \\ 0.9 & 0.3 & 0.4 & 0.5 & 1 \\ \hline 0.7 & 1 & 0.3 & 1 & 0.8 \\ 0.8 & 0.1 & 0 & 1 & 0.1 \\ \hline 0.4 & 0.6 & 1 & 0.7 & 0.5 \\ 0 & 0.5 & 0.2 & 1 & 0.3 \\ 1 & 0 & 1 & 0.4 & 1 \end{array}\right].$$

M is a super fuzzy matrix we see each $M_{ij}$ is a fuzzy matrix; $1 \le i \le 3$ and $1 \le j \le 2$. We see some of the matrices in M are fuzzy square matrices.

***Example 1.1.20:*** Let us consider

$$M = \left[\begin{array}{c|cc|ccc|c} 0.1 & 1 & 0.3 & 0 & 0.7 & 0.2 & 0.5 \\ 0.6 & 0.7 & 1 & 0.4 & 1 & 0.3 & 1 \\ 1 & 0.5 & 0.2 & 1 & 0 & 0.5 & 0.7 \\ \hline 0 & 1 & 0 & 0.3 & 0.4 & 1 & 1 \\ 0.3 & 0.2 & 1 & 0.5 & 0.7 & 0.9 & 0.4 \\ \hline 0.9 & 0.1 & 0.3 & 1 & 0 & 0.8 & 0.2 \\ 1 & 0.3 & 0 & 0.9 & 0.7 & 1 & 1 \\ 0.5 & 1 & 0.2 & 0.8 & 0.6 & 0.4 & 0 \\ \hline 0.4 & 0 & 0.3 & 1 & 0.5 & 0 & 0.2 \\ 0 & 1 & 0.4 & 0 & 1 & 0.3 & 1 \end{array}\right]$$

$$= \begin{bmatrix} M_{11} & M_{12} & M_{13} & M_{14} \\ M_{21} & M_{22} & M_{23} & M_{24} \\ M_{31} & M_{32} & M_{33} & M_{34} \\ M_{41} & M_{42} & M_{43} & M_{41} \end{bmatrix}.$$

We see some of the $M_{ij}$ are square fuzzy matrices some of them are fuzzy row vectors and some just fuzzy column vectors.

We see ($M_{11}$ $M_{12}$ $M_{13}$ $M_{14}$) is the mixed super fuzzy row matrix.

$$\begin{bmatrix} M_{11} \\ M_{12} \\ M_{13} \\ M_{14} \end{bmatrix}$$

is just a mixed fuzzy super column vector.

Likewise ($M_{41}$ $M_{42}$ $M_{43}$ $M_{41}$) is only a mixed fuzzy super row vector. We can say ($M_{31}$ $M_{32}$ $M_{33}$ $M_{34}$) is a mixed fuzzy super row vector and each of the matrices have only three rows but different number of columns.

***Example 1.1.21****:* Let us consider the fuzzy supermatrix

$$M = \begin{bmatrix} M_{11} & M_{12} & M_{13} \\ M_{21} & M_{22} & M_{23} \\ M_{31} & M_{32} & M_{33} \end{bmatrix}$$

where each $M_{ij}$ is a fuzzy square matrix. Then we say this M is a perfect square fuzzy supermatrix.

M is given by

$$
M = \left[\begin{array}{ccc|ccc|ccc}
0.3 & 1 & 0 & 1 & 0.7 & 0.9 & 0.2 & 1 & 0.1 \\
0.2 & 1 & 1 & 0.6 & 1 & 0 & 0.4 & 0 & 0.7 \\
0.6 & 0 & 0.5 & 1 & 0.7 & 0.2 & 0.7 & 1 & 0.4 \\
\hline
1 & 0.7 & 0.9 & 0.3 & 1 & 0 & 1 & 0.6 & 0.3 \\
0.8 & 1 & 0.5 & 0.5 & 0 & 0.7 & 0.4 & 0.7 & 0.9 \\
0.6 & 0 & 0.4 & 0 & 1 & 0.2 & 0.2 & 0.1 & 0 \\
\hline
0.6 & 1 & 0.7 & 1 & 0 & 0.2 & 1 & 1 & 1 \\
1 & 0 & 0.3 & 0 & 1 & 0.3 & 0 & 0.7 & 0 \\
0 & 0.2 & 1 & 1 & 1 & 1 & 1 & 0 & 1
\end{array}\right].
$$

***Example 1.1.22:*** Consider the fuzzy supermatrix given by

$$
M = \left[\begin{array}{cc|cc|cc}
1 & 0.3 & 0.2 & 0.7 & 0.7 & 0.4 \\
0 & 0.1 & 1 & 1 & 0.2 & 1 \\
\hline
0.7 & 0.4 & 0.5 & 0 & 0.1 & 0 \\
0.3 & 0.5 & 1 & 0.6 & 0.5 & 1 \\
\hline
0.2 & 1 & 0 & 1 & 0.6 & 0.9 \\
1 & 0.2 & 0.7 & 0.4 & 0 & 1
\end{array}\right]
$$

$$
\begin{bmatrix}
M_{11} & M_{12} & M_{13} \\
M_{21} & M_{22} & M_{23} \\
M_{31} & M_{32} & M_{33}
\end{bmatrix}.
$$

M is a perfect fuzzy square supermatrix, each $M_{ij}$ is a $2 \times 2$ fuzzy square matrix.

***Example 1.1.23:*** Let M be a super fuzzy matrix

$$
M = \begin{bmatrix}
M_{11} & M_{12} & M_{13} & M_{14} \\
M_{21} & M_{22} & M_{23} & M_{24} \\
M_{31} & M_{32} & M_{33} & M_{34}
\end{bmatrix} =
$$

$$\left[\begin{array}{ccc|ccc|ccc|ccc}
0.3 & 1 & 0.7 & 1 & 0 & 0.9 & 0.7 & 0.6 & 0.2 & 0.9 & 1 & 1 \\
1 & 0.6 & 1 & 0.6 & 0.2 & 1 & 1 & 0 & 1 & 0 & 0 & 0.7 \\
0.5 & 0.7 & 0.8 & 0.5 & 1 & 0.3 & 1 & 1 & 1 & 0.6 & 0.5 & 1 \\
\hline
0.1 & 0.2 & 0.3 & 0.4 & 0.5 & 0.6 & 0.7 & 0.8 & 0.9 & 1 & 0 & 0 \\
1 & 1 & 0 & 0.3 & 0.1 & 0.2 & 0.4 & 0.5 & 0.7 & 0.6 & 0.8 & 0.9 \\
0 & 0.6 & 1 & 1 & 0 & 1 & 1 & 1 & 0 & 0.2 & 1 & 1 \\
\hline
1 & 0.2 & 0.8 & 0 & 0.9 & 0.8 & 0.2 & 0.4 & 0.1 & 0.7 & 0.6 & 0.3 \\
0.9 & 1 & 0.7 & 0.9 & 0 & 1 & 0.3 & 0.7 & 1 & 1 & 0 & 0.5 \\
0.6 & 0.3 & 1 & 0.8 & 1 & 0 & 1 & 0 & 0.3 & 0.8 & 1 & 1
\end{array}\right]$$

M is a fuzzy super square matrix but it is clearly not a perfect fuzzy super square matrix

Now as in case of matrices we can also partition fuzzy matrices and make it into a fuzzy supermatrix.

We just illustrate this by an example.

***Example 1.1.24:*** Let us consider a $6 \times 7$ fuzzy matrix

$$M = \begin{bmatrix}
0.7 & 0.2 & 0.7 & 0.6 & 0.5 & 0.4 & 0.3 \\
0.6 & 0.1 & 0.6 & 0.2 & 0.1 & 1 & 0 \\
0.5 & 0 & 0.5 & 0.3 & 0.2 & 0.1 & 0.9 \\
0.4 & 1 & 0.1 & 1 & 0.5 & 0 & 0.4 \\
0.3 & 0.9 & 1 & 0.3 & 1 & 0.2 & 1 \\
1 & 0.8 & 0.2 & 0 & 0.1 & 1 & 0
\end{bmatrix}.$$

Suppose we draw a line between the columns 2 and 3 and 6 and 7 respectively. Let us draw lines between the rows 1 and 2 and between 4 and 5 respectively.

The resultant M is as follows. We denote it by $M_p$ i.e., M after partitioning.

$$M_P = \left[\begin{array}{cc|cccc|c} 0.7 & 0.2 & 0.7 & 0.6 & 0.5 & 0.4 & 0.3 \\ \hline 0.6 & 0.1 & 0.6 & 0.2 & 0.1 & 1 & 0 \\ 0.5 & 0 & 0.5 & 0.3 & 0.2 & 0.1 & 0.9 \\ 0.4 & 1 & 0.1 & 1 & 0.5 & 0 & 0.4 \\ \hline 0.3 & 0.9 & 1 & 0.3 & 1 & 0.2 & 1 \\ 1 & 0.8 & 0.2 & 0 & 0.1 & 1 & 0 \end{array}\right]$$

$$= \begin{bmatrix} M_{11} & M_{12} & M_{13} \\ M_{21} & M_{22} & M_{23} \\ M_{31} & M_{32} & M_{33} \end{bmatrix}$$

where ($M_{11}$ $M_{12}$ $M_{13}$) is a mixed fuzzy super row vector where as ($M_{21}$ $M_{22}$ $M_{23}$) and ($M_{31}$ $M_{32}$ $M_{33}$) are fuzzy super row matrices.

Further

$$\begin{bmatrix} M_{11} \\ M_{21} \\ M_{31} \end{bmatrix} \text{ and } \begin{bmatrix} M_{12} \\ M_{22} \\ M_{32} \end{bmatrix}$$

are fuzzy super column matrices where as

$$\begin{bmatrix} M_{13} \\ M_{23} \\ M_{33} \end{bmatrix}$$

is only just a mixed fuzzy super column vector.

Now having seen some special types of fuzzy supermatrices, we now proceed on to recall how the transpose of a supermatrix is defined.

**DEFINITION 1.1.13:** *Let X = ($X_1$ $X_2$ … $X_n$) be a mixed super row vector the transpose of X denoted by*

$$X^t = (X_1\ X_2\ \ldots\ X_n)^t$$
$$= \left(X_1^t\ X_2^t \cdots X_n^t\right)^t$$

*is a mixed super column vector. Likewise if*

$$Y = \begin{bmatrix} Y_1 \\ Y_2 \\ \vdots \\ Y_m \end{bmatrix}$$

*is a mixed super column vector then*

$$Y_t = \begin{bmatrix} Y_1 \\ Y_2 \\ \vdots \\ Y_m \end{bmatrix}^t = \begin{bmatrix} Y_1^t \\ Y_2^t \\ \vdots \\ Y_m^t \end{bmatrix}^t$$

*is a mixed super row vector.*

We illustrate this by the following examples.

***Example 1.1.25:*** Let

$$X = (X_1\ X_2\ X_3\ X_4)$$
$$= (3\ 1\ 4\ 5\ |\ 0\ –\ 3\ 2\ 1\ –\ 4\ 7\ |\ 7\ 0\ 3\ 1\ 2\ |\ 0\ 3\ –\ 5)$$

be the given mixed super row vector.

Now

$$X^t = (X_1\ X_2\ X_3\ X_4)^t$$
$$= \left(X_1^t\ X_2^t\ X_3^t\ X_4^t\right)^t$$

$$= \left[\begin{array}{c} 3 \\ 1 \\ 4 \\ 5 \\ \hline 0 \\ -3 \\ 2 \\ 1 \\ -4 \\ 7 \\ \hline 7 \\ 0 \\ 3 \\ 1 \\ 2 \\ \hline 0 \\ 3 \\ -5 \end{array}\right].$$

***Example 1.1.26:*** Let

$$Y = \left[\begin{array}{c} Y_1 \\ Y_2 \\ Y_3 \\ Y_4 \\ Y_5 \end{array}\right]$$

be the given mixed super column vector where

$$Y = \begin{bmatrix} 1 \\ 0 \\ -3 \\ 2 \\ \hline -3 \\ 5 \\ \hline 7 \\ 6 \\ 2 \\ \hline 1 \\ 2 \\ 3 \\ 4 \\ 5 \\ \hline -7 \\ 1 \\ 2 \end{bmatrix}.$$

Now

$$Y^t = \begin{bmatrix} Y_1 \\ Y_2 \\ Y_3 \\ Y_4 \\ Y_5 \end{bmatrix}^t = \begin{bmatrix} Y_1^t \\ Y_2^t \\ Y_3^t \\ Y_4^t \\ Y_5^t \end{bmatrix}^t = \left[1\ 0\ \text{-}3\ 2 \mid \text{-}3\ 5 \mid 7\ 6\ 2 \mid 1\ 2\ 3\ 4\ 5 \mid \text{-}7\ \ 1\ 2\right]$$

is clearly the mixed super row vector.

Now we proceed on to define the transpose of super row matrix and super column matrix.

**DEFINITION 1.1.14:** *Let $A = [A_1 \mid A_2 \mid \ldots \mid A_n]$ be a super row matrix where $A_i$ are $p \times t_i$ matrices; $1 \leq i \leq n$; denoted by $A^t = [A_1 \mid A_2 \mid A_3 \mid \ldots \mid A_n]^t = [A_1^t \mid A_2^t \mid \ldots \mid A_n^t]^t$*

$$= \left[\begin{array}{c} A_1^t \\ \hline A_2^t \\ \hline \vdots \\ \hline A_n^t \end{array}\right].$$

*$A^t$ is a super column matrix. Thus the transpose of a super row matrix is a super column matrix. Similarly the transpose of a super column matrix is a super row matrix given by*

$$Y = \left[\begin{array}{c} Y_1 \\ \hline Y_2 \\ \hline \vdots \\ \hline Y_m \end{array}\right],$$

*the column matrix where $Y_i$'s are $s_i \times t$ matrices $i = 1, 2, \ldots, m$. The transpose of Y is denoted by*

$$Y^t = \left[\begin{array}{c} Y_1 \\ \hline Y_2 \\ \hline \vdots \\ \hline Y_m \end{array}\right]^t = \left[\begin{array}{c} Y_1^t \\ \hline Y_2^t \\ \hline \vdots \\ \hline Y_m^t \end{array}\right]^t = \left[Y_1^t \mid Y_2^t \mid \ldots \mid Y_m^t\right]$$

*is the super row matrix.*

Now we proceed onto illustrate this by the following examples.

***Example 1.1.27:*** Let X = [$X_1$ | $X_2$ | $X_3$ | $X_4$] be a super row matrix where

$$X = \left[\begin{array}{cccc|cc|ccccc|c} 2 & 1 & 0 & 5 & 3 & 1 & 2 & 3 & 4 & 5 & 6 & 1 \\ 0 & 7 & 2 & -1 & 1 & 4 & 1 & 1 & 0 & 2 & 1 & 7 \\ 1 & 0 & 1 & 2 & -1 & 0 & 0 & 7 & 8 & 9 & -1 & 8 \\ -1 & 5 & 0 & 3 & 0 & 1 & 2 & 3 & 4 & 5 & 6 & -4 \end{array}\right].$$

$$X^t = [X_1 \mid X_2 \mid X_3 \mid X_4]^t = [X_1^t \mid X_2^t \mid X_3^t \mid X_4^t]^t$$

$$= \left[\begin{array}{cccc} 2 & 0 & 1 & -1 \\ 1 & 7 & 0 & 5 \\ 0 & 2 & 1 & 0 \\ 5 & -1 & 2 & 3 \\ \hline 3 & 1 & -1 & 0 \\ 1 & 4 & 0 & 1 \\ \hline 2 & 1 & 0 & 2 \\ 3 & 1 & 7 & 3 \\ 4 & 0 & 8 & 4 \\ 5 & 2 & 9 & 5 \\ 6 & 1 & -1 & 6 \\ \hline 1 & 7 & 8 & -4 \end{array}\right].$$

Clearly $X^t$ is a super column matrix.

***Example 1.1.28:*** Let

$$Y = \begin{bmatrix} Y_1 \\ Y_2 \\ Y_3 \\ Y_4 \\ Y_5 \end{bmatrix} = \left[\begin{array}{ccccc} 3 & 1 & 0 & 1 & 5 \\ 0 & 1 & 7 & 0 & -2 \\ \hline 1 & 4 & 1 & 2 & 1 \\ 9 & 5 & 8 & 3 & 0 \\ 0 & 1 & 2 & 3 & 4 \\ 5 & 6 & 7 & 8 & 9 \\ \hline 1 & 0 & 2 & 3 & 7 \\ \hline 1 & 0 & 1 & 0 & 1 \\ 0 & 1 & 0 & 1 & 0 \\ 1 & 1 & 1 & 1 & 1 \\ \hline 8 & 0 & -1 & 0 & 1 \\ 1 & 1 & 0 & 2 & 0 \end{array}\right]$$

be the given super column matrix. Now

$$Y^t = \begin{bmatrix} Y_1 \\ Y_2 \\ \vdots \\ Y_5 \end{bmatrix}^t = \begin{bmatrix} Y_1^t \\ Y_2^t \\ \vdots \\ Y_5^t \end{bmatrix}^t = [Y_1^t \mid Y_2^4 \mid ... \mid Y_5^t] =$$

$$\left[\begin{array}{cc|cccc|c|ccc|cc} 3 & 0 & 1 & 9 & 0 & 5 & 1 & 1 & 0 & 1 & 8 & 1 \\ 1 & 1 & 4 & 5 & 4 & 6 & 0 & 0 & 1 & 1 & 0 & 1 \\ 0 & 7 & 1 & 8 & 2 & 7 & 2 & 1 & 0 & 1 & -1 & 0 \\ 1 & 0 & 2 & 8 & 3 & 8 & 3 & 0 & 1 & 1 & 0 & 2 \\ 5 & -2 & 1 & 0 & 4 & 9 & 7 & 1 & 0 & 1 & 1 & 0 \end{array}\right],$$

clearly $Y^t$ is a super row matrix.

Now we proceed on to define the transpose of a supermatrix M.

**DEFINITION 1.1.15:** *Let M be a supermatrix given by*

$$M = \left[\begin{array}{c|c|c|c} M_{11} & M_{12} & \cdots & M_{1n} \\ \hline M_{21} & M_{22} & \cdots & M_{2n} \\ \hline \vdots & \vdots & & \vdots \\ \hline M_{m1} & M_{m2} & \cdots & M_{mn} \end{array}\right]$$

*where each $M_{ij}$ is a $s_i \times t_j$ matrix $1 \leq i \leq m$ and $1 \leq j \leq n$. Now*

$$M^t = \left[\begin{array}{c|c|c|c} M_{11} & M_{12} & \cdots & M_{1n} \\ \hline M_{21} & M_{22} & \cdots & M_{2n} \\ \hline \vdots & \vdots & & \vdots \\ \hline M_{m1} & M_{m2} & \cdots & M_{mn} \end{array}\right]^t = \left[\begin{array}{c|c|c|c} M_{11}^t & M_{12}^t & \cdots & M_{1n}^t \\ \hline M_{21}^t & M_{22}^t & \cdots & M_{2n}^t \\ \hline \vdots & \vdots & & \vdots \\ \hline M_{m1}^t & M_{m2}^t & \cdots & M_{mn}^t \end{array}\right]^t$$

$$
= \left[\begin{array}{c|c|c|c}
M_{11}^{t} & M_{21}^{t} & \cdots & M_{m1}^{t} \\ \hline
M_{12}^{t} & M_{22}^{t} & \cdots & M_{m2}^{t} \\ \hline
\vdots & \vdots & & \vdots \\ \hline
M_{1n}^{t} & M_{2n}^{t} & \cdots & M_{mn}^{t}
\end{array}\right].
$$

*$M^t$ is again a supermatrix.*

We illustrate this situation by some examples.

***Example 1.1.29:*** Let M be a supermatrix

$$
M = \begin{bmatrix}
M_{11} & M_{12} & M_{13} & M_{14} \\
M_{21} & M_{22} & M_{23} & M_{24} \\
M_{31} & M_{32} & M_{33} & M_{34}
\end{bmatrix}
$$

$$
= \left[\begin{array}{ccc|ccccc|cc|cccc}
3 & 0 & 1 & 7 & 3 & 4 & 5 & 1 & 3 & 4 & 1 & 0 & 3 & 4 \\
5 & 7 & 8 & 0 & 1 & 3 & 4 & 0 & 1 & 2 & 0 & 1 & 5 & 6 \\
8 & 1 & 0 & 1 & 1 & 0 & 0 & 1 & 0 & 1 & 9 & 0 & 2 & 0 \\ \hline
2 & 0 & 5 & 1 & 1 & 0 & 3 & 4 & 1 & 1 & 1 & 1 & 1 & 1 \\
7 & 8 & 0 & 0 & 1 & 5 & 7 & 2 & 0 & 1 & 0 & 1 & 1 & 0 \\
1 & 0 & 1 & -1 & 0 & 0 & 1 & 0 & 1 & 0 & 2 & 0 & 0 & 1 \\
9 & 8 & 4 & -7 & 8 & 9 & 0 & 1 & 0 & 0 & 1 & 2 & 3 & 4 \\ \hline
0 & 1 & 1 & 1 & 2 & 3 & 4 & 5 & 7 & 2 & 1 & 1 & 0 & 7 \\
1 & 0 & 1 & 3 & 0 & 5 & 0 & 7 & 1 & 0 & 7 & 0 & 1 & 1
\end{array}\right].
$$

Now

$$
M^t = \left[\begin{array}{c|c|c|c}
M_{11} & M_{12} & M_{13} & M_{14} \\ \hline
M_{21} & M_{22} & M_{23} & M_{24} \\ \hline
M_{31} & M_{32} & M_{33} & M_{34}
\end{array}\right]^t
$$

$$= \left[\begin{array}{c|c|c|c} M_{11}^t & M_{12}^t & M_{13}^t & M_{14}^t \\ \hline M_{21}^t & M_{22}^t & M_{23}^t & M_{24}^t \\ \hline M_{31}^t & M_{32}^t & M_{33}^t & M_{34}^t \end{array}\right]^t$$

$$= \left[\begin{array}{c|c|c} M_{11}^t & M_{21}^t & M_{31}^t \\ \hline M_{12}^t & M_{22}^t & M_{32}^t \\ \hline M_{13}^t & M_{23}^t & M_{33}^t \\ \hline M_{14}^t & M_{24}^t & M_{34}^t \end{array}\right].$$

i.e.,

$$M^t = \left[\begin{array}{ccc|cccc|cc} 3 & 5 & 8 & 2 & 7 & 1 & 9 & 0 & 1 \\ 0 & 7 & 1 & 0 & 8 & 0 & 8 & 1 & 0 \\ 1 & 8 & 0 & 5 & 0 & 1 & 4 & 1 & 1 \\ \hline 7 & 0 & 1 & 1 & 0 & -1 & -7 & 1 & 3 \\ 3 & 1 & 1 & 1 & 1 & 0 & 8 & 2 & 0 \\ 4 & 3 & 0 & 0 & 5 & 0 & 9 & 3 & 5 \\ 5 & 4 & 0 & 3 & 7 & 1 & 0 & 4 & 0 \\ 1 & 0 & 1 & 4 & 2 & 0 & 1 & 5 & 7 \\ \hline 3 & 1 & 0 & 1 & 0 & 1 & 0 & 7 & 1 \\ 4 & 2 & 1 & 1 & 1 & 0 & 0 & 2 & 0 \\ \hline 1 & 0 & 9 & 1 & 0 & 2 & 1 & 1 & 7 \\ 0 & 1 & 0 & 1 & 1 & 0 & 2 & 1 & 0 \\ 3 & 5 & 2 & 1 & 1 & 0 & 3 & 0 & 1 \\ 4 & 6 & 0 & 1 & 0 & 1 & 4 & 7 & 1 \end{array}\right].$$

Now as in case of simple or ordinary matrices we see transpose of $M^t$ i.e., $(M^t)^t = M$.

Now we illustrate and define some simple rules of supermatrix multiplication. For more information about supermatrices please refer [125, 316]. We just see $M = [M_1 \mid M_2 \mid \ldots \mid M_n]$ is a super row matrix where each $M_i$ is a $t \times p_i$. Let $X = (x_1, x_2, \ldots, x_t)$ be a ordinary or simple row vector. Then we define $X \circ M$ the product of row vector with the super row matrix is defined as follows $X \circ M = (x_1, x_2, \ldots , x_5) \circ [M_1 \mid M_2 \mid \ldots \mid M_n] = [N_1 \mid N_2 \mid \ldots \mid N_n]$ where $N_i$ is a $1 \times p_i$ row matrix and $[N_1 \mid N_2 \mid \ldots \mid N_n]$ is a super row vector. Thus the product of the row vector with a super row matrix yields a super row vector.

We illustrate this situation by the following example.

***Example 1.1.30:*** Let $M = [M_1 \mid M_2 \mid M_3 \mid M_4]$ be a super row matrix given by

$$M = \left[\begin{array}{ccccc|cc|cccc|c} 4 & 1 & 0 & 2 & 1 & 3 & 1 & 1 & 2 & 3 & 1 & 2 \\ 1 & 0 & 1 & 5 & 0 & 0 & 5 & 3 & 1 & 0 & 1 & 3 \\ -1 & 6 & 0 & 1 & 2 & 7 & 2 & 1 & 0 & 1 & 1 & 1 \\ 3 & 0 & 1 & 2 & 3 & 0 & 4 & 0 & 1 & 0 & 1 & 0 \end{array}\right].$$

Let X = (1 0 1 0) be the given row vector. Now

$$X \circ M = [1\ 0\ 1\ 0] \circ \left[\begin{array}{ccccc|cc|cccc|c} 4 & 1 & 0 & 2 & 1 & 3 & 1 & 1 & 2 & 3 & 1 & 2 \\ 1 & 0 & 1 & 5 & 0 & 0 & 5 & 3 & 1 & 0 & 1 & 3 \\ -1 & 6 & 0 & 1 & 2 & 7 & 2 & 1 & 0 & 1 & 1 & 1 \\ 3 & 0 & 1 & 2 & 3 & 0 & 4 & 0 & 1 & 0 & 1 & 0 \end{array}\right]$$

$$= \left[ (1\ 0\ 1\ 0)\begin{bmatrix} 4 & 1 & 0 & 2 & 1 \\ 1 & 0 & 1 & 5 & 0 \\ -1 & 6 & 0 & 1 & 2 \\ 3 & 0 & 1 & 2 & 3 \end{bmatrix} \middle| (1\ 0\ 1\ 0)\begin{bmatrix} 3 & 1 \\ 0 & 5 \\ 7 & 2 \\ 0 & 4 \end{bmatrix} \middle| \right.$$

$$\left. (1\ 0\ 1\ 0)\begin{bmatrix} 1 & 2 & 3 & 1 \\ 3 & 1 & 0 & 1 \\ 1 & 0 & 1 & 1 \\ 0 & 1 & 0 & 1 \end{bmatrix} \middle| (1\ 0\ 1\ 0)\begin{bmatrix} 2 \\ 3 \\ 1 \\ 0 \end{bmatrix} \right]$$

$$= [\ 3\ 7\ 0\ 3\ 3 \mid 10\ 3 \mid 2\ 2\ 4\ 2 \mid 3\ ]$$

is a mixed row vector.

Now suppose we have a super column matrix

$$T = \left[\begin{array}{c} T_1 \\ \hline T_2 \\ \hline \vdots \\ \hline T_m \end{array}\right]$$

where each $T_i$ is a $p_i \times t$ matrix, i = 1, 2, …, m, i.e., each $T_i$ has only t columns. Suppose $X = [X_1 \mid X_2 \mid \ldots \mid X_m]$ is a super row vector where each $X_i$ has $p_i$ elements in it then we can define the special product denoted by $X \circ_s T$ which gives a row vector. That is

$$X = \left[ x_1^1 \ldots x_{p_1}^1 \mid x_1^2 \ldots x_{p_2}^2 \mid x_1^m \ldots x_{p_m}^m \right] \times$$

$$\left[\begin{array}{cccc} t_{11}^1 & t_{12}^1 & \cdots & t_{1t}^1 \\ \vdots & \vdots & & \vdots \\ t_{p_1 1}^1 & t_{p_2 1}^1 & \cdots & t_{p_1 t}^1 \\ \hline t_{11}^2 & t_{12}^2 & \cdots & t_{1t}^2 \\ \vdots & \vdots & \cdots & \vdots \\ t_{p_2 1}^2 & t_{p_2 2}^2 & \cdots & t_{p_2 t}^2 \\ \hline \vdots & \vdots & & \vdots \\ \hline t_{11}^m & t_{12}^m & \cdots & t_{1t}^m \\ t_{pm_1}^m & t_{pm_2}^m & \cdots & t_{pm_t}^m \end{array}\right]$$

$$= \left[ \sum_{i=1}^{m} x_j^i \; t_{j1}^i \; , \sum_{i=1}^{m} x_j^i \; t_{j2}^i \; ,\ldots, \sum_{i=1}^{m} x_j^i \; t_{jt}^i \right]$$

$$= (y_1, y_2, \ldots, y_t)$$

is the row vector.

We illustrate this situation by the following example.

***Example 1.1.31:*** Let

$$M = \left[\begin{array}{c} M_1 \\ \hline M_2 \\ \hline M_3 \\ \hline M_4 \end{array}\right]$$

$$= \left[\begin{array}{ccc} 3 & 0 & 1 \\ 1 & 1 & 2 \\ 0 & 1 & 1 \\ 0 & 1 & 0 \\ \hline 1 & 0 & 2 \\ 0 & 1 & 0 \\ \hline 1 & 1 & 1 \\ 1 & 0 & 1 \\ 1 & 1 & 0 \\ 0 & 1 & 1 \\ 0 & 0 & 1 \\ 1 & 1 & 1 \\ \hline 2 & 1 & 3 \end{array}\right].$$

Let X = (0 1 0 1 | 1 0 | 1 0 1 0 0 0 | 1) be the mixed super row vector (given). $X \circ_s T = [6\ 5\ 8]$, i.e. we just ignore the partition but straight multiply it as if we are multiplying the usual $1 \times 13$ row vector with $13 \times 3$ matrix. We do ignore the partition for the compatibility of the resultant.

We give yet another example for the reader to have a clear understanding of this product.

***Example 1.1.32:*** Let us consider a super column vector

$$M = \begin{bmatrix} M_1 \\ M_2 \\ M_3 \\ M_4 \\ M_5 \end{bmatrix} = \left[\begin{array}{ccccc} 1 & 0 & 3 & 4 & 5 \\ 0 & 1 & 0 & 1 & 0 \\ \hline 1 & 2 & 0 & 0 & 3 \\ 5 & 6 & 0 & 1 & 0 \\ 7 & 1 & 0 & 0 & 1 \\ 1 & 0 & 1 & 1 & 0 \\ \hline 7 & 1 & 0 & 0 & 0 \\ 1 & 0 & 1 & 1 & 1 \\ 0 & 1 & 0 & 0 & 1 \\ 1 & 1 & 0 & 1 & 1 \\ 1 & 1 & 1 & 0 & 1 \\ \hline 1 & 0 & 1 & 0 & 1 \\ \hline 2 & 1 & 0 & 2 & 0 \\ 0 & 2 & 0 & 0 & 1 \\ 1 & 0 & 1 & 0 & 0 \end{array}\right].$$

Suppose X = (1 0 | 1 0 0 1 | 0 1 0 0 1 | 1 | 0 1 1) be the mixed super vector given, we find

$$X \circ_s T = [7\ 5\ 8\ 6\ 12]$$

is a simple row vector i.e., is a $1 \times 5$ row matrix.

Here it is pertinent to mention that this type of product is not defined in [125]; we have specially done this mainly to suit our super fuzzy models which we will be defining in the next section and using them in chapter two of this book.

It is still important to mention here that we are not very much interested in the way in which supermatrix products are defined. When we use them in fuzzy super models, we need compatibility further the super row vectors action in several situations give only a simple row vector so we use this method of product and call it as the special product of super matrices. We may in case of fuzzy supermatrices replace this special product by special min max functions or operator or max min operator according to our need.

Now we proceed on to define the special product of super matrices with super row vector. It is once again pertinent to mention that we are only interested in studying these special types of product, which always produces a super column matrix or super row matrix or a supermatrix only with a simple row vector or a super row vector.

For none of our fuzzy models which we will be defining in this book need product of a super row matrix with super column matrix or a product of a supermatrix with another supermatrix and so on.

**DEFINITION 1.1.16:** *Let*

$$M = \begin{bmatrix} M_{11} & M_{12} & \dots & M_{1n} \\ M_{21} & M_{22} & \cdots & M_{2n} \\ \vdots & \vdots & & \vdots \\ M_{m1} & M_{m2} & \cdots & M_{mn} \end{bmatrix}$$

*be a supermatrix where each $M_{ij}$ is a $p_i \times q_i$ matrix, $1 \le i \le m$ and $1 \le j \le n$. Consider*

$$X = \left( x_1^1 \ldots x_{p_1}^1 \;\middle|\; x_1^2 \; x_2^2 \ldots x_{p_2}^2 \;\middle|\; \ldots \;\middle|\; x_1^m \ldots x_{p_m}^m \right)$$

*a mixed super row vector.*

$$X \circ M = \left( y_1^1 \ldots y_{q_1}^1 \;\middle|\; \ldots \;\middle|\; y_1^n \ldots y_{q_n}^n \right)$$

*the product done as in case of a super row vector with a super column matrix exactly ($q_1 + q_2 + \ldots + q_n$) times, this product is called the super special product of the mixed super vector with a supermatrix.*

We illustrate this situation by two examples.

***Example 1.1.33:*** Consider the supermatrix

$$M = \begin{bmatrix} M_{11} & M_{12} & M_{13} \\ M_{21} & M_{22} & M_{23} \\ M_{31} & M_{32} & M_{33} \\ M_{41} & M_{42} & M_{43} \end{bmatrix}$$

$$= \left[\begin{array}{ccc|cccc|ccccc} 0 & 1 & 2 & 0 & 1 & 0 & 1 & 1 & 0 & 1 & 0 & 1 \\ 1 & 1 & 0 & 1 & 1 & 0 & 0 & 0 & 1 & 0 & 1 & 0 \\ \hline 1 & 0 & 1 & 1 & 0 & 0 & 0 & 1 & 0 & 1 & 1 & 0 \\ 0 & 1 & 0 & 0 & 1 & 0 & 0 & 0 & 1 & 0 & 0 & 1 \\ 2 & 1 & 2 & 0 & 0 & 0 & 1 & 0 & 0 & 1 & 0 & 1 \\ 3 & 1 & 2 & 0 & 0 & 1 & 0 & 1 & 1 & 0 & 1 & 0 \\ \hline 1 & 0 & 4 & 1 & 1 & 0 & 3 & 0 & 1 & 0 & 2 & 1 \\ \hline 1 & 1 & 1 & 1 & 0 & 0 & 1 & 0 & 0 & 1 & 1 & 1 \\ 0 & 1 & 1 & 1 & 1 & 0 & 0 & 0 & 0 & 0 & 1 & 1 \\ 0 & 0 & 1 & 0 & 0 & 1 & 1 & 0 & 0 & 0 & 0 & 1 \end{array}\right].$$

Let X = [1 0 | 0 1 1 0 | 1 | 0 1 0] be the given mixed super row vector. Now

$$X \circ_s M \quad = \quad [3\ 4\ 9 \mid 2\ 4\ 0\ 5 \mid 1\ 2\ 2\ 3\ 5].$$

Clearly $X \circ_s M$ is a mixed super row vector.

Now we give yet another example.

***Example 1.1.34:*** Let M be a supermatrix given by

$$M = \begin{bmatrix} M_{11} & M_{12} & M_{13} & M_{14} \\ M_{21} & M_{22} & M_{23} & M_{24} \\ M_{31} & M_{32} & M_{33} & M_{34} \end{bmatrix}$$

$$= \left[\begin{array}{ccc|cccc|ccccc|c}
1 & 0 & 1 & -1 & 1 & 0 & 1 & 1 & 0 & 1 & 0 & 1 & 1 \\
0 & 1 & 0 & 0 & 0 & 1 & 0 & 0 & 1 & 0 & 1 & 0 & 0 \\
\hline
1 & 1 & 0 & 0 & 0 & 1 & 0 & 1 & 0 & 0 & 0 & 0 & 0 \\
1 & 0 & 1 & 1 & 0 & 0 & 0 & 0 & 1 & 0 & 0 & 0 & 1 \\
0 & 0 & 1 & 0 & 1 & 0 & 0 & 0 & 0 & 1 & 0 & 0 & 0 \\
1 & 0 & 0 & 0 & 0 & 0 & 1 & 0 & 0 & 1 & 0 & 0 & 0 \\
\hline
1 & 0 & 1 & 1 & 0 & 0 & 1 & 0 & 0 & 1 & 0 & 0 & 0 \\
1 & 0 & 0 & 0 & 1 & 0 & 0 & 1 & 0 & 0 & 0 & 1 & 1 \\
1 & 1 & 0 & 0 & 0 & 1 & 0 & 0 & 0 & 0 & 1 & 0 & 0
\end{array}\right].$$

Suppose X = [0 1 | 0 0 1 1 | 0 1 0] be the given mixed super row vector to find $X \circ_s M$. $X \circ_s M$ = [2 1 1 | 0 2 1 1 | 1 1 2 1 1 | 1] = Y. Y is again a mixed super row vector. We can also find $Y \circ_s M^T$ = [11 4 | 5 5 5 5 | 6 7 5] = $X_1$. $X_1$ is again a mixed super row vector.

**Remark:** In the super fuzzy models to find the hidden pattern of the dynamical system we need to find the effect of a mixed super row vector X on a supermatrix M using the special product and suppose Y is the resultant we find $Y \circ_s M^T$, say the resultant is $X_1$, we then calculate $X_1 \circ_s M$ and so on. That is why we have defined in these supermatrices special product of them. We further wish to mention we would be using fuzzy matrices in case of fuzzy matrices, we will be using some other special operator just as $\circ_s$, the special product.

However we illustrate this situation also by a simple example.

***Example 1.1.35:*** Let

$$M = \begin{bmatrix} M_1 \\ M_2 \\ M_3 \\ M_4 \end{bmatrix}$$

be a fuzzy super column matrix given by

$$M = \left[\begin{array}{ccccc} 0 & 1 & 0 & 1 & 0 \\ 1 & 1 & 0 & 0 & 1 \\ 0 & 1 & 0 & 1 & 1 \\ \hline -1 & 1 & 1 & 1 & 0 \\ 0 & -1 & 0 & 1 & 1 \\ 1 & 0 & 1 & 1 & 1 \\ 1 & 0 & 1 & 1 & 1 \\ \hline -1 & 0 & 1 & 1 & 1 \\ 0 & 1 & 1 & 1 & 0 \\ \hline 0 & 1 & 0 & 1 & 1 \\ 1 & 1 & 1 & 1 & 1 \\ 0 & 1 & 1 & 0 & 1 \end{array}\right].$$

Let X = [1 0 0 | 0 0 0 1 | 0 1 | 0 1 0] be the mixed fuzzy row vector. We define a new operation which we call as super fuzzy special product denoted by '$\odot_s$'. $X \odot_s M$ = [2 3 3 4 2].

Now we in case of $\odot_s$ threshold the resultant $X \odot_s M$ as follows;

If

$$X \odot_s M = [a_1\ a_2\ a_3\ a_4]$$

put $a_i = 1$ if $a_i \geq 0$; put $a_i = 0$ if $a_i \leq 0$.

So that

$$X \odot_s M = [1\ 1\ 1\ 1\ 1] = Y \text{ (say)}.$$

Now we find

$$\begin{aligned} Y \odot_s M^T &= [2\ 3\ 3 \mid 2\ 1\ 4\ 4 \mid 2\ 3 \mid 3\ 5\ 3] \\ &\hookrightarrow [1\ 1\ 1 \mid 1\ 1\ 1\ 1 \mid 1\ 1 \mid 1\ 1\ 1]. \end{aligned}$$

'$\hookrightarrow$' denotes the row vector has been thresholded and updated.

## 1.2 Super Fuzzy Relational Maps

In this section for the first time we define the notion of super column fuzzy relational maps, mixed super row fuzzy relational maps and super fuzzy relational maps. For the definition and properties of FRMs refer [337-8,367].

The three new models have been introduced mainly to study the problem relating media and reservation for the first time. To this end we have just in the section 1.1 introduced the new notion of fuzzy supermatrices, fuzzy super column matrices, fuzzy super row matrices, fuzzy super row vector and fuzzy super column vector. Throughout this book by simple or ordinary row vector or column vector we just mean the usual row vector or column vector respectively.

Further in this section we will be using only fuzzy super matrices, fuzzy super row matrices and so on. Further we call a supermatrix, which has their entries from the set {–1, 0, 1} to be a fuzzy supermatrix. This convention is being used by several of the fuzzy theorists when they construct fuzzy models, like Fuzzy Cognitive Maps (FCMs) and Fuzzy Relational Maps (FRMs).

Now we proceed on to define and describe super column fuzzy relational maps model.

**DEFINITION 1.2.1:** *Suppose we have some n sets of experts. These n sets of experts form some n distinct category of groups may be based on age or profession or education or gender or any other common factors. This model can be described as a multi set of experts model i.e. we have sets of experts i.e.; not a multi expert model but multi set of experts model. Thus we have n sets of experts each set may contain different special features. But the only common factor is that they all agree to work upon the same problem with a same set of attributes.*

*The model which we would now be defining would be called as the super column fuzzy relational maps for this model is constructed using the fuzzy super column matrix hence we to*

*specify the use of fuzzy super column matrix call this model as super column Fuzzy Relational Maps model (super column FRM model).*

*We have just n sets of experts working on the problem with some m sets of attributes say $A_1, A_2, \ldots, A_m$. Now suppose the n sets of experts are such that the first set has $n_1$ experts, the second set has $n_2$ experts and so on and the nth set has $n_n$ number of experts.*

*Let*

$$M = \left[\frac{M_1}{\frac{M_2}{\frac{\vdots}{M_n}}}\right]$$

*be a fuzzy super column matrix which is the super column dynamical system of the super column FRM model; having m number of columns and the number of rows of $M_1$ will be $n_1$, that of $M_2$ will be $n_2$ and so on and that of $M_n$ will be $n_n$. Now $M_i$ will be the connection fuzzy matrix of the $n_i^{th}$ set of expert related with the fuzzy relational maps model given by the $n_i^{th}$ set of experts using the m-set of attributes $A_1, \ldots, A_m$, having $n_i$ number of experts. In this way $M_1$ will be the associated connection fuzzy matrix of the first set of experts which has $n_1$ number of experts and so on. Thus*

$$M = \left[\frac{M_1}{\frac{M_2}{\frac{\vdots}{M_n}}}\right]$$

*will be the fuzzy super column matrix associated with the n sets of experts $n_1, n_2, \ldots, n_n$ called as the super column fuzzy relational maps model, defined using FRMs.*

*Now we just give what are the domain and the range spaces of this super column fuzzy relational maps model. Clearly the domain space is a fuzzy super mixed row vector relating all the n sets of experts who have worked with the model. Any element X of the domain multi space D would be of the form*

$$\left(x_1^1\ x_2^1\ \ldots\ x_n^1\ |\ x_1^2\ x_2^2\ \ldots\ x_{n_2}^2\ \ |\ \ldots\ |\ x_1^n\ x_2^n\ \ldots\ x_{n_n}^n\right).$$

*If we consider X = [1 0 0 … 0 | 0 0 1 … | … | 0 0 … 0 1 0] this implies we have the opinion of the first expert from the $n_1^{th}$ set of experts, the third experts opinion from the $n_2^{th}$ set of experts and so on and $n_{(n-1)}^{th}$ experts opinion from the $n_n^{th}$ set of experts. Further all other experts at that moment are remaining silent for that special state vector X. Further the state vector in the domain space in this case always is a mixed fuzzy super row vector, the state vectors take only values from the set {0, 1}.*

*The range space of this model has state vectors which are simple row vectors taking its entries from the set {0, 1}. Now we just indicate how this model works. Suppose M is the mixed fuzzy super column matrix associated with the super column fuzzy relational maps model. We have the domain space D to be a mixed fuzzy super row vector with entries only from the set {0, 1}. Thus if $X \in D$*

$$X = \left(x_1^1\ x_2^1 \ldots x_{n_1}^1\ \middle|\ x_1^1\ x_2^2 \ldots x_{n_2}^2\ \middle|\ \cdots\ \middle|\ x_1^n\ x_2^n \ldots x_{n_n}^n\right)$$

*where $x_i^j \in \{0,1\}$, j = 1, 2, …, n and $1 \leq i \leq n_t$; t = 1, 2, …, n. Now any element Y in the range space of R will be $(m_1, \ldots, m_m)$ where $m_i \in \{0,1\}$. Having defined the model we will just sketch the functioning of the model. Let $X \in D$. $X \odot_s M$ ($\odot_s$ denotes the special product of X with M); be Y, then we find the product $\odot_s$ using the usual product of X with M, thus the resultant is only a row vector in R. But while multiplying Y with $M^T$ we partition the resultant row vector just in between the elements $x_{n_1}^1$ and $x_1^2$, $x_{n_2}^2$ and $x_1^3$ and so on i.e., between the elements $x_{n_{n-1}}^{n-1}$ and $x_1^n$. So if $Y \odot_s M^T = Z$, then $Z \in D$.*

*This process is repeated. In fact the updating and thresholding of the mixed fuzzy super row vectors and the fuzzy row vectors are carried out at each stage for the following reasons.*

1. *The super column dynamical system M can recognise only those fuzzy row vectors or mixed super fuzzy row vectors only when its entries are from {0,1}, which implies the on or off state of the nodes / attributes.*

2. *This alone guarantees that this operation of special product* $\odot_s$ *terminates after a finite stage i.e., it either becomes a fixed point or it repeats itself following a particular pattern i.e., this resultant is supposed to give the super hidden pattern of the super column dynamical system.*

Next we proceed on to define and describe the new row super fuzzy relational maps.

**DEFINITION 1.2.2:** *Suppose we have n sets of attributes related with a problem which is divided into different sets and some n experts view about it and give their opinion. Each of these n-sets, view the problem in a different angle. Then to construct a model, which gives the consolidated view of the problem. That is at each stage the problem is viewed in a very different way. We construct a single model so that, the hidden pattern is obtained.*

*Let* $E_1, E_2, \ldots, E_n$ *be the n experts who study the problem, they all study the problem and give opinion on the n sets of attributes,* $A_1, A_2, \ldots, A_n$ *where* $A_1$ *has* $(a_1^1, \ldots, a_{n_1}^1)$ *number of attributes* $A_2$ *has* $(a_1^2, \ldots, a_{n_n}^2)$ *number of attributes and so on. Thus* $A_n$ *has* $(a_1^n, \ldots, a_{n_n}^n)$ *number of attributes. Using the set* $E_1, E_2, \ldots, E_n$ *as the domain space and* $(a_1^1, \ldots, a_{n_1}^1)$ *as the range space, let* $M_1$ *be the related connection matrix of the FRM-model.*

*Let $M_2$ be the related connection matrix of the FRM model using $E_1, E_2, ..., E_n$ as the domain space and $(a_1^2, a_2^2, ..., a_{n_2}^2)$ as the range space and so on.*

*Thus if we consider the fuzzy super row matrix,*

$$M = [M_1 \mid M_2 \mid ... \mid M_n]$$

$$= \begin{array}{c} \\ E_1 \\ E_2 \\ \vdots \\ E_n \end{array} \begin{array}{c} \begin{array}{ccccc} a_1^1 a_2^1 \cdots a_{n_1}^1 & a_1^2 a_2^2 \cdots a_{n_2}^2 & & a_1^n a_2^n \cdots a_{n_n}^n \end{array} \\ \left[ \begin{array}{c|c|c|c} \quad & \quad & \cdots & \quad \\ & & & \\ & & & \end{array} \right] \end{array}$$

*M will be the related super fuzzy row FRM matrix which is the super row dynamical system for the super fuzzy relational map model. This super fuzzy FRM model uses $(E_1, ..., E_n)$ to be the domain space; which has its state vectors to be usual simple fuzzy row vectors where as the range space state vectors are mixed super fuzzy row vectors given by elements of the form*

$$A = \left(a_1^1\ a_2^1 ... a_{n_1}^1 \mid a_1^2\ a_2^2 ... a_{n_2}^2 \mid ... \mid a_1^n\ a_2^n ... a_{n_2}^n\right)$$

*where $a_i^j \in \{0, 1\}$; $j = 1, 2, ..., n$ and $1 \le i \le n_i$, $i = 1, 2, ..., n$. Thus when we have $X \in D$; $X \odot_s M = Y$. $\odot_s$ is the special product of X and M where $Y \in R$ and*

$$Y = \left(y_1^1\ y_2^1\ ...\ y_{n_1}^1 \mid y_1^2\ y_2^2\ ...\ y_{n_2}^2 \mid \cdots \mid y_1^n\ y_2^n\ ...\ y_{n_n}^n\right),$$

*the elements of Y are from the set {0, 1}.*

*Now $Y \odot_s M^T \hookrightarrow X_1 \in D$ and so on until we arrive at a fixed point or a limit cycle to be the super hidden pattern of the super row dynamical system M.*

Now we illustrate first the two models by concrete examples. These examples are just only an indication of how the model works and is not a real world problem.

***Example 1.2.1:*** Let us have some 3 sets of attributes related with the domain space given by 3 sets of experts say $\left(E_1^1\ E_2^1\ E_3^1 \mid E_1^2\ E_2^2\ E_3^2\ E_4^2 \mid E_1^3\ E_2^3\ E_3^3\ E_4^3\ E_5^3\right)$ and suppose they choose to work with some 5 attributes ($a_1$ $a_2$ $a_3$ $a_4$ $a_5$) for a given problem P. Let D and R denote the domain and range space respectively.

The related fuzzy super relational row matrix is given by

$$\begin{array}{c} \\ a_1 \\ a_2 \\ a_3 \\ a_4 \\ a_5 \end{array}\begin{array}{c} \begin{array}{ccccccccccccc} E_1^1 & E_2^1 & E_3^1 & E_1^2 & E_2^2 & E_3^2 & E_4^2 & E_1^3 & E_2^3 & E_3^3 & E_4^3 & E_5^3 \end{array} \\ \left[\begin{array}{ccc|cccc|ccccc} 1 & 0 & 0 & 1 & 0 & 0 & 1 & 1 & 1 & 0 & 0 & 0 \\ 1 & 1 & 0 & 0 & 1 & 1 & 0 & 0 & 1 & 0 & 1 & 0 \\ 0 & 0 & 1 & 0 & 0 & 0 & 0 & 0 & 0 & 1 & 1 & 0 \\ 1 & 0 & 0 & 1 & 0 & 1 & 0 & 1 & 0 & 1 & 0 & 0 \\ 0 & 1 & 1 & 1 & 0 & 0 & 1 & 0 & 0 & 0 & 1 & 1 \end{array}\right] \end{array}.$$

Suppose X = (1 0 0 0 0) with only the node $a_1$ in the on state and all other states are in the off state to find the effect of X on the super dynamical system M;

$$\begin{aligned} X \odot_s M &\hookrightarrow (1\ 0\ 0 \mid 1\ 0\ 0\ 1 \mid 1\ 1\ 0\ 0\ 0) \\ &= Y \in R. \end{aligned}$$

Now

$$\begin{aligned} Y \odot_s M^T &\hookrightarrow (1\ 1\ 0\ 1\ 1) \\ &= Z \in D. \end{aligned}$$

$$\begin{aligned} Z \odot_s M &\hookrightarrow (1\ 1\ 1 \mid 1\ 1\ 1\ 1 \mid 1\ 1\ 1\ 1\ 1) \\ &= T. \end{aligned}$$

We see the super hidden pattern given by X is a fixed binary pair {(1 1 1 1 1), (1 1 1 | 1 1 1 1 | 1 1 1 1 1 1)}.

Now we illustrate super column FRM model by a simple example.

***Example 1.2.2:*** Let us consider a set of fixed attributes ($a_1$, $a_2$, … , $a_6$) which is associated with a problem P. Suppose we have 3 sets of experts given by

$$\left(e_1^1\ e_2^1\ e_3^1\ e_4^1\ e_5^1 \mid e_1^2\ e_2^2\ e_3^2 \mid e_1^3\ e_2^3\ e_3^3\ e_4^3\right),$$

$e_j^i$ represents an expert or a set of experts satisfying the same set of criteria. Now taking the 3 sets of experts along the row and the attributes along the column we assign to the domain space the 3 sets of criteria. Now taking the 3 sets of experts along the row and the attributes along the column we assign to the domain space the 3 sets of experts given by the mixed fuzzy super row vector and the range space consists of the ordinary fuzzy row vector ($a_1$… $a_6$) where $a_i \in \{0,1\}$; i = 1, 2, 3, …, 6. We would wish to state that we are giving only an illustration and this model has nothing to do with any of the real world problems or models.

The associated connection relational fuzzy super column matrix T is given by

$$T = \begin{array}{c} \\ e_1^1 \\ e_2^1 \\ e_3^1 \\ e_4^1 \\ e_5^1 \\ e_1^2 \\ e_2^2 \\ e_3^2 \\ e_1^3 \\ e_2^3 \\ e_3^3 \\ e_4^3 \end{array} \begin{array}{c} \begin{array}{cccccc} a_1 & a_2 & a_3 & a_4 & a_5 & a_6 \end{array} \\ \left[\begin{array}{cccccc} 1 & 0 & 0 & 0 & 0 & 0 \\ 0 & 1 & 0 & 0 & 0 & 0 \\ 0 & 0 & 1 & 0 & 0 & 0 \\ 0 & 0 & 0 & 1 & 1 & 0 \\ 0 & 0 & 0 & 0 & 0 & 1 \\ \hline 1 & 0 & 0 & 1 & 0 & 0 \\ 0 & 1 & 0 & 0 & 1 & 0 \\ 0 & 0 & 1 & 0 & 0 & 1 \\ \hline 0 & 0 & 0 & 0 & 0 & 1 \\ 0 & 0 & 0 & 1 & 0 & 0 \\ 0 & 0 & 1 & 0 & 0 & 0 \\ 1 & 0 & 0 & 0 & 0 & 0 \end{array}\right] \end{array}.$$

Now we want to find the super hidden pattern of the mixed fuzzy super row vector.

$$X = (0\ 0\ 0\ 1\ 0 \mid 0\ 1\ 0 \mid 0\ 0\ 0\ 1)$$

on the fuzzy super dynamical system T.

$$X \odot_s M \hookrightarrow (1\ 1\ 0\ 1\ 1\ 0) = Y \text{ (say)}.$$

$$Y \odot_s T^t \hookrightarrow (1\ 1\ 0\ 1\ 0 \mid 1\ 1\ 0 \mid 0\ 1\ 0\ 1) = Z \text{ (say)}.$$

$$Z \odot_s T \hookrightarrow (1\ 1\ 0\ 1\ 1\ 0) = Y_1 = (Y).$$

Thus the super hidden pattern of the super dynamical system is a fixed binary pair given by {(1 1 0 1 1 0), (1 1 0 1 0 | 1 1 0 | 0 1 0 1)}.

Now we proceed on to describe the super FRM model. This model comes handy when several sets of experts work with different sets of attributes.

**DEFINITION 1.2.3:** *Suppose we have some problem P at hand and we have n sets of experts $N_1, N_2, N_3, \ldots, N_n$ where each $N_i$ is a set of experts, $i = 1, 2, \ldots, n$. Further we have some p sets of attributes; i.e., we have $M_1, \ldots, M_p$ sets of attributes. We have (say) some experts work on some sets say $M_i, M_k, \ldots, M_t$, $1 \leq i, k, \ldots, t \leq n$. Likewise some other set of experts want to work with $M_s, M_r, M_l, \ldots, M_m$, $1 \leq s, r, l, m \leq n$ where we may have some of the set of attributes $M_i, M_k, \ldots, M_t$ may be coincident with the set of attributes $M_s, M_r, \ldots, M_m$.*

*Now we cannot have any of the fuzzy models to apply to this. Then now we get a new model by combining the two models, which we have described in the definitions 1.2.1 and 1.2.2. We set the fuzzy super matrix as follows. Let the $N_i^{th}$ set of experts*

*give their opinion using the $M_j^{th}$ set of attributes then, let $P_{ij}$ denote the connection FRM matrix with the $N_i$ set of attributes forming the part of the domain space and $M_j$ attributes forming the range space. This is true for $1 \leq i \leq n$ and $1 \leq j \leq p$.*
*Thus we have a supermatrix V formed with these np fuzzy matrices*

$$\begin{array}{c} \\ N_1 \\ N_2 \\ \vdots \\ N_n \end{array}\begin{array}{c} \begin{array}{ccc} M_1 & \cdots & M_p \end{array} \\ \left[\begin{array}{ccc} P_{11} & \cdots & P_{1p} \\ P_{21} & \ldots & P_{2p} \\ \vdots & & \vdots \\ P_{n1} & \ldots & P_{np} \end{array}\right] \end{array}.$$

*It may so happen that some of $P_{ts}$ may be just $t \times s$ zero matrices, $1 \leq t \leq n$ and $1 \leq s \leq p$. Clearly $P_{is}$ are fuzzy matrices with entries from the set {–1 0 1}. It may so happen that some set of $N_t$ experts may not want to use the set of s attributes and give their opinion in which case we have that associated matrix $P_{ts}$ to be a zero matrix. Clearly this model can be visualized as a combination of the two models described in the definitions 1.2.1 and 1.2.2.*

*This model is defined as the super FRM model. That is if we consider any of the rows of this fuzzy supermatrix V then we have that row is again a mixed fuzzy super row matrix given by $[P_{j1} \mid P_{j2} \mid \ldots \mid P_{jp}]$ now $1 \leq j \leq n$. As j varies over n we get the total of V.*

*Similarly if we consider any column of V we see it corresponds to*

$$\left[\begin{array}{c} P_{1t} \\ \hline P_{2t} \\ \hline \vdots \\ \hline P_{nt} \end{array}\right]$$

*where $1 \leq t \leq p$. Thus by taking all the super fuzzy columns we get the super fuzzy matrix. Thus each column of V is a mixed fuzzy super column matrix. Now this super dynamical system V performs the work of the both these models described in definitions 1.2.1 and 1.2.2, simultaneously.*

Now how does this model function?

We have both the domain and range space of this fuzzy super FRM model to be only mixed super row vectors given by $[N_1 \mid N_2 \mid \ldots \mid N_n] \in D$ and $[M_1 \mid M_2 \mid \ldots \mid M_p]$ in R.

Thus $X \in D$ would be of the form

$$X = \left(x_1^1 \; x_2^1 \ldots x_{n_1}^1 \;\middle|\; x_1^2 \; x_2^2 \ldots x_{n_2}^2 \;\middle|\; \cdots \;\middle|\; x_1^n \; x_2^n \ldots x_{n_n}^n\right)$$

like wise and $Y \in R$ would be of the form

$$Y = \left(y_1^1 \; y_2^1 \ldots y_{p_1}^1 \;\middle|\; y_1^2 \; y_2^2 \ldots y_{p_2}^2 \;\middle|\; \cdots \;\middle|\; y_1^p \; y_2^p \ldots y_{p_p}^p\right) \in R.$$

Clearly the $x_{ts}^j$ and $y_{rs}^k$ are from the set {0, 1}; $1 \leq j \leq n$, $1 \leq k \leq p$, $1 \leq t \leq n_t$, t = 1, 2, …, n and $1 \leq r \leq p_r$; r = 1, 2, …, p.

Now we illustrate the functioning by an example; this example is however not a real world problem only an example constructed to show the working of the model described in definition 1.2.3.

***Example 1.2.3:*** Let us suppose we have at hand a problem P which is worked out by a set of 3 sets of experts and they give their views on 4 sets of attributes given by the fuzzy super matrix V.

$$
V = \begin{array}{c} \\ E_1 \\ \\ \\ \\ E_2 \\ \\ \\ \\ E_3 \\ \\ \end{array}
\overset{\begin{array}{cccc} \quad A_1 \quad & \quad\quad A_2 \quad\quad & \quad A_3 \quad & \quad A_4 \end{array}}{\left[\begin{array}{cccc|ccccc|ccc|ccc}
1 & 0 & 0 & 0 & 0 & 0 & 0 & 0 & 0 & 1 & 0 & 0 & 0 & 1 & 0 \\
0 & 1 & 1 & 0 & 0 & 0 & 0 & 0 & 0 & 0 & 1 & 1 & 0 & 0 & 1 \\
0 & 0 & 0 & 1 & 0 & 0 & 0 & 0 & 0 & 0 & 0 & 1 & 1 & 0 & 1 \\
0 & 0 & 0 & 0 & 0 & 0 & 0 & 0 & 0 & 0 & 0 & 0 & 0 & 0 & 0 \\
\hline
0 & 0 & 0 & 0 & 1 & 0 & 1 & 0 & 1 & 0 & 0 & 0 & 0 & 0 & 0 \\
0 & 0 & 0 & 0 & 0 & 1 & 0 & 1 & 0 & 0 & 0 & 0 & 0 & 0 & 0 \\
0 & 0 & 0 & 0 & 1 & 0 & 0 & 0 & 1 & 0 & 0 & 0 & 0 & 0 & 0 \\
\hline
1 & 0 & 0 & 0 & 1 & 0 & 0 & 0 & 0 & 1 & 0 & 0 & 0 & 1 & 0 \\
0 & 1 & 1 & 0 & 0 & 1 & 0 & 0 & 1 & 0 & 1 & 1 & 1 & 1 & 0 \\
0 & 0 & 1 & 1 & 0 & 0 & 1 & 1 & 0 & 0 & 1 & 0 & 0 & 0 & 1 \\
1 & 1 & 0 & 0 & 1 & 1 & 0 & 0 & 0 & 0 & 0 & 1 & 1 & 0 & 0 \\
1 & 0 & 1 & 0 & 1 & 0 & 0 & 1 & 1 & 1 & 0 & 1 & 1 & 1 & 0
\end{array}\right]}
$$

$$
= \begin{array}{c} E_1 \\ E_2 \\ E_3 \end{array}
\overset{\begin{array}{cccc} A_1 & A_2 & A_3 & A_4 \end{array}}{\left[\begin{array}{c|c|c|c}
M_{11} & M_{12} & M_{13} & M_{14} \\ \hline
M_{21} & M_{22} & M_{23} & M_{24} \\ \hline
M_{31} & M_{32} & M_{33} & M_{34}
\end{array}\right]}
$$

where $M_{ij}$'s are fuzzy matrices which correspond to the connection matrices of the FRM. We see the set of $E_1$ experts do not wish to give their opinion on the set of attributes $A_2$. Likewise the set of experts $E_2$ do not wish to give their opinion on the set of attributes, $A_1$, $A_3$ and $A_4$; only they give their views on $A_2$ alone.

However the set of experts $E_3$ have given their views on all the four sets of attributes $A_1$, $A_2$, $A_3$ and $A_4$. Thus we have $M_{12}$ = (0) and $M_{22} = M_{23} = M_{24} = (0)$.

Now we just illustrate how this dynamical system V functions. V will be known as the super dynamical system of the fuzzy super FRM maps. Suppose we wish to study the effect of the state vector

$$X = (1\ 0\ 0\ 0 \mid 0\ 1\ 0 \mid 0\ 0\ 0\ 1\ 0) \in D$$

on the dynamical system V. We get X $\odot_s$ V where $\odot_s$ is the super special product.

$$X \odot_s V \hookrightarrow (1\ 1\ 0\ 0 \mid 1\ 1\ 0\ 1\ 0 \mid 1\ 0\ 1 \mid 1\ 1\ 0) = Y \in R.$$

Now

$$Y \odot_s V^T \hookrightarrow (1\ 1\ 1\ 0 \mid 1\ 1\ 1 \mid 1\ 1\ 1\ 1\ 1) \in D$$

and so on.

We proceed on until we arrive at the fixed point or a limit cycle. This fixed point or a limit cycle which will form a binary pair will be known as the super hidden pattern of the dynamical super system. Now having seen the three types of new super fuzzy models we apply them to the problem of media's role in OBC reservations in the next chapter.

**Chapter Two**

# ANALYSIS OF THE ROLE OF MEDIA ON RESERVATION FOR OBC USING SUPER FUZZY MODELS

In this chapter we use the new super fuzzy models constructed in chapter I, to analyze the role of media on 27 percent reservation for the OBC in institutions like IITs, IIMs and AIIMS. This chapter has six sections. Section one gives a brief description of the attributes given by the experts. Section two uses the super FRM model to study the interrelations between the social problems and the reservations for OBC. In section three the role of media is analyzed by the different category of people who served as experts using super fuzzy mixed FRM models

The role played by the media in OBC reservations and the government authorities concerned is analyzed using super fuzzy column FRM models in section four. Section five analyzes the role played by the media in OBC reservations using FCM models. The final section gives an analysis of reservations for OBC, and media's responsibility in the education system as observed by the British and the role of present day caste division based on birth in India.

## 2.1 Brief Description of the Attributes given by the Experts

We in this section analyze the reservation problem in the context of media and the related issues of rights, social reformation, social justice, social discrimination, social tolerance and so on. First we list the set of attributes given by the experts and describe them in a line or two before we use the relevant fuzzy model. The set of attributes supplied by a collection of experts is as follows.

1. Social Rights: By social right the experts some of whom were socio scientists, psychologists and educationists said they use this term in two ways, one the social right a person has towards the society immaterial of caste or creed so that the normal functioning of society remains undisturbed. Only when the social rights are misused we find disturbance in the society agitation or protest or riots or finally revolution taking place in the nation. For example, in the context of reservation it is the social right of every OBC/SC/ST to make use of the rule put by the parliament to have 27 % of seats in all educational institutions run by the central govt. But it is the misuse of social rights by the anti reservation protestors for their protest demands against the social right of the majority, the social right given or declared by the govt to improve or uplift the majority of the society. These improvement programmes cannot be taken as discrimination for otherwise the sons of the soil will loose everything and someday become extinct species of the nation. The other side is if for long their social rights are denied they may end up as terrorists or dacoits for their very living. Then who can be held responsible, the govt and more so are the anti reservation protestors. As long as caste is in the nation the social rights will also vary from caste to caste said the socio scientists.

2. Social Reformation: In India any one with the social out look or awareness makes an observation between the comforts enjoyed by a few upper castes and Brahmins and the poor social conditions under which the OBC/SC/ST, thrive even worst than animals, it would be the natural instinct of any human (a social creation i.e. a man/woman) to put forth some social reformation so that not only the sharp gap is narrowed but at the same time the working poor get some benefit from the nation. So it is the social right of the parliament to make some social reformation to suit the circumstance and situations. Unless such timely social reformations are made, the nation may be forced to face a riot or a revolution. Once again it has become pertinent to mention that to make social reformation the courts have no power the parliament is the policy maker of it, in India as the parliament functions by the peoples representatives. So any social reformation can be thought of as a people reformation. If a few people protest or function against such social reformation because they cannot use it they are liable for punishment as they are against the very development of the nation. For only by social reformations the govt can make steps towards the development of the nations. A few people belonging to a particular society living in a comfortable super state is by no means a scale to say the nation is developed. The real development lies in the last mans development. As of today in India, the rural poor living conditions are questionably very disappointing. He has no access to health centres, no nutritious food, no toilet facilities, no clean water (a recent survey shows the ground water in many of the states in India are polluted by fluoride and chemicals which is hazardous to health). At least in the time of British rule he could get clean water and unpolluted food grains and unpolluted air which in spite of poverty helped him to live without diseases for a long spans, say even over 90 years or so, but in India now a days a rural poor becomes invalid and cannot work even at the age of 45. Who is responsible for all this? Is this real development of the nation or otherwise. So any social reformation made to upgrade the living conditions of the majority cannot be protested or questioned says the expert.

3. Social Justice: Social justice is a very special term used by the experts. They are in many cases rules or laws not made by the socially well of Chief Justices and Judges of the nation who belong to the upper castes or Brahmins; they are the justice contemplated by the socio politicians(i.e., political leaders or politicians or reformers) who are concerned mainly about the majority poor suffering, any form of socio discrimination these leaders have the power to make the politicians pass such order which the experts term as Social Justice; for instance Periyar was a social reformer who could do anything politically in T.N. V.P Singh is both a social reformer as well as a political leader for he was our P.M and held office, unlike the social reformer Periyar. Like wise Ambedkar, Sahul Maharaj, Vasant Moon and so on. So the social justice need not be a law appearing in the courts of law. It can be done by means of amendment in the constitution or addition to constitutions passed in both the houses. Thus the reservation made for OBC in the educational institutions run by the central government can be termed as social justice to be more precise experts said. The social justice passed by parliament cannot be even questioned by the Courts of Law for the democracy is nothing but the rule of the people. When people approve something as social justice to respect the majority of the people when passed with over whelming majority should be accepted unquestioned and cannot be challenged.
4. Social Discrimination: When the experts said social discrimination they went on to describe in pages so it was very difficult for us to put in a few lines. The experts were asked to define their ideology of “Social Discrimination” mainly in the context of OBC reservation or in general “reservations “in India. So good old topics of discrimination by the Aryans from the time of their entry till date cannot be mentioned for it is by itself a research. So we mean taking only into account the “caste” which is based on birth as given in the “laws of Manu”, which even rules the nation which is unfortunately the laws practised in the courts. By social discrimination the experts mean only those problems faced by an individual or a society or a group due to their only reason being based on the caste by birth i.e. largely 'non Brahmins'. A non Brahmin cannot sit with a Brahmin in a dining hall and dine together in the temple! This is social discrimination. Unless this is annulled how will a Sudra and a SC/ST live in peace in. For the society always takes up things easily which are evil. It takes even centuries to adopt the good and healthy policies, by people at large for every man wants to have some one inferior to him to satisfy his ego! Caste in India is based only in the mind set up. Thus clearly the history of social discrimination entered India only through the kyber pass by the Aryans. So any rule or imposition of policies by the govt to annul social discrimination should not be debated by the upper castes enlisted in the Manu. For clearly reservation for OBC in the educational institutions run by the central government is to help the OBC from being socially discriminated by the Brahmins claims the experts.

5. Social Tolerance: By the term social tolerance the experts imply the following. We the OBC /SC/ST have a very high degree of social tolerance that is why we have so long tolerated all the atrocities done by the Brahmins on us like denying us education, denial of right to walk in the public streets, denial of good posts, denial of certain good food (for it is mentioned in the Manu that ghee which they use in the Yagas and their staple food should not be consumed by the Sudra). These social rights like education is still denied and we like fools have been just tolerating them. The protest made by them (Brahmins) against reservation for OBC in central institutions of higher learning can be termed as social intolerance. We became socially intolerant when we were denied the right to walk in public streets near and around the temples. The hero of Vaikam, Periyar led a very massive protest against the practise by the Brahmins, when even the pigs, cows and dogs could freely use these roads. So our social tolerance was put to test and it was not in any way wrong on the path of the OBC/SC/ST to claim their right to use the roads near by the temples. Our tolerance in such cases would only portray us more inferior to pigs, dogs and cows. Now when the govt is implementing thro' an constitutional amendment the OBC reservation that too in the educational institutions run by the central govt there is nothing for the Brahmins/upper castes, to protest and the media to blow it up. They should have the social tolerance and above all the broad mind to welcome that change; but unfortunately they have become socially intolerant and made all indecent protest against reservation for OBC.

6. Social Intolerance: This is just the opposite concept of social tolerance.
7. Social Balance: If the nation has to function smoothly or if the nation should progress for the betterment of the majority it is very much essential that the society is socially balanced to some extent. For if social balance is lost nation has to face internal chaos or revolution. After 60 years of independence with the suicide of farmers and acute unemployment problems among the rural uneducated poor the govt. felt a need for reservation for the OBC to bring in some social balance between the OBC and the Brahmins and upper caste. The socio scientists and politicians felt it a dire need to implement reservations for the OBC in educational institutions run by central government. If the social imbalance persists we will have riots and at times even terrorism and so on. Now itself we see mainly in the festive seasons in T.N we hear more of chain snatching, pick pockets and so on. This study is done according to the socio scientist who feel to bridge the balance such hazardous thefts take place. The reservation for OBC has been contemplated by the socio scientists and politicians to balance the society so as to avoid further social imbalance which may cause every form problem questioning the very peace of the nation.
8. Social Imbalance: This is the term opposite to the term social balance when the difference between the rich and the poor is increasing day by day; govt. must take some reformative steps to lessen the gap. We see even the urban poor living conditions are distinctly very different from that of the urban rich. This form of social set up in any nation cannot hold on for a very long time. Further the major difference is that the richness is associated with specific castes and the poverty is associated with the OBC/SC/ST. i.e., the very economic condition in general in India is based on caste that is why the govt had taken punitive steps to implement reservation which is of course an action which cannot immediately relieve the poor OBC/SC/ST but at least give some mental satisfaction claims the experts. Thus one has to settle the social imbalance, otherwise it will lead to the disturbance in the social equilibrium of the nation.

9. Social Equilibrium: When the experts said social equilibrium they said that it is highly dependent on the social balance of the nation. If it takes even one day for the social balance to be struck than that is a threat to social equilibrium. Even without a days threat social imbalance can be caused by a terror attack or a mass suicide (for instance mass suicide by farmers) and so on. It is up to the government to take immediate action to set things in the right track. This setting up of the right track means only making any rule or policy to the majorities welfare not for the comfort of the socially and economically powerful, for such appeasement of social and economically powerful of the minority may in a long run result in a chronic displacement of the social equilibrium resulting in the disintegration of the very nation.

   The majority of the experts who were socio scientists very clearly made this comment and wanted the authors to record it without any form of dilution or distortion. All of them said the OBC reservation is an educational policy laid down for the socially discriminated majority of the society. Without the majorities social development it is impossible for India to remain as developed. One may think their (upper caste and Brahmins) advancement both in economy and learning is a development, this will not last long for the majority discontentment and despair may route into the tilting of the nations equilibrium permanently and disintegrating it. For the economy only in the hands of the socially better off society who form the minority of the nation may not for a long period of time make the majority, socially deprived people to be tolerant towards them, a day may come when riots or revolution is certain to take place. Now as signals of discontent and despair one may see stray incidents taking place in India in the states like A.P, Manipur, Assam, Rajasthan etc. Any such act of despair is very contagious for in no time they can plague the whole nation. What will the integrated India do? The seriousness of the problem in India is that the money amassed is not even the sons of the soil. They make talk and tall claims of making an integrated Hindu India! What is this Hindu India? A nation cannot have one religion be it Africa or China or Russia or U.S or U.K! In the world all countries have people forming a heterogeneous group of religions, but to speak of even reservation for only one religious group in India is a threat to peace. For the people who follow other religions are all sons of this soil unlike the most Brahmins claim themselves to be the minority Hindus. If history is traced using linguistics it is very clear that even the very word 'Hindu' is a contemptuous term used by Muslims and Christians to represent Indians. So it is all the more very doubtful about any particular religion followed by the sons of the soil of India. If there was a religion followed by them they should have a substantiating religious book. Till date historically the Indians (sons of the soil) don't have any religious script which they followed. At this stage it is pertinent to mention that the Vedas written by the Brahmins is not and cannot be claimed as the religious book of Indians for the two fundamental reasons.

1. They are not the sons of the soil for they should not read Vedas per laws of Manu of the Brahmins.
2. They are just 3% of the nation population who have settled in India for their own personal and selfish benefit. It is still important to make a mention that unlike the British or the Muslims they are unconcerned about the development of the majority who form the rural poor and who are also socially and economically very backward mainly because all the economy and education is in the hands of the 3% minority foreigners. This inequality is certain to tilt the social equilibrium of the nation unless the govt takes strong positive steps to mainly educate the majority; for giving them money or land cannot solve the problem it can be thought of only as an interim step or an immediate step. It is no true solution claims the experts.

3. Education for all that too in institutions run by the govt; be it central or state govt claimed the experts. British favoured our education, but Brahmins can give us anything but not education. If comforts are more enjoyed by upper castes denying the even very basic living of the backward and SC/ST we see the social equilibrium would be lost, leading to social inequilibrium which means the society has the danger to be dismantled.

10. Social Inequilibrium: By this term the experts who were mainly socio scientists said they say the nation lost its social equilibrium or is in social inequilibrium state, it means the nation once disintegrated into pieces it cannot again become whole or integrated. They said (socio scientists said) when Russia attained its social in equilibrium it was forcefully divided. The great USSR was sadly disintegrated into bits. It was a hard blow to see USSR dismantled. Socio scientists claim the socio in equilibrium was very mild in USSR in spite of its division. This may not be the case of all nations which lose their social equilibrium. Thus if reservation for OBC is not given a permanent divide of the nation based on the caste is certain; this change will not occur under peaceful conditions but only under blood shed. So the socio scientists who were the experts said it is high time all the politicians put their head together and save the nation. After all the Aryans who have come to India are aliens, they are least bothered about the real progress of the nation for they are very much bothered only about their comforts and their getting the best out of the nation. If the nation is disintegrated they will leave the nation seeking for better prosperity said the experts. Due to the modernization, the OBC/SC/ST, agriculturist have failed to get even the very living from their agricultural profession which had lead to mass suicide by the farmers. This form of reservation is the only solution to avert any form of riot in the nation claimed the experts. Further the experts said only after 60 years of independence the govt felt that OBC must be given some percentage of reservation in the educational institutions run by it(govt).This is not a accidental plan or for vote. This has been done only to lessen the gap mainly socially between the upper castes and the OBC. For economy cannot and will not lessen the gap when it comes on to social status. What ever be the economic status of a Sudra he remains low before an upper caste but if the Sudra is educated that too from the institutes like IIT or IIM or AIIMS has social status in par with an upper caste; that is the way to lessen the social gap it is essential that reservation is given to OBC in these institution. For if the social gap increases the nation is sure to attain the social inequilibrium state warned the experts.

11. Social Conscience: Most of the experts said if the public especially the upper castes ever had at least to some degree a social conscience they would not have made any anti reservation protests. Further they added that the anti reservation protests were made only by a select few backed by the Hindutva force who had no guts to say when the same was tabled in the two houses. They were very stealthy, religiously a fanatic Hindutva backed youth calling them under the banner “youth for equality”, the direct brain behind them was Dr. Venugopal, the director of AIIMS and their mission was supported by some IIT-ians’ in Chennai, but could not agitate in the Periyar land said the experts. So these set of comfort and power loving selfish people were the one who made almost all the anti reservation strikes, said the experts. If further one traces the lawyer who appeared for Dr. Venugopal, he is the council for Sri Sankaracharya of Kanchi Mutt, Chennai and he is also counsel for BJP/VHP and RSS. Now one can understand the strong social bondage between them said the experts. How can they have any social conscience when it is a matter related with the majority who are the sons of the soil said the experts. So it is only the lack of social conscience among the upper castes that has made them to protest. People with even a little of social conscience would never be against any of the policies which is helpful to the majority. “The OBC reservation” said the experts, was basically to help the uneducated majority towards education in institutions run by govt, so any person against such development lack social conscience. It took 60 long years for the govt to have some social conscience to think about the majority who happened to be the sons of the soil but even before the implementation or even introduction in phases or otherwise, for the people who are just foreigners not even the sons of the soil made such protest which was well projected by the media is a clear act of selfishness and they have no basic social conscience but only selfishness and self centeredness which has made them to act or protest in an arrogant way making use of the full support of the media and the court which are in their monopoly said the experts. The experts further said there is always a limit to everything and their selfishness had no bound which would lead to their permanent fall. If at all they had even an atom of social conscience, they would not have acted in this senseless manner, claimed the experts.

“Social conscience” is an important and a selfless trait which each and every real/true citizen must possess for making the nation a better living place for everyone. If the social conscience is lost in the individual or by a section of people, it is always a danger to the nation. Further in India, the media has no social conscience. It functions in a way which never reflects in any sphere that it has social conscience but only it is dominated by the caste conscience said the experts.

12. Social Understanding: This is a term used by the experts who said social understanding is something more than the social conscience. All people with social conscience may not and need not have the social understanding but a person in general with a social understanding will have at least some percentage of social conscience claimed the experts. They said caste feelings may make a person act without social conscience but if he/she has a social understanding i.e. the better understanding of the social set up of India and the social problems faced by the majority. For the policy makers in almost all cases are never the representatives of the majority but only from the dominant castes, no policy made by them ever caters to the real needs of the majorities development said the experts. They were always an eye wash or to be more critical the policy makers have never made any policy towards positive development towards the true sense of progress. This is very much evident from the fact that millions in this nation are illiterates, for in India if a person can draw his/her name with the correct spelling he/she is declared as an literate with this meaning of literacy so many millions remain illiterates. If the media or the protesting upper caste ever had any social understanding and really are interested in steering the nation towards progress will they ever protest for OBC reservation. With this, the reader is requested to know the concept of social understanding. Here the experts wanted us to make the note the media never functions with any social understanding for they are blind about it.

13. Social Awareness: One may not have social understanding but may be aware of the society's wants needs for better development and progress of the nation. For instance the politicians only after 60 years of independence became aware of the fact that the majority need to get good education at least form the educational institutions run by the central government. To make this a policy or rule they passed an act in both the houses which was passed with an over whelming majority. The act was to give 27% of reservation for the OBC in central govt. run institutions. The political leaders were very much aware of the fact that the neglect of the majority in the nation would not only make the nation retard in progress may in a long run tilt the equilibrium or the very peace of the nation. So more with fear the political leaders decided to have the right for the majority to be educated in these institutions, for the gap between the OBC and the upper castes were very widening socially. Only one section of the people that too not even the sons of the soil had control over economy, media, law and education hence they exploited the labour of the OBC and paid them meagrely which was not even enough to have a square meal, leave alone to meet both the ends. Thus the globalization or modernity or information age never had any form of impact on the lives of these people. So the only way to raise some of them was to give good education and the govt cannot ask the private institution to do that is why they put first the rule 27% reservation for the OBC in the central govt run educational institutions. Any socio scientist with a little of social awareness will welcome this change or policy. Any person with some social awareness will accept that reservation for OBC is essential. Thus the experts want to emphasise that even if social awareness is felt by them still with selfish motives will not accept the reservation. For the upper castes are the people who are first aware of all problems yet they think they can deny the majority as long as they can, so that they can enjoy the maximum benefit of the nations prosperity. But lack of awareness to accept the dire need of the hour of the majority may create all forms of other social problems like terrorism, dacoits, plundering the rich and so on thus marring the peace of their (upper caste) very living. If they do not understand and act accordingly the nation is forced to face dire consequences said the experts. The media by its act has miserably failed to make the upper castes and the public socially aware of the situation.

14. Social Commitment: Every citizen has a social commitment both individually and collectively. It may not be shocking to see some terrorist feel it is their social commitment to help the poor and the needy so they do not mind acting in an erratic way to achieve this: for they are so sensitive to the sufferings of the people of their surroundings, they do murder or any other base act like looting the rich or so on to supply the needy and poor. When the nation fails to cater to the needs of the majority it leads to such type of terrorism they may be called by any name but they have a leader and the leader is the brain behind the followers in fact the benefactors of this group will be very faithful to their leader so when the experts mentioned the term social commitment it has varied meaning. In a sense the upper castes and the Brahmins think it is their social commitment to fight against the reservation for OBC for they think they may loose their opportunity if 27% OBC become their competitors. Here it is worth while to mention that this only shows their real commitment, it is their selfish attitude to usurp the basic needs of the majority. It is here that the media has utterly failed to be socially committed, for the responsibility of the media is to project the nation commitment towards the majority of the nation to see they lead a life with some social dignity, equality and comfort. So this reservation for OBC may lead to social equality which may otherwise be a threat to the nations peace, media has acted against its proper norms of social commitment towards the nation. It is the social commitment of the very nation to provide good education for the majority who are the sons of the soil of this nation. Thus we mean and use this word in different shades depending on the need of it.

15. Social Arrogance: When the nation has declared that 27 percent of reservation is a need of the hour for OBC in the institutions run by the central govt, it is the arrogance of the Brahmins and upper castes to protest. We call this the social arrogance because socially they are high, as only claimed by them by the laws of Manu so they have the social arrogance to fight for the OBC reservation. For the SC/ST who are socially below the OBC according to the laws of Manu. Here it is pertinent to mention that even the OBC did not display this sort of social arrogance for the reservations of SC/ST. But the social arrogance of the Brahmins and the upper castes is very much fondled and supported both by the media and the courts, claims the experts. In fact they went on further to say that the very media is socially arrogant to project such false propaganda against OBC reservations. So it is the irresponsible behaviour of the media to project such false propaganda to the public regarding reservations.
16. Social Responsibility: The experts define this term as some shades stronger than the concepts social awareness or social conscience. For when we say some one is responsible, it means some what their social right or social responsibility. So when the experts say social responsibility it is a collective responsibility as well as a individual responsibility (one may ask towards what or whom once again) towards the majority or towards the nation or more so towards their own fellowman. It is all the more the duty of the media to behave socially responsible; for the media in the opinion of the experts have behaved socially very irresponsible by projecting beyond proportions a few anti reservation protests that too mainly by the same set of upper caste/ Brahmin women students and doctors. It is the media under the instructions of the caste dominated parties like BJP/VHP/RSS blew up these small uncommitted protests(for the protestors were playing cards, reading novels in the reclined posture or sleeping all under the shamiana with air coolers that too in the AIIMS campus)done with unconcern as if forced by others to do. Except for the false propaganda of the media these protests are nothing. It is very unfortunate they could pull of like this for over 40 days unconcerned of the poor patient problems. The media was so irresponsible for it did not give any importance to the patient problems and the deaths of patients due to these protests. The climax of the protests was the court ordering the govt to pay the doctors in their period of absence when they abstained from their work wilfully to protest against the very govt, who was paying them and who were their employers. Thus in this situation both the media and the court have failed to be responsible said the experts. Thus when the govt or politicians after 60 years of independence acted socially responsible the court and the media acted socially irresponsible over a serious issue. The doctors who were socially responsible for their patients failed to be so which resulted in the death of a few of the patients.

17. Socially Irresponsible: This attribute will be used just in an opposite way to the concept socially responsible or social responsibility.
18. Social Maturity: When any policy or an act is made to benefit the majority of people the persons or public who are

not benefited by it must also have the social maturity to welcome such social changes. For instance the reservation for OBC in the educational institutions run by the central govt. is in no way affecting any caste or creed for they (OBC) too should be given an opportunity to be educated in the govt institution which is claimed to be of good standard. Since it is in the first place govt property no one should even protest it for all the policy is also laid only by the govt for the welfare of the nation. But the Brahmin and the upper castes do not have the social maturity to accept it though they may not have the mind to welcome it said the experts. This order is made to extend some sort of social uplift for the majority who are socially backward. Social backwardness cannot be changed by the economic upliftment. If the social backwardness is to be changed the one and only way is to educate people, give them the best of education in a good educational institute. This alone has the capacity to change their social state. Unless the social state is changed it is impossible to declare the nation as the developed one that is why with all the adoption of globalization, modernization and information technology, ruling the nation, the nation continues to be under developed, for the social status of the majority are very backward, the additional criteria is the majority of these people who live in rural areas are also economically backward. So the reservation for OBC is to be welcomed and accepted by one and all, who long for the development of the nation. Only socially immature people are protesting against it. The added pity is the media by supporting them, also has proved itself to be socially immature. So we can use the attributes social maturity and social immaturity as opposite terms as claimed by the experts.

19. Social Evil: OBC reservation is not a social evil but the media and the Brahmins are projecting OBC reservation as a social evil. So the experts wanted to take this attribute also for discussion. The Brahmins and the upper castes have the social arrogance to say that the reservation would retard the nation development hence according to them it can be termed as the social evil.

20. Social Equality: The pity is that the OBC are socially inferior to the Brahmins propagated by the laws of Manu and by the constitution that castes were accepted. So now social equilibrium may be lost if they are going to be treated socially unequal. Once social equilibrium is lost no one can once again make India an integrated one. It would become like Russia. So only the political leaders and the scientists after lots of hesitation and discussions brought back the issue of reservation for OBC in higher education; for the only way to wipe of the social inequality is to give them good higher education. No amount of good school education can make them socially equal. So only with lots of research the political leaders found the only way to attain social equality is to give reservation for OBC in higher education.

    If social equality is not attained even after a century the nation would only continue to remain as a developing one. So the govt. made reservation for OBC only to attain social equality. Social equality can only be made by reservations for the social problem or the social evil in India is caste. Caste that is based on birth is the first and the main threat to social equality. Caste based on birth firmly denied education for OBC so only after the study of the problems the socio scientists and the concerned political leaders strongly felt to achieve social equality the one and only means was to give reservation for OBC in the educational institutions of higher learning at least run by the central govt. The main and only purpose of this move is to establish social equality.
21. Social Harmony: A nation if it has to maintain at least some degree of social harmony it has to make rules or laws so that the majority of the nations population are not socially strained and do not plan to become anti social elements. For social harmony economic equality alone will not work as India is a uniquely caste plagued nation. So for social harmony, the equality could be attained only through education. Thus the only tool to have social harmony in India is by giving higher education to one and all in institutions at least run by the central and state govt. If social harmony is disturbed the nation is certain to face social unrest and revolution. For no nation can go on with social unrest for anything untold even can happen at anytime. So people cannot live in peace i.e. social peace will also be at stake.

22. Social Welfare: The terms social harmony and social welfare are not the same attributes. We hope the clever reader will be in a position to understand the subtle difference between these two concepts. We often have come across the term or the phrase "Social Welfare State" or just Welfare State". What really one means by it. The experts said they make use of this phrase in the following way. They say in the context of the caste ridden India when they say "Social Welfare" of the state they mean that the discriminated majority which is based on caste should enjoy the laws or orders of the govt which would be to their welfare socially. All economic welfare schemes can only make some changes in economy but surely not socially said the experts. So social welfare in something like right of them to study in good educational institutions. The female children are given some support to get more higher education and govt takes steps to monitor they do not become school drop outs and so on. Here also the experts want to make the readers understand that economic equality alone cannot make the state a socially welfare state. Higher education is the main criteria to make India a Social Welfare State, may be with some economic aid for the poor OBC/SC/ST. To strike social harmony in the nation, the main criterion is the nation must be a social welfare state, claimed the experts.
23. Social Peace: The experts by this attribute meant peace in the society i.e. in the nation peace among all caste people. Unless peaceful atmosphere or peace exists in the nation, the nation cannot make any step towards progress. For instance 'social peace' means peace among people, different castes and creeds which alone can lead to social harmony i.e. social peace alone can lead to social unity and above all social stability claims the experts. When will social peace prevail in any nation, only when one and all have an access to good higher education without discrimination. It is a pity in India even after 60 years of independence we lack social peace for the access to good higher education for the majority of the OBC/SC/ST is impossible say not even a dream. The experts ask who is responsible for this; the nation has failed till date to make provisions for the OBC/SC/ST for good higher education even in the educational institutions run by government. It is a pity when this topic is touched or discussed the Brahmins/upper castes become hysteric and the media projects well the cries of hysteria and their by sealing the opportunity for the majority of OBC/SC/ST to get good higher education even from the institutions run by the govt. Economic betterment alone cannot lead to social peace, only social betterment coupled with improved economy can pave way for the better understanding and guaranteeing social peace said the experts. It is a futile work for the govt to help only economically for in their opinion all economic plans coupled only with education have borne good fruits.

24. Social Unity: The experts define this concept in a very unique way. By “Social Unity” they mean the society enjoys a unity; i.e. all castes people welcome any venture of the nation united. Presently they say the nation has no social unity for when reservations in higher learning for OBC was announced the upper castes and Brahmins made indecent “symbolic” protests. OBC however got a complete support from SC/ST. So social unity does not exist in India. We wish to bring to light that suppose the govt announces increase in the Dearness Allowance for central and state govt employees; it is pertinent to mention that a daily wager, a mason, a painter or a sweeper who don’t have leave salary or medical allowance or city allowance or House Rent Allowance (HRA) etc do not benefit in any way. These class of workers who form the majority of the population has no say over this or cannot make a protest for they have no education to understand the benefit enjoyed by the educated and at the same time does not lead to any form of betterment of the socially deprived for majority of these labourers are only from the socially deprived classes i.e.

form the OBC/SC/ST. It is a pity they cannot voice against the govt to scrap, the increase of DA for the employees. Here we cannot say because of social unity the bill was passed. It is only the ignorance of these set of people the bill was passed. But when govt announces a small percentage for them in the field of higher learning the social unity is lost. Here it is pertinent to mention the nation had never enjoyed anything as social unity. That is it is not an exaggeration if we say that the nation till date has never had anything as social unity claims the experts. They further say till the caste system is abolished or till 90% of the OBC/SC/ST have higher education in institution of excellence the nation will never have anything as social unity.

25. Social Stability: Now the experts discuss the concept/attribute social stability i.e. when can we say the nation is socially instable? When they say socially stable they mean that the society is stable and will not become easily displaced for instance even the displacement of a few castes will not upset the stability of the society then we say the nation is socially stable but displacement of the majority will certainly dismantle the stability; so is the issues of reservation for OBC. If OBC reservation is not given after passing it as a bill in both the houses certainly it will affect the stability of the nation for if they join and make protest, no compromise can ever be achieved said the experts (Gujjars Agitation).

Now the experts have given 25 attributes when we said we are going to analyse the problem using super fuzzy models. When we argued about these attributes they also described them in a line or two and said they mean different shades of meaning and they said they have given them as the main tool we are going to use is fuzzy theory. Most of these attributes have exact opposite attributes/concept like social equilibrium; the exact opposite of it being social inequilibrium and so on. Only in few places we have given the opposite attributes also.

## 2.2 Super Row FRM Model to Study the Role of Media on OBC Reservation

In this section we use super FRM model to study the effects of social problems in relation with the media on OBC reservations.

We now list the attributes before we use them in the fuzzy super models for getting the experts opinion.

$S_1$ – Social Rights
$S_2$ – Social Reformation
$S_3$ – Social Justice
$S_4$ – Social Discrimination
$S_5$ – Social Tolerance
$S_6$ – Social Intolerance
$S_7$ – Social Balance
$S_8$ – Social Imbalance
$S_9$ – Social Equilibrium
$S_{10}$ – Social Inequilibrium
$S_{11}$ – Social Conscience
$S_{12}$ – Social Understanding
$S_{13}$ – Social Awareness
$S_{14}$ – Social Commitment
$S_{15}$ – Social Arrogance
$S_{16}$ – Social Responsibility
$S_{17}$ – Social Irresponsibility
$S_{18}$ – Social Maturity
$S_{19}$ – Social Evil/Havoc
$S_{20}$ – Social Equality
$S_{21}$ – Social Harmony
$S_{22}$ – Social Welfare
$S_{23}$ – Social Peace
$S_{24}$ – Social Unity
$S_{25}$ – Social Stability

The attributes $S_1$, $S_2$, $S_3$, $S_4$, $S_7$, $S_{10}$, $S_{11}$, $S_{13}$, $S_{14}$, $S_{16}$, $S_{20}$, $S_{21}$, $S_{22}$ and $S_{25}$ are taken as the attributes of the domain space.

The following attributes related with reservation and media was given by the set of experts A, B, C and D given below are taken as the nodes of the range space of the super FRM.

$R_1$ – Reservation for OBC is …
$R_2$ – Absence of reservation for OBC would lead to …
$R_3$ – Medias dominated by Brahmins/upper caste has worked against reservation causing …
$R_4$ – Political leaders have given reservation to OBC …
$R_5$ – The verdict of court on reservation for OBC would create in the Nation …
$R_6$ – Public feel reservation for OBC would make …

The experts have taken only these 6 major attributes. Four groups of socio scientists and educationalists were taken. The experts from SC/ST formed the group A, OBC formed group B, minorities formed the group C, the upper caste/Brahmins, socio scientists and educationalists formed the group D.

We used the super row fuzzy model with four rows.

| | | $S_1$ | $S_2$ | $S_3$ | $S_4$ | $S_7$ | $S_{10}$ | $S_{11}$ | $S_{13}$ | $S_{14}$ | $S_{16}$ | $S_{20}$ | $S_{21}$ | $S_{22}$ | $S_{25}$ |
|---|---|---|---|---|---|---|---|---|---|---|---|---|---|---|---|
| | $R_1$ | 1 | 1 | 1 | 0 | 1 | 0 | 1 | 1 | 1 | 1 | 1 | 1 | 1 | 1 |
| | $R_2$ | 0 | 0 | 0 | 1 | 0 | 1 | 0 | 0 | 0 | 0 | 0 | -1 | 0 | 0 |
| **A** | $R_3$ | 0 | 0 | 0 | 1 | 0 | 1 | 0 | 0 | 0 | 0 | -1 | 0 | 0 | -1 |
| | $R_4$ | 1 | 1 | 1 | 0 | 1 | 0 | 1 | 1 | 1 | 1 | 1 | 1 | 1 | 1 |
| | $R_5$ | 0 | 0 | -1 | -1 | -1 | 1 | 0 | 0 | 0 | 0 | -1 | 0 | 0 | -1 |
| | $R_6$ | 0 | 1 | 1 | 0 | 1 | 0 | 1 | 1 | 1 | 0 | 0 | 1 | 1 | 1 |
| | $R_1$ | 1 | 1 | 1 | 0 | 1 | 0 | 1 | 1 | 1 | 1 | 1 | 1 | 1 | 1 |
| | $R_2$ | 0 | 0 | 0 | 1 | 0 | 1 | 0 | 0 | 0 | -1 | 0 | -1 | 0 | 0 |
| **B** | $R_3$ | 0 | 0 | 0 | 1 | 0 | 1 | -1 | 0 | 0 | 0 | -1 | -1 | 0 | -1 |
| | $R_4$ | 1 | 1 | 1 | 0 | 1 | 0 | 1 | 1 | 1 | 1 | 1 | 1 | 1 | 1 |
| | $R_5$ | 0 | 0 | -1 | -1 | 0 | 1 | 0 | 0 | 0 | -1 | 0 | 0 | 0 | -1 |
| | $R_6$ | 0 | 1 | 0 | 0 | 1 | 0 | 1 | 1 | 1 | 0 | 0 | 1 | 1 | 1 |
| **C** | $R_1$ | 1 | 1 | 1 | 0 | 1 | 0 | 1 | 1 | 1 | 1 | 1 | 1 | 0 | 1 |
| | $R_2$ | 0 | 0 | -1 | 1 | 0 | 1 | 0 | 0 | -1 | 0 | 0 | -1 | -1 | 0 |
| | $R_3$ | 0 | 0 | 0 | 1 | 0 | 1 | -1 | 0 | 0 | -1 | 0 | 0 | 0 | -1 |
| | $R_4$ | 1 | 1 | 0 | 0 | 1 | 0 | 1 | 1 | 1 | 1 | 1 | 1 | 0 | 1 |
| | **$R_5$** | 0 | -1 | 0 | -1 | 0 | 1 | -1 | 0 | 0 | -1 | 0 | 0 | 0 | -1 |
| | **$R_6$** | 0 | 1 | 1 | 0 | 1 | 0 | 1 | 1 | 0 | 0 | 1 | 1 | 1 | 1 |
| **D** | **$R_1$** | 0 | 0 | 0 | 1 | 0 | 1 | 0 | 0 | 0 | 0 | 0 | 0 | 0 | 1 |
| | **$R_2$** | 1 | 1 | 1 | 0 | 1 | 0 | 0 | 0 | 1 | 0 | 1 | 1 | 1 | 1 |
| | **$R_3$** | 0 | 0 | 0 | 1 | 0 | 1 | 0 | 0 | 0 | 0 | 0 | 0 | 0 | 0 |
| | **$R_4$** | -1 | 0 | -1 | 1 | 0 | 1 | -1 | -1 | 0 | -1 | 0 | 0 | -1 | 0 |
| | **$R_5$** | 1 | 1 | 1 | 0 | 0 | 0 | 1 | 1 | 1 | 1 | 0 | 1 | 1 | 1 |
| | **$R_6$** | 0 | -1 | 0 | 0 | -1 | 0 | -1 | 0 | 1 | 0 | -1 | 0 | 0 | -1 |

Hence it is important to mention that when -1 is used it means the negation of the statement, i.e. if one expert gives value -1 to social right it implies it is not a social right or just like 0 which means nothing related to social right. The use of -1 only gives more emphasis in the negative sense.

Now the super row fuzzy relation maps model M is used in the analysis of 4 different categories of socio scientists about the impact of reservation in social angles.

Suppose the node social right is in the on state we study the effect of the on state on the node $S_1$ on the dynamical super system when all other nodes are in the off state that is the effect of the 4 groups of people about reservations for OBC in educational institutions of higher learning run by the central government.

Thus the given state vector

$$X = (1\ 0\ 0\ 0\ 0\ 0\ 0\ 0\ 0\ 0\ 0\ 0\ 0\ 0).$$

Now the effect of X on the dynamical system is given by

$$XM^T = (1\ 0\ 0\ 0\ 0\ 0\ 0\ 0\ 0\ 0\ 0\ 0\ 0\ 0) \times M^T$$

$$
= \left[\begin{array}{cccccc|cccccc}
1 & 0 & 0 & 1 & 0 & 0 & 1 & 0 & 0 & 1 & 0 & 0 \\
1 & 0 & 0 & 1 & 0 & 1 & 1 & 0 & 0 & 1 & 0 & 1 \\
1 & 0 & 0 & 1 & -1 & 1 & 1 & 0 & 0 & 1 & -1 & 0 \\
0 & 1 & 1 & 0 & -1 & 0 & 0 & 1 & 1 & 0 & -1 & 0 \\
1 & 0 & 0 & 1 & -1 & 1 & 1 & 0 & 0 & 1 & 0 & 1 \\
0 & 1 & 1 & 0 & 1 & 0 & 0 & 1 & 1 & 0 & 1 & 0 \\
1 & 0 & 0 & 1 & 0 & 1 & 1 & 0 & -1 & 1 & 0 & 1 \\
1 & 0 & 0 & 1 & 0 & 1 & 1 & 0 & 0 & 1 & 0 & 1 \\
1 & 0 & 0 & 1 & 0 & 1 & 1 & 0 & 0 & 1 & 0 & 1 \\
1 & 0 & 0 & 1 & 0 & 0 & 1 & -1 & 0 & 1 & -1 & 0 \\
1 & 0 & -1 & 1 & -1 & 0 & 1 & 0 & -1 & 1 & 0 & 0 \\
1 & -1 & 0 & 1 & 0 & 1 & 1 & -1 & -1 & 1 & 0 & 1 \\
1 & 0 & 0 & 1 & 0 & 1 & 1 & 0 & 0 & 1 & 0 & 1 \\
1 & 0 & -1 & 1 & -1 & 1 & 1 & 0 & -1 & 1 & -1 & 1
\end{array}\right.
$$

$$
\left|\begin{array}{cccccc|cccccc}
1 & 0 & 0 & 1 & 0 & 0 & 0 & 1 & 0 & -1 & 1 & 0 \\
1 & 0 & 0 & 1 & -1 & 1 & 0 & 1 & 0 & 0 & 1 & -1 \\
1 & -1 & 0 & 0 & 0 & 1 & 0 & 1 & 0 & -1 & 1 & 0 \\
0 & 1 & 1 & 0 & -1 & 0 & 1 & 0 & 1 & 1 & 0 & 0 \\
1 & 0 & 0 & 1 & 0 & 1 & 0 & 1 & 0 & 0 & 0 & -1 \\
0 & 1 & 1 & 0 & 1 & 0 & 1 & 0 & 1 & 1 & 0 & 0 \\
1 & 0 & -1 & 1 & -1 & 1 & 0 & 0 & 0 & -1 & 1 & -1 \\
1 & 0 & 0 & 1 & 0 & 1 & 0 & 0 & 0 & -1 & 1 & 0 \\
1 & -1 & 0 & 1 & 0 & 0 & 0 & 1 & 0 & 0 & 1 & 1 \\
1 & 0 & -1 & 1 & -1 & 0 & 0 & 0 & 0 & -1 & 1 & 0 \\
1 & 0 & 0 & 1 & 0 & 1 & 0 & 1 & 0 & 0 & 0 & -1 \\
1 & -1 & 0 & 1 & 0 & 1 & 0 & 1 & 0 & 0 & 1 & 0 \\
0 & -1 & 0 & 0 & 0 & 1 & 0 & 1 & 0 & -1 & 1 & 0 \\
1 & 0 & -1 & 1 & -1 & 1 & 1 & 1 & 0 & 0 & 1 & -1
\end{array}\right]
$$

$$= \quad (1\ 0\ 0\ 1\ 0\ 0 \mid 1\ 0\ 0\ 1\ 0\ 0 \mid 1\ 0\ 0\ 1\ 0\ 0 \mid 0\ 1\ 0\ {-1}\ 1\ 0)$$

= Y'.

Y' is thresholded to get the super fuzzy row vector Y= (1 0 0 1 0 0 | 1 0 0 1 0 0 | 1 0 0 1 0 0 | 0 1 0 0 1 0).
Now we study the effect of Y on M;

$$YM \hookrightarrow (1\ 1\ 1\ 0\ 1\ 0\ 1\ 1\ 1\ 1\ 1\ 1\ 1\ 1) = Z.$$

(Here '$\hookrightarrow$' denotes the resultant vector has been updated and thresholded). Thus only the nodes social discrimination and social inequality is in the off state. Thus we see this model gives the final majority opinion but the super fuzzy vector T shows group D is however acting just opposite or distinctly different from the other 3 groups for them social right is not reservation but social right is established only but not following reservation and for them the verdict of the S.C on reservation is a social right. Now we study the effect of Z on the super fuzzy dynamical system.

$$ZM^T \hookrightarrow (1\ 0\ 0\ 1\ 0\ 1\ |\ 1\ 0\ 0\ 1\ 0\ 1\ |\ 1\ 0\ 0\ 1\ 0\ 1\ |\ 0\ 1\ 0\ 0\ 1\ 0) = W;$$

the fuzzy super row vector after thresholding it. We see from this resultant super fuzzy row vector that when the social right node alone was in the on state, the nodes $R_2$, $R_3$ and $R_5$ were in the off state for all the 3 group of experts A, B and C and for the group D it was very surprising to see that according to them absence of reservation for OBC would lead to the social right in the society and the verdict of the court on reservation for OBC would create social right in the nation. Thus we see from this analysis using the mathematical model viz; fuzzy super FRM model we see that the hidden pattern clearly professes the mind set of the Brahmins and the upper castes. It is very surprising to see that the Brahmins / upper castes backed by the Hindutva forces like BJP/RSS/VHP have the guts to say things right opposite to the majority more so against the natural justice or against a natural phenomenon. The analysis only shows that these minority Brahmins who are not even the sons of the soil have the guts to speak against the social reforms for the majority who are not only denied education but all forms of comfort in their life. How to view this selfish attitude? Most of the experts say the Hindutva political leaders of BJP/RSS/VHP are much against the very thought of educating the Sudras or the OBC for it is against the laws of Manu, further it is higher education alone that can make them think so, the Brahmins do not want the majority to think for if they start to 'think' they will not be doing slavish work for them so they are frightened not only they will become their competitors but also they will not have someone to work for them for very paltry pay. Money is the only goal coupled with protecting their identity through the Hindutva cover. That is why they hold to the religion which was never pro founded properly except which gave the laws by Manu which can be viewed as an anti social anti justice, discriminative selfish code to promote Brahmins and Brahmins alone. It is still unfortunate to note that they wrote some laws to flourish in an alien land and how can this be in any way suitable for the toiling patient mass; sons of the soil condemns the experts. Further as they were the policy makers from the good old time till date the laws pronounced is followed blindly with no rhyme or reason by our honourable courts claimed the experts. Who are our law givers, Brahmins or those who protect Hindutva for by doing so they can climb the ladder of success said the experts. Until the Hindu religion is replaced by original Indian culture, India can never dream to become a self sufficient nation. The developments, which it boasts of, is only the development of the paltry 3% foreigners, Brahmins said the experts. The status of the last man has only deteriorated in these 60 years. This is very evident from the mass suicide of the farmers! Farmers once were adorned by all the poets for they were worshiped for their selfless service. They never exploited the land and they were worshippers of nature, which gave prosperity to their profession. The Brahmins cheat the common man with "dharba grass" with no knowledge in any science or technology they have been cheating the nation and have taken the highest position of policy makers to win the nation! Is this not clear from the acts of NKC asked the experts?

Suppose we take a state super vector from the domain space to find the effect of it on the super dynamical system M.

Let

$$S = (1\ 0\ 0\ 0\ 0\ 0 \mid 1\ 0\ 0\ 0\ 0\ 0 \mid 1\ 0\ 0\ 0\ 0\ 0 \mid 1\ 0\ 0\ 0\ 0\ 0)$$

belong to the domain space of the fuzzy super FRM.

$$SM \hookrightarrow (1\ 1\ 1\ 0\ 1\ 0\ 1\ 1\ 1\ 1\ 1\ 1\ 1\ 1).$$

(Here it has become important to mention that since we are working with the same set of attributes in the domain space for all the four groups of experts A, B, C and D while thresholding make in the resultant vector

$$T = (a_1, \ldots, a_{14});$$

if $a_i \geq 2$ put 1, if $a_i \leq 2$, put 0) "$\hookrightarrow$ symbol denotes the state vector has been updated and thresholded".

Let

$$S.M \hookrightarrow (1\ 1\ 1\ 0\ 1\ 0\ 1\ 1\ 1\ 1\ 1\ 1\ 1\ 1) = Y.$$

Now

$$Y \cdot M^t \hookrightarrow (1\ 0\ 0\ 1\ 0\ 1 \mid 1\ 0\ 0\ 1\ 0\ 1 \mid 1\ 0\ 0\ 1\ 0\ 1 \mid 1\ 1\ 0\ 0\ 1\ 0).$$

Thus when the node "$R_1$ – Reservation for OBC is" alone in the on state we see according to group D only and mainly the nodes $R_2$ and $R_5$ has come to on state and this makes them feel that only absence of reservation for OBC and the S.C verdict to stay OBC reservation are ones that will give social rights, social reformation, social justice and so on. It is once again important and pertinent to mention from the analysis using this super model leads us to understand that the Brahmins has that height of arrogance to talk totally against reservations for OBC as if it is a social crime. The study of super fuzzy mathematical model projects their true character says the experts. Their views are just opposite to that of the other three groups. How very fanatic, the Brahmins can behave when it comes to reservations one can easily understand. Even while getting their opinions as experts. They were very much against it but they know us that we do not at any place dissuade any form of opinion from any one who wishes to give, for we are just set to work only after truth and truth alone. So any opinion for or against does not matter to us. Only during and after analysis they were the group interested to know the results and we just showed them the analysis also. We just mention the effect of yet another state vector from the range space of the super fuzzy model. Let

$$X = (0\ 0\ 0\ 1\ 0\ 0\ 0\ 0\ 0\ 0\ 0\ 0\ 0\ 0)$$

be the given state vector where only the node $S_4$ from the range space is in the on state; i.e. only the node social discrimination is in the on state and all other nodes are in the off state. The effect of X on the fuzzy super dynamical system M is given by

$$XM^T \hookrightarrow (0\ 1\ 1\ 0\ 0\ 0 \mid 0\ 1\ 1\ 0\ 0\ 0 \mid 0\ 1\ 1\ 0\ 0\ 0 \mid 1\ 1\ 0\ 1\ 1\ 0\ 0) = Y.$$

Let Y denote the effect of X on $M^T$ where Y is a super fuzzy row vector. Now we see

$$YM \hookrightarrow (0\ 0\ 0\ 1\ 0\ 0\ 0\ 0\ 0\ 0\ 0\ 0\ 0\ 0) = X' = X \text{ (say)}.$$

Thus we see the super hidden pattern of the super fuzzy dynamical system is given by the binary pair {(0 0 0 1 0 0 0 0 0 0 0 0 0 0), (0 1 1 0 0 0 | 0 1 1 0 0 0 | 0 1 1 0 0 0 | 1 0 1 1 0 0)}. We see when $S_4$ alone is the on state in the resultant, all the three groups A, B and C, the resultant reads as if the node social discrimination alone is the on state the nodes $R_2$ and $R_3$ alone come to on state there by indicating that absence of reservation for OBC would lead to social discrimination and media's dominated by Brahmins/upper castes has projected reservation as social discrimination. The group D alone gives that social discrimination is caused by reservation for OBC. $R_1$, the media dominated by the Brahmins has projected reservation as a social discrimination – $R_3$ and $R_4$ which says the giving of reservation by the political leaders to OBC is a social discrimination. Thus we see the way Brahmins who are not the sons of the soil view the very reservation for OBC in the central govt institutions of higher learning as a social discrimination.

Further it is interesting to note that the resultant of X gave a super hidden pattern which was a fixed point who also from this solution says that Brahmins have no mind to accept the reservation even for a short period.

Here having studied the situation using super fuzzy FRM model, we give here the exact number of experts involved in each of the groups, A, B, C and D. The group A contained 241 experts opinion though we targeted for 250. The group B comprised of 492 experts opinion though we used 500 experts. Group C contained 91 experts opinion though we used 100 experts. Finally the group D had 150 experts and 142 of them alone gave their opinion. From the critical study we used the super fuzzy FRM model to give the super hidden pattern. All the resultants give only super hidden pattern which was only a fixed point.

Mathematical observations made from the super fuzzy FRM model:

1. The resultant of every state vector gave way to the super hidden pattern which was always a fixed point. In any dynamical system if we get the resultant for every state vector to be a fixed point it implies that the problem under study does not yield any changes in due course of time i.e. the problems are time independent says the experts. Further the views are not flexible.

2. The on state of the only node social discrimination does not yield or make any other node of the range space of the super fuzzy FRM to comes to on state which clearly implies the reservation for OBC in the educational institutions of higher learning is not a social discrimination only the denial of the reservation for OBC in these institutions is a social discrimination. Further resultant hidden pattern shows that this state vector from the range space is a fixed point and it had nothing to do with the other nodes of the range space.

   On the contrary the on state of this node $S_4$ made the nodes $R_2$ and $R_3$ of the domain space to be in the on state in the three groups A B and C. Further the on state of $R_2$ and $R_3$ only meant that the absence of reservation for OBC is a social discrimination for these institutions are after all run by the central govt. and so at least in these institutions one needs reservations for these OBC to make some improvement in their social set up. Also the act of media which is dominated by the Brahmins was always projecting the reservation for the OBC in higher learning in the govt run institution as a social discrimination. Now the on state of the only node $S_4$, i.e., social discrimination gives the on state of the nodes $R_1$, $R_3$ and $R_4$ by the group D people; which according to them imply that reservation given to the OBC is a social discrimination; media also have projected reservation as social discrimination and the act of the political leaders to announce reservation for OBC is also a social discrimination!

3. The on state of the node $S_{10}$ in the range space of the super fuzzy dynamical system results in the on state of $S_4$, i.e., the social discrimination. Thus social inequilibrium is caused by the social discrimination and no other node comes to on state in the range space. The on state of social inequilibrium results from the absence or denial of reservation for OBC in educational institutions of higher learning run by the central govt. Also the social inequilibrium has been caused because of the false propaganda about reservation of the Brahmin dominated media. The social inequilibrium is caused by the verdict of the S.C on staying the reservation for OBC have come to on state in the fuzzy super row vector of the domain space in the group A, B and C and all other nodes in A B and C remain in the off state. Further from group D it is shown that social inequilibrium is because of reservation for OBC in educational institutions of higher learning. The media is also a cause for social inequilibrium and political leaders who have supported the cause for passing the order about reservation for OBC in higher learning. Thus these Brahmin socio scientist have given just the opposite of what has been given by these OBC socio scientist group of experts, socio scientist group formed by SC/ST and the minorities. Thus we see from this study; specially for the Brahmin, when it comes to nations prosperity versus their selfish interests, comforts and personal profits they ruthlessly and shamelessly stand for their selfish interests comforts and personal profits alone. This is clear from the results derived using the super fuzzy FRM model.
4. It is still surprising to see that our study reveals that when only the node $R_3$ of the domain space in the on state in all

the four groups A, B, C and D given by the fuzzy super row vector

Y = (0 0 1 0 0 0 | 0 0 1 0 0 0 | 0 0 1 0 0 0 | 0 0 1 0 0 0).

We see

$$YM \hookrightarrow (0\ 0\ 0\ 1\ 0\ 1\ 0\ 0\ 0\ 0\ 0\ 0\ 0\ 0).$$

Thus in this case also for the state fuzzy super vector Y once again the super hidden pattern is a fixed pair given by

{(0 0 0 1 0 1 0 0 0 0 0 0 0 0),
(0 1 1 0 0 0 | 0 1 1 0 0 0 | 0 1 1 0 0 0 | 1 0 1 1 0 0)}.

Thus we see the act of media would make the nation suffer from social discrimination and social inequilibrium. Thus the analysis clearly shows the media has absolutely failed to make any form of positive impact on the nation when it made propaganda against reservation. Further this social discrimination and social inequilibrium has resulted from the anti reservation trend set up the media.

5. The present social status of the nation is that the denial of reservation to the OBC in the institutions of higher learning is a threat to social unity, denial of social rights of majority, it is against the principal of social reformation and social justice. Their act of staying the OBC reservations amounts to social discrimination, the nation will be forced to act against the very norms of social tolerance, soon the social balance would be lost resulting in riots and revolution by the majority calm OBC resulting in social imbalance of the nation. Once the social equilibrium of the nation is lost the nation would become dismantled and disintegrated like Russia claims the socio scientists. The very denial of reservation to the OBC shows the nation is acting without any social conscience or social commitment or social understanding of the dire consequences without any social awareness because the few who talk against the reservation for OBC are the upper caste / Brahmins who are a socially arrogant do not know their social responsibility for they are foreigners to our nation so only act in this irresponsible way! They mainly lack social maturity. They want to grab the maximum economy, comfort and education from the son of the soil. These acts of the Brahmins will certainly ruin the social harmony, social welfare, social peace, social unity and finally the social stability of the nation resulting in permanent disaster!

### 2.3 Super Fuzzy Mixed FRM Model to Study the Role of Media in Falsely Blaming the Government and Supporting Dr. Venugopal

They have given the following attributes to work with this very extra ordinary situation. Each concept or attribute given by them is described in a line or two for the reader to understand in what direction or meaning these concepts / attributes are used in this super fuzzy model.

$M_1$ - The Brahmin dominated media projected the government (Minister) as folly in the administrative work.

$M_2$ - The Brahmin dominated media projected Dr. Venugopal as a world level cardiologist and a man of repute and honesty.

$M_3$ - The Brahmin dominated media was neutral to both.

$M_4$ - The casteist media projected only upper caste as perfectionist.

$M_5$ - The authority Sudra was insulted beyond limits by the media.

$M_6$ - The subordinate Brahmin was portrayed as a man of virtue and renowned doctor.

$M_7$ - The media took sides on this issue by writing the support of BJP to Dr. Venugopal.

$M_8$ - The media failed to see the true problems in AIIMS as pointed out by the Health Minister who is the president of AIIMS.

$M_9$ - The media failed to project the protests carried out in the AIIMS campus which was backed by Dr.Venugopal and the Hindutva fanatics.

$M_{10}$ - The media failed to investigate the truth behind the issue with ulterior motives and selfish intensions prevailing in the admissions in AIIMS.

$M_{11}$ - The protest carried out by youth for equality in the AIIMS amounting to disobedience of the High Court order was not properly projected by the Brahmin media.

$M_{12}$ - Till date the wrong and the unlawful ways in which the P.G medical seats are filled in AIIMS is not discussed by the media.

$M_{13}$ - The media failed to bring out the patients death caused by the striking doctors of AIIMS. No investigation about the death was carried out.

$M_{14}$ - Court gave the marvellous order that they should be paid for the period of strike which spread over 40 days. Their absence of work and that too agitating against the govt under which they were working was a blow to the norms of natural justice! The experts ask will ever any group from OBC get this sort of judgement. Even daily wagers when they go on strike they do not get their pay. The only reason is backed by Hindutva and religious leaders!

$M_{15}$ - The administrative frauds carried out by Dr. Venugopal is not brought out by the media or by the court openly.

These are the 15 attributes mainly agreed by all experts. Now as groups they will be used to construct the fuzzy super model. The four categories of experts are doctors, politicians, public and medical students. This is taken as the domain space attributes. Both the domain and range space attributes are super fuzzy row vectors.

The nodes of the domain space given by the doctors are

$D_1$ - doctors belonging to SC/ST castes

$D_2$ - doctors belonging to OBC

$D_3$ - doctors belonging to Brahmin communities/upper castes

$D_4$ - doctors from minority religions

$D_5$ - doctors who are neither against nor favour reservation.

The nodes related with politicians. We had to take the opinion of the politicians for this case is the Union Cabinet Central Minister versus the Brahmin director and the move to give reservations for the OBC was only contemplated by the political leaders to socially uplift the OBC. So we are forced to interview and take their opinion also in this regard.

$P_1$ - Political leaders belonging to OBC
$P_2$ - Political leaders belonging to SC/ST
$P_3$ - Political leaders from upper castes and Brahmins
$P_4$ - Political leaders from minority communities

We do not accept passive political leaders for they cannot stay neutral on a social issue. Because they have voted in both the houses no body ever said he/she was neutral in both the houses.

Next we give the list of public, who are experts to give their views

$C_1$ - Educated public from the OBC
$C_2$ - Uneducated public from the OBC
$C_3$ - Educated public from SC/ST
$C_4$ - Uneducated public from SC/ST
$C_5$ - Public from the upper castes/Brahmins
$C_6$ - Educated public from the minorities
$C_7$ - Uneducated public from the minorities

Now the list of students who given their views about the occurrences in the AIIMS under director, Dr. Venugopal.

$S_1$ - Medical OBC students from AIIMS and other medical institutes
$S_2$ - Medical upper caste and Brahmin students from AIIMS / other medical colleges
$S_3$ - OBC students other than medicine
$S_4$ - Upper caste/Brahmin students from other faculties other than medicine
$S_5$ - SC/ST students from all faculties
$S_6$ - Minority students from all faculties

Now we give choice for the four groups to select any number of attributes from the given 15 attributes. These 15 attributes were also formed using the experts from all the categories. Now we use fuzzy model in general and super fuzzy mixed FRM model in particular to analyze the problems which is given by the super fuzzy matrix M

| | $M_1$ | $M_2$ | $M_3$ | $M_4$ | $M_5$ | $M_6$ | $M_7$ | $M_8$ | $M_9$ | $M_{10}$ | $M_{11}$ | $M_{12}$ | $M_{13}$ | $M_{14}$ | $M_{15}$ |
|---|---|---|---|---|---|---|---|---|---|---|---|---|---|---|---|
| $D_1$ | 1 | 1 | 0 | 1 | 1 | 1 | 0 | 1 | 0 | 1 | 1 | 1 | 1 | 1 | 1 |
| $D_2$ | 1 | 1 | 0 | 1 | 1 | 1 | 0 | 1 | 0 | 1 | 1 | 1 | 1 | 1 | 1 |
| $D_3$ | 0 | 0 | 1 | 0 | 0 | 0 | 1 | 0 | 1 | 0 | 0 | -1 | 0 | 0 | -1 |
| $D_4$ | 1 | 1 | 0 | 1 | 1 | 1 | 0 | 1 | 0 | 1 | 1 | 1 | 1 | 1 | 0 |
| $D_5$ | 1 | 1 | 0 | 1 | 1 | 1 | 0 | 1 | 0 | 1 | 0 | 0 | 1 | 1 | 1 |
| $P_1$ | 1 | 1 | 0 | 0 | 1 | 1 | 0 | 1 | 0 | 1 | 1 | 1 | 1 | 1 | 1 |
| $P_2$ | 1 | 1 | -1 | 0 | 1 | 1 | 0 | 1 | -1 | 1 | 1 | 1 | 1 | 1 | 1 |
| $P_3$ | 0 | 0 | 1 | 1 | 0 | 0 | 1 | 0 | 1 | -1 | 0 | 0 | 0 | 0 | -1 |
| $P_4$ | 1 | 1 | 0 | 0 | 1 | 1 | -1 | 0 | 0 | 1 | 1 | 1 | 1 | 1 | 0 |
| $C_1$ | 1 | 1 | 0 | 0 | 1 | 1 | 0 | 1 | 0 | 1 | 1 | 1 | 1 | 1 | 1 |
| $C_2$ | 1 | 1 | -1 | 0 | 1 | 1 | 0 | 1 | -1 | 1 | 1 | 0 | 1 | 1 | 0 |
| $C_3$ | 1 | 1 | 0 | -1 | 1 | 1 | -1 | 1 | 0 | 1 | 1 | 1 | 1 | 1 | 1 |
| $C_4$ | 1 | 1 | 0 | -1 | 1 | 1 | 0 | 1 | -1 | 0 | -1 | 1 | 1 | 0 | 0 |
| $C_5$ | 0 | 0 | 1 | 1 | 0 | 0 | 1 | 0 | 1 | 0 | 0 | 0 | -1 | -1 | -1 |
| $C_6$ | 1 | 1 | 0 | 0 | 1 | 1 | 0 | 1 | 0 | 1 | 1 | 0 | 1 | 1 | 1 |
| $C_7$ | 1 | 0 | 0 | 1 | 1 | 1 | 0 | 1 | 0 | 0 | 1 | 1 | 1 | 1 | 1 |
| $S_1$ | 1 | 1 | 0 | 1 | 1 | 1 | 0 | 1 | 0 | 0 | 0 | 0 | 0 | 0 | 0 |
| $S_2$ | 0 | 0 | 1 | 0 | 0 | 0 | 1 | 0 | 1 | 0 | 0 | 0 | 0 | 0 | 0 |
| $S_3$ | 1 | 1 | 0 | 0 | 1 | 1 | 0 | 1 | -1 | 0 | 0 | 0 | 0 | 0 | 0 |
| $S_4$ | -1 | 0 | 1 | 0 | -1 | 0 | 1 | -1 | 1 | 0 | 0 | 0 | 0 | 0 | 0 |
| $S_5$ | 1 | 1 | 0 | 1 | 1 | 1 | -1 | 1 | -1 | 0 | 0 | 0 | 0 | 0 | 0 |
| $S_6$ | 1 | 1 | -1 | 0 | 1 | 0 | 0 | 1 | 0 | 0 | 0 | 0 | 0 | 0 | 0 |

Now having given the experts opinion in the form of a super fuzzy FRM connection matrix we find the super hidden pattern for one or two vectors given by the experts, here. However we have worked with several to derive the conclusion. Suppose one is interested in studying the situation when the node $M_3$ alone of the range space is in the on state and all other nodes are in the off state. To study the effect of the statement "The Brahmin dominated media was neutral to both". Thus let

$$X = (0\ 0\ 1\ 0\ 0\ 0\ 0\ 0\ 0\ 0\ 0\ 0\ 0\ 0\ 0)$$

be the given state vector; to study the effect of X on the super fuzzy dynamical system.

$$XM^T \hookrightarrow (0\ 0\ 1\ 0\ 0 \mid 0\ 0\ 1\ 0 \mid 0\ 0\ 0\ 0\ 1\ 0\ 0 \mid 0\ 1\ 0\ 1\ 0\ 0) = Y.$$

We find

$$Y \cdot M \hookrightarrow (0\ 0\ 1\ 0\ 0\ 0\ 1\ 0\ 1\ 0\ 0\ 0\ 0\ 0\ 0) = X_1 \text{ (say)}.$$

$$X_1 \cdot M^T \hookrightarrow (0\ 0\ 1\ 0\ 0 \mid 0\ 0\ 1\ 0 \mid 0\ 0\ 0\ 0\ 1\ 0\ 0 \mid 0\ 1\ 0\ 1\ 0\ 0)$$

after appropriate thresholding of the resultant we get $Y_1$ which is nothing, but Y leading to the fixed binary pair {(0 0 1 0 0 | 0 0 1 0 | 0 0 0 0 1 0 0 | 0 1 0 1 0 0), (0 0 1 0 0 0 1 01 0 0 0 0 0 0)}. We see the on state of the node $M_3$ – The media is neutral to both alone lead to a fixed point which makes only the nodes $M_7$ and $M_9$ to the on state in the range space i.e. the media was not neutral for it took sides on this issue by writing the support of BJP to Venugopal for the BJP asked Anbumani to resign for the post which shows the media was biased and was acting openly against Dr. Anbumani Ramadoss. Further the media was not neutral for it did not project the true fact that the protests carried out in the campus was backed by Dr.Venugopal and Dr.Venugopal was in turn protected and strongly backed by the Hindutva fanatics. On the other hand the resultant fixed point from the domain super space made only the nodes $D_3$, $P_3$, $C_5$, $S_2$ and $S_4$ to the on state which clearly shows only the doctors belonging to upper caste or Brahmins dire say the media is neutral in both the cases. Also the political leaders from the upper castes and Brahmins claim the media has acted neutral in both the cases. The public from these sections also blindly say the media has acted neutral. Further only the upper castes and Brahmin medical students as well as the students from the other faculties accept that the medical has acted neutral in both the cases. Consequently all these have made on the state that BJP is supporting Dr. Venugopal activities and the media has completely black out the fact that all protests made in the AIIMS is instigated at large but the director Dr. Venugopal and the Hindutva fanatics and certain Hindu religious heads form the south.

Suppose only the node $M_5$ alone in the range space is in the on state in the given state vector and all other nodes are in the off state. To find the super hidden pattern of the state vector using this super dynamical system. Let

$$T = (0\ 0\ 0\ 0\ 1\ 0\ 0\ 0\ 0\ 0\ 0\ 0\ 0\ 0\ 0)$$

be the given state vector.

The effect of T on $M^T$ is given by

$$T \cdot M^T \hookrightarrow (1\ 1\ 0\ 1\ 1\ |\ 1\ 1\ 0\ 1\ |\ 1\ 1\ 1\ 1\ 0\ 1\ 0\ |\ 1\ 0\ 1\ 0\ 1\ 1) = R;$$

where R is a super fuzzy mixed row vector.

Now we calculate

$$R \cdot M \hookrightarrow (1\ 1\ 0\ 1\ 1\ 1\ 0\ 1\ 0\ 1\ 1\ 1\ 1\ 1\ 1) = T_1.$$

$$T_1 \cdot M^T \hookrightarrow (1\ 1\ 0\ 1\ 1\ |\ 1\ 1\ 0\ 1\ |\ 1\ 1\ 1\ 1\ 0\ 1\ 0\ |\ 1\ 0\ 1\ 0\ 1\ 1) = R_1\ (= R).$$

Thus super hidden pattern is given by the binary pair {(1 1 0 1 1 1 0 1 0 1 1 1 1 1 1), (1 1 0 1 1 | 1 1 0 1 1 | 1 1 0 1 | 1 1 1 1 0 1 1 | 1 0 1 0 1 1)}.

Thus we see that when only the input the Brahmin media has treated Dr. Anbumani Ramadoss with bias by insulting him in every possible way all the nodes in the range space come to on state except $M_3$ , $M_7$ and $M_9$ which implies that the Brahmin media was not neutral in the issue, media was partial by taking sides and reporting in support of Dr. Venugopal and above all the media failed to project the truth that Dr. Venugopal was behind all the anti reservation protests carried out in the AIIMS campus which was against the order of the high court, that such type of protests or agitations cannot be held with in the campus.

On the contrary the director had helped the anti reservations protestors by putting a Shamiana in the premises of AIIMS and the facility of air coolers that is why they could lie down in the shade reading novel / sleeping, playing cards, can carry the protests in a novel way ! Finally from the super fuzzy row vector given as the hidden pattern we see that only the Brahmin groups in all the four categories did not agree to the fact that the media had insulted the Sudra authority limitlessly. The authors have worked with several on state of the attributes and have given the consolidated analysis as observations.

Observation derived from the super fuzzy mixed FRM model is as follows:

1. The media failed to be neutral while reporting the misunderstandings between Dr. Venugopal and Dr. Anbumani Ramadoss. The Sudra authority was portrayed as an embodiment of faults and on the contrary the sub ordinate Brahmin Dr. Venugopal was pictured as a dutiful renowned professional. The media failed to show the nasty side of Dr. Venugopal who could instigate the students to protest and abstain from their duty for a long period of over 40 days. Is this a dutiful character of Dr.Venugopal or a selfish Brahmin casteist personality? The media has not only biased and blocked the facts but has given false perspective of the whole situation.

2. The media which should have acted in a way to educate people of the true situation has lied and not only made false propaganda against reservation but also against the president Dr. Anbumani Ramadoss of AIIMS. For instance only after 20 days of stir Dr. Anbumani Ramadoss went to AIIMS and was checking patients, this was reported in the first page of the 'The New Indian Express' dated 1-6-2003. This was very indecently brought out in the media saying he wanted cheap popularity so only he came to work in AIIMS. By chance a Brahmin Union Minister they would have said the simple nature of the Union Minister with a resolute duty conscience have come to serve the suffering people and the director would have been not only sacked from the post but sent home and the courts would justify the acts. What is law? What is the power of Brahmins?

3. When AIIMS sacks Venugopal they said it was a blatant slap on medical profession – it is reported in the first page. Sacked AIIMS chief moves to HC and in that he challenged Dr. Anbumani Ramadoss to be holding an office of profit as AIIMS President. With this sacking of Dr.Venugopal once again he instigated them to continue the strike which was never brought out in the Brahmin media. On 6th, he approaches the court and on 8th he returns to AIIMS as director. Immediately BJP advised Dr. Anbumani Ramadoss to stop meddling with AIIMS. So a sudra in power who is an authority not only for AIIMS but for Dr. Venugopal as a President of AIIMS cannot question Venugopal says the Brahmin monopolized BJP and is supported by the Brahmin media. This sub ordinate Brahmin is making a plea in the Delhi high court to disqualify Dr. Anbumani and his authority. The reader here is requested to contemplate if at all can ever a non Brahmin subordinate go to court against the Brahmin authority. Will ever the court accept such a plea? The paper projected the lasting pain of Dr.Venugopal who said to the reporters that he joined as a 18 year old from a village in A.P Rajamundhry and his credentials were boosted to a very high degree. BJP openly sought the resignation for Dr. Anbumani. The H.C, Delhi on Aug 24 has warned the Union Minister Dr.Anbumani from acting as president of the AIIMS governing body if the centre failed to file a reply on his suitability for the post. So a Brahmin subordinate questions the suitability of the Sudra authority using the H.C, for always from good old days law had been their monopoly. Can a non Brahmin do all these acts?

4. It is an important factor to observe that when the opinions are taken from all categories of people, especially we see that the Brahmins and some of the upper castes always hold a view just in a diametrically opposite way. This only clearly shows that the Brahmins are highly concerned of themselves and their community. A small inconvenience or to be more critical any form of encouragement especially in education given to the OBC (Sudra) they will revolt as if the govt and the nation are at heading towards disaster. This is mainly done by them for it is very clearly given in the laws of Manu that the Sudra, SC/ST should not be given education, so any act passed or practised in India which is against the laws of Manu will not only be condemned by them and they will also hand in hand be supported by high court and supreme court. So will they give stay by setting their own people to file a PIL (public interest litigation) and the H.C | SC will immediately admit their case and judgement would by all means be against the majority and against all norms of justice especially for the OBC/SC/ST. The important factor to be observed here is the Brahmins are very much scared that if the OBC/SC/ST by chance become educated they will certainly outshine them and this will not only make them shun also they (OBC/SC/ST) would start to think and question their atrocities. That is why all their opinions in case of reservations for OBC in education particularly in higher learning, is right opposite to the non Brahmins view. This has been observed by the facts spoken by them, but also this has been established mathematically using fuzzy super models from the data collected from them.

5. It is still important to mention that even while getting the data from them several of them told the authors that even if their views were against the public or majority or in particular non Brahmins it should be given as it is. To the best of authors' knowledge none of their views have been ever touched or polished or formatted but given as it is only used in these models for solution. In fact we authors have in no way made any change in any views given by anyone for we want to get a true conclusion, a true vision by this scientific analysis. It is still important for the authors to record that their view or their opinion has not been given in any place. We authors have worked like computers to record analyse and give the solution. Most of the expert opinion if they were educated was obtained by writing only for the uneducated technically trained experts / students filled in the questionnaire for them.

6. Only while discussing and getting the results about the stand of the media every one including movement and several non Brahmin Central Cabinet Ministers made it very clear that their views about OBC reservations in the higher education, that too the institutions run basically by the central govt was never recorded. Thus the media mainly dominated by the Brahmins have never done their duty or failed to do their duty by bias reporting about OBC reservations. The media which should have acted very responsible in the social issues have acted in a very irresponsible way. This they have done on selfish stand forgetting about the real development of the nation and above all about the true problems of the majority. How can the nation deny for over 60 years to the majority mass good higher education? If the Brahmin mass which is a paltry Hindu minority with 3% population that too who are the foreigners of this nation have been not only denying the majority of their rights to get good education from the nation but on the top of it, they are continuously and ceaselessly exploited by the non Brahmins for their labour which is so poorly paid that they work for over 8 hours a day and yet cannot afford to have a square meal. Then why is the media supporting them and the court helping them as if they are the victims of reservation. This shows that the laws of Manu is still ruling the nation and nothing more that is why even now the court, the media and all the Brahmins are dead against reservation in higher education for OBC. If Sudras become educated how will the legal dharma expounded in the laws of Manu be followed. That is why they are perfect in their dealing for a Sudra should not be educated! So almost all the models very clearly exhibit that the Brahmins and the Brahmin media are against reservation for the OBC in the higher learning that is why the Supreme Court are giving stay in keeping with the laws of Manu.

7. The media has been very biased in protecting Dr. Venugopal who had done all acts against reservation, he is the root cause for all the protests in the north by the medicos. Without his instigation and the support the students would not have stepped out of the class room or the OP dept to carry out these types of strikes. Not only he instigated them but also gave them the essential protection to their leave, pay or brake in services. That is why the court (laws of Manu) gave them the support and under dealing work from the legal side. Can any section of the nation will pass as a special law and pay anyone who have fully abstained from work that too working in the very govt who have employed them and get pay for over a period of one month? Who is financing them? Whose tax money is it? They can behave and protest against the nation act as anti nationals, yet the nation should pay them with all benefit for over a month? What is this law? No law except laws of Manu alone can give this unjustifiable justice! So now the study makes clear that how the Brahmins can act when it is to their interest! The stand united to disintegrate the nation and ruin the majority. In this country even if a daily wager takes leave on account of his bad health or problems, he cannot get his pay; 'no work no pay' is the only rule in such a case how can the doctors get pay that too when they make protest against the govt under which they are employed for 40 day or so. What are the laws of Manu? Because Dr.Venugopal happens to be a Brahmin and the strike is also carried out by them (Brahmin) the court has given a new law? Did ever this Brahmin media make any criticism about the doctors paid in their period of absence from work, and they were agitating / protesting against the govt? Did the media analyse the move they were getting support from the courts for their striking period? Thus the media has failed miserably to be neutral even in this issue. The media should have written articles. When the media has been writing so many articles against reservation as if the move is a threat to nation where as in actuality it is the only means which can prevent nation from disintegration. But the media in reality knows the pros and cons of everything yet they are projecting this false propaganda only out of selfish personal interest. So the media in the mater of reservation for OBC are only giving wrong signals and false propaganda. This act of media will in a long run ruin the nation and leave the nation bereft of any real progress.

8. The media has failed in its duty for it has not investigated the real admission process in the P.G medical seats in AIIMS. On the contrary projected the stray protest with so many colourful photos.

9. All the non Brahmin experts uniformly agree on the basic fact that media has not acted neutral in the issues between Dr. Anbumani Ramadoss and Dr. Venugopal. Further the study confirms that the BJP and the Hindutva force have backed Dr. Venugopal. They were against the Sudra Dr. Anbumani. It is also important at this stage to mention that even the court has unnecessarily harassed Dr. Anbumani for it should have asked only the govt for clarification and not Dr. Anbumani who held the office of the President of the AIIMS on the basis, he was the Central Union Health Minister. If at all a Brahmin would have held this post and if he was a qualified doctor both the media and the court would have written pages and pages about his performance as they (the Brahmin media) have done now in case of Dr. Venugopal. For the only reason Dr. Anbumani is a Sudra and Dr.Venugopal a Brahmin who has more voice with them they are taking and projecting in this way. Here also the media has acted as a casteist, which fault should not be done by any good media with a balance mind and a true mission to contribute or work for the nation development.

10. The media has failed to be neutral in this issue for it has not given any thing about the seriousness of the death of the patients due to their strike. This is a very cruel neglect on the part of the striking doctors for they have not followed their very professional ethics so the media if it was not biased should have flayed the striking doctors, but on the other hand only their anti reservation protests were given undue propaganda their failure in their professional ethics were not even mentioned. When Dr. Anbumani urged the doctors to get their work. Dr. Anbumani was portrayed as an inconsiderate person. Media should have given enough coverage and the courts should have ordered for some compensation for the deaths of the patients for the death was only due to the strike and nothing more. Why was the death of patients never found properly mentioned and investigation is only due to three factors :

a. Patients are poor from non Brahmin castes.
b. No one is backing them like Hindutva force and media backing Dr.Venugopal.
c. If the patients death was projected, the striking mass will loose sympathy on the other hand would be against their irresponsible and anti social activities. So the media tried its level best to mention properly about the death. Here only if at all the media tried to educate the public of the real issue no one would ever support the anti reservation activities. So the media tried to keep ignorant the public of the true facts and falsely projected as if, the doctors were affected by reservations for OBC even when not even a set of OBC students have joined under reservations in these institutions. This only clearly shows the extra careful nature of the comfort loving and selfish Brahmins. They are clearly unconcerned about the millions of people who are deprived of good higher education due to the social problems.

11. Once again the media has failed to be neutral for it said two have been consumed by fire of reservations. The news was given as if both were medical students. But on the proper enquiry it was clearly revealed that one who had some minor burns was only a gutka seller and the other did not set himself ablaze only was admitted for … Thus why a gutka seller was set ablaze? Was it done on his own? Was it for some propaganda that reservation is protested by student?

No proper report and the real brains behind the activities of the gutka seller was given by any media! Is this media dharma or Brahmin media policy to hide all true facts related to reservation and project only those false and cooked up anti reservation protests! Like a gutka seller portrayed as a medical student setting himself against reservation and so on. When will the Brahmin media speak the truth? The nation will become developed only when the media gives unbiased information to the people on any social issue!

12. Finally another biggest bias the Brahmin media did was to totally hide or block out all the pro reservation protests and project only the anti reservation protests. Even the pro reservation protests carried under the banner of medicos' forum for equal opportunities in the AIIMS campus was totally blocked out by the media. None of the protests were ever projected or any importance or even a mere coverage was made. This is very unfortunate for the media failed to act justly in covering the issues of reservation, especially when the agitation or protests was organized in support of reservation. Only those anti reservation protests are covered with all great colours.

13. Another bias displayed by the media in the issue of Dr.Venugopal versus Dr.Anbumani is that when the President of AIIMS, Dr. Anbumani asked the people working in AIIMS to remove the shamiana put up by the anti reservation protestors no one removed it for this was the first visible, open disobedience of the subordinate Dr. Venugopal against Dr. Anbumani. By chance one of our non Brahmin Director had ever disobeyed the Brahmin president, the spot dismissal would have taken place. The media would have written pages about the arrogance of the Sudra subordinate. In this case as the subordinate happened to be Brahmin and the authority was a Sudra only reporting to show the powerlessness of Dr.Anbumani was clearly put out by the Brahmin media. No mention of insubordination of Dr. Venugopal was mentioned as if Dr. Venugopal was Dr.Anbumani's authority. Everything projected by the media showed Dr. Anbumani was ruining and bringing down the repute and the standard of the AIIMS. Several articles about AIIMS autonomy came! How can any institution fully financially supported by the Central Govt become autonomous under the only control of Dr.Venugopal. Only after Dr.Anbumani took the post of President as he was the Central Union Law Minister, he as a doctor; found from some of the non Brahmin colleagues of Dr.Venugopal, the bad functioning of the institution. Only after Dr.Venugopal became the director of the AIIMS several short coming in admissions in research etc. took place. As a shrewd doctor and an able administrator, Dr. Anbumani found them. This reformative actions of Dr. Anbumani was seen wrongly by Dr.Venugopal in a casteist angle and with the support of the Hindutva forces and with the help of media, he (Dr. Venugopal) found it an appropriate movement to carry out the protest under the pretext of anti reservation. That is why he supported all the anti reservations casteist protests held by them. Further when he took up the personal issue as an anti reservation protest the Hindutva backed him; not only BJP has openly supported him by asking Dr. Anbumani Ramadoss to resign of his post as the president of AIIMS. Court also indirectly supported him by using order in two days. Can any non Brahmin get such fast reinstatement or stay form his very dismissal as director that too by a Union Cabinet Minister. This is the mantra of caste! The media was in full swing to support him.

In all we observe that in the issue of AIIMS problems and the protests made by the doctors against reservation for the OBC and the tug of war between Dr. Anbumani, the AIIMS, President and his subordinate Dr.Venugopal of AIIMS, the media was never neutral. In fact the court was also in all the legal dealings supported only the Brahmins. So it was only laws of Manu prevailing at the court when it came to reservation in higher education for OBC

## 2.4 Use of Super Column FRM Model to study the Interrelation between the Government and Public

MHRD minister was very clear in his mission for when he tabled Quota bill no mention of creamy layer among OBC. He clearly mentioned the Bill was brought in order to benefit millions of students belonging to socially and economically weaker sections of society. The nodes given by the experts are as follows:

$T_1$ - The media was playing a very low level cheap politics when it made statements like Arjun Singh was implementing the quota for OBC to become better than the P.M Manmohan Singh. The media also said in a sneaky way that to bring down P.M in the eyes of Sonia, the Congress chief, he was bringing quota. In truth quota for OBC in the institutions of higher learning, run by the central govt was a social reformation aimed to develop the nation for thousands of years only one caste have developed in leaps and bounds, educationally and economically the question of their social status need not be discussed for even as they came as beggars yet by the laws of Manu, they claimed them selves to be the highest caste twice born; born out of the mouth (or head of Brahma) (Pl. do not ask for any scientific explanation of these for they are found in the laws of Manu [43]). So their social conditions always remained superior to the very sons of the soils. Only now the govt has thought over the majority, to uplift their social status it has suggested a reservation for them that too a paltry 27 percent in higher learning. This cannot be accepted by the nomads and are protesting it by deviating it from all reality using the media, by speaking and portraying Arjun Singh who made the suggestions to be a person acting for his personal benefits.

$T_2$ - For social upliftment or development, good higher education is the only way. So a rationalistic person

working for the nation development would first uplift socially the majority so the govt. move can never be called as personal profit but as one for the national profit (cause)

$T_3$ - The media projected Arjun Singh's acts was questionable by E.C. So he had to clarify before E.C, though he had actually not created any flaws. This propaganda are made only to dilute the true issues regarding reservations for OBC.

$T_4$ - Several lengthy articles came in the editorial and print media as news that Arjun Singh's implementation of quota would ruin congress …

$T_5$ - The only way to make the majority whose resources are misused and exploited by the minority Hindu Brahmins by exploiting their labour very badly can be changed and a better social and economic conditions can be developed by using their resources properly was only by giving them good higher education. So that too after 60 years of independence and thousands of years of exploiting them by the Brahmins who sat as policy makers or advisors to the king, only now the govt. had the mind to give reservations for this mass only in 2006.

$T_6$ - The govt. only suggested the reservation for OBC in Higher Learning in educational institutions that too run by the central govt. Only this was passed by both the houses with a over whelming majority. That is why Arjun has sharply asked "Why do we have parliament if we cannot honour our commitment on quota? "Media if it had been little responsible it should have first made a nice propaganda for reservation in support of reservation for the very amendment and the act that 27 percent reservation for OBC in higher learning in educational institutions was passed in over whelming majority, hence everyone must respect the people's right and hence the policy cannot be flouted or even debated; but the Brahmin dominated selfish media by all its acts has projected reservation as an anti social act but in reality it is socially a reformative act for the actual progress of the nation. The media is miserably acted against people and above all against public and nation.

$T_7$ - The print Brahmin dominated media has written several essays decrying reservation as an act of nation's development and merit will suffer. This act of Brahmin media is very questionable and what is the merit of this media to publish articles which is anti national and anti people.

$T_8$ - The medicos and the IIT-ians made protests only against Arjun Singh! This was construed by experts as only a drama! Arjun Singh likes to make this policy without the support of both the houses he can do nothing! So if both the houses had not recommended that too in this case with a over whelming majority nothing would have happened. So the media joining hands with the anti reservation protestors ridiculing the political leaders by names and burning their effigies is wrong. Above all why media too has lost its sanity by writing articles against reservation. That is why it is not wrong if we say that the media has become blinded with selfishness and in haste has lost its reasoning. The media should have educated the protesting medicos and the IIT-ians that it is not Arjun Singh alone we have got reservation but it is govt and above all the MPs from both the houses; but with a school boy mentality media has joined hands with the students in decrying the HRD minister.

On the other hand in the pro reservation protests of course thanked the UPA govt.

$T_9$ - Media should have talked ill of and ridiculed each and everyone in both the houses for favouring reservation.

This act only makes one strongly suspect BJP Hindutva and Brahmin MPs who were in mind against reservations for OBC, for the fear of their position in public supported reservation for OBC. What does this imply if they do not now give reservation they would be shown out in the coming elections and nothing. What does this show the majority of the population are OBC who decide who should rule and nothing more!

$T_{10}$ - The cunning act of media was, it was giving all sorts of wrong feed back to the public that centre was looking for exit route on reservation row and Manmohan Singh is against it and so on. But how can then the 104th Constitution Amendment Bill was passed in 379 votes in Lok Sabha on 21-12-05 and 172 votes in favour of it in Rajya Sabha only one against in both the houses and one abstaining from the Lok Sabha i.e., only 3 did not vote in favour of it of the total of 554, 551 (379+172) favoured it! Then why create such wrong waves about reservations that too by a very few upper castes and Brahmin medicos.

$T_{11}$ - This reflect that pro reservation protestors are naturally much more intelligent than the anti reservation protestors. They understand and grasp the whole situation unlike the hasty selfish Brahmins. They (OBC) pro reservation protestors are not backed by anyone and they are socially and above all economically from a very poor strata yet much more intelligent than the anti reservation protestors.

$T_{12}$ - Media is bias to project the govt, as that one that has committed a horrendous mistake by implementing reservations!

$T_{13}$ - Media has projected Manmohan Singh as the anti reservation person there by trying to bring a split with in the congress. Thus media is trying to divide the congress.

We give the views of the public, political leaders and educationalists taking them as three distinct groups, P, L and E where P is the group of public. L is the set of political leaders and E the category of educationalists.

The following are the experts forming the public group P

$P_1$ - OBC people
$P_2$ - SC/ST people
$P_3$ - Minorities
$P_4$ - Brahmins

L – Experts from Political Leaders

$L_1$ - BJP/Pro Hindutva / Pro Brahmin / Political Leaders
$L_2$ - Dravidian Political Leaders like MDMK excluding AIDMK
$L_3$ - Minority Party Leaders
$L_4$ - Casteless Political Leaders
$L_5$ - Non Congress, Non communist, Non BJP parties like PMK, Telugu Desam etc.

We list the sets experts from educationalists E

$E_1$ - OBC
$E_2$ - Brahmins
$E_3$ - SC/ST
$E_4$ - Minorities
$E_5$ - OBC, NRIs.

Now we make use of the super fuzzy column FRM model to analyze the problem.

In this super fuzzy column FRM model P, the attributes related with the domain space are usual or simple row vectors given by ($T_1$ $T_2$ … $T_{13}$) the Range space of this model is a mixed fuzzy super row vector given by ($P_1$ $P_2$ $P_3$ $P_4$ | $L_1$ $L_2$ … $L_5$ | $E_1$ $E_2$ … $E_5$).

The super row matrix associated with the super fuzzy column FRM model given by the expert is as follows.

| | $P_1$ | $P_2$ | $P_3$ | $P_4$ | $L_1$ | $L_2$ | $L_3$ | $L_4$ | $L_5$ | $E_1$ | $E_2$ | $E_3$ | $E_4$ | $E_5$ |
|---|---|---|---|---|---|---|---|---|---|---|---|---|---|---|
| $T_1$ | 1 | 1 | 1 | 0 | 0 | 1 | 1 | 1 | 1 | 1 | 0 | 1 | 1 | 1 |
| $T_2$ | 1 | 1 | 0 | 0 | 0 | 1 | 1 | 0 | 1 | 1 | 0 | 1 | 1 | 1 |
| $T_3$ | 1 | 1 | 1 | 1 | 0 | 1 | 1 | 1 | 1 | 1 | 0 | 1 | 1 | 1 |
| $T_4$ | 1 | 1 | 1 | 0 | 1 | 1 | 1 | 1 | 1 | 1 | -1 | 1 | 1 | 1 |
| $T_5$ | 1 | 1 | 1 | 0 | 0 | 1 | 1 | 1 | 1 | 1 | 0 | 1 | 1 | 1 |
| $T_6$ | 1 | 1 | 1 | 0 | -1 | 1 | 1 | 1 | 1 | 1 | 0 | 1 | 1 | 1 |
| $T_7$ | 1 | 1 | 1 | 0 | 0 | 1 | 1 | 1 | 1 | 1 | 0 | 1 | 1 | 1 |
| $T_8$ | 1 | 1 | 1 | 0 | 0 | 1 | 1 | 1 | 1 | 1 | 1 | 1 | 1 | 1 |
| $T_9$ | 1 | 1 | 1 | 0 | 1 | 1 | 1 | 1 | 1 | 1 | 0 | 1 | 1 | 1 |
| $T_{10}$ | 1 | 1 | 1 | 1 | 0 | 1 | 1 | 1 | 1 | 1 | 0 | 1 | 1 | 1 |
| $T_{11}$ | 1 | 1 | 1 | 0 | 1 | 1 | 1 | 1 | 1 | 1 | -1 | 1 | 1 | 1 |
| $T_{12}$ | 1 | 1 | 1 | 0 | 0 | 1 | 1 | 1 | 1 | 1 | 0 | 1 | 1 | 1 |
| $T_{13}$ | 1 | 1 | 1 | 1 | 0 | 1 | 1 | 1 | 1 | 1 | 0 | 1 | 1 | 1 |

Now the experts want to study the effect of the on state of the node $T_7$ alone and all other nodes are in the off state; i.e.,

$$X = (0\ 0\ 0\ 0\ 0\ 0\ 1\ 0\ 0\ 0\ 0\ 0\ 0)$$

belongs to the domain space; to study the effect of X on P.

$$X \circ P \hookrightarrow (1\ 1\ 1\ 0 \mid 0\ 1\ 1\ 1\ 1 \mid 1\ 0\ 1\ 1\ 1) = M_1 \text{ (say)}$$

($\hookrightarrow$ denotes the resultant that has been thresholded as, if $a_i > 6$ put 1 $a_i < 6$ put 0 in the resultant.). $T_1$ is a fuzzy super row mixed vector. Now we study the effort of $T_1$ on P.

$$M_1 \circ P^T \hookrightarrow (1\ 1\ 1\ 1\ 1\ 1\ 1\ 1\ 1\ 1\ 1\ 1\ 1) = X_1 \text{ (say)}$$

i.e., all the nodes come to on state in the domain space

$$X_1 \circ P \hookrightarrow (1\ 1\ 1\ 0 \mid 0\ 1\ 1\ 1\ 1 \mid 1\ 0\ 1\ 1\ 1)$$
$$= T_2\ (=T_1);$$

which is again a super fuzzy mixed row vector. Thus we arrive at the hidden pattern which is a fixed super vector given by the binary pair {(1 1 1 1 1 1 1 1 1 1 1 1 1), (1 1 1 0 | 0 1 1 1 1 | 1 0 1 1 1)}. We see when only the node $T_7$ i.e., the Brahmin dominated media has written several essays decrying reservation as an act that would ruin merit and the progress of the nation all the nodes in the domain space come to on state showing very strongly that the media has been very biased in projecting the reservation as an act of social discrimination. In fact in the fuzzy super vector from the range space we see only the nodes $P_4$, $L_1$ and $E_2$ did not come to on state and all other nodes in the range space also came to on state their by making one know clearly that Brahmin public, the Brahmin political leaders and the Brahmin educated alone do not agree on the factors showing that the Brahmin dominated media has been acting continuously against reservation for OBC in educational institutions of higher learning run by the central govt.

Further this model show mathematically in general and fuzzy theoretically in particular that he Brahmins, be it they belong to educated class or just as a public and more so as a Brahmin party (BJP RSS …) they are dead against the nation giving reservations for OBC. As suggested by the experts if we take the node $L_1$ alone in the on state in the range space i.e. BJP/pro Hindutva pro. Brahmin political leaders and all other nodes are in the off state. The effect of

$$M = (0\ 0\ 0\ 0 \mid 1\ 0\ 0\ 0\ 0 \mid 0\ 0\ 0\ 0\ 0),$$

on the fuzzy super mixed row vector on P is given by

$$M_1 \circ P^T \hookrightarrow (0\ 0\ 0\ 0\ 0\ 0\ 0\ 0\ 0\ 0\ 0\ 0\ 0)$$

which is the zero vector which says the BJP/RSS/VHP or any of the Brahmin dominated parties do not agree upon any form of bias on the part of the Brahmin dominated media. They say only we are viewing it wrongly. But the article carried out in the print media clearly shows the venom in the vision. We have worked out with several other state vectors in the on state and several other experts and arrived at the conclusions and observations. Observations obtained from the super column FRM model.

1. This fuzzy super column FRM model has acted in a very sensitive way for we see when certain nodes from the range space are in the on state in the resultant binary pair one gets the domain space to be just the zero row vector i.e. none of the nodes in the domain space come to on state which has become as one of the unique features of this model.
2. In this model we see some of the state of the certain nodes in the domain makes all the Brahmin media related nodes of the fuzzy super vectors of range space to remain in the off state, which is also another special feature of this particular model.
3. It is interesting to see that when any one of the nodes $T_i$ in the domain space vector is in the on state all the other nodes of the domain space come to on state signifying that all the nodes are interrelated in a very strong way and the media has been acted bias in case of reservation for OBC and so also against the govt.
4. The media has falsely trying to bring a split with in the Congress. This sort of false propaganda is unnecessary for the media which is condemned by the OBC educationalist, public and political parties. The same view is also held by the minorities and SC/ST belonging to public, educationalists and political parties.
5. The false propaganda of the Brahmin media tarnishing the image of govt. that it was acting against the progress and the development of the nation was also well brought out by this super fuzzy column FRM model.
6. The media has failed to express as well failed to understand that even if Arjun Singh wishes to bring reservation for OBC in educational institutions of higher learning or wishes to make any amendment he cannot do anything. It is only the full support of the both the houses i.e. out of 554 strength 551 favoured the $104^{th}$ constitution amendment. It is ignorance and prejudice of the Brahmin media to even understand this simple situation. As the immature medicos under the strong instigation of the Director of AIIMS Dr.Venugopal with no proper knowledge of anything carried placards against the govt., the Brahmin media instead of correcting them with school boys maturity joined hands with them by making propaganda against the govt. Unfortunately the Brahmin print media has wasted its precious time and economy by penning several articles against reservations. Only by giving education to these people alone i.e., this social reformation alone can lead to the true progress of the nation more so of the majority who happen to be the son of the soil.

7. It is unwarranted on the part of the media to publish the burning of effigies of political leaders in the medicos protest or calling him by names like ‘monster’ etc for if the media has courage, let it write against all the 541 political leaders who said yes to the amendment. It is very indecent on the part of the Brahmin dominated political parties to support the reservation in houses but instigate people or build new banners like “Equality for Youth” and so on and make those persons to protest under this banner. This is an act of shame on the part of these Hindu fanatic politicians to act on double standard. They not only give these protesting medicos money but legal support for they are not terminated from their courses of study they are paid for the period of over 40 days of strike. Above all they made protests only under the shade of the shamiana with air coolers in the very campus. What a drama? This was the drama staged by the shameless Hindutva politicians said hundreds of experts who have given views about the anti reservation activities in the AIIMS campus which was instigated by Dr. Venugopal.
8. If the non Brahmin medicos had staged this type of protests even outside the campus in scorching sun even with no water to drink, against a Brahmin Union Cabinet Central Minister, then what would have been the plight of these non Brahmin medicos they would have been dismissed from their courses of study and sent home or to prison. Their future would have been doomed. Will the PM or the President come to talk with them? This is the greatest discrimination the non Brahmins face in the modernized century. So only to wipe of this social discrimination they need reservation in the institutions of “Higher Learning”, demanded the experts. The experts asked who would even back these non Brahmin protestors who form the majority so only they need reservation for social equality!

9. Because the govt. happened support the non Brahmin and secondly he has taken up an issue which supports the non Brahmin, the govt. was very cheaply treated and projected by the Brahmin media. The experts very strongly protested to the cheap, devious acts of the media.
10. One cannot forget that from day one i.e., 21-12-05, when the bill was passed till date 12-03-07 he had been single minded very serious about his getting reservation for OBC in higher learning for that alone can economically lift the majority of the OBC and save them from social discrimination. Arjun Singh even after the SC stay was very firm that in a day or two (12-04-07) the govt should remove the stay on reservations. On March 30-02-2007, Arjun said that govt will try to vacate the stay by SC.

    Arjun Singh dismissed as a “very subjective view” the charge that the quota issue has divided the nation. The HRD minister rejected as a “canard” that the Congress was not enthusiastic in pursuing the issue. Here also the media always tries to split the party or bring confusion within the party. The very statement that the quota issue was dividing the nation was a very false propaganda, the experts claimed already we are divided based on the caste. The first broad division is Brahmins and non Brahmins, will they accept caste is dividing the nation, so from today on wards no mention of caste will be made and no symbols of castes be used!

    If quota is dividing the nation, caste has already divided the nation! What steps are they going to take? Why is not the media writing even a single article in this direction they asked? All their acts are selfish and in no place they had worked for the development of the nation or for the majority’s good; they only strive for their over pampered caste, comfort and live at the cost of others labour.
11. Arjun had said that the congress would ensure the validity of the law. The party’s national secretary D.Raja said, “It is a retrograde judgement … Parliament has to intervene effectively. Time has come for Parliament to have an effective say on this. It cannot be left to the judiciary". The left front was also critical of the government handling of the issue. CPM General secretary Prakash Karat said, "We want the central government to take necessary steps in parliament to ensure the implementation of the law".

Reacting to the judgement CPI general secretary AB Bandhan hoped that the govt. would take effective steps in the wake of the judgement. Manmohan Singh has clearly said that the courts should not cross their limits in the conference attended by the Chief Ministers of the states and the judges of the courts. It is pertinent to take a note that Manmohan Singh had made this comment in the presence of the Chief Justice of the Supreme Court.

He clearly said that the parliament, legal department and the administration of the nation are three pillars and the legal dept has acted in a way to separate the thin line between the three pillars. The law dept must be very careful to see how far it can use its power over the other department. The legal dept crossing its limits only shows that there is difference between these pillars. Taking the rights of other department is a clear cut evidence of exceeding its limits…

He further said if there is not a perfect understanding between the three pillars, would only hinder the services to the public and the nation. It should respect the action and rights of the other dept then only a good relationship can be created. Each and every dept has their rights and responsibilities. Such responsibilities should be carried out unitedly. As the nation has over 2.5 crores cases pending in the Honourable courts they should do it fast.

Nowadays to protect their property using their political power they go to courts. This has become a very bad custom or tradition. It is very essential before admitting such cases the courts should analyse the cases fully. It is important to mention here that the P.M has talked daringly in this way. The P.M. deserves appreciation for his bold comments about the courts.

Here also the media which was only a Brahmins media voicing only their comforts and never bothered about any common mans and majority social welfare has very clearly now shown its colour. Even 'The Hindu' which had been in several issues of reservation was acting not that biased has shown its colour. The Brahmins had been in the ocean of happiness when the S.C put a stay on the 27 percent reservations to OBC in the higher education in the educational institutions run by the central government. They greeted each other by distributing sweets and celebrated the occasion when the majority of the nation was in despair and anger. The Hindu, the Brahmins, the casteist daily satisfied its expression by writing headlines like in meritorious institutions like IITs, IIMs reservations cannot enter.

The minister for HRD has very clearly said that till the next judgement comes these institutions should stop admitting any students. He further said this order is applicable to medical agricultural and engineering institutions. .

However Veerappa Moiley says that the S.C in its interim stay has not given any stay over the 27 percent reservation for OBC. In this situation these institutions of higher learning; have been very hasty in showing the 'stay by the S.C.' as a cause, has planned to admit students without reservations. Manmohan Singh's very bold statement, that the courts should not cross their limits was very much supported by the experts.

Now the experts wanted to record the following:

a. The media was projecting as if Arjun Singh and Manmohan Singh had dispute and misunderstanding over the 27 percent reservation for the OBC in the educational institutions of higher learning run by the central govt. Now it has been proved beyond doubt that there was no misunderstanding between them on this issue, but Arjun Singh was trying to get better of Manmohan Singh, the reasons best known to him. Thus how the Brahmin media is playing such nasty of bringing the reservation issue in between to mar the govt. commitment to OBC.
b. Once again the true face of the Brahmin media is very clear when the Hindu blabbered out that reservation for OBC cannot enter the institutions of merit on meritorious ground like IITs and IIMs. All the while the Hindu was trying to behave neutral only when the Supreme Court ordered a stay on the 27 percent OBC reservation, the Brahmin daily Hindu lost its false cloak and unknowingly said reservations for OBC cannot enter the meritorious institutions like IITs and IIMs as it will make the nation loose the merit. But when HRD Minister took the strong step of withholding the admissions until further order, the Brahmin media has continuously behaved against reservations by projecting reservations as a socially discriminative move.

c. All Brahmin dominated media projected the celebrations of the stay by S.C with such big photos.

d. On the same day the IIMs and IITs said they would go by the directive of the government on providing 27 percent reservation for other backward classes. If the media was unbiased it would have also analysed the pros and cons of the SC stay.

e. Once again the experts pointed out that no one other than the youth for equality, can challenge this issue, for they are the complete brain child of the hindutva parties like BJP/RSS /VHP and above all the chiefs of the Brahmin religious mutts like Kanchi Sankaracharya mutt, Swami Sivananda etc. So the youth for equality is itself a very strong Hindutva backed set up only because these set of Brahmins were behind this in the first place the petition got admitted and in the second place the verdict was in their favour. The experts argued will the laws of Manu ever give any ruling affecting the interests and comforts of the Brahmins. Even in remote of dreams such things will not occur. Or can one ever imagine that a non Brahmins youth can ever in any court get a stay on any issue passed by the parliament which is in favour of Brahmins. If this is the case, how can we ever say the courts does justice to non Brahmins bombarded many non Brahmins who were layers themselves. In fact Dr. Veeramani has highly condemned the indisciplinary in the Judiciary. Dr. K. Veeramani himself was a reigning advocate! These views he has published in the daily 'Viduthalai' for which he is the editor in chief of that daily. All these cannot take place with the situation reversed a group of

> non Brahmin students filing case against their Brahmin ministers cannot even take place if granted such a situation arises the case will be rejected at the very initial stage. Even by chance the petition gets admitted it is for certain that the order will be never in favour and the case would be dismissed. Thus now one can understand how correctly said by our Prime Minister Manmohan Singh that parliament, judiciary and administration are the three pillars always in India, the pillar of judiciary is against the other two when any positive socially reformative move to up lift the majority who are socially discriminated by caste is aimed at.

This has been in practice for the past 60 years. Dr. Veeramani, the DK Leader who presided over a seminar went a step ahead and accused the two judges for lacking in judicial honesty and discipline.

The judges had failed to mention the name of Dravidar Kazhagam, the only organization which got itself impleaded in the case and submitted arguments in the 36 page interim order. He said this only showed the judicial dishonesty of the judges. It is still surprising to see that a nine judge bench of the Supreme Court has passed as orders on the basis of the recommendations of the Mandal Commission report in 1994. However the two judge bench had reversed the order. This is judicial indiscipline, Veeramani said and added “Let us be prepared to enter the (battle) field”.

Senior advocate and professor Ravi Varma Kumar also former chairman of the Backward Classes Board of Karnataka said that the nine reasons given by the judges for the interim stay were absurd and perverse. The derogatory language used by the judges in their judgement was a big insult and should be expunged, he added.

CPI state unit general secretary D.Pandian claimed that the parliament was the supreme body in a democracy. The Supreme Court had no right to interfere in its policy decisions.

Viduthalai Chiruthargal Katchi general secretary Thol Thirumavalavan assured that all necessary assistance from the Dalit's to ensure that social justice right of the OBC. Tomorrow it might fall on the SC/ST.

It is an unfortunate occurrence that within a year the Brahmin media had to face their shameful and disgraceful way how their conjecture was proved to be a blatant lie for instance a few of them like; Manmohan Singh did not support quota. Arjun Singh was not getting along well with Manmohan Singh; there was split in congress over the reservation issues. Thus the experts wished that the Brahmin media will at least from now on wards correct its mode of activities while reporting. They advised it to be neutral and never selfish and act as a mouth piece of Brahmins and Hindutva forces.

## 2.5 Analysis of Role of Media on Reservation for OBC using Fuzzy Cognitive Map

Now we use the fuzzy cognitive maps FCM model to study the 27 percent OBC reservation in central institutions of higher learning. The main features like merit, quality, social reformation caste etc. how media has failed to analyse will be seen are dealt with

$R_1$ - Reservation for OBC in institutions of higher learning is the need of the hour to have stability of the nation.

$R_2$ - Reservation for OBC is as social right of a majority and nation responsibility.

$R_3$ - Views of educated non Brahmins on reservation.

$R_4$ - Reservation will ruin the standard of these institutions.

$R_5$ - Reservations will divide the nation.

$R_6$ - Reservation for OBC alone can keep the nation integrated.

$R_7$ - Media has given a false propaganda (as social evil) to the world and public about reservation.

$R_8$ - Media has miserably failed to be neutral on the topic of reservation for OBC it has shown the caste fanaticism.

$R_9$ - Ban on caste would give no scope for reservation. So the need of the hour is banish caste or give reservation.

$R_{10}$ - Caste is a socially discriminative negative factor in India based on birth; so reservation is needed to annual social discrimination which is a positive action.

$R_{11}$ - Brahmin view reservation as a negative force which has lead to media discriminative action.

$R_{12}$ - Without reservations the nation would become a place of terrorism for majority of the OBC are very poor living in rural areas and also socially discriminated

$R_{13}$ - Education alone can lift the social barrier of social discrimination based on caste.

$R_{14}$ - Economy alone can by no means lift the social barrier of social discrimination.

$R_{15}$ - The development of the nation is based only on the development of the majority that too sons of the soil. All other talks are only superficial!

$R_{16}$ - The economy and the major wealth of the nation is in the hands of the Brahmins who are nomads, who came to India for green pastures as beggars and now have usurped the total wealth and economy of the nation in the name of Hindu religion. Why the religious centres and temples under the Brahmin control have amassed the wealth of India and they passed on laws and orders in the BJP time itself that the wealth of these Hindu (Brahmin) controlled mutts cannot be taxed or they do not pay income tax. What a special law for the special people.

$R_{17}$ - School education for the OBC in these 60 years of India’s independence has not developed the rural OBC, socially. If all of them migrate to the city, majority of them only work as daily wagers doing the menial jobs of masons, sweepers etc. So only reservation in higher education alone can socially uplift them who are so much ruined by the superstitions of the Hindu caste and religion which has been forcefully practised for thousands of years by the Brahmins only for of their personal income

$R_{18}$ - Views of SC/ST on OBC reservation and its relation to other attributes.

Now using these 18 attributes related to the OBC reservations given by the majority of the experts, we give the opinion of the experts using the tool of Fuzzy Cognitive Maps (FCM). We have been suggested by the experts to use this model for the primary reason this model along is easily understandable by the majority of the experts, secondly it gives the hidden pattern of the system, which is very essential for socio scientists to make predictions.

The connection fuzzy matrix V is given below.

| | $R_1$ | $R_2$ | $R_3$ | $R_4$ | $R_5$ | $R_6$ | $R_7$ | $R_8$ | $R_9$ | $R_{10}$ | $R_{11}$ | $R_{12}$ | $R_{13}$ | $R_{14}$ | $R_{15}$ | $R_{16}$ | $R_{17}$ | $R_{18}$ |
|---|---|---|---|---|---|---|---|---|---|---|---|---|---|---|---|---|---|---|
| $R_1$ | 0 | 1 | 1 | 0 | 0 | 1 | 1 | 1 | 1 | 1 | 0 | 1 | 1 | 0 | 1 | 1 | 1 | 1 |
| $R_2$ | 1 | 0 | 1 | 0 | 0 | 1 | 1 | 1 | 1 | 1 | 0 | 1 | 1 | 0 | 1 | 1 | 1 | 1 |
| $R_3$ | 1 | 1 | 0 | 0 | 0 | 1 | 1 | 1 | 1 | 1 | 0 | 1 | 1 | 0 | 1 | 1 | 1 | 1 |
| $R_4$ | 0 | 0 | 0 | 0 | 0 | 0 | 0 | 0 | 0 | 0 | 1 | 0 | 0 | 0 | 0 | 0 | 0 | 0 |
| $R_5$ | 0 | 0 | 0 | 0 | 0 | 0 | 0 | 0 | 0 | 0 | 1 | 0 | 0 | 0 | 0 | 0 | 0 | 0 |
| $R_6$ | 1 | 1 | 1 | 0 | -1 | 0 | 1 | 1 | 1 | 1 | 0 | 1 | 1 | -1 | 1 | 1 | 1 | 1 |
| $R_7$ | 0 | 1 | 0 | 1 | 1 | 0 | 0 | 1 | 0 | 0 | 1 | 0 | 0 | 1 | 0 | 0 | 0 | 0 |
| $R_8$ | 0 | 0 | 1 | 1 | 1 | 0 | 1 | 0 | 0 | 0 | 1 | 0 | 0 | 1 | 0 | 0 | 0 | 0 |
| $R_9$ | 1 | 1 | 1 | 0 | 0 | 1 | 1 | 1 | 0 | 1 | 0 | 1 | 1 | 0 | 1 | 1 | 1 | 1 |
| $R_{10}$ | 1 | 1 | 1 | 0 | 0 | 1 | 1 | 1 | 1 | 0 | 0 | 1 | 1 | 0 | 1 | 1 | 1 | 1 |
| $R_{11}$ | 0 | 0 | 0 | 1 | 1 | 0 | 0 | 0 | 0 | 0 | 0 | 0 | 0 | 1 | 0 | 0 | 0 | 0 |
| $R_{12}$ | 1 | 1 | 1 | 0 | 0 | 1 | 1 | 1 | 1 | 1 | 0 | 1 | 1 | 0 | 1 | 1 | 1 | 1 |
| $R_{13}$ | 1 | 1 | 1 | 0 | 0 | 1 | 1 | 1 | 1 | 1 | 0 | 1 | 1 | -1 | 1 | 1 | 1 | 1 |
| $R_{14}$ | 0 | 0 | 0 | 1 | 1 | 0 | 1 | 1 | 0 | 0 | 1 | 0 | 0 | 1 | 0 | 0 | 0 | 0 |
| $R_{15}$ | 1 | 1 | 1 | 0 | 0 | 1 | 1 | 1 | 1 | 1 | 0 | 1 | 1 | 0 | 0 | 1 | 1 | 1 |
| $R_{16}$ | 1 | 1 | 1 | 0 | 0 | 1 | 1 | 1 | 1 | 1 | 0 | 1 | 1 | 0 | 1 | 0 | 1 | 1 |
| $R_{17}$ | 1 | 1 | 1 | 0 | 0 | 1 | 1 | 1 | 1 | 1 | 0 | 1 | 1 | 0 | 1 | 1 | 0 | 1 |
| $R_{18}$ | 1 | 1 | 1 | 0 | 0 | 1 | 1 | 1 | 1 | 1 | 0 | 1 | 1 | 0 | 1 | 1 | 1 | 0 |

Now the experts supply certain state vector the hidden pattern of them we find using the dynamical system V. Let

$$X = (0\ 0\ 0\ 0\ 0\ 0\ 0\ 0\ 1\ 0\ 0\ 0\ 0\ 0\ 0\ 0\ 0\ 0)$$

be the given state vector where only the node $R_9$ i.e. Ban on caste gives no scope for reservation so the need of the hour is ban caste or give reservation for OBC is in the on state and all other nodes are in the off state to find the effect of X and V i.e. to obtain the hidden pattern of X using the dynamical system V

$$X \cdot V \hookrightarrow (1\ 1\ 1\ 0\ 0\ 1\ 1\ 1\ 1\ 1\ 0\ 1\ 1\ 0\ 1\ 1\ 1\ 1) = X_1 \text{ (say)}$$

($a_i > 9$ put 1 if $a_i < 9$ put 0 as a thresholding function) ('$\hookrightarrow$' as usual denotes the resultant is updated and thresholded).

$$X_1 \cdot V \hookrightarrow (1\ 1\ 1\ 0\ 0\ 1\ 1\ 1\ 1\ 1\ 0\ 1\ 1\ 0\ 1\ 1\ 1\ 1) = X_2 \text{ (say) } (\therefore X_2 = X_1).$$

Thus the hidden pattern is a fixed point and only the nodes $R_4$, $R_5$, $R_{11}$ and $R_{14}$ are in the off state, i.e. Reservation will ruin the standard of these institutions have no effect on the system and Reservation will divide the nation is unaffected by the on state of $R_9$. Further the Brahmins are of the view that Ban on caste itself is an impossibility (view got from discussions) how come the other alternative and above all $R_{14}$ is for clearly reservation in higher education has nothing to do with the economy of the people will come to on state.

Now suppose the experts want to work with the node of $R_{15}$ alone and all other nodes remain in the off state. To find the hidden pattern of the state vector

$$Y = (0\ 0\ 0\ 0\ 0\ 0\ 0\ 0\ 0\ 0\ 0\ 0\ 0\ 0\ 1\ 0\ 0\ 0)$$

on the dynamical system V.

$$Y \cdot V \hookrightarrow (1\ 1\ 1\ 0\ 0\ 1\ 1\ 1\ 1\ 1\ 0\ 1\ 1\ 0\ 0\ 1\ 1\ 1) = Y_1 \text{ (say)},$$

$$Y_1 \cdot V \hookrightarrow Y_2 = (= Y_1).$$

Thus the resultant of Y is a fixed point and only the nodes $R_4$, $R_5$, $R_{11}$ and $R_{14}$ are in the off state and all other nodes come to on state thus the development of the nation means and is only the development of the majority, sons of the soil nothing more is intricately connected with all other nodes except the view reservation will ruin these institutions standard and divide the nation, the views of Brahmins who believe that development depends on the economy forgetting the very fact that economy alone would not make the nation a developed one, only by giving them higher education alone, one can develop the nation.

Now we proceed on to work with the on state of the node $R_{11}$ alone and all other nodes of the state vector remain in the off state. Let Z = (0 0 0 0 0 0 0 0 0 0 1 0 0 0 0 0 0 0) be the given state vector to find the effect of Z on V. The hidden pattern of Z using the dynamical system V is given by

$$Z \cdot V \hookrightarrow (0\ 0\ 0\ 1\ 1\ 0\ 0\ 0\ 0\ 0\ 1\ 0\ 0\ 1\ 0\ 0\ 0\ 0) = Z_1 \text{ (say)},$$

$$Z_1 \cdot V \hookrightarrow (0\ 0\ 0\ 1\ 1\ 0\ 0\ 0\ 0\ 0\ 1\ 0\ 0\ 1\ 0\ 0\ 0\ 0) = Z_2\ (=Z_1).$$

Thus the hidden pattern is a fixed point and all nodes in favour of reservation remains in the off state and only the nodes $R_4$, $R_5$ and $R_{14}$ come to on state there by indicating according to Brahmins reservation will ruin the standard of these institutions of higher learning and the implementation of reservation will divide the nation. They also view economy plays a vital role in giving reservations for they say reservation based on economy can be given.

Now we have worked with several other state vectors in the on state using the C-program given in [327, 346]

## 2.6 Observations based on this Analysis by Students and Experts through Seminars and Discussions

I. The analysis made it clear that the media was acting in a very biased way against the very welfare of the nation and propagating all false and selfish views about reservations like the "merit" of these institutions are at stake if one implements reservation. Further the reservation should be based on economy. The experts bombarded the Brahmins by asking how much percentage of the Brahmins are living below poverty line if we take the OBC majority of them living in the rural areas are living very much below the poverty line. All Brahmins live in luxury and pomp, which a rich man in the village cannot obtain. Not only that the majority of the nations wealth is in the hands of the Brahmins religious trusts and ex IIT-ians of India who are Brahmins. If this is the case the readers are requested to ponder over the fact how reservation should be based on economy. Only fact derived from their claim is that if it is said reservation based on economy they will give false pay certificates or income certificates and get all the seats in the higher education run by the central govt institution like IIT and IIM and AIIMS.

II. Even as of today if anyone takes the statistics of the Brahmins studying in these institutions they would very easily cross 50 percent of the seats. Hence the govt is requested instead of taking the census of OBC, let it take the percentage of Brahmins studying in these institutions for the past 10 or 15 years and based on this statistics and based on the percentage from them serving in western countries now give the OBC at least 50 percent reservations say just for 10 years so that they who have been deprived of their rights can at least be compensated. For justice delayed is justice denied!

III. If the nation is facing any form of economic crisis let it now demand the money back which was spent on educating these upper caste students in these institutions of higher learning from all those who gone abroad either for study or for earning. As over 20% of the nations wealth is in their hands India (nation) can run these institutions in a better way for OBC.

IV. Let the govt pass on a very strong rule that whoever goes abroad after obtaining the degrees from the central govt institutions of higher learning must pay back the money which the govt has spent on them to give that particular degree. If such a rule is put none of the Brahmins would come to study in these institutions. It is to be well understood that they want the degree from these institutions not for real knowledge but mainly as a good stamp or a passport to leave India. This should be curtailed for the tax money of the honest OBC should not be spent on these nomad Brahmins who go to other countries to serve them for money. What is the benefit for the nation for educating them except a permanent economic loss to nation. If by chance this practice was among the OBC/SC/ST, the jealously and the selfish Brahmins would not have remained spectators like us but have made such strong rules that none of the OBC/SC/ST leave the nation after getting their degrees from these institutions. So now the experts suggest that those who leave abroad should pay back the amount the govt. has spent on them giving their degree.

V. The experts strongly suggest that only reservation in the higher education alone can give social equality to the OBC and also wipe out the social discrimination of sudrahood which imply as per Manu, a sudha means born to a prostitute. The only alternative is ban all caste and caste system or give reservations as per their strength. Further making economically better the status of sudras would by no means wipe out their sudrahood only educating them will make them free of their sudrahood. Here once again it has become pertinent to mention that it is clearly said in the laws of Manu that sudra should not learn even they should not hear to the oral recitation and it is further said in the Manu that if they try to listen i.e., hear or learn they will be punished by pouring hot liquid lead into their ears. If that is the rule in laws of Manu how will Brahmins ever agree to give education for OBC. That is why the SC which practices laws of Manu had given a stay on the 27 percent reservation for OBC. Thus the only solution to uplift and make OBC socially equal is by giving education. Another reason for the Brahmins denying them education is once they become educated i.e. they get higher education they would start to think and reason so by no means they would become a slave to the religious practices which gives a means of living to the Brahmins.

VI. Finally they fear if OBC are given higher education, they would out shine the Brahmins intellectually. The issue that education alone can lift the social barrier of social discrimination based on caste was established in the analysis of the problem using FCM model.

The experts were ashamed of the persons who boast that the nation is developing rapidly, in reality the nation is developing for the past 60 years and the beauty is the rich Brahmins have become richer and the poor OBC have become poorer and they suffer the stigma, social discrimination based on birth. So only providing reservations in higher education they can become better and there by the nation can boast of its development. It is unfortunate that the nations power (low) economy and education are only in the hands of the Brahmins and in fact at their disposal so only the poor in these 60 years have not only become poorer but victims of health hazards and other social problems.

It is still very shameful to state that these Brahmins have the guts to protest against Reservation for the OBC and the media is playing the catalyst role by projecting and propagating all sorts of false views and news about OBC reservation in institutions of higher learning. Further it is a pity to note that media and the legal system have joined hands to give false and biased news about reservation. The very protests were well schemed and instigated by the BJP/hindutva fanatic forces. Finally the parliament is over ruled by the court. Is this democracy or bureaucracy? Is this the rule of brahmins or the rule of the people elected persons?

In India we are undergoing a period in which all things are taking place in a topsy-turvy way. The experts sternly warned one and all this way of functioning will not go long. The nation will be dominated by terrorism. Once a nation which stood for peace and tranquillity would be changed into a nation of war and riots, all these can be attributed only to the acts of the selfish Brahmins who have gone very fat in the economy. Such form of imbalance will not continue for long and it is certain, the nation may face disintegration if the political leaders do not take appropriate actions at an appropriate time! It is high time the politicians over look the courts and act according to the views of the parliament. The courts have failed to act justly they have lost the point because of selfishness and self centeredness. Media has already lost its neutrality. In such a chaotic situation it is better the political leaders take actions so that the nation does not become divided due to the non implementation of reservation for OBC in higher education in the central govt institution. Brahmins cannot claim IITs, IIMs or AIIMS as their property for the sole reason they have been selfishly enjoying it as a whole. It is high time they become rationalistic or they will be forced to loose!

The experts went on giving points in favour of reservation without which the nation may have to face such untold problems which may wipe out the race of Brahmins (nomads) from the Indian soil. Majorities patience cannot be taken for granted. They want both caste to stay and reservation to go how can such a thing be justified. They are even against the very norms of natural justice, said the experts. Because these arrogant Brahmin have media in their monopoly so spread all things which would lead to nations divide and permanent disintegration; assert the experts.

VII. It is unfortunate to note that the nation's economy is in the hands of foreigners who have nothing to do with the sons of the soil or the natives of India. This is a big risk and threat to the nation. Since they were the policy makers of the nation and both media and legal dept is in their hands or to be more open only at their disposal and they happen and continue to be the Rajagurus of the administration or executive be it a king rule or a democratic PM's rule. The nation is facing a threat to peace and also a threat from these nomad minorities. It is high time the politicians take a objective look at those facts; who is having the nations wealth and important departments of the nation? When British was only serving India but did not hoard the economy of the nation as the Brahmin foreigners denying the natives, specially higher education in institutions that too run by the govt. The British foreigners did not deny the sons of the soil education like the Brahmins. The truth is only after they (British) came to India they started the rules by which the OBC/SC/ST were permitted to get education. Imagine the guts of these nomads to make laws and rules and say we should not be given higher education in these institutions. Why is the SC questioning their right, but denying OBC/SC/ST rights? The experts said, suppose we were not the sons of the soil and they happen to be the sons of the soil what would be the position of the OBC/SC/ST as in times before independence. OBC would not be allowed to put on dress, except tie a loin cloth around their waist, which the OBC/SC/ST can have from the dead bodies or the dead persons. For clearly their law (Laws of Manu does not permit any good dress or jewellery for the OBC/SC/ST). When they have to be like animals, where comes' the question of education, now the so called higher education.

The only suggestion given by the experts is that now all committees and all policy makers must be only OBC/SC/ST for Brahmins being foreigners will not and by no means are going to work for the welfare of the OBC/SC/ST. The govt should see to it that the judges of all the courts are OBC/SC/ST; unless such revolutionary changes are slowly brought out by our political leaders the nation may face a revolution resulting in the disintegration of the nation or they (Brahmins) would be forcefully made to evict India.

For instance take the IT dept. only the Bangalore Brahmin is hoarding both money and power.

Now look at these IITs, IIMs and AIIMS though they are the ones run by the Central Govt who are the directors of all these institutions only Brahmins? What is the percentage of their total population less than 5%? But they dominate all professions of authority by forcefully denying good posts to the OBC/SC/ST intellectuals! Now they have crossed their limits by saying OBC/SC/ST cannot have reservation in these institution? What are they leading the nation to? CRISIS. It is high time they are put in their place for otherwise India will lose piece!

Another instance take media, it is in their monopoly. Take the courts in India, they monopolize it! So unless they were curtailed, they will once again preach laws of Manu and say OBC/SC/ST are only to serve us and should not expect any form of benefit or should not possess money, for the property or luxury enjoyed by OBC/SC/ST pains the Brahmins (Manu). Till date their mentality is not changed a little. What is their right to protest the govt rule passed in the parliament that 27 percent reservations must be given to OBC in higher education in higher institutions run by central govt? This is the fool proof of their selfishness and self centeredness.

It is high time they are stopped from such anti national activities! It is still a pity when the medicos made protests against reservations the PM and president met to plead them to stop the protests. If at that time itself these medicos who abstained from their duties resulting in the death of innocent patients were dismissed from their studies or posts; the nation would not have faced a SC ruling, a stay on reservations! The readers can understand what sort of arrogance and selfishness and heartlessness are displayed by these nomad beggars!

Why did the court admit the petition in the first place when it was against the nation policy and anti nationalistic in view? Is it not an anti democratic act? Who are these Brahmins to protest for 27 percent OBC reservation? What is their personal loss? So the readers and public must be educated of the threat they are facing within the India. Brahmins are the biggest threat to India when it comes to reservation in higher education for OBC! One may wonder why the authors have used the word threat to India! We were advised by the experts to use it! Because of their anti national anti social and anti democratic acts certainly if reservation are not implemented as per the amendment of the constitution passed with 541 out 545 majority, the nation is certain to face a revolution or terrorism were the peace of the nation will become a question mark. Already we see naxals attack in AP. It is the dissatisfied and discriminated lots becoming united with some form of ideals and goals before them. When the political leaders do not do any justice to the socially discriminated majority who are the sons of the soil (not the foreigners of India) certainly some unrest is sure to take place in the nation which will continue until some form of social reformation takes place or in a more positive way some social reformation is carried out in a forceful way to benefit the affected. It is therefore advisable said the experts, the political leaders take the relook at the situation and implement the said reservation before the social unrest takes place. It will be wise to deal with the reservation which has been delayed for 60 years and it is nearly 2 years after the amendment was passed with over whelming majority in Dec. 2005.

OBC are usually very sensitive but calm people if their calmness is disturbed it would be like a Tsunami in India. So let the so called Brahmins measure the pros and cons and not use media or courts which are under their control to postpone the reservations. They (Brahmins) with minority strength, have the guts to publish in the print media if reservation for OBC in higher education is implemented, they would become terrorist. Did they ever think what would be their status if the majority OBC who are denied their social justice and social rights by them (Brahmins) carry out their protests? The experts warned them of dire consequences. In fact they have to run out of India through the same kyber pass. Let the govt. not delay this very sensitive issue for it pertains to the social rights of the majority the experts say .They (experts) are wonder struck at the acts of these Brahmins. The only question again and again they (experts) pose to the reader is who are they to protest the govt which is implementing some programs to support the socially discriminated based on birth? They at this point welcome the socio scientist to help the nation out of this problem. Why have not the Brahmins in these 60 years never thought of the merit of the nation? What is merit, based on caste or based on real achievements? If based on real achievements what have the achieved to help the poor rural man from social discrimination and poverty? Even in British time they were not so poor and this much socially discriminated. They could have a square meal a day. But with the advent to Brahmin achievements the last man toils over 8 hours a day and denied of a square meal? Is this nation progress? What is the use of rich becoming rich? Why the political leaders who visit them only for votes are not taking appropriate actions? The Brahmins know very well to play the delaying tactics in everything, only the govt have to establish their power and implement reservation for OBC, which is implementation of social equality. First the govt. court, media and the public should teach them (Brahmins) their limits about OBC reservations and their voice in this. Because in India, the courts and media are of their monopoly or to be more factual functions for them and under it is essential for the govt to intervene and render justice. For if the public take this issue in their hands the nation is sure to face a disaster said the experts. Can we ever imagine the reverse situation? Suppose tomorrow the govt passes an act about the limit of property movable or otherwise to the foreigners (Brahmins) in India! What will they do! The experts said such days are close for the exploitation of the Brahmins over the non Brahmins have reached a breaking point.

The minute the limit goes it will be beyond any compromise for them to live in peace in India! The experts further warned that the Brahmins were always only enslaving the non Brahmins by superstitions, when they (Brahmins) themselves do not believe and practice it. During Periyar time it was little better. Now once again after the advent of BJP rule they have grown up. They by these acts, systematically hinder the progress of the non Brahmins; there by the nations progress is only retarded!

The readers may be very much surprised to see why the authors are speaking of superstition in the nation. To know the peak of superstition the authors are forced to give a small instance that has taken place in a village in Tamil Nadu. During the rule of BJP there was a sudden increase in the number of monkey statues more so in Tamil Nadu. Monkey statues were made at very big heights say up to 90 or more on the sides of the main roads. These BJP and RSS are monkey worshippers for they feel and preach that by worshipping monkey they gain strength and wealth. More so the Hindutva, VHP and RSS worship these monkeys and this practice have been passed on to the sons of the soil also. The monkey statues are more worshipped by bachelors. Now on 14-04-07 was the Tamil New Years day. In Vizhupuram district in Tamil Nadu recently the BJP/RSS/Hindutva Brahmin forces have installed a monkey statue of 90 feet few years ago. On 14-04-07, this 90 feet monkey was worshipped by anointing 10000(ten thousand) litres of milk. We are forced to record this information on the following grounds

1. The leading Tamil Daily had wasted A4 size of space to put the colour photo of the monkey statue with the information!
2. In T.N there is a rise in the price of milk for which there is protest made by political leaders and public, they why waste 10,000 litres of milks.
3. The Vizhupuram village is a poor rural area if instead of anointing the 90 feet monkey statue they could have distributed 1 litre of free milk to 10,000 poor families on the New Years Day at least they would have celebrated the year more happily.

4. Did the monkey statue drink the milk? What is the significance in wasting 10,000 litres of milk on a mere stone?
5. It is a waste of time, economy and energy.

The Brahmin nomads only taught such form of superstition to the poor socially backward and uneducated. That is why by chance they (OBC) become educated certainly they will question and would become rationalistic! Are these foreign Brahmins doing any good to the nation? Certainly the British or the Muslims did not steer the nation to such stupidities and superstitions. That is why the experts blast these Brahmins as anti nationals only their progress is important to them.

1. Is not this act in the information technology age of modernization very degrading and disgusting to any one with little reason?
2. After some 200 to 500 years will not other countries criticize India more particularly Tamil Nadu as a barbaric nation of Monkey Worshippers.
3. Is 90 feet statue of monkey a necessity in a poor village of Vizhupuram?
4. What is the monkey doing to help the poor people?
5. The Brahmins think by installing such big monkey statue they are trying to keep the villagers as Hindus. So it is an important Hindutva agenda said the experts.

Now just having given small instance of how they spread in the innocent rustic minds of the villagers such superstition and stupidity which certainly is of no use to these people or to the nation at large, the Brahmin by these act cheat the natives and keep them under their control using the ignorance of the natives. The Brahmins fully know the monkey statue is not going to do anything to anyone but they spoil the very outlook of the socially deprived non Brahmins by which they follow the Hindu religion. Even if one of two protest about it, they will be terrorized by saying the stone monkey would punish them. But in reality the Brahmins know the stone monkey statue can do nothing.

Above all the media is trying its level best to keep them under control by these superstitious acts. Is media doing any service to the nation or really hindering the real development of the nation? It is up to the reader to know the facts! Is the media really educating the public or ruining them? That is why the media dominated by Brahmins have always been propagating information which will never make the public / commoner any progress, but retard all their development claims the experts.

6. It is pertinent to mention here that a Jesuit Father Robert de Noble from Italy writes how the education was denied even as early as (during the rule of Thirumalai Nayakkar) (1623-1659). He says during his rule 10,000 students were studying, not even a single student was a non-brahmin; all the 10,000 were only brahmins; he has recorded this in his comparative literature book. Further he says that in other castes, non Brahmins were not given the right to study. Further in all posts the Brahmins alone were there and the majority were always denied all types of good posts. Arrogant Brahmins found happiness by making the non Brahmins as slaves.

Further Robert says these Brahmins dominated the govt and using their influence made policies and rules which were beneficial and were towards development of them(Brahmins).
This shows that even as early as the 17th century, the Brahmins were very clear in their ideology not to give education to the non Brahmins. Now the world have become advanced it is a computer age but the Brahmin mind still has not become advanced, or modernized still he continues his good hold ideas which is clear from the fact that they are protesting for 27 percent OBC reservation in institutions of higher learning and the SC (court of Manu) has approved it by granting a stay. It is a surprise only after 60 years of independence, the political leaders thought of to wipe out the social discrimination faced by the OBC, but this issue has been in the process for over 15 long months and finally a stay is given by the Honourable Supreme Court. Thus we see the Brahmins have never developed as regarding caste. That is why in this machine age they do all nasty things. Second the political leaders should quote such instances in public meetings and other meetings and educate the mass; then only the non Brahmins will not fall in the trap of the Brahmins. We record with the deep gratitude to the westerners especially the Britishers and the Jesuit Fathers who have given us not only the right to learn but also given us good education in the higher educational institutions.

Except for them (westerners) the non Brahmins would only be leading a barbaric life worshipping cows and monkeys and working as slaves to the Brahmins.

7. The Brahmins said that if reservations are given to OBC the standard of these institutions would go. What is the solution they have given say at least after these 60 years of independence. Even today over 33 crores of people in India go without a square meal a day. What is the use of merit which the media is falsely propagating everyday that too only after the govt opened the topic of reservations for OBC in institutes of higher learning. Did the media ever open about the 33 crores people in India not having a square meal a day? What these meritorious have achieved? They eat more than 3 times a day with all nutrition food rich in every thing. If they have really achieved something first should not they eradicate poverty! atleast this starvation. If even today the Hindu religious mutts come open to feed them, they can so do for a century with the amount of wealth they are hoarding. So it is high time the govt makes a policy for at least 60 more years to give education in these institutions of higher learning only to OBC/SC/ST and see how merit vanishes in the first place and monitor every 5 years, how many crores of people in India go without a meal a day. Unless the govt takes drastic and stern steps in this issue, the nation is forced to face a revolution. Think for a minute if these 33 crore people in frustration turn to be terrorists to rebels what will be the status of the nation? Can these 3% population make tall orders against implementation of 27 percent OBC reservation in institutions of higher learning.

In the first place will these institutions be functioning; after it has made 33 crores terrorists or rebels.

What will happen to its so-called merit? What is the nations peace? How will the media function? Media will then be forced to report terrorism! So it will be forced to report about the activities of the OBC/SC/ST? Who is destroying and disintegrating the nation? Who is hindering the nation progress? Who is destroying the nation peace? Who are they to say that reservation cannot be given to OBC? What is their local stand to say so? Is it from their mutt? Is it from their wealth? It is from the govt tax payers the non Brahmins. For here it is pertinent to mention that the Hindu religions do not show their income or any thing about their movable and immovable property! So it is high time they think for a while and act wisely!

VIII Tolerance has a limit. 33 crores people starvation is never projected in the Brahmin media. Only rich pomp and their achievements alone find their place in these Brahmin media! The first signal has come from the starving. The noble, the intellectuals and the politicians can no more ignore this. So let them not delay to do atleast the partial justice of giving the OBC at least 27 percent of reservations in the institutions of higher learning.

Since education and the employment go hand in hand we see in employment also all good posts which pay them well and which have more power are only dominated by them till date. That is even after 60 years of independence the situation is not changed. When it comes to posts with authority like directors of govt institutions, chiefs of Sires' etc are only held by them. They say they alone are qualified and meritorious. It is pertinent to mention in most of these cases, several of their subordinates would be much more qualified! They get these posts based only on the following facts:

1. They know whom to approach and they have many god fathers to do this for them, on the contrary even a meritorious OBC/SC/ST because of not having god fathers would not be even considered for these posts.
2. The media would exploit the situation to see to it Brahmins are posted.

3. As most of the administrators who have more power are only Brahmins they easily stall all moves by govt. Now did they not get a stay on the reservation for OBC?

What is the mantra behind it? Who gave the stay? The Supreme Court gave the stay. We ask in the first place how was the stay petition filed and got admitted? Once again can any such move by govt to benefit the Brahmins be stayed in the SC by a non Brahmin? Even in the filing level it would be rejected? These are facts and they are presented by the OBC/SC/ST intellectuals who had served as experts in various categories of analysing the role of media in the OBC reservations! When the intellectuals have started to feel restless over the issue, the issue cannot be taken by the govt lightly. Most of them feel the denial of reservation to OBC that too after 60 years that is after supported by 351 out of 355. If denied the reactions of them will be more on the revolt will make a peaceful nation into a turbulent one.

The court orders have wounded the feelings of politicians apart from the OBC. For it says, “State is empowered to enact affirmative action to help the backward classes, but it should not be unduly adverse to those who are left out. Now did the left out in these 60 years ever think they were enjoying someone’s shares and did the court ever think they were adverse to OBC till date. Are they asking for the back log of the 60 years? It is damaging to say “help” OBC for it is their “right”. is never a help. When right of a person is taken it amounts to a grievous mistake first it should be punished!

The bench further said, “reservation cannot be permanent and appears to perpetrate backwardness”. This is given by a learned why don’t they think even after 60 years of independence 33 crores in India living in the backwardness even for a square meal!

What can be called as forwardness? Even before giving or implementing if they make a statement of this kind what does this show asked the experts; only bias and nothing more! Did in any place the statement reservation for OBC is permanent was ever mentioned, asked the experts.

Indicting the government for enacting such a law, the bench said, no where in the world castes queue to be branded as backward. Nowhere is there a competition to become backward.

With this act the subject of the equality is unduly put under strain. “The experts wanted the authors to record their views. So if the authors do not include, next time they would not serve as experts. This group is a heterogeneous mixture of intellectuals from doctors, senior councils, doctorates IAS and committed socio scientist. So authors had no other choice to include it for when the book comes in print they would not loose faith in us or so. They said only in India, we say majority of the people as Sudras. Sudra according to Manu means born to a prostitute and a list of castes are given according as when the prostitute conceives with a Brahmin, Vaishiya, Sudra and so on. So if this caste even after many thousands of years cannot be abolished but practised as if practised in very olden days which is existing no where in the world!

As long as caste system exists in India, the experts said reservation would be in India. On 24th Dec. 2006, the Brahmins have asked reservations for them! What would the learned bench comment on this? Will this also perpetrate the backwardness or perpetrate forwardness and development?

On the contention of the centre that it had taken the 1931 census as basis for fixing 27 percent quota for OBC, the bench said what may have been said the data in 1931 census cannot be a determinative factor now.

Here it is pertinent to mention the views given by the experts. They (experts) said that the learned judges had failed to learn that the past history about the census from 1951 on ward.

“Sardar Vallabhbhai Patel; independent India’s first home minister had decided in 1950 that there will be no caste based census from 1951 on wards, when the first census took place in free India. “B.P Mandal had himself written to three successive home ministers and had repeatedly requested them to conduct a caste based census. But he was refused saying that Sardar Patel independent India’s first Home Minister had decided in 1950 that … free India. P.S Krishnan advisor to the Human Resource Development Ministry for this supreme court said Patel’s reasoning was that the country needed unity and such a census could have caused divisiveness. Krishnan said that Mandal wrote three Do letters to three home ministers namely H.M

Patel in 1978-79, Y.B. Chavan in 1979 and Gynani Zail Singh in 1980.

Krishnan who was Social Welfare Secretary in the 1990 V.P Singh government when the Mandal Commission judgement came and was also the founding member secretary in the national commission for backward classes from 1993 to 2000 was brought in by Human Resource Development Minister Arjun Singh last year to help the government defend its position in the supreme court as the OBC quota legislation was challenged. While the UPA government to satisfy the Supreme Court will have to overturn Sardar Patel's decision of 1950 for having a caste based census. Krishnan says, 'There is adequate data for implementing reservation for OBC in education. Data is not a problem". He says every survey ranging from Census to National Family Health Survey To NSSO survey – puts the percentage of OBC more than 27 percent. So 27 percent reservation is a very safe figure, that cannot be questioned, "said Krishnan who has been instrumental in preparation of all counter – affidavits filed by the centre.

He also said that it is incorrect to say that 1931 census was the sole basis of projecting an OBC population to 52 percent where as the Mandal Commission followed a four-pronged strategy. "It considered the state lists collated evidence presented before it, toured the country made random surveys, before coming to a conclusion of the 52 percent OBC population.

When everything was so nicely presented before the SC by the HRD, the Brahmin bastion BJP has the guts to say that the presentation of the case from the govt side was not done properly. Now the Registrar General of India had commented that to take the census we need 2,00, 000 enumerators and the caste of the nation would only be reflected by the castes of these two lakh enumerators. So even now if we take a census it would be only fully biased he said. Further the only question the experts wish to place before the reader is, "Why does the SC imagine that the OBC population would be even less than 27 percent?" Has any genocide occurred on the OBC so that the OBC population has reduced to 27 percent they questioned the authorities. The bench said the concept of Creamy layer is not

prima – facie relevant as contented by the centre. The bench further said there is no explanation as to why there is no firm data for determining backwardness” As the result unequals are treated as equals. For this solution is given in the last chapter of this book. We do not want repetition so we have restrained from giving it here.

The petitions had challenged the law contending that it would have wide ramification and divide the country on caste basis, resulting in anarchy and affecting communal harmony as well as jeopardising the constitutional rights to equality. Vote banks should be replaced with talent banks and not with caste politics.

IX Whether they (govt) give reservation or not, no one can ever dispute the fact that even today the nation is divided on castes. So reservation has nothing to contribute to the division, for the division is existing. Even today that upper caste especially Brahmins don’t dine with OBC/SC/ST in the premises of the Kanchi Mutt. It is still not surprising to record even in AIIMS, SC/ST students cannot dine in the hostel premises with upper caste. We ask one of the reader ‘Is this a caste divide or a reservation divide’? Only caste divide as reservation is yet to be implemented. Secondly the divide was so shockingly noted even in the $16^{th}$ century by Robert Novel. No one could be educated only Brahmins students of strength 10000 were educated and well fed. Still it would be shocking to the reader in T.N we have a village in Coimbatore named Singanallore. This incident has taken place in 1930. The MLA of Singanallore was one Mr. R.K Shanmugam. He later became the finance minister of India. In the Singanallore Village there was a primary govt school (i.e. only upto $5^{th}$ std.) for SC/ST cannot learn with others, that is the caste divide in India. The students can reach the school only by going through a road in which there were only houses of Brahmins. So the Brahmins never allowed the SC/ST to go to school. Whenever the inspector came to inspect the school(once a year) the Brahmins would take these SC/ST students in a van to the school and only for one day they (SC/ST) will stay in the school premises till the inspector leaves. For the British in those days had put the rule that every school run by the govt must admit SC/ST students to get the aid to run the school with total govt. support. Thus one day these children (SC/ST) would attend in a year. This is a real happening. Is it caste divide or reservation divide?

Caste divide has been ages since the nomads entered India. Only after their entry it was said that caste is based on birth and a large is based on birth and a large class of OBC people became sudras, born to a prostitute and another group untouchables. It is still not clear from the laws of Manu how they made the two castes as distinct Periyar in fact has very clearly said that the Sudra given to the majority of the non Brahmins is worst than the notion of untouchables. Now is the nation divided based on castes or is it one?

Already the nation has divisions due to language. So it is a false allegation for everyone knows very well that the nation is already divided based on caste and the implementation of 27 percent reservation to the OBC is by no means going to divide the nation. At least by educating these Sudras, the nation in due course of time will make some amendments for neglecting the majority of the sons of the soil who have been socially discriminated ever since the entry of the Aryans. It is pertinent to mention the vanniyars who were rulers of the nation have become now the most backward and their literacy rate is also very low. The shocking information is that they are poorly employed in spite of being OBC (MBC). Likewise Mallars who were also rulers have become untouchables. It is a disgrace the Aryans entered as nomads and have infected a peaceful united integrated nation with castes, since their entry the nation remain divided.

Until the laws of Manu are totally banned and Vedas are banned in India, there in no chance of seeing India as an integrated one. The Brahmin arrogance and the supremacy they practice on the sons of the soil has been proved beyond doubt, when they had the courage to protest against reservation and using the principles of Manu have stayed the rule passed by the govt that too in both the houses. One of the reasons given by them is reservation will divide India what a hippocratic and anti democratic statement is it! The nomads in reality are anti nationals that is why they are protesting some social reformation. They are least bothered about the majority welfare.

Only they are most bothered about their welfare or nation welfare. They have not changed they want only the Manu rule that is do not give education to OBC.

The Brahmin media will not give any true information on any matter if it concerns OBC/SC/ST. Only their anti protests were blown out of proportion by the media. This view was given by Periyar as early as 1937 even before India’s independence. So the media has been functioning to promote only the welfare of Brahmins and has always blocked out the problems of progress or of achievements of the non Brahmins.

The Brahmin media to block out or totally hide the reservation issue will talk only about other things like; now the marriage of Amitabachan’s son, Sethusamudram project etc. They are well motivated and planned in doing these acts. Suddenly why monkey statues are built or consecrated and this news given the maximum publicity? Even on April 18, 2007, Wednesday, “The New Indian Express” has given in page 9 nearly a quarter page photo of the demon being stamped by the monkey statue. The mouth of the demon is so large that people can go through its mouth. Further they are proud to keep on record the monkey statue in a newly consecrated in a temple in New Delhi is 32.9 metres in height is perhaps the tallest of its kind in Delhi. It took 12 years to construct the temple! When 33 crores of people are starving without even a meal a day! What a waste? How the Brahmins are unnecessarily wasting the money on the monkey statue when nearly one fourth of the nations population are starving without one meal a day. Why is the starvation and the starvation deaths that is occurring every day is not reported.

The readers are expected to ponder over the monkey statue and the money spent on building it and its utilitarian purpose to the nation. Are the one fourth population say 33 crores important or the building of a monkey statue for over 12 years and consecrating it important! Will the Brahmins do anything good for the native, let the readers ponder over this!

They have usurped the wealth, education and all form of comforts, social identity from the natives and left the natives in poverty, ignorance and superstition. Now they have lost the rationality and have started to work for the making of the monkey statue in the $21^{st}$ century. It is high time reformers and people with rationalist ideas take the control over the situation and educate them. Once they are educated naturally they would start to develop and take the lead otherwise they will die as slaves. It is the dire need of the hour. The natives curse their past karma for their poverty forgetting their wealth and land which was once in their possession have been usurped by the Brahmins. They should be made rationalistic then only they can come up in life; the only way to make them so is by educating them. The Brahmins are frightened for once the OBC/SC/ST become educated they will be rationalistic they will not worship the monkey statue and they will not also help the Brahmins to build such huge statues and they will ridicule the Brahmin for worshiping a monkey. If the 97% fail to visit temples, clean and maintain it what will be the Brahmin priests plight? Who will feed them? So only they are frightened their power will go! Above all they will ask these, OBC/SC/ST to visit the monkeys temple or some other temple every day and worship it so that they can have wealth or comfort and so on. How can the OBC/SC/ST get wealth or comfort when the Brahmins are the sole custodians of them? This is only to get crowd in the temples and enough money and materials to run the temple as well as fill their stomach. So only higher education is the solution to uplift the natives and this can be achieved only through reservation in higher education. That is why when govt. puts up the order of reservation for OBC they have stalled it using the SC. The judgement is only a blow to social justice or socially equality. The judgement is supporting social discrimination for equal opportunity is not provided for a large section of natives on the false ground it affects the social equality.

Any order or rule cannot be based on promoting only a few and denying a majority. All these 60 years only the Brahmins have grown beyond proportion in wealth, power and education at the total cost of the non Brahmins which amounts to social discrimination based on caste. Unless such situations are changed the nation is deemed to face disaster!

Now we in this section take the opinion of over 1000 students from arts colleges from various districts of TN, about 27%

reservation for OBC in institutes of higher learning run by the central government.

X The students formed a heterogeneous group. This study we felt was important for mostly the courses in the institutions like IITs, IIMs and AIIMS are professional and usually these arts college students are looked down by these students of these institutions like IITs, IIMs and AIIMS. So we wanted to find whether these arts college students supported 27% reservation for OBC or not. These students belonged to all castes, all religions and were both from poor and rich social strata. For these students who do certain special courses as evening college students are rich and also have a better educational background than the majority of the day college students. A brief introduction about need for OBC reservation in higher learning, the govt stand and the oppositions and favours was discussed with them for 10 to 15 min and they were given over 40 min to write their views critically about OBC reservations in higher learning. We made this approach for in the first place we wanted them to think and react; for if we provided any set of questionnaire they would tick or write mechanically. Their true involvement would not be there so only we adopted a new method of approach to get their views. For being in the age group 17 to 25 the future generations of the nation their feeling and opinions were very much important for our study. Majority of them did air their views in a very sincere way. This is given in the last chapter of this book. Some the striking suggestions in support of and also a few against reservations which surprised us is recorded in the following in a line or two.

a. Some of them said when the govt is starting new educational institution it should in all two categories; employment and student admissions they should follow the "rule"; do not call it by reservation call it as "equal opportunities for all" or equality for all in govt by giving 3% to Brahmins, 97% percent to non Brahmins shared by Muslims, Hindus, Christians as per their strength. This should be done for all posts teaching and administrative and student admissions. If the govt had been doing like this now the question of reservation will not raise and all economy of the nation will not be in the hands of the Brahmins, foreigners to our nation. Our elders have failed miserably, at least now let them wake. Also the caste which is now so much dominating and dividing the nation would have been settled long before. Since a particular caste Brahmins had occupied all good posts they have used their power to admit only Brahmins in these institutions in these 60 years. They have been enjoying solely and selfishly all the benefits of the nation at the cost of the non Brahmin sufferings. That is why now when asked to part with a paltry percentage of 27% they are staging mega protests. The carelessness and lack of foresightedness on the part non Brahmins have now made, the Brahmins dominate the non Brahmins and rule them. This is the case with courts, all administrative posts which have power in the centre and all posts in ISRO DRDO etc. At least now let the govt in centre or state, when they introduce a new office or a educational institution or a department or any other work place. Let it recruit people from all castes by giving only 3% to Brahmins and 97% to non Brahmins. The govt has once burnt its fingers atleast let it be very careful. Let it rectify this by not giving any posts to Brahmins in these institutions still their strength comes down to 3%. Unless govt takes such steps no one can even curb the arrogance of the Brahmins said this group of students. If in the nation social equality should prevail and social discrimination to be annulled this is the only way!

b. The students strongly protested the acts of the Brahmin medical students who protested against reservation for OBC. They said only the amendment has been passed with a overwhelming majority in both the houses (i.e., 351 out 355 has favoured the bill). What is the right of these foolish, selfish Brahmin students to protest? The first day of protest itself the govt should have dismissed them. If that was done will any other student ever come to road and talk or protest against reservation. Here also they said govt failed to take necessary steps to punish the Brahmin students when their protests had no local stand; for who are they to say no to reservation for OBC? Suppose the OBC protests what will they do? So the acts of Brahmin student is very mean to the core they condemned. Any act of social reformation by the govt should not be questioned or protested that too when it is done to favour a suffering majority who face social discrimination they said. Some students of arts college said the anti reservation protests of the Brahmins confirm them only as anti social and anti national elements.

c. Students said the politicians with Hindutva ideals were the brain behind all these protests that is why they made all arrangements to these medicos like Shamiana to give complete shade and they were provided with air coolers and they protested lying on the bed sheets under the shade either reading novels or sleeping and made a novel protest. Only Brahmins alone can remain so comfortable and enjoy all things even during their agitation. To crown all this they were paid their usual pay for over 40 days and were given attendance for it would affect their course. Who in India can ever get all support? Asked some students.
d. The Brahmins media was happy to cover in colour and then flash it with headlines and photos in the first page. We cannot get all these even in our dreams they said. The Hindutva politicians if they had something like guts or courage they should have protested in the parliament and not instigate the students. They do not want to make any issue in the parliament for they are sure they will lose their MP post next time. They want power as well as votes and still protest. This is exactly the character of Brahmins, they said.
e. They all criticised Dr. Venugopal and said when he was a student he could not protest that is why like a coward he was supporting the whole show from behind for if he openly protests reservation for OBC he is protesting against his own authority claimed the students. He made the maximum use of the situation by providing every thing to the protesting students including the campus. Can any student protest in the campus by erecting a shamiana and a air cooler asked some of the boys without being arrested for over 40 days. We would be chased out of the road when we make protests on the road in the scorching sun, they said.

f. Even when the protests are carried out by students of arts, colleges we never make such under democratic protests like the Brahmin students. They had proved themselves anti nationals and anti socials for what they are saying is nations merit will go if reservation is given when they have no feeling for the nation, an height of cunningness! What is their answer, merit or food at least once a day for the 33 crore people which is important? What have they driven the nation to, in these 60 years of independence? Even now the night soil is carried on the head by OBC/SC/ST. What is the advancement towards social discrimination? How many Brahmins have carried the night soil on their head? When govt does something to establish social equality and put an end to social discrimination based on caste (Here they wanted to say now various and some MBC do the work of scavengers and sanitary workers) who are they to protest. Let them (Brahmins) for at least 60 more years carry the night soil on their heads said the students of arts colleges. Some of them said our parents also do the same job. Even in class room Brahmin boys do not sit near us we are only MBC. We said if you don't sit our side it is good for you can never clean the shit in your heart which you are carrying for thousands of years. We ignored them said a group of boys. They said our language, our dress, our skin makes them to run away even if our parents like them (Brahmins) had been in shade, AC and ate fruit juice, milk, ghee and curds; we too would be like you, unfortunately we do not get and some of us have not seen them yet they said. The interaction of these students was very much educating us and gave the authors lot of surprise and shock! These Brahmin boys would reserve the first bench in front of the teacher. If we sit in that place they would wipe it and then sit so to teach them a lesson we sit in their place till the teacher comes and they every day wipe the desk and table and sit. As they have not complained about us we too maintain a stone silence with that group. Even our assignment sheets they don't touch! These students said we have thousands of instances like this. How we are treated in the hostel in the dining room etc. they said!

## Chapter Three

# EXCERPT OF NEWS FROM PRINT MEDIA AND SUGGESTIONS AND COMMENTS BY THE EXPERTS

This chapter gives the views and suggestions of the experts about reservation for OBC using the print media.

We had a heterogeneous group of experts from socio scientists, economists, journalists, educationalists, politicians, public, students, farmers, govt. and private employees, doctors, engineers, psychologists, administrative officers and teachers. They were divided as groups basically depending on the place they belonged to so that they could meet regularly and carry out discussions everyday or on weekends. We covered 19 places stretching all over Tamilnadu apart from a group, which worked from all over the world carrying discussions over the Internet. Periodically we collected the information from them with the views and comments. They not only comment on the news, which appear in the print media but also give suggestions and views. They have very openly commented about the short comings of the media as well as the selfish attitude of the upper-caste which was using media to stall the 27% reservations for the OBC in the institutions like IITs, IIMs, AIIMS, etc run by the central govt. of India. We have mostly given the views and suggestions verbatim. The authors have acted only as machines mechanically transferring the information of the experts.

We have taken the excerpts from The New Indian Express, The Hindu, Times of India, Deccan Chronicle, etc. We have also not given the views chronologically as the groups happened to be a very heterogeneous one, some were responding much more faster than others and some were little slow as they could not spare time on weekdays and only on week ends they gave their views. However we have remained impartial and unbiased in giving their views in this chapter.

The New Indian Express

The first direct talks between anti reservation medicos and Prime Minister Manmohan Singh tonight failed and the students and doctors decided to continue their strike, which has thrown basic health services out of gear in several cities. Singh had an hour long meeting with the medicos during that he told them that he fore sees a massive expansion of the higher educational system which will see a huge growth in educational opportunities available to all classes and categories of students. The prime minister assured the medicos that the road map laid down in the decision of the UPA coordination committee addressed the concerns of all categories of students.

Immediately after the meeting with the prime minister the general body of resident doctors and medical students met for two hours and decided to continue their strike.

In a related development two IITs termed the govt. move to increase the number of seats ‘disastrous’ saying it will adversely affect the academic quality of the premier institutes and asked for its reconsiderations. It would be most disastrous to impose a 27.5% quota on IITs in an ostensibly fair way by increasing the number of seats.

The experts viewed the news of our prime minister meeting the striking medicos was itself not amiable to any rhyme or reason. For they feel if the majority of the OBC made a collective strike that their right to study in a govt. institution of their choice was violated, by the few Aryans who were agitating, what would the govt. do? The experts warned the nation if the quite sufferings and rights of the majority is not given due consideration even after 6 decades of independence, the consequences would be dire and irrevocable!

They criticized the IITs tall claim that the reservation will affect disastrously the academic quality of primer institutes; for he says certain IITs had no guts even to place their Ph.D and M. Tech thesis in open. They are kept in lock and key for their plagiarism would come to light.

Sunday Express

Sensing strong reactions against the union HRD minister's move for OBC reservations in educational institutions the center is mulling over a strategy for a retreat. It has now taken the stand that states could implement their own reservation policies … while Manmohan Singh and his close associates including Kapil Sibal, Hansraj Bharadwaj and Panab Mukherjee are against reservations, others such that Sharad Pawar, Lalu Prasad Yadav, Ram Vilas Paswan besides Arjun Singh are for it.

It is learnt that Arjun Singh was asked by the prime minister to remain silent on the issue and air his views only after the Assembly poles…

More over Arjun Singh's hand was forced by the near unanimous political consensus that the constitution should be amended in the light of the verdict of the Supreme Court.

But later Arjun Singh sought to play the Mandal Card to the advantage by hinting at quota in IITs and IIMs invoking the enabling legislations. Arjun Singh piloted the $104^{th}$ constitution Amendment Bill which was passed in the Lok Sabha on Dec 21, 2005 with 379 votes in favour and one against and one abstaining. On Dec. 22, 172 of them voted in favour of the amendment and only two against it. However after passage it became the constitution $93^{rd}$ amendment.

The amended article 15(5) states noting in this article or sub-clause (G) of clause (1) of Article 19 shall prevent the state from making any provision by law, for the advancement of any educationally backward classes of citizens or for the scheduled caste or scheduled tribes insofar as such special provisions relate to their admission to educational institutions, including private educational institutions whether aided or unaided by the state other than the minority educational institutions referred to in clause (1) of article (30).

By the constitutional Amendment the government has provided for reservations for OBC's in all educational institutions including private whether aided or unaided, excepting minority educational institutions.

From the report by Anita Saliya "Centre looking for exit route on reservation row" the reader is made to understand the hard truth that out of 381 vote 379 voted in favour in Lok Sabha and the Bill was passed in the Rajya Sabha 172 voted in favour and only two were against it on Dec. 22, 2006. Thus OBC reservation was welcomed by the elected political leaders. What does all this mean? Every one including the Aryan Brahmins are fully aware of the fact that reservations for the OBC in institute of so called quality is essential in order to uplift the majority on one side and secondly to see that the investment done on IITs and IIMs are not exported for service. If majority of them is educated certainly they will not seek to fly to other nations for money; after getting the good foundation in India which amounts to day light swindling of India and behaving like traitors by serving other nations. However a group of experts maintained that Arjun Singh was not for 27% reservations for OBC, his acts had some ulterior motives. They said he was against the Prime Minister Manmohan Singh.

The New Indian Express

"Strikes in medical colleges and hospitals are illegal according to a Supreme Court order and patients cannot be taken to ransom. Hospital strikes means playing with people's lives"; said Union Health Minister. The Minister however expressed the hope that the situation will return to normal soon "We don't want to take such drastic steps and we condemn any brutal action taken by police we hope the students and doctors realize the gravity of situation and return to work" said Ambumani. He added the agitation is uncalled for" as the Human Resource Development minister is ready to talk to them"

The experts felt even the govt. hinting the raising the number of seats in premier educational institutions to accommodate the new 27% OBC quota is itself very unfair and in human.

When a bill is passed on with over 99.5% support what courage the few Aryans have to protest.

Any reform is made only to help the less privileged. It is but unnatural on the part of the selfish Aryans to protest, when the govt. tries to uplift the majority who seek good education. It is very cruel on the part of the nomad Brahmins to have stalled mandal. This time if they try to repeat history the expert says govt. should take firm steps to show them their proper place in the society. Any govt. elected by the democratic way should only toil for the majority and for the minority Brahmins who are already in the realm of life needs to know the truth says the experts. The denial of admissions to OBC in central govt. institutions amount to discrimination by caste and social status. They are not realizing the harm they have done to the nation i.e., to the majority for so many years.

Express news service

Extending support to the proposed 27% reservation for OBC in institutions of higher education, Bharathiya Janashakthi Party founder Uma Bahrain called for a national debate to ensure that the benefits of reservation would really reach the deserving ones among the OBC's.

Addressing a press conference here on Sunday Uma Bharati … said "under the present education system in India, there are a lot of disparities between the urban elite schools, which charge high fees and the rural government schools where the poor children study under leaking roofs using broken pencils".

**"To know and to defer false choice between merit and social justice**

The continuation of vitriolic attacks on the knowledge commission by prominent political leaders speaks poorly of the reverence for free and informed discussion on the country. It betrays in addition a refusal to read the inclusive nature of India's edge in the global knowledge economy. At the heart of the current problem is much more than the specific fate of the commission, possibly the most crucial policy initiative of the UPA government. The proposal to reserve seats of OBC's in educational institutions derives from the need to make access to opportunity socially inclusive. Without social inclusiveness no society can be truly democratic or be confident of meeting its real potential in political or economic term. The trouble with politicians currently agitating the commissions's 6-2 divide on adopting quotas in this; they place allegiance to reservations above that larger objective of social inclusiveness.

The debate began by the knowledge commission is too important to be allowed to be stifled by the political correctness or fear of allowing non governmental inputs in education policy… The ludicrous charge that the commission's members are transgressing some constitutional norm by expressing disagreement with parliament's reservation legislation is not just undemocratic – it is an assault on the prerogative to harness informed opinion. The commission derives its existence from that express will of the PM to draw expertise to sharpen India's knowledge edge.

Infact indications from the Prime Minister's office are that he along with a group of ministers constituted for the purpose is intent on finding innovative ways of escaping from the current zero sum mind set. This is important. The mindset had already set off a very dangerous confrontation based on a false choice between merit and social justice. Induction of intellectual representatives of industry and politician in undertaking reform in the education sector in order to achieve both merit and social justice is the key to removing anxiety and a sense of grievance. And it cannot be done as long as MP's and ministers remain intolerant of disagreement.

Very many experts reacted to this article. Over 50 of them said they had written a response letter to this article. The intolerant Brahmins monopolized press did not publish even a single letter of disagreement says over 50 elite intellectuals and educationalists who wants really to uplift the nation. What is merit to them? What is social justice? It is long years they have usurped over merit and social justice under that pretext, and left the majority miserable.

The media is in their hands so none of our protests even come out to public. When the term "reservation" comes they boil, burst and burn. What does this show they want to lead a monopolized comfortable life least bothered about the nation.

Who is the monopoly of our industries? Who are pinnacles of knowledge and intelligence ask them? They will say only they are meritorious so only they are.

They have been enjoying and living in luxury all these years, denying the rights of the majority. The first and the tall claim is that, including socially backward is not truly democratic. So does it imply all these years with no reservation for the OBC the nation was functioning democratically! What is democracy? According to Brahmins; is denial of all things including education and economy to the majority is democracy! What is the average age in which a rural poor agriculture labour dies in these years? Has science and technology reached them or done anything good to them?

The majority are under the tyranny… tyranny of the so-called Brahmin policy makers majority do not enjoy true democracy. They do not even have a square meal a day but still continue to work over 10 hours a day? Is this democracy to them? Are they not the maximum sufferers after 60 years of independence? Why no one is ever bothered about them? What is the intellectual support played by the knowledge commission for their sake? How many farmers suicide took place in India? How many are Brahmins? What is social justice the knowledge commission rendered to them? Except for the wrong agricultural techniques which is headed by a non agriculturist Brahmin MS Swaminathan, had ruined them.

The New Indian Express

"The union ministry for health and family welfare today issued "termination notices" to senior residents in all medical colleges in Delhi. They have been given 24 hours to return to work after which the termination will come into effect. The ministry had issued the senior residents notices on Monday calling their move to join the striking doctors as "contempt of court". Following this some doctors resumed work but many remained on strike.

The ministry had to take some view if they don't come back. They are not affected with the proposed policy and they should come back rather than agitate without any reason that is why termination notices have been issued to them. Health secretary P.K. Hola said representatives of the agitating doctors presented a memorandum to union health minister Anbumani Ramadoss, and Pranab Mukherjee.

The New Indian Express

The students and resident doctors today decided to continue their agitation after another meeting with the President APJ Abdul Kalam, which failed to break the deadlock over the proposed 27% quotas for Other Backward Classes (OBC).

Hospital services were affected for the 13th day in succession. Clinics and chemists shut shop till noon as part of the IMA medicial bandth. Faculties was mainly affection the left AIIMS took mass casual leave to show their support for the students affecting the main Out Patient Department (OPD) as well as the parallel one’s set up in the early days of strike.

The experts who spoke about their 13 continuous days, the life’s of the poor patients where some of them may need medical attention immediately. They are least bothered about the ethics of their profession but so much bothered about the reservation! They have no right over Government decision to uplift a section of people who form the majority population of the nation.

The New Indian Express

It was mentioned that bar dancers will perform stage shows here next month to raise funds for youth for equality to fight against the centre’s reservation policy.

Express news service.

Representatives of the Indian Medical Association (IMA) from the four southern states have come out strongly in support of the 27 percent reservations in Central Government run institutions.

At the IMA Tamil Nadu’s 60th Annual Conference TAMCON which began on Saturday; representatives from Tamil Nadu, Kerala, Karnataka and Andhra Pradesh adopted a resolution to his effect.

Already there is 69 percent reservation in Tamil Nadu, 60 percent in Andhra Pradesh and 50 percent in Kerala and Karnataka…

Express news service

Condemning the agitation against reservation in the northern part of the country Dr Ramadoss said “it has been instigated by elements with vested interests and a section of the media. "He urged the center not to yield to pressure and withdraw the move to introduce 27 percent reservation in Central Government educational institutions. Several experts said that Dr Ramadoss was very correct in pointing out it is the menace of the media which is functioning to the full advantage of the anti reservation policy. The News Paper media is not doing properly its duty by favouring a section of the people. Any media's policy must stand unbiased and report both sides. It is unfortunate that the majority of the media in India is run by the Brahmins and is very strong in upholding their views in a biased way.

Express news.

Termination notices were today issued to striking medicos in state-run hospitals in the capital even as the Delhi government today came out with ads. to recruit junior and senior doctors to cope up with the situation.

Express news

While the campaign against the proposed reservations for OBC's in institutions of higher education spreading in other metros, Chennai continues to campaign for reservation. As part of the campaign an SMS sent in favour of reservation gives the data regarding the dominance of forward castes among the faculty of IIT, Madras. The SMS reads as follows: "Chennai IIT Professor 400 among them Brahmins 282 (70%) other FC's – 40% (10%), BC 57 (14%), SC/ST-3 (0.75%), Christians 15 (3%), Jains 3 (0.75%) and Muslims nil.

The total population of Brahmins in India is 2.1% Is it a justice? Will the reservation in Higher education be the solution to this injustice or not? Raise your voice by sending SMS to atleast 10 members", the SMS reads. It is not sure whether the statistics mentioned in the SMS is correct or not.

Meanwhile the Doctors association for social Equality has appealed to the doctors, medical students party leaders, mass organization and people to participate in the human chains to be formed on May 23, 2006 throughout the state by the Campaign Committee for Reservation to support the Central Government initiatives to implement 27% reservation for OBC's in higher educational institutions.

The general secretary of the association Dr. G.R. Ravindranath said it was regrettable that some organization was demanding exclusion of creamy layer from benefiting from the 27 percent reservations and demanding it on the basis of economic criteria.

These demands were unwanted at this juncture and they were diverting the main core issue, he said adding that reservation for the OBC's should not be delayed and diluted at any cost.

The association urged center to put an end to the confusion and implement the 27 percent reservations in to for OBC in higher educational institutions without any delay from the current academic year.

Peter Alphonse Congress whip said, "I suggest to all those who claim that reservation will go against merit and qualification to look at Tamil Nadu where the policy has been in existence since 1927. The Tamil Nadu experience shows that reservation is not against merit" he said emphasizing that reservation would provide adequate opportunities to the OBC.

The article by Seena Menon unreserved comments

In 2004 Rajeev Goswami died a lonely death in the Safdarjung Hospital. Two years later the imaginary demon he aspired to burn along with his self immolation bid in 1989 seems to be raising its head again….

An engineering student, I know how much one has to do to get into such elite institute like IIT. Implement the reservation scheme and a second anti mandal protest will start. That is for sure. What is really sad about this is that all you guys who oppose are Brahmins and who have for 1000 years oppressed and kept the highly intelligent minds out of the focus… now it's time for pay back… just sit back and enjoy the show.

What is the Pandora box that the not so unsuspecting HRD minister has opened? He uttered what was recommended in the Mandal Commission almost twenty years before. It means that the ministry is proposing to reserve seats for other backward classes (OBC) in educational institutions including central institutes such as IITs and IIMs. These institutes are not just learning centers but brand names. When companies come to recruit students from IIT they will not understand the reforms. They will not change their recruitment patterns says a student of IIT Madras.

Reservation in prestigious institutions like IITs and IIMs is highly observed as anybody can get admission. The future of the talented and the deserving students are at stake says a second year medical student from the Kottayam Medical College.

This brings us to the one major contention – merit.

"Reservations do not rob institutions of their excellence, instead they create new circles of excellence", believes Dr. Jayarajan of Institute of Development Alternatives, Chennai. Moreover most of the professional institutes have eligibility criteria. For instance only those who clear the screening test such as GATE are selected to IITs.

The reason why central governments have taken their time to implement the basic recommendations, however is clear to most, "It is the lack of political will" say Dr. Jeyarajan. Despite the passionate protests there are certain voices that see the statement for what it is. First of all I feel the issue has been overblown and does not deserve this much attention say Roshan (name changed) of IIT - Madras.

"Reservations for the OBC's are devised to divide the society vertically and horizontally and to exploit a large number of people – poor, deprived and ignorant – so that a few on top may continue to enjoy power and privilege" says Tony a medical student at Manipal.

The views of few experts who gave their opinion about this article

They all first pointed out and objected to the statement given by a student of IIT, Madras which is as follows:

"When companies come to recruit students from IIT, they will not understand the reforms. They will not change their recruitment patterns".

They (experts) protested that it is not IIT alone which has campus interviews made by companies and recruitment taking place, Anna University has MIT Bharathidasan University and so many institutions why even colleges carry out such interviews. The expert says, only Brahmin owned companies will face problems!

For they have to make change in their pattern of recruitment to get Brahmin students! Our OBC students are not that dull let these arrogant Brahmin fellows first understand it said some experts! Because they happen to be brighter than the Brahmins they do not want the government to give some opportunities by way of reservations, so that some justice is done to the majority.

A student of Kottayam Medical College says "Reservation is prestigious institutions like the IITs and IIMs is highly absurd as anybody can get admissions. The future of the talented and the deserving students are at stake". The experts asks what is the meaning of "talented and deserving" that too just after school final exams. For marks say, just pass in the final exams got by a rural poor boy shows 100 fold merit than a city bread Brahmin who gets over 90% for the circumstances and the encouragement and above all the atmosphere both get ought to be compared. Who is talented? Who deserves a seat? Certainly not the Brahmin fellow says these experts. All people with little reasoning will accept our argument says the experts. Place the problem for analysis before any socio scientist all will agree only on our conclusions.

All these years that is 6 decades after independence did deserving students with real talent given any opportunity? Why no one has ever bothered about the future of the really talented and deserving students all these days?

The basic claim is that even those who are going to be selected are persons who should clear the screening or entrance tests as said by Dr. Jeyarajan of Institute of Development Alternatives, Chennai. The experts feel that unnecessarily, the few leading media is blowing up the issue, falsely to get a maximum feeling of anti reservation wave among the public. It is only more a media myth than any reality feels the experts.

The New Sunday Express

Sure the OBC did not face untouchability and most of them did not suffer from the worst oppression of caste system. But they have suffered from systematic disadvantage in accessing education and middle – class jobs. Look at its effect today/ according to the National Sample Survey out of 1000 upper caste Hindus in Urban area, 253 are graduates. Among the Hindu OBC the figure was only 86 per 1000. The picture gets worse if we look at post graduate and professional degrees. Caste-wise breakup from another study shows that access to higher education still reflects the traditional caste hierarchy; the rate of highly educated is 78 per 1000 among the Hindu Brahmins around 50 or plus for other twice born caste Hindu's Christians and Sikhs (with the exception of Rajputs who now include many upwardly mobile non dwijas) but only 18 for the OBC and even less for SC and ST. The inequalities in the level of educational attainment of different caste groups are still unacceptable large. The situation is not an outcome of any natural differences in IQ of different caste groups or uneven desire to peruse higher education. These differences are principally an outcome of unequal opportunities. That is why the government needs to step into this; OBC or other Backward Classes or Backward Communities other than Scheduled castes and Scheduled Tribes. These were most 'Shudras' in the traditional Varna hierarchy: below the dwija or the twice born but above the untouchable communities. But not all sudras are recognized as OBC's…

All upper castes (all those form any religion who do not qualify currently for SC, ST or OBC quota) are about 33 percent of population. Even after the OBC reservation, at least 50 percent of the seats in higher education will still be open to them. So strictly in the narrow caste share calculation it is not clear how the upper caste are being deprived of their due. That is the problem is not their share of cake but the real problem is that of the very small size of the cake…

Caste is a very useful criterion for several reasons. One, the original discrimination is access to education took place on the basis of caste; the same criterion needs to be used for reversing that discrimination. Two, caste is still a very good proxy for various kinds of social and educational disadvantages and the single best predictor of educational opportunities. Three, caste and economic hierarchy tend to fuse at the upper and lower end; the poor are likely to be lower caste and the upper caste likely to be well to do. Finally caste certificates tend to be more reliable than other proofs of disadvantage especially and notoriously unreliable certificates of income…

Several other arguments as that by the Brahmins are projected. We give the views of the expert.

This article is a welcome one, for at least some percentage of truth is brought out. We first want to ask did any of the OBC's wage a war like the Brahmins when the reservation for SC/ST was given. When we the majority that too sons of the soil are so tolerant and are taking up the policy of and encouraging the policy the Brahmins are fighting for the 27% percentage reservation which is only a paltry amount in comparison with their population strength argues the experts.

Apart from this nearly 33% of these upper castes are having even after giving 27% of reservation, 50% to enjoy. Our politicians and socio scientists and policy makers are to be lawful and proper they must have given 33% reservation to these upper castes and at least 44% to OBC's. This 44% will be only an injustice in comparison with their population, for majority of the Muslim converts are only from the OBC and SC's. Strong laws must be laid down by the government to punish these antisocial elements. These striking medicos are like terrorist or anti-nationals for they are denying medical facility to poor patients most of whom are non Brahmins. For no one knows the seriousness of the patient health conditions. Denial of medical attention to patient is against their medical ethics. Stern action must be taken for their acts of boycotting their work for 19 days continuously.

The point Yogendra Yadav brings out in this article that, the original discrimination in access to education took place on the basis of caste, the same criterion needs to be used for reversing that discrimination. India can be a nation without reservation based on caste if there was no caste. But this cannot be even a dream in India by the twice born. Thus reservation is the only means that can render some justice to the majority feels the experts. Their protests are anti social, anti democratic and anti nation for they are not for making India a developed nation!

May 27, two youths on Saturday attempted self – immolation as medicos opposing the 27 percent quota for OBC

in higher education intensified their stir, staging massive rallies in the national capital and other places.

Rishi a Gutka vendor in East Delhi and Surendra Mohanty a post graduate in SCB medical college hospital in cuttack tried to immolate themselves as medicos here decided to continue their, over two week long stir. Rishi 23 suffered 30 to 40 percent burns on his arms and the chest when he set himself afire during the Delhi aao desh bachao (come to Delhi save the nation) rally at the Ram Lila grounds in Delhi. He was admitted to the Lok Nayak Jayaprakash Narain Hospital.

Mohanty attempted self immolation but did not suffer burn injuries. He was shifted to the ICU after his condition appeared to deteriorate; authorities were investigating whether he had consumed any poisonous substance….

"Anti – quota fire rages on, 2 try to burn in it". They said first govt. must keep an enquiry about how the gutka vendor Rishi was burnt or went for self immolation on his own? Mohanta did not suffer burn injuries thro! Was Rishi given a large sum of money to act so?

If these things would have taken place by the OBC trying to burn a Brahmin man or cajoling him to self immolate what would be the role of media in projecting this? Why is the media avoiding to discuss and to publish news about this matter. An open enquiry must be made and the culprit behind this must be punished appropriately.

They said the Brahmins only do not want "Caste based reservation" this very clearly shows still they want to cling on to the caste mentioned in the laws of Manu. Experts wondered when caste was so much important to them and they do not want to say, "No caste system in India but wants no caste reservation"!

Express news service

Members of the National Knowledge Commission (NKC) the country's top notch think tank are divided on the thorny issue of reservation in central educational institutions with six of the eight members favouring the status quo. Two members are however in favour of reservations. The commission headed by Sam Pitroda met over the last three days to discuss several issues pertaining to education including reservation for Other

Backward Communities (OBC) in central educational Institutions proposed by Union Human Resources Development minister. In a statement issued to the media Pitroda said, “The NKC firmly believes that a knowledge society must be a socially inclusive society and this social inclusion must be reflected in the educational institutions”….

The experts are very much pained by the media’s expression the NKC is divided on the “thorny issue of reservation in central educational institutions”. Is the govt. helping the majority who are lying as helpless last men to study in central educational institutions a thorny issue! Is it a thorny issue for the media or for the members of the NKC asks the experts? So the social equality of the last man by giving reservation in central educational institution has become a thorny issue!

The term “socially inclusive society” was highly condemned by many experts. Few openly said, the major wealth of the nation is in the hands of Sankara Mutt’s Narayana moorthys, TTK’s, and Ambanies. How can OBC’s reservation be based on any form, the NKC’s seek? Majority of the OBC’s are not only deprived of good education and good or well paid jobs above all they are socially very backward they added.

**An editorial by BibekDebroy** Bidebroy@phdccimail.com of The New Indian Express …

When we think of higher education, we must recognize that more than 80,000 Indian students have headed for this US and many more else where. For obvious reasons these students come from richer sections of society. We can’t insist that US embassies have an OBC quota for student visa. Indian institutions are setting up shop abroad and not all of them require permission from HRD.

We can’t insist they have an OBC quota when they function in Dubai… Indian Institutions of higher education have a comparative advantage. Do we want them to set up shop abroad instead of in India? Through WTO negotiations, higher education will eventually be opened up. Don’t we want domestic providers to be equipped to handle the competition … By no stretch of imagination higher education, a public good … Consider for example data on entry into IIT’s which should not be too difficult to process, since it already exist. Let us look at cross classification of entrants based on their caste backgrounds and as functions of whether their schooling was in government or in private schools. I am prepared to bet that deprivation and inequality of access will be higher for government school. If we are rational that should tell us something not only about correlates of backwardness but also about how public policy should be formulated. Instead we will be preoccupied with quotas in public institutions of higher education and eventually when we realize this has not worked we will want to extend such quotas to private institutions says Bibek Debroy.

We made a set of experts who were socio scientist and educationalist to comment on this article. They put the following points in response to Bibek.

The importation of 80,000 Indian students to US and many more elsewhere is waste of investment on the part of India, so if reservations are done in higher education that too to the backward classes it would certainly be a benefit to the nation for they would by no means seek to leave to US or other countries for they would more be interested to settle in their own states, there by economy and brain drain of Indians can be totally saved to a very large extent.

No one demanded OBC quota for student visa. This sort of statements only show the arrogance and their basic superiority based on caste and always a thought of treating other castes inferior to them and nothing more. This sort of statement show their true colour and caste arrogance, dominated by intolerant selfish attitude towards the nation development.

They are trying to dilute the issue by taking about the standards of government schools which is wrong. For B.Tech courses in IITs, is clearly not a higher education, but just a graduate professional course only. If the avenues are opened to them certainly India can produce better quality trained people than before says the experts. For the IIT's entry point is nothing but training and not intelligence or any form of quality the experts claim. They (Aryans) have money to pay the coaching institutions and get entry.

Express news service

"Consult IIT, IIM directors says Narayana Murthy".

This is unwanted for who are the directions of IITs and IIMs asks the experts, after all they are appointed only by the govt. so why should they be consulted; as a govt. servant they should obey the orders of the govt.

As he has earned fat economy in the name of Infosys Chief Mentor and Chairman of IIM Ahmedabad, he has to be true to his caste!

The Govt. must constitute a committee to study the work done by all the IITs. This alone can give the true picture of the work done by them.

He also tried to dilute the issue by asking the government to improve the school standards.

These experts ask why only as soon as the government spoke of ‘OBC quota’ i.e. 27 percent reservation for OBC’s they are speaking about the improvement of the school standard. Why they have not in these over 60 years after independence talked about improving the standards of the school? What does this really imply? This is yet another mischief of the media!

The Chairman of IIM (A) talks about consulting IIT and IIM directors about increase in the number of seats from Bangalore on 19 May; telepathically on the same day from Ahmedabad 19 May director of IIM (A) Dholakia says “increase seats at this juncture as its plate is full”. “IIM-A in a massive expansion plan has increased its seats from 180 in 2002 to 400 by 2007. At this juncture our plate is totally full” he told press reporters here.

The New Indian Express

There are signs that the reservation crisis is getting defused with the formula mooted by Pranab Mukherjee of 27 percent reservation for OBC in premier institutions in a phased manner and at the same time increasing the number of seats to placate the agitating students.

If one is honest the whole episode has been triggered not out of love for OBC or a desire for their welfare via reservation route, for which a case could be made…

Mandal II is a story of Arjun Singh getting the better of Manmohan Singh. As a political ploy it could be called masterly. He moved once his Rajya Sabha seat was secured. The story goes that Sonia Gandhi had expressed dissatisfaction at the way certain education schemes were being run. In a preemptive strike he has played a card which has made it very difficult for Manmohan Singh to dump him or change his portfolio without inviting the charge of being anti OBC.

What ever be the rhetoric around, Mandal II to justify to what OBC's have seen has been a cynical play of power politics at work, without little concern for consequences for the education system or for students, OBC or others. It is one more example of how little our politicians care for anyone else, or the extent to which they can go to serve their own ends. It also shows a further hardening of public cynicism about politicians as a breed, which is needless to say dangerous for democracy.

The news analysis by Neerja Chowdhury is one which blames the politicians for implementing or trying to implement 27 percent reservation for OBC in premier institutions. For this she blames Arjun Singh and projects as if Manmohan Singh is helpless. The author has clearly brought the stand of Arjun Singh, he being an upper caste man never wanted reservations to OBC, in fact wanted not only save his position but aspired to become the P.M. In the analysis the author criticizes the mean acts of Arjun Singh for his double stand on OBC reservation which is a danger to democracy. Is anti reservation slogan more dangerous to the nation for it affect a majority's right of getting educated in an institute run by their elected govt. asks the experts?

The news analysis hints at the poor state of affairs in higher education under Arjun Singh, claimed a group of the experts. Reservation to OBC according to the experts is more pro-democratic than antidemocratic. Thus atleast after 6 decades of independence the experts feel that Congress has acted in a little democratic way and clearly Arjun Singh cannot be given any credit for it, said some experts.

Express services

Taking an unscheduled suo motto notice of reports of patients stranded in hospitals across the country, the Supreme Court today sent a stern warning to doctors who are on strike, protesting the quota bill saying the agitation might be construed as "Contempt of Court". "Patients are suffering; protests could be contempt of court" – warning to doctors.

The case came up for its first hearing yesterday with two PIL's are by supreme Court Advocate Asoka Kumar Thakur and the other by Shiv Khera on which the bench issued notices to the centre saying that the "(OBC quota) policy is capable of dividing the society on caste basis and as a matter of fact it has to be adjudicated upon."

It gave centre eight weeks to reply and another six weeks for the petitioners to file rejoinders if any".

The experts reacted to this in a very different way. They (Brahmins) may feel OBC quota policy is capable of dividing the society but we warn the nation, the non implementation may leave the nation in constant terror and fear!

The experts condemned the court for admitting the very case filed by them for it amounted as anti social, anti democratic and antidevelopment of the nation.

Only, when the majority gets benefit by any policy put forth by the nation and if a particular minority section is against it, it cannot be even relooked. The only way to solve the problem is punish sternly the people who rise anti reservation quota slogans; feels the section of OBC lawyers and educationalist from Chennai.

The acts of these striking doctors are unrationalistic; they must be dismissed from their services for being anti govt. and anti nation. For when they are employed by the govt., they have no right to fight against a govt. policy which was passed in both the houses with over 99.5% (votes) as an amendment.

Some experts said the employee union leaders of IIT Madras one Ramadoss and Gajendran were dismissed from their services for striking against the administration of IIT (M). Their appeal in the court was also dismissed. So if law is the same for one and all i.e. govt. employees, these doctors must be dismissed as they are protesting against the very govt. who is their employees. Thus the suspension of the students and dismissal of the striking doctors is the only means to solve the problem says the experts, if law is the same for one and all.

The New Indian Express

Vice-Chairman of the National knowledge Commission (NKC) P.M., Bhargava has asked PM. Manmohan Singh to implement 27 percent OBC reservations from the next academic year 2007-2008 to give time and opportunity to the educational institution to upgrade, their infrastructure.

Bhargava, a biotechnologist was one of the two members to support reservation the high profile NKC… The experts feel the media is doing its best role in the spread of its anti reservation policy. It has tried its level best to project all possible points in all possible angles and shades. Not only that through different personalities. It has tried hard to take us the social status as the next big issue. It is also trying hard to win the legal side so that its defeat in the political side can be compensated feels the experts.

Express news.

Health Minister Anbumani Ramadoss was checking patients at the All India Institute of medical Sciences in New Delhi. However he (Anbumani Ramadoss) clarified that the decision to reserve 27 percent seats for OBC in higher educational institutions had been taken by parliament and the UPA government, was committed to it.

End strike, there would not be any punitive action-SC. Observing that the damage done to the patients was 'irretrievable' the court ordered the medicos to call off forth with all strikes protests demonstrations or any form of dissent.

"If in the meantime service of any doctor has been terminated for participating in the strike, it should be recalled by the government the bench said, adding that the concerned doctors should be granted three days time to rejoin.

Several of the experts said they were very much surprised to see the soft approach of the SC. If the strike was for pro-quota by the OBC what would have been our plight-will ever the SC be so soft; questioned several experts. We would have been ruthlessly dismissed from services and an enquiry would have been set up against us added the experts. This was no imagination, such things have happened in IIT madras.

The New Indian Express

The state Assembly today passed a bill to provide reservation in private educational institutions in Tamil Nadu for BC, MBC, SC and ST.

Hither to 69 percent is in vogue in government and aided institutions only. Now it has been extended to the self financing institutions and deemed universities.

The Bill was introduced by higher Education Minister K Ponmudi in the House.

The Bill said after careful consideration of the population in various categories and the present state of their advancement in education the state government has taken a policy decision to extend the existing level of 69 percent reservation in admissions to educational institutions other than the minority educational institutions for ensuring advancement of minorities in the state. The Bill was passed by voice vote without discussion!

Most of the experts were struck with wonder when the Brahmins are not striking against this 69% quota in private educational institutions. They said the Brahmins know fully well that their institutions would be permanently closed in T.N.

Further not even a single protest was made. The Brahmin Association of T.N is lost in this situation. Thus when the govt. makes a more to up lift the nation majority all these protests and strikes are unwarranted added the experts. A few experts were very much upset by the way the central govt. was tolerating the anti quota protestors. In this the role of media is very significant said the experts.

They wondered why only majority of the protest was made only by the docs and students of AIIMS. They all felt the anti-quota protestors must be treated and punished as antisocial elements.

The New Sunday Express

Human Resources development minister Arjun Singh rubbished suggestions that his drive on reservations for OBC was intended by a desire to embarrass P.M Manmohan Singh of whom he was “jealous”. This is a canard which is below contempt.

Only that person who does not know what kind of respect and regard I hold for Sonia Gandhi could say. She is a leader and whatever she decides is acceptable for me”, he told Karan Thapar on Devil’s Advocate Programme on CNN-IBN.

He was replying to a question whether he was jealous of the PM as the congress president had picked up Manmohan Singh for the top job.

Replying to a question whether the P.M was in the know of things he said, “I think a very motivated propaganda was on this issue. Providing reservation to OBC was in the public domain right from December 2005 when the Parliament passed the enabling resolution”.

The P.M he said was fully aware of it and he had met the OBC leaders to convince, them to support this legislation. How can you say he was unaware; he asked.

Express news.

In the face of stepped up protests by students from across the country against the proposed quota for OBC in top educational institutions, P.M. Manmohan Singh on Saturday urged striking medicos to call off them agitation and to have faith in the govt.’s “sincerity”.

The P.M made the request within hours of saying that his government was committed to the betterment of the backward classes and minorities.

The striking doctors demanded no penal action against the striking doctors and publication of a white paper by the government making its stand on reservation.

Toughening its stand the health ministry warned that stringent action would be taken against medicos who continue to abstain from work. Advertisement appeared in the press for recruitment of doctors in Central Govt. hospitals Officials said walk in interviews would begin on Monday.

Meanwhile the group of ministers created to deal with the issue on Saturday forwarded its report to the P.M favouring implementation of the quota for OBC’s as soon as possible and recommending an increase in seats in elite educational institutions – Agencies.

Experts wonder why only medical students are raising their voice; infact anti quota protests are fully absent from other stream of education. Most of the experts felt some sort of closed fraud must be practiced in AIIMS. Only to hide it and safe guard it they are making such protests they claimed!

Further youth for equality is more based by BJP/RSS Brahmins and it is a very small group trying to fight a very undemocratic issue. This itself will show how socially aware is the youth for equality group. The experts very strongly came upon the group and asked them first to change their name for they are not youth for equality but only youth for the inequality with BJP/RSS brain behind and functioning for the protection of Brahmins, they claimed.

These experts said they are going to launch a organization called youth for equal representation and are going to work for reservation for OBC not for 27% but for 55%. They (experts) all said a few women who claimed themselves as medicos alone made anti-quota slogans and the print media especially the leading news paper never failed to project it in the first page with photos and big captions. This mainly makes one feel that is only medicos women students be it Delhi or Gawahti or Jalandhar or AIIMs protests alone is shown! What is the myth behind this!

Express news

Both IIM (C) and IIT Kharagpur have decided to increase the number of seats. IIM-c have been increased the seats from 270 to 330, while IIT-Kharagpur would take 2505 students in UG and PG courses this year as compared to 1862 last year.

When these two institutions have the mindset what hinders others is the question proposed by some of the experts.

The New Indian Express

"Under pressure from allies like the NCP, the Maharashtra Government today decided to restore reservation and would be 13 percent for SC seven percent for ST and 30 percent for OBC.

Apprehending that the govt. may dilute the provision of reservations for SC, ST and OBC, a pro-quota of medicos demanded the implementation of reservation in strict terms and the P.M make the governments stand clear on the matter. Through the course of anti reservation stir the Congress party has been shifting its stance frequently with different ministers adopting different stances in public. The medicos Forum for Equal opportunities said;

"threatening an all India coordinated movement of the pro quota doctors and progressive forces in the coming week the forum demanded implementation of quota for SC, ST and OBC in the strictest possible sense of the word" and equal opportunities for all in quality education PTI.

This shows their Aryan character.

Express news service

Even as anti reservation movement is gathering steam in some parts of the country, the Tamil Nadu Indian Medical Association (IMA) and the "Doctors Association for Social Equality"(DASE) have come together to launch a 'Campaign Committee for Reservation' to support the centre's move to provide 27 percent reservation for other Backward classes (OBC) in Central Educational Institutions.

Apart from the doctors organizations which support the controversial move of the union Human Resources Development Ministry, all other forces and organizations supporting social justice would be brought into the committee G.R. Ravindranath, General Secretary of DASE said.

Strongly condemning the strike call on the reservation issue, IMA Tamil Nadu, DASE, Physician for Peace Health India Foundation (PPHIF) and other supporting organization alleged that IMA Delhi had taken an autocratic decision and acted against social justice. Whatever be the view of IMA Delhi on reservation quota issue, we do not subscribe to that. Tamil Nadu IMA whole-heartedly supports Union HRD Minister.

The strike call was instigated by the upper caste medicos but a majority of doctors, else where in the country were in favour of the new reservation quota proposals they said.

C.N. Deivanayagam, President of PPHIF said the claims of anti reservation lobby regarding deterioration of quality of education due to reservation were absolutely unfounded. The concept of high ranking students becoming better professionals is a myth when we go through the statistics we can get enough and more examples to prove this theory wrong. The upper caste medicos in Delhi are trying to dictate terms and they do not have any right to stall the work. "They are violating the basic code of conduct for doctors and also the medical ethics. They should be put behind bars" he said.

He also hailed the decision of the Mumbai doctors who were arrested by the police during the protest against the Centre’s move to go for a silent protest rally. “They have the right to register their protests but not at the cost of the patients he said.

The experts observed in this article the terms like “… support the controversial move of the Union Human Resources Development Ministry….” “… anti reservation steam is gathering stream” are needless for the media which should be unbiased in presentation of news and views. Further most of the experts welcome the strong view of Prof. C.N. Deivanayagam, President of PPHIF that the protesting medicos must be put behind the bars. The experts blamed the centre for their soft dealing with the striking medicos.

Express news

“Why have Parliament if we can’t honour our commitment on quota?” A day after the cabinet committee on political affairs struck a conciliatory note on the contentions quota issue HRD minister Arjun Singh tied the Manmohan Singh Government in another knot by questioning the very existence of Parliament if it could not honour its “commitment” over reservation.

It stark contrast to the debate raging outside members of the Lok Sabha demonstrated complete unanimity in favour of reservation for OBC in the course of a discussion during question hour today.

Arjun Singh said that the commitment to OBC reservation is not of government alone. This commitment has been expressed by our Parliament. What is the justification of parliament’s existence if we are not able to honour that commitment?

Arjun Singh said the commitment made to SC, ST and OBC by the House for reservations in private unaided institution, is irrevocable and this govt. and parliament are committed to that. A proposal to implement that in a …. I, as a servant of this House have been entrusted to see that it is implemented…

Express news

The University Grants Commission said there was a provision for 10 percent hike in seats every year under the 10th Five year plan. “Going by the funds sanctioned to colleges and universities, we have a provision for 10 percent increase in seats annual in all educational institutions" UGC chairperson Sukhdev Thorat told reporters here.

The New Indian Express

Express news. CII against reservation. Apex industry chambers confederation of Indian Industry (CII) today strongly objected to any mandatory reservation for socially underprivileged in the private sector even as it forecast 8.3 percent GDP growth target for this financial year.

Mandatory reservation in any form is not conducive to competitiveness of the industry. It is not acceptable R. Seshasayee the new president of CII said in his inaugural press conference.

However even as he objected to any mandatory reservation policy for S.C and S.T he said industry needed to take positive action to empower the backward classes to join the mainstream.

Experts condemned the new president of CII R.Seshasayee views about reservation as inhuman and unfit to be a president. They further added he should be sacked and a non-brahmin must be made as president.

The New Indian Express

Even as the protest against reservation at premier institutions is continuing in North India the pro-reservation movement is gaining strength in Chennai.

At a human chain organized by Doctors Association for social equality, political parties cutting across ideological policies came together near the memorial hall to join hands demanding the promulgation of an ordinance for reservation. Senior leaders including CPI secretary, D Pandian and senior leader Nallakannu, Dravidar Kazhagam President K. Veeramani, DPI General Secretary Thol Thirumanavalan, Puthia Thamilzhagam leader K. Krishnasamy, Tamil Nationalist leader P. Nedumaran, General Secretary of Doctors Association for Social equality G.R. Ravindranath and cadres of various political parties, student bodies and organizations including medicos participated in the agitation. The gathering demanded that the centre implement the reservation without any dilution and refused compromise on this. Interestingly, patented enemy of the left and other minorities, US secretary of State Condleeza

Rice, received all round applause when it was mentioned that she had come forward in life with the help of affirmative action that the United States follows.

Unions extended support; the BSNL, other backward class and postal / RMS / MNS. OBC employee welfare association also insisted that the central government implemented 27 percent reservation to OBC in all central education institutions.

In a demonstration at the main entrance of General Post Office on Rajaji Salai on Tuesday, demonstrators warned that if the centre failed to implement it, demonstrations would spread all over the country.

The ICF/OBC welfare Association praised Union HRD Minister Arjun Singh and former Prime Minister V.P. Singh for taking the initiative. They were holding the demonstration flaying those who opposed reservations. The experts still maintained that the media has played the major havoc. They are noted with pride no politician waged his tongue against reservation for the fear of his very existence as a politician depended on the OBC and not on the 3% Brahmins who were agitating and getting the Brahmin media 100% support; as if reservation for the majority of the people of the nation as anti democratic.

The New Indian Express

Offering to talk to striking medicos, Prime Minister Manmohan Singh to night said government was working towards an “amicable settlement” on the reservation issue and asked them to call off their agitation.

“I have no hesitation in talking with any one, so many students are in agony something that pains me” he told reporters on the side lines of the function to mark the completion of two years of UPA in office. The government, the Prime minister said was convinced of the need to have a fair, just and inclusive education system. Noting that the country’s youth was the greatest asset, the P.M. said” We will pay special attention to the needs of the socially backward classes where ensuring that no deserving student is denied an opportunity to secure education”. The striking medicos responded by saying they have taken the appeal in a positive manner … (Picture caption A medical student reads a book even as she takes part in the on going hunger strike, held as part of their anti quota protest in the AIIMS campus in New Delhi).

Most of the experts felt hurt for the P.M pleading them to stop their strike. The P.M ought not to beg them when suggesting an implementation for the OBC's which is their birth right to get education or atleast get a chance to read/educate in a institution that is run by the govt. When these institutions are run by the tax payer money denial of opportunity to the majority who form the OBC class is undemocratic. More so is the strike done by the 3% minority have been enjoying such good education all these 60 years! Who is going to account for the millions of students who have been denied opportunity for many years.

Express news

"Never has the media been so one sided in the post – independent era Ramadoss charged".

The party's founder leader Dr. S. Ramadoss who had already met leaders of National Parties in New Delhi to press for immediate implementation of 27 percent reservation for OBC told journalists here on Monday that he would take up the issue again at the meeting of UPA leaders in the capital on May 23. Besides the PMK leader would visit Maharashtra, Karnataka and Andhra Pradesh to garner more support for reservation.

Ramadoss said, he met Karunanidhi again on Monday in this regard, "Karunanidhi should do the role of Periyar in this issue", he added.

Claiming that 2 percent of the population was funding 70 percent reservation in Central institutions; he said, "This is one battle we will win.

Criticizing English news channels for repeatedly telecasting the same, footage, Ramadoss said that the channels were trying to make people believe that the entire nation was revolting over the issue. "Never has the media been so one sided in the post independence era" he charged.

Unwilling to buy the theory that merit was at stake, Ramadoss armed with statistics of MBBS admissions in the state in 2004, said that students of backward classes, Most Backward Classes and Scheduled Castes were the rank holders and not students belonging to Forward Castes.

As many as 77 percent admitted to 12 medical colleges in the state were from BC, MBC and SC/ST categories, he said. Among rank holders only 6 were from the forward community on the first 10. There is no truth in the claim that meritocracy will suffer”.

Several of the experts praised Dr. Ramadoss for bringing out the media’s false colour in the reservation policy. They all said the media was overdoing anti reservation protests and no reporting or giving single line news about pro reservation activities or even meetings. They feel media which is to be unbiased is full of not only bias but is falsely propagating the anti democratic acts. Except for a handful of women medicos that too only from the north no one in India is protesting against reservation.

If this is the case they (experts) are not able to comprehend why the media is reporting falsely.

Express news

The government is putting together material that would form the basis for arriving at the decision of 27 percent reservation for OBC in higher education, even as it ruled out plans to review the quota policy.

The ministry has to put together all the material available, to reach conclusions of 27 percent … The material will be put together; wait for the material” Finance Minister P. Chidambaram said in an interview to Karan Thapar on the programme “Devils advocate” to be telecast on CNN-IBN on Sunday night.

The Court on May 29 gave the Centre eight weeks to reply to a series of queries on reservation.

When repeatedly asked that his answer suggested the government had announced a decision in favour of reservation with out evaluating the need he said, “Government will answer in a proper forum in Parliament in the Supreme Court (but) not in an interview”.

He however said the government would compile its answers on the basis of various state reports of Other Backward Classes (OBC), the Mandal report and NSSO surveys. Chidambaram who is a member of the group of ministers looking into reservations for OBC in higher education said the government has no intention of reviewing the caste-based reservation policy.

"As I understand … There is no ground to review whether there should be reservation or not. There is no ground at all… If a review means questioning the justification of reservation I say no" he said, Union Minister for Statistics and Programme Implementation.

Experts could feel some thing fishy on the part of the students demanding to set up a non political commission to review the existing reservation policy.

This very clearly proves one and all that no political party or no politician has courage in India to say no to reservation or make any amendment on the bills which was passed in over whelming majority in both the houses. That is why they wanted a non political commission to be set up.

As usual two line information was projected about a protest made in Mumbai on June 10 by a group of pro-quota activists on Saturday.

The experts feel that when it is a pro-quota strike the media never wish to cover it properly. They just reluctantly report it in a line or two as if their protest is very frivolous!

Express news

The indefinite strike on the controversial OBC reservation continued for the tenth day today with the Indian Medical Association (IMA) splitting as pro-reservationists formed a parallel Indian National Medical Association (INMA).

The Indian Association of Tour Operators (IATO) sought the president intervention to break the impasse.

The INMA comprising OBC, SC/ST, IMA members alleged that the Indian Medical Association with the fat backing of a political party (BJP/RSSVHP) were managing the medical strike for political mileage.

The INMA alleged that the IMA had been hijacked by upper caste people and before pledging support to the anti reservation agitation they had not sought the opinion of the backward class members.

The PMK a constituent of the UPA today staged a demonstration at Jantar Mantar here protesting against the anti reservationists.

The protestors led by former minister and MP of Chidambaram E Ponnuswamy urged the centre to immediately promulgate an ordinance implementing 27 percent reservation for OBC in the higher educational institutions – ENS.

When in Delhi itself the PMK a constituent of the UPA stage a demonstration in the prime location like Janatar Mantar the media coverage is so poor it reports in a line or too. Participants were not simple students or doctors but a former minister and M.P still the media suppresses that news. The experts heavily condemned the partial or biased attitude of the Brahmin monopolized media for their undemocratic and arrogant way of functioning.

The New Indian Express

Commission needs to focus on quality. The six-two split of the knowledge commission gives holistic depth to the current debate on mandal. II. In Bangalore this week a majority of the members of the commission came out in sharp misgiving about the centres proposal to extend the benefits of reservation in central educational institution to OBC. The commission was set up recently by the prime minister to advise him directly on educational reform and challenges in meeting its knowledge potential. The first measure of its institutional relevance will be taken now. The debate on reserving Seats for OBC is much larger than the question of quotas. It is a pointer, whether Arjun Singh intended to be so or not to India's failure in making education socially inclusive. To its credit the knowledge commission has lashed out this fundamental agenda from the bare bones of regrettably narrow and obviously political agenda set in motion by the human resource development minister.

It is in this context that the minority views in the commission gains significance. Even as its vice chairman P.M. Bhargava affirmed commitment to quotas he too emphasized the need to broad base access. Thus he saw as, a necessary condition in extending reservations for OBC in higher education that would require the government to set up 4,00,000 high schools and desists from cutting back on the number of candidates chosen on the basis of 'merit'.

In the transparency of the discussions the knowledge commission has been wholesomely responsive. Given its remit the onus would entirely be on the Central Government and more specifically, the prime minister to whose office the commission directly reports-to take forward the debate. But commission is gauged not by the wise words they record, but by the effect of their deliberations. A good example of this is the kelkar task force on fiscal responsibility; in budget management. We therefore await both imaginative and pragmatic suggestions by the knowledge commission on giving higher education in the country; an excellence edge and in making it more accessible.

Experts were shocked by the Brahmin arrogance of the knowledge commission, to say falsely that it would affect merit they condemned the P.M for making such a caste fanatic knowledge commission for a country where the representation of the majority was nil. They stamped the findings of the knowledge commission as akin to laws of Manu and whole argument had no merit but only biased rubbish from the sewage tank. The experts were full of anger and sharp criticism on the Brahmin commission throughout the discussion.

At the first out set they want to educate both the Brahmins and the politicians that it is not higher education just a graduate degree that too when both, the so called meritorious well placed Brahmin students and the poor rural students take up the same syllabus and the same question paper. This meritorious getting over 90% is not any merit or quality contribution to the nation as seen from the past.

They challenged if the same comfort and encouragements is given to our students they would easily get hundred percent. Why they are the blood succors of our nation. The experts further they said an ex-director of the IITs, Dr. (Prof) Anantha Krishnan also a veteran educationalist says the merit or quality as the Brahmins say is never at stake. The experts ask what have these Brahmin done in these 60 years of independence except have made, the majority poor. One can see the steady increase in the majority problem even to get basic comforts. On the contrary all Brahmins are not only holding good/well paid power wielding office and more than 60% of the nations economy is in the hands of the 3% Brahmin.

The expert asks, why now they feel the need for good high schools? What did they do to the OBC all these 60 years? They using their power over the politicians have ruined us because we were submissive! Now denying seats in govt. run institution in the name of quality will by no means be accepted. If they feel India would loose quality let them quit India and go and settle themselves in nations with high quality. They also sternly warned the media for its biased reporting! They have the guts to say OBC reservation would lead to India's failure in making education socially inclusive. Some of the experts said that they are contented to live in a failed India if reservation for the OBC is really going to be so; and they would see to it the govt. will implement it soon.

They even after 60 years of independence using words like merit excellence and so on, try to fool the majority. Even the very thought they would be denied some seats in the IITs, IIMs and AIIMS makes them act like anti socials and terrorists, claims the experts. They think we are fools because we have not refuted them and we have tolerated them in spite of our strength and merit.

The New Indian Express

The Government would be forced to bring –in a legislation if private sectors voluntarily who do not provide reservation for weaker sections in society, minister for Social Justice and Empowerment Meira kumar said here on Monday.

She referred to the two year time frame in the UPA common minimum Programme for providing reservations in private sectors and some more time could be given to them for taking steps in this direction.

But at the end if they don't do it voluntarily the govt. would be compelled to bring a legislation to enforce reservation in private sectors"

Guest Column – P. Radhkrishnan:

The writer is a professor of Sociology at the Madras Institute of Development Studies.

Given the constitutional provisions for the uplift of the weaker section since the 1950s and the chaotic proliferation of literature on these provisions, the ignorance of most of these debating, the merits of reservation is appalling.

It is also appalling that among the self appointed educationalists are some retired Vice-chancellors who for the sake of remaining in the news project a different image before the media. They behave as if they are the constitutional pundits and saviours of the socially deprived but practice exactly the opposite in the institutions in which they have their sinecure. Here they do not follow even elementary rules of professional etiquette and democratic governance, yet sections of the media suck up to their unscrupulous characters and stream their photographs.

The constitutional provisions for job reservations for socially deprived groups is qualitatively different from the constitutional mandate that the state shall make special provisions for the advancement of the socially and educationally backward classes in particular the scheduled castes and scheduled Tribes. The latter is a comprehensive package of affirmative action.

The best example of an appraisal of this package is in the recent report of the constitution Review Committee where in the context of job reservations, attention is drawn to making use of Article 46.

The cacophony and confusion caused by vested interests of particular irresponsible politicians playing with the emotions and sentiments of the gullible masses have left hardly any scope for the return of reason and fairness in the public realm. Only this can help assess the causes, contexts and performances of any public policy.

Nobody has asked why the ministry of human resource development (MHRD) and the state governments have not come out with a status paper on affirmative action in general, and reservation in particular – say (a) the educational levels of the different social groups such as the scheduled tribes, OBC and others (b) how many applications the central universities and so called elite institutions IITs and IIMs have been getting every year (for analytical purposes take only the last two or three years) from different social groups (c) how many of the applicants get admission and how (d) whether there is any overt or covert attempt by any of these institutions to exclude any section of society (e) what kind of socio economic profile these institutions have in terms of their student teacher in take and so on.

No doubt Indian society is notorious for its discriminations and exclusionary practices. While the British administration for long used these practices as part of its patronage politics, in the second half of the 19th century it admitted, in the context of the then depressed classes that what existed was abstract justice. When its attempts admit the depressed classes children in public schools failed, following protests and boycotts by caste Hindus, it started separate schools for the depressed classes. They still exist in some states.

When the victims of past and persisting discrimination are not able to get even minimalist state support to equip themselves through education at least upto the pre university level, it is the state that should give an explanation. After all it is the states failure for foster talent and institutions and built new institutions and expand the education system commensurate with demand for education.

Reserving a few slots at the top without allowing and enabling the deprived groups to climb up even the lower rungs of the education ladder is political chicanery. If data available are any indication as of now India’s demand for higher education is by 35 percent of the relevant age group. Contrast this with the present enrolment of 9 to 11 percent compared to 45 to 85 percent in developed countries. Contrast this also with the fact that India’s present outturn of degree of holders is just about 7.5 percent.

As Muslims have only more or less the same population is close to the SC population, it is strange that nobody has said anything about their plight so far as the state have failed for 55 years in its first charge of providing basic education to India's unwashed millions and is not doing so even now, the distortion and dilution of a large and comprehensive affirmative action package to the political largesse of reservation in educational institutions cannot be anything but diversionary and a fig leaf to cover its own sloth and shoddiness.

So the issue for debate is not whether one is for or against reservation, but how to make the state accountable to the citizens and ensure quality education to the entire rising generation with special provision for the socially deprived.

Now the experts view this article in the following angles;

1) The point raised by Prof.P.Radhakrishnan about no body talking or helping a hand about the plight of Muslims is very sad. The main reason they attribute for it is that the Hindutva agenda. So they are anti-Muslims and anti-Christians hence any way only two years have passed since the congress which is also according to the expert semi hindutva. The reservation for the minority Muslims will have to take place only in due course of time. Only now slowly they are realizing that they are denied all opportunities in central educational institutions also. The truth about that is to be unraveled.

Other reason for Muslims not getting any reservation in India is after all they are mainly the converts of SC/ST or OBC/MBC who are the sons of the Indian soil. Like the Aryans the Muslims have not come to India from Arabia or Pakistan they are the original sons of the soil most of them from the deprived castes or Sudras so they are doubly stigmatized in the Indian land one by religion and other by their caste in India!

They are forgetting the hard truth that Muslims were the one who shouldered in the quit India movement to win our freedom. Now they form the neglected lot where as the 3% Aryans are dictating terms to the govt. Only now the OBC should fight for their reservation atleast as an act of loyalty to the sons of the soil feels the experts.

The experts felt P. Radhakrishnan seems to be biased over the reservation issue for he says the state has failed for 55 year in its first charge of … who is the root cause for the government's failure. This can only be attributed the 3% Brahmins who are acting as policy makers for India and who are very keen on stalling any progress in the part of the SC/ST/OBC in keeping with the laws of Manu. It is a pity that none come forward to say eradicate caste, based on birth. Do not follow Varnashra dharma, we will create a casteless India, but on the contrary they want caste based on birth and they (Aryan nomads) want to remain in all high posts and deny any form of reservation even in govt. run institutions. First let Radhakrishnan criticize himself before he speaks of the inefficiency of the state. Does it not speak in truth that the Brahmins who are the policy makers of India has been so biased to deny all forms of reservations to OBC. The brahmins at the first place should be ashamed to carry on protests against reservations. The experts blamed the media for publishing articles on anti-reservations and denying the publication of all articles which are pro reservation. Several of the experts said how the articles on pro reservation hardly could find any place in media. They also heavily blamed the media as it was only run by the Brahmins and it was only for their propaganda. So having the media in their hands they fish out for anti reservation articles and publish them surplusly there by trying to establish that the very notion of reservation is anti national and anti democratic.

The experts added IITs IIMs, AIIMS and some of the centrally run central educational institutions in India have been serving as brahmins strong forts so when any place / reservation is asked even by the govt. who is their authority they have the guts to strike, protest and ask govt. to review the very policy of reservation. The experts portray this as the limitless arrogance of the brahmins. They have the cunningness to state the very reservation policy as the vested interest of the irresponsible politicians. Who is irresponsible? The brahmins or the politicians let them answer? Who has vested interest the so called 3% Aryans have the vested interest on their development alone which has made them blind when they protest against reservation for the majority, deprived of good education and basic needs in life.

Some of the experts said they are vexed with the same way (method) they protest for reservation. They after 57 years of independence are side tracking the major issue of reservation in central govt. institutions by asking the govt. to built good schools or by saying the govt. schools are of poor standard and they should be improved. Why suddenly they are so much concerned about govt. schools and their standard now!

Express news service

Bharatiya Jan Shakthi (BJS) president and former Chief Minister of Madhya Pradesh Uma Bharati, on Tuesday called

for a national debate on the proposals of the centre to bring in reservations to the OBC in higher educational institutions.

Bharati who arrived here from Rameswaram said that she was for reservations in elite educational institutions as disparities between the rich and the poor and unemployment problem were posing a great threat to nation. "The educational system is poor in general across the country. As the disparities are intense, students cannot compete on merit" Bharati said and called upon the protesting medicos to withdraw the anti quota stir.

The New Indian Express

With the debate of reservation of jobs in private sector simmering, commerce minister Kamal Nath has said that having quotas in private jobs would not impact FDI flows into the country, which grew over 40 percent last fiscal to 7.5 billion dollars.

There would be no impact on FDI flows he told reporters on sidelines of CII annual session on Wednesday when asked about the effect of quotas for the socially underprivileged in private sector...

Express news

Agitating medical students have said they will launch nation wide protests next week in continuation with their on going indefinite strike even as the centre pleaded before the Supreme Court that it is not inclined to apply reservation for admission to the post graduate courses in All India Institute of Medical Science (AIIMS) during the current academic session 2006-2007...

The agitation against the proposal to provide 27 percent reservation to OBC in central Government educational institutions was getting support from students of several institutions.

The experts disposed of their claims that several students from other institutions were supporting them. The experts claim only AIIMS was agitating with the major exposition to their strike was well projected in the major leading dailies like The New Indian Express, The Hindu and the D.C. All of them are monopolized and owned by Brahmins, so naturally the anti-reservation agitation by a repeated handful of students have been blown up out of proportion by them added the experts. Further the experts said in truth why the AIIMS students were against reservations was that in their institution they had a very queer type of reservation viz. if 30 students did their graduate study in medical and if 70 students were to be admitted in P.G. courses all the 30 students were guaranteed of their seat to P.G. courses that is the right of the sons of the soil in P.G seats of AIIMS! This has already come to light! When they have enjoyed such type of nice reservation based on the right as the alumni will they ever accept OBC reservation! Now the sons of the Indian soil the OBC cannot claim any reservation in India but because they have studied in AIIMS even if they just get 60% in U.G. they are sure to get PG course in AIIMS.

Express news service

International trainer and world's largest selling, "you can win" book author Shiv Khera has strongly opposed caste based reservation system and favoured reservation on economic basis.

In this special address at the Rotary international district assembly held at Vellalar Arts College here on Sunday he dubbed the reservation issue a sensitive one and referred to the agitations in North India. People do not like casteism. But, the politicians nurture it for their own benefit and thus divide people he charged....; But the present caste-based quota would partition each and every village", he warned.

The experts said the Brahmin Shiv Khera Opposes caste based quota but he is not against the existing caste based society. Why every brahmin accepts caste based society as per Laws of Manu, the cruelty of caste, based on birth but no brahmin wants caste based quota?, asks the experts.

The New Indian Express

"Unfazed by protests over the issue of OBC quota Human Resources Development Minister Arjun Singh today rejected suggestions for a 'relook' into the proposals for reservations in the elite educational institutions and said it was up to Prime Minister Manmohan Singh to decide when to bring the matter to the cabinet. Singh projecting himself as the cause however left it to the Union Cabinet and Parliament to take a decision on the matter.

"This is entirely for the cabinet to decide; as soon as the Prime Minister allows it to come in the cabinet a decision will be taken "he told reporters who wanted to know a bill in this regard would come up in the ongoing brief session of parliament.

Slamming the National knowledge commission headed by Sam Pitroda majority of whose members had rejected the quota proposal, Singh said, well with all due respect to the Great knowledge Commission I must point out to them that they are not above the constitution".

The students submitted a memorandum to the governor saying the proposed 49.5 percent reservations in IIMS, IITs and AIIMS would "completely destroy these institutes that are centres of excellence".

Minorities front warns and medicos to call off agitation.

Medicos forum for equal opportunities opposes agitation.

AIIMS cancels leave for all senior doctors.

The experts however praised the views of Arjun for his sternness and rejects, relook. They also added that he has put the great knowledge commission in its place by asserting that they are not above the constitution, however blamed the media for projecting the protests of handful of medical students and junior doctors. They wanted the govt. to probe into the striking junior doctors and medical students to find out the true brain behind these activities. They also suggested that these medical students must be suspended for a year for such heinous antisocial activities. Such punishments alone can make the other students community not to venture in such anti social acts against the govt. and the constitution. The experts, majority of whom were socio scientists claimed that when the students are irrational at such young age how would they be as grown ups! What kind of contribution to the society will they make? Can it ever be social equality? Even a child of 10 years age, can spell out the differences of the working class and the high class and would say the working class must be given some sort of help! Why no rhyme or reason plays any role in these medicos? Only their caste and their security is upper most in their mind and they absolutely lack any sense of responsibility claims the experts!.

They all uniformly claimed, the very undemocratic role played by the brahmin media is the root cause of all confusion, they ask who are the students to say and submit a memorandum that reservations in IIMS, IITs and AIIMS would completely destroy these institutes which are centres of excellence"?

First let them come up on open and say the achievements they have made in these 60 years!

What do these students know about merit? Is getting good marks by studying in the good school in final exams, enjoying all types of comforts and extra tuitions an indication of merit? Let them reverse their situation with a poor rural boy studying in a govt. school. The pass by the poor boy is more than 90% by these so called meritorious students claims, a group of socio scientist for both of them answer the same set of question papers. Thus the merit of them claimed by them is farce!

By denying OBC/SC/ST admission even after 60 years of independence they have destroyed the democracy in India. All these years, they have done injustice to the majority under the unscrupulous terms like, merit, excellence and so on. This time the experts have challenged that they would not be fooled. We would not rather compromise for 27 percent reservation. We basically need at least 52 percent said some of the experts as 60 years is a very long period for which they have been cheated added a few experts.

The New Indian Express

Is this the beginning or the end? With the UPA government now deciding to introduce a 27% quota for OBC in all higher educational institutions, starting 2007, this atleast is what those in the forefront of the anti-quota stir are likely to believe. The rapidity with which the government came to this decision without engaging with the many counter-arguments advanced of escalating quotas may impel such a view.....

Nevertheless one might ask why the same state over so many years has failed to make a great success of its quota regime. After all if state supported institutions fail to fill their quotas, and over decades and face no penal action for not meeting their constitutional obligations surely some of the blame must fall on the state. Similarly what has come in the way of institutions running quality remedial programmes to help those from less privileged backgrounds compete and perform well? Or is it that once having instituted a regime of quotas interest in affirmative action wanes. Worse why is it that over all these years we still do not have authoritative studies on the working of different affirmative action schemes including quotas? The debate is not about whether caste matters or even the reality of social exclusion.

This move for dalits than the OBC's is undeniable.

The experts feel that the policy makers for dalits and OBC are only Brahmins and upper castes so how can any form of proper reformative steps be taken by them. We see even after 60 years of independence that the 27% OBC quotas in education institution run by govt. as per constitution amendment is ruthlessly protested by the selfish Brahmins! So how can any form of reformation to uplift the majority who are unfortunately non-brahmins can ever take place. Unless the policy makers for dalits and Sudras (OBC) are dalits and OBC we cannot succeed in any big reformation in their status. This is very clearly seen in by the way the National Knowledge Commission acted against reservation!

The New Indian Express

The government has spoken: 27 percent OBC reservations plus a commensurate rise in seats in institutions of higher learning. This is the facile quick fix response of a government in a hurry to create a lasting constituency for itself.

To expert an indifferent education bureaucracy to deliver on the promise of expanding the number of seats almost over night strains credulity…, The five top IIT's which take about 1700 students at the graduate level every year now have to suddenly increase their seats by 27 percent without even a blue print of how they are to go about doing so. Remember there will have to be a cumulative increase in the number of seats over the subsequent four or five years, because each student who finishes the first year's course will have to be accommodated in the second year course and so on … Even if we are to suspend reason and imagine that this mammoth task will be achieved it is sure to result in an overall deterioration of standards. It took nearly 30 years for the IITs to become institutions of excellence that they are today. Unless adequate time and resources are spent in expanding the existing infrastructure and recruiting the right teachers and staff, increase of seats will only lead to a sharp decline in the quality of higher learning in the country.

Ironically this cavalier attempt at social engineering is taking place at a time when the country can ill afford its inevitable consequences. A dilution of excellence a deterioration of standards, a lowering of commitment to educational outcomes would severely undermine India's prospects as a serious contender in the knowledge economy. It is the same old story: legislate in haste and repent at leisure.

This article was well tackled by our educationist and human work force experts as follows:

Once again the experts were concerned on the loose usage of terms like dilution of excellence, a deterioration of standards, a lowering of commitment to educational out comes, would severely undermine Indian's prospects as a serious contender in the knowledge economy.

The experts ask, where do the IIT-ians go after getting their degree? Is it not only brain drain and a shear waste of investment on these graduates? Atleast OBC will remain in India to serve India as they have more love for their soil, than the Brahmins who are just aliens to our land and who always seek prosperous place.

The New Indian Express

The Supreme Court today observed that the proposed 27 percent quota in all central educational institution is capable of dividing this great nation on caste basis and issued notices to the centre and six ministers asking them to answer what was the 'rationale' behind the classification of a class of citizens as "OBC".

Beginning the hearings on a PIL, a vacation Bench of Justice Arijit Pasayat and Justice Lokeshwar Singh Panta made it clear that "as a matter of fact it (the case) requires judicial review" and set three crucial questions on the union of India reply:

- What is the basis of the norms for fixing OBC quota.
- What is the rationale for this

- If the proposal is implemented what are the modalities of such implementation and the basis for the modalities?

The centre ministries of social welfare, HRD and others were given eight weeks to reply to the question and six weeks " there after" for the petitioners to give their" rejoinders if any".

The petitioners Supreme Court senior advocate Ashoka Kumar Thakur and Shiv Khera sought a direction questioning of the amendment to the constitution, Article 15(5) paving the way for enactment of a law for a further revision of 27 percent apart from the constitutionally sanctioned SC/ST reservation of 22.5 percent.

In a brief order, the judges said the "policy, if adopted will divide this great country on caste basis", and has severe social and political ramifications" the judges said there need not be a stay at this juncture as we have already put them (the respondents mainly the centre) on enough alter. In view of today's development as the Supreme Court is seized of the matter, the judges orally observed that there be no more agitations but we are not making it (the observation) a part of the order. Additional solicitor General Gopal Subramaniam present in the court accepted the notices on behalf of the centre and told the court that the constitutionality of the issue (the amendment) has to be debated to which the judges agreed. It's likely that a constitution bench will be formed to determine the issue after the notice and rejoinder periods. Both Thakur and Khera contended that the constitutional amendment should be struck down first as it has fathered and mothered the law.

"The SC/ST and other backward classes (Reservation of seats) in educational institutions act".

The experts were asked to comment on the statement of the SC; "OBC quota in all central educational institutions" is capable of dividing this great nation on caste basis:

There were nearly 40 experts who were socio scientists, psychologists, educationalists, lawyers, doctors, technically trained people, students and teachers'. The following were the major 4 questions put by the SC to govt.

- OBC quota can divide nation, tell us how you decided?

- What is the basis of the norms for fixing OBC quotas?
- What is the rationale for this?
- If the proposal is implemented what are the modalities of such implementation and the basis for the modalities.

At the first set, all of the experts refuted the statement OBC, quota can divide the nation! Already people are divided because of caste system by Varna of manu, so when the people are divided already "what does the statement OBC quota can divide nation"? mean. Any common OBC /SC man always feels the pangs of discrimination and of the Brahmins monopoly in all fields except in scavenging, they added. They are the policy makers for us feels the OBC/SC. The majority of the economy of India is in the hands of them!

All seats for any form of professional education in good colleges/ university are only occupied by them, by the fictitious words of merit, quality or excellence! When we look at the legal side not even 10% are OBC Judges. All are only Brahmins! Do not all this has already divided the people and hence the nation? Atleast if govt. by some amendments patches this and gives some form of representation / reservations for OBC the problem of insecurity of the majority can be averted.

They have been passive observers and never receivers of any benefit even after nearly 60 years of independence. It is high time some corrective action is carried out really before the nation faces any trouble. No logic and no scientist will ever accept in India a 3% Brahmins who have invaded India are occupying or to be more precise usurped all fields of importance treating the 97% population who are natives as people with no merit or talent. High time they give up all form of protests using the media which is doing the wrong job.

We experts ask, what is the power of the Brahmins to question the constitution? Is rule of SC powerful or is the constitution which passed the rule in over whelming majority superior? First let them be rational? We appeal to the blind folded irrational functioning of the SC's when it comes to reservation! It amounts to loose talking for "Thakur and Khera, contented that the constitutional amendment should be struck down first as it has fathered and mothered the law, "The SC/ST and other Backward classes (Reservation of Seats) in educational institutions Act".

The experts feel they are careless in their talks and above all they are so much agitated when one talks of reservation for OBC or quota for OBC so in their agitation they totally loose their reason and mind. That is why; they are least bothered about the majority and are so much worried about them.

It is high time that the OBC are given their due before it becomes a bigger divide of nation. For according to the experts only the 3% Brahmins will be the most affected lot for the OBC as a national policy accepted the 22.5% reservation for SC/ST. Now if the Brahmins protest against the OBC reservation what prevents the OBC to protest for their right! So it is high time the politicians do not go back on their statement.

All these protests by the Aryan brahmins was only more like an hysterical act in which, in majority of the cases a handful of brahmin women lead the protest. Most of the protests only came from the medical institutions, in the north recorded the experts and the Brahmin women protesting photos were captured in large sizes by the Brahmin media as the media is only in their monopoly! The Brahmin media must be first educated to give importance to all forms of news irrespective of caste and creed, should be the basic policy of any media. It is unfortunate in India only the news wished by the Brahmins alone is projected and not any other news or views. The media was so harsh and cunning to hide a pro-reservation protest for OBC.

Anti reservation protests alone was projected by media. It is a media myth which is a playing a havoc in the nation! They may say they are bothered about the great nation, but none of them ever say about the majority of the people right or equality! This is very unfortunate claims the experts.

Express news

With a debate on having quotas for the socially underprivileged raging across the country, the government is considering suggestions for protecting identity of SC/ST and OBC candidates appearing for Public Service Commission exams including UPSC.

The planning commission in a paper sent to the PMO, has suggested doing away with the current practice of issuing distinct roll numbers to SC/ST and OBC candidates appearing for such exams as it exposes them to discrimination by interviewers. These numbers would be replaced with a serial number common up to the interview level to guard against any discrimination.

Express news

The All India Minorities Front (AIMF) has threatened to take on the striking medicos on the streets, unless they called off their agitation against the quota hike.

Terming the protesting doctors and students as misguided the AIMF appealed to them to leave the path of confrontation and with draw their agitation unconditionally; otherwise the members of AIMF would be forced to take on them on the streets a statement made by the front on Sunday PTI.

The New Indian Express

**Dethrone royalty in academia,** by Anant Sudarshan:

HRD minister Arjun Singh's announcement on reservations in the IITs and IIMs to almost 50 percent has set in motion a debate as inconclusive as it is heated. Amidst the cacophony of voices, it is important to ask the question – how real is the problem for which reservation must be a solution? Or are they an instance of terrible policy born of populist polices?

To answer this let us examine some common criticisms or reservations. The first, most often heard is that reservations destroy merit images of incompetent doctors, useless engineers and inept doctors are trotted out as a dire warning of what might happen once anything other than an examination is allowed to determine admissions. A second argument is that quotas do not reach the truly poor and instead are gobbled up by a creamy layer. Again there is no dearth of anecdotal evidence about SC/ST student who drive up to college in fancy cars. Third objection made by a few, is that it is not necessary by reservations that are bad but reservation based on caste, that seats could be reserved for economically backward sections instead.

There are serious arguments are so persistent because to a large extent they seem individually true. As a recent student from one of the IITs in my personal experience it is true that SC/ST students (there are a few exceptions) tend to have done much worse in the JEE within IIT as well many of them have low grade point averages. It is also true that students getting in through quotas are often from the middle class families yet hidden within these common observations is an inherent contradictions and one that hints at a real problem in our society. If the academic performance of students belonging to certain sections of our society (the so called backward caste) remains poor even though they are not economically badly off and even though they are a subset selected through intense competition (for even the quota seats are very hard fought prizes) then surely there is a serious issue that needs to be addressed. The very fact that quota students are not being held back by economics alone points to caste as a constraining factor.

Sadly it seems true that for certain groups in our society, some professions and institutions remain unwelcoming and difficult not because of a lack of money but due to deeper more subtle biases. This is an entirely different problem from that of economic backwardness and it is this difference that lies at the heart of a rationale for reservations. In its essence this is a problem similar to the low representation of girls in these institutions. Across caste lines families provide more support to boys studying to get through hugely demanding entrance examinations. As a result women are a minority in many areas of study and often perform worse in examinations, not because they are less intelligent but because they have a very different support structure.

Reservations in our “institutions of excellence” will not uplift India’s poor or those in the lower castes (a thousand seats are a drop in the ocean). But they might lessen this artificial segregation. The argument that quotas dilute merit is seductive but flawed. Students benefiting from quotas dilute are often middle class with a reasonable school education. Also quota seats are fiercely contested and successful candidates are highly motivated. Unfortunately grades and marks are strongly dependent on social and family backgrounds and attitudes and do not necessarily reflect inherent ability. American universities have a less demanding style of examination than IIT's yet some of the world's best engineers and scientists emerge from them.

Indeed while a number of remarks are made about the poor performance of "backward" castes. It is interesting to note that IIT students admitted from outside India also do badly in internal examinations. Unless we plan to argue that SC/ST students US and French exchange students and Indian expatriates of different castes are all some how intrinsically lacking in merit, it seems for more likely that students sometimes have lower grades not because they have no ability but because strong cultural and social difference still divide society. It is these differences that a policy of reservation seeks to address.

There remain serious questions about whether quotas are the best approach. But caste divisions cannot be confused with general economic poverty. We need comprehensive studies and clear indices to measure the effectiveness of quotas. However as long as facts suggest an inherent class of royalty in India's knowledge economy, we need to work towards change. That process is painful and tumultuous but it is neither unnecessary nor avoidable.

The writer is at Stanford University in the U. S anants@stanford.edu

It is very important and interesting to note that anant Sudarshan is an ex-IIT-tian, of Delhi. So it is important to listen to his views as it can be a true picture about reservations. They also confirmed that it is beyond economy which holds them backward. That can be only caste. For caste has played and even today plays a major role in the social structure of India. If students are from backward or dalit they very often become the victims of ill treatment from their Brahmin and upper caste teachers. If they are dark in complexion says the expert the case of finished. If he / she comes from an uneducated poor family the case is to be permanently buried. So it is according to the expert only reservation that can save them. Unless reservations in educational institutions run by central govt. is made the nation can never be called a developed one for the majority would continue to be poor and good higher education would be only a dream for them.

The New Indian Express

The weekend's disturbances in medical colleges and government hospital across India, which is likely to intensify in the days to come, provide yet another example of middle – class activism taken to the extreme. The students are junior doctors may have the conviction of their cause through chain SMSes no doubt helped stroke those fires – but there can be no justification for stopping medical services. Not only does it go against the oath they have taken or will take but it will also lose the public sympathy and support. The strikes are mostly at government run hospitals, which are the last resort for the poorer sections of society, to deny them treatment is an act of callousness. They also hint strongly a knee-jerk reaction on Mandal II because however flawed the governments stand policy on reservations in higher education there is no disputing the fact that social inequality does exist.

This is also time for the politicians and political parties to take their responsibilities rather more seriously. The main opposition parties have not criticized Arjun Singh's proposals strangely neither have they criticized the protests. One can understand their tacit support for quotas – it is borne out of political compulsions. But they cannot bat for both sides, which is the suspicion that will accrue the longer they keep silent. This is also no time for the HRD minister to be scoring political points in attacking the National Knowledge Commission and its chairman. He has more pressing matters at hand and he could begin by talking to the students.

The strikes are not justified but the protests are given the scale and spread, at least a symbol of how the student community feels. They deserve a fair hearing. But first please tend to your patients.

This article was very much criticized by the experts mainly in two angles. One if the doctors and students are so selfish and do not think of others that is patients, can they ever think of the society and the majority of the deprived OBC's, their sufferings and their backwardness. If they are a little humane will they ever carry such protest which mainly has affected the poor patients? It is an ardent bitter truth the majority of the people are OBC are they live below poverty line with no education, or any form of comfort. Can we say they are lazy? They are much more hard working than the Brahmins and these striking medicos. Majority of them work more than 8 hours a day still cannot expect or practically have a square meal a day. They are paid from Rs 40 to Rs. 150 maximum for their labour over 8 hours labour and only when they work they get this pay with no pension or provident fund or medical allowance or HRA or city compensatory allowance. Why no upper caste persons or no politician is ever thinking about this, leave alone helping them? So if govt. just says reservations for them they join together make a beastly protest and the media is very handy in helping them to blow up the protest in multi dimensions.

Why do these politicians do not have courage to say we may only increase the percentage if you protest, we will dismiss or suspend you (protesting medicos) from your services or studentship? Why are they budging to their unlawful selfish protests? Suppose these OBC make a counter protest can these upper caste live in peace and luxury? Can the politicians find their rule in peace? Please let the Brahmins, upper caste and the politicians not forget that it may not be long when the silent OBC's take up to protest the nation at that time would come to stand still and the Brahmin must choose to leave the nation. Already the south is in fire. Let the govt. and the arrogant Brahmins choose to live in peace or leave their position and nation respectively once for all! Why the media is reporting they need a fair hearing when their protest is not fair?

Express news

NKC's member convener Pratap Bhanu Mehta who heads the Centre for Policy Research said "Arjun Singh's claim that the knowledge commission and the chairman did not know the constitutional amendment is absolutely preposterous. In fact most members have written on the subject".

Mehta among the first to protest quota proposal said. He is obfuscating the issue. My view is that even the constitutional amendment is an enabling legislation. And it does not require the govt. to extend (OBC) quota in (academic) institutions.

NKC's vice chairman PM Bhargava top bio technologist and one of the two members to support the HRD minister's quota proposal said "His (Arjun Singh) casting aspersions on the entire knowledge commission, is totally unwarranted, unbecoming and incorrect".

Unfortunately all policies are only favouring the uppercaste both economically, educationally and in employment. See IT is ruled by Narayana murthy, a Bangalore Brahmin who dictates to the nation. Agriculture by M.S. Swaminathan who owns in crores of fund from India and abroad, Industries are only in the pockets of a few Brahmins, Ambanis and Iyengars.

That we see majority OBC are striving for a square meal a day inspite of working over 8 hours. Unless politicians take onto ponder over the problems and take appropriate steps the nation would face problems. The experts of the group included socio scientists public and post graduates from different fields.

Some said NKC as usual did not have even one non Brahmin, OBC. All the policy makers unfortunately happen to be only from the arya Brahmins! How can OBC have anything even as lawful? The very fact the striking Brahmin medicos were least bothered about patients such is the mental status of the selfish arrogant Aryans. Further the experts said until govt. extends reservations to all posts especially in the field of teaching in institutes run by central govt. the reservation in studies may not be of much benefit to the students. For it is high time the govt. to learn that if the arya Brahmin can be so rude to the patients how will they ever be supportive to the OBC and SC/ST students in the field of education as teachers and research guides. The atrocities of them over the OBC/SC/ST students in certain institutions are very well known.

The New Indian Express

While medicos in different parts of the country are waging war against the proposed 27 percent reservation for OBC in higher educational institutions, Tamil Nadu has dared to take an opposite view favouring reservations.

Doctors in state have been expressing their support for reservations staging demonstrations for speedy implementations of union HRD minister proposal to introduce reservations. The Tamil Nadu Medical Council too has favoured reservations.

“The entire medical fraternity of Tamil Nadu is supporting the central government reservation policy. It is really unfortunate that medical students and doctors are protesting against the quota system in Delhi”. Prakasam, president Tamil Nadu medical council told express today.

When compared to Tamil Nadu where there is a higher percentage of reservation for most backward classes and backward classes, there is only reservation for OBC in the northern states…

Experts accused the BJP of organizing the protests strikes secretly and welcoming the OBC, quota openly they said that this could be evidenced from the fact that the media gives undue importance to every anti-reservation protest and blows it out of proportion. The experts allege that the BJP does not dare to oppose the OBC quota because it is afraid of losing the vote bank of the majority. Therefore they welcome it externally for the sake of political survival but also secretly dance behind the screens and lend all financial and logistic support to conduct the strike against reservations. If they oppose the reservation policy, they are scared that none of the regional parties will ever enter into an alliance with the BJP.

Express news service

The BJP on Wednesday welcomed the UPA govt. decision to implement reservations for the OBC in higher educational institutions but criticized it for failing to build a consensus on the issue.

“We welcome the govt. decision. But the govt. has miserably failed to take into consideration, opinions of all sections of society”. BJP general Secretary Arun Jaitley told reporters here.

The experts criticized the double standards of the BJP. Their second charge is that the “… govt. has miserably failed to take into consideration, opinions of all sections of society”. Experts say that; in the language of the BJP, “all sections of the society” refers exclusively to the Brahmin; the “upper-crust”, who form about 3% of the nation’s population.

Express news service

Even as the government has decided to extend 27 percent reservation in institutions of higher education for Other

Backward Classes (OBC) from June 2007 it is yet to ascertain the precise population of this community. Yogendra Yadav noted social scientist and professor at the centre for the study of developing studies here… Yadav said "In my opinion the OBC population in our country is between 40-44 percent. I am basing this on the national election studies conducted over the past two decades by our centre". But BP Mandal commission… in 1990 has estimated the OBC population in the country at 52 percent….

The experts came strongly on the anti reservation people who wanted to do away with reservation, were asking questions about OBC population. (1) A tight slap to them can be that certainly OBC are greater than 27% of the Indian population (2) Will the 3% Aryan Brahmins have the mind or heart to regulate themselves only for 3% of seats in all educational institutions 3% on all high posts. We warn them a day would come when such laws will be implemented if they do not learn now to accept and adopt the life and let live policy.

One of the experts pointed out the news item from PTI in The New Indian Express in the form photos of 25-5-06 which had the photos of the effigies of Sonia, Manmohan Singh and Arjun Sign with the writings on the effigy of Sonia Gandhi.

It was a big shock to find the photo in the leading daily as PTI news which shows "A man sweeping near the effigies of Congress president Sonia Gandhi, Prime minister Manmohan Singh and HRD minister Arjun Singh before they were burnt by the medical students protesting against the reservation policy on the AIIMS campus, New Delhi on Wednesday.

The experts want to put forth to the public, politicians and the OBC the following questions.

Can they misuse their rights like this? For AIIMS is a govt. institution and our nations PM, HRD minister and Sonia Gandhi, all of them are non Brahmins so only they burnt their effigy who are govt. employees!

Can a non Brahmin group burn the effigy of the a Brahmin leader. What will they do even such act was contemplated?

1. The govt. should have to unconditionally terminate or dismiss the medicos on the spot.

2. The experts said the AIIMS medicos are fully backed by Hindutva fanatics like BJP, VHP and RSS, that is why they do such base acts with this courage.
3. Above all who are they to demean sweeper jobs? If reservation was in existence the upper caste and Brahmin would have been sweepers!

Express news

President APJ Abdul Kalam, today asked doctors and students agitating against the proposed 27 percent OBC quotas to end their strike even as the agitating students called for "civil disobedience" tomorrow. While telling two groups of students who met him today that the quotas would be introduced next year, he gave an assurance that the number of seats in central educational institutions would be increased.

Can one say, in this nation all are given equal opportunities? When will the majority get at least some basic comforts and above all social rights? Only the reader should find an answer. Can we ever claim our nation to be a developed one? It is at this moment this fact was brought to light by one of the experts – It is reported from the July 2006 issue monthly magazine "Makkal Kalam." That is a hostel for 100 girl students. These students were provided with only one tube light. The minister for Adidravida welfare T.N, M.S. Tamilarasi when went on an inspection to these hostels found this shocking information.

Now can these students ever come to the status of the striking medicos? Thus if these rural poor from non brahmins ought to become developed not only reservations must be given to them but for them all education in these institutions must be made free; if the parents are not graduates and are from rural areas with paltry employment, demanded the experts. The experts said what answer can these striking medicos give for the status of these students? Caste alone is the answer. What is Brahmins responsibility towards them and above all to the nation? Is it fair on their part to protest reservation? Infact the reservation which govt. has announced is not sufficient. It is unnecessary for our president to see them and assure anything. Let them throw open the problems of the poor OBC and strictly punish the striking medicos. The media which should have done this job has failed to do it as it is only a Brahmin media.

The writer is a professor of economics, Jawaharlal Nehru University pradipla@mail.jnu.ac.in.

Caste is being used as the sole criterion in public politically oriented, towards positive discrimination; categories like OBC and SC are being treated as essentially homogeneous groups. Almost the entire population to each of these categories, it is alleged suffer from uniformly a high degree of deprivation. But what is the evident? The link between caste and deprivation should be examined at the macro level, small sample studies do not suffice in the large diverse country…

In 1911 the average work participation rate was 42 percent for high castes, 54 percent for the OBC and 58 percent for the SC. The average literacy rate was about 11 percent for high caste 1 percent for OBC and 13 percent for SC. Thus both the literacy rate and the work participation rate would appear to confirm that caste is a very good indicator of deprivation.

Such aggregate statistics are usually provided in support of caste based public policy…

Thus the use of caste as a criterion in public policy is not justified by the available evidence.

The article was given before a section of experts. Some of whom were professors teachers and socio scientists for opinion. They gave the following points (1) it is to be noted that the major portion of the nation economy is in the hands of the upper caste. When we take the caste of OBC/SC/ST they say work participation rate in some villages; in case of these non brahmins would be greater than 80% but if we take the work participations rate of an upper caste and arya brahmins it can be less than 20%. The reasons attributes for this are (1) The land and majority of that villages property would be in the hands of the upper caste and the arya Brahmins (2) The most disgraceful thing to be noted is that as an economist they should not only study the work participation rate but the relative pay they get must be first studied.

Can any where in the nation the sons of the soil toil for over 8 hours yet be denied a square meal with no medical facility or HRA or PF or pension to the majority of these workers are daily wagers. The whole family works and the child labour in at its best. Does it mean the discrimination is not caste based? Do you see any child labour practice among arya Brahmins? How many of them are daily labourer? How many are scavengers? How many plough the land? Statistics for all these must be taken. All this will prove and establish beyond any doubt that only caste is the factor. All well paid white coloured jobs are done by them. They receive not only pay something extra the tax of the poor with least labour if not money. The divide in India is only caste and if the arya Brahmins do not want reservation first, “let them come out and throw away the caste and become for sometime scavengers, sweepers and household servants” advised the experts.

The New Indian Express

In a fresh bid to end the anti reservation medicos strike the govt. tonight proposed to examine their key demand for experts committee on quota in central educational institutes, but the students and the doctors rejected the proposal and decided to intensify the agitation by calling a bandh in New Delhi on 31 May.

The government also said it was committed to expansions of medical educational facilities and institutions to ensure that students from all sections of the society have adequate opportunities to meet their aspiration for higher education.

The run-up to implement 27 percent quota from June next year will be utilized in creating and expanding necessary infrastructure and other facilities (including teaching faculty and other staff) corresponding with over all increase in the number of seats so that there is no cut of in the total seats available for general category students a govt. statement said.

For ensuring the quality of teaching appropriate measures including relaxation in the age of retirement of professors will be taken up so that the services of qualified faculty remain available, he said….

It is the decision of Sonia, PM, Arjun. Rejecting allegations that he was doing hard politics from a soft ministry HRD minister Arjun Singh today said that the decision for reservation in higher educational institutions was not his personal agenda but it has the endorsement of UPA chairperson Sonia, PM and the UPA govt.

Experts came very heavily on the fact, all these days when OBC the majority of the nations populations was denied or never had any opportunity for higher education in institution run by central govt. no one ever bothered about it. When they are little awake and the govt. is giving a paltry percentage of reservations the arya Brahmins are so much protesting and saying their opportunity is lost. Are we asking for backlog reservations? Even if all the seats are declared by the govt. as reserved seats in all govt. educational institutions experts said the arya brahmins have no voice to dispute over the issue; for over 60 years they have solely enjoyed every thing.

Suppose the OBC/SC/ST students protest saying that we do not receive any quality education from the teachers in the govt. run institutions; what will the govt. do?, asked a set of experts. The immediate action would be to suspend them or ask their parents to meet the authorities or punish them in different ways. But when an arya Brahmins puts forth the same question with his arrogance it is unfortunate our PM explains to them and promises to keep up quality. We are unable to understand the soft stand taken by the PM and our President.

Express news

Opposition leader L.K. Advani, today demanded setting up of a specialist committee to study the impact of reservation and examine the alternative remedies for affirmative action which would not undermine merit or excellence.

The experts said ages since the arya Brahmins have been using useless terms like 'merit' 'excellence' 'quality', 'standard' and so on. The only indirect claim is that only their production have these qualities. They continue to treat us only as second class citizens, said the experts. Our only advice is let them do quality meritorious work in any other nation and leave behind our nation for peace and development! The article of V.T.R in Dalit Voice "Merit is my foot is sufficient reply to them.

Privatize affirmative action by Jaithirth Rao, The writer is a Chairman and CEO, Mphasis, jerryrao@expressindia.com

Our great govt. has decided that since they cannot improve primary and secondary education they will have quotas in the higher education both in elite government institutions (IITs,

IIMs, etc) and in private institutions (except those run by the minorities who can do what ever they want) in accordance with constitutional amendments which are passed by parliament with the explicit objective of subverting considered supreme court judgments…,

So my friends wrote to me asking me why I was "silent" on the reservations issue. I told them that I was silent because I am convinced that we are faced with a fact accompli. Like it or not this battle is lost, reservations will not go away. Their strength hold will increase. No group which gets concessions is going to let go of it. More and more groups will spend time energy and money getting themselves classified as "backward". Otherwise buy fake certificates.

The decisions of our honorable central cabinet (populated doubtless by honorable men and women), the diktats of our feudal – socialist HRD minister the directives of our much amended constitution do not apply in.., Sri Lanka, the Maldives, Bhutan, Nepal, Bangladesh, Singapore, Dubai and Pakistan. They are all exempt and some have enlightened rulers. Eureka! Let every educational entrepreneur worth his salt open up colleges in these unreserved countries. All who can afford and some who cannot (but who will borrow an opportunity for our banks!) will become willing students in these colleges. Given the most engineering medical and the management college involve hostel living anyway, the additional burden may not be much. And trust me, the holders of the "unreserved" foreign degrees will command higher wages than holders of reserved domestic degrees. The immutable laws of the labour market shall prevail and once again our cabinet of central planning aficionados will be made fools of, albeit after adding transaction costs to the economy and unbearable and avoidable angst to thousands of deserving "non backward" students who are not rich enough to cross the borders to the new colleges that spring up.

I had completely under estimated Indian entrepreneurship. Indian organizations have already opened up medical colleges in Nepal, in Dubai. They obviously had good astrologers for consultants. They knew that the UPA government bored with nothing to do (upgrading primary schools or fast forwarding the judicial process being too lowly tasks for them) would announce the great reservations mela. Having anticipated this, our entrepreneurs have already moved. Now all that they need to do is expand and grow. Even the much persecuted (uniformly and impartially by the NDA and the UPA) IIM directors must have been prescient, hence their desire to expand to Singapore and so on.

On a different and serious note I believe that corporate India has another business opportunity. In a country where there are tens of millions of unemployed reserving jobs either in public or private sector will achieve very little. If we lap into the entrepreneurial skills of our people who revel in self employment (we are a nation of hawkers and of farmers both classic self employed categories) we can create an explosion of productive capacity corporate India should introduce an affirmative action programme in the management of its vendors, suppliers, dealers and distributors.

My friend Vaidhanathan from IIM Bangalore says that the need of the hour is the "Vaishya – isation of all of India's castes. Let us go for it. A cabinet or parliamentary decree will result in one more law and the corrupt will side step and which while enriching numerous govt. inspectors will leave India's poor (of all castes) pretty much where they are.

The experts welcomed Jaithirth Rao aim of leaving India for education. Atleast if they leave India, OBC and SC/ST will have some place. They further said all the while the Brahmins were getting education from India from institutes of national importance from govt. financed by tax payer and after getting the best of education left India and served other nations for money. Infact not only they grew fat in economy it is said that 10% of India's wealth is in the hands of IIT aluminums. He further added that they are very much against reservations and do a lot of spade work to see that reservation is not implemented. They have a very strong network based on their community. It is better if they get education from countries like Sri Lanka, Maldives, Bhutan, Nepal, Bangladesh, Pakistan then certainly India will not lose money and energy over them to serve other nations. Atleast the poor rural OBC, SC and ST will get a better share of seats from these good institutions.

The experts further said the loose way of talking by our friend only projects his frustration, it is like the "fox saying the grapes are sour". If they had faith in what our friend has talked why they have not practiced it. See our friend makes a blunder by stating India's poor would remain as the same poor. This is also good but the Brahmins in these 60 years have only made the poor still poorer. This is what has made the experts feel sad. Further these experts still ridiculed Jaithirth Rao for saying we are the nation of hawkers and farmers, clearly how many arya Brahmins are hawkers and farmers? What maximum an uneducated Brahmin can do is beg or work in a temple as a priest. Other than this the lazy uneducated arya Brahmins do not know and do not want to do any job. We suggest why not for a while for a change these arya brahmins take up the job of scavenging, coolies and sweepers and let the OBC, SC and ST read for some time. Thus, if they are willing to exchange certainly we will view them as patriots and has a true feeling towards our India. But we the experts doubt very much now why they want is get educated from universities in Pakistan, Dubai, Nepal, Bangladesh, Bhutan, etc. Some of experts said they are awaiting the good movement when they will quit India. For they have done more harm to the Indian soil than any good! Now if govt. wants to banish the backwardness based on caste they are protesting as if all the government run institution are owned by them! Let them remember caste by birth is not practiced in any nation in the world.

Express news

I write to resign as member convener of the National Knowledge Commission. I believe the commissions mandate is extremely important, and I am deeply grateful that you gave me the opportunity to serve on it. But many of recent announcements made by your government with respect to higher education lead me to the conclusion that my continuation on the commission will serve no useful purpose. The knowledge commission was given an ambitious mandate to strengthen India's knowledge potential at all levels. We had agreed that if all sections of Indian Society were to participate in and make use of the knowledge economy, we would need a radical paradigm shift in the way we thought of the production dissemination and use of knowledge. If some ways this paradigm shift would have to be at least as radical as the economic reforms you helped usher in more than a decade ago. The sense of intellectual excitement that the commission generated stemmed form the fact that it represented an opportunity to think boldly honestly and with an eye of posterity. But the government's recent decision (announced by Honourable Minister for Human Resource Development on the floor of parliament) to extend the quotas for OBC in central institutions the palliative measures the govt. is contemplating to defuse the resulting agitation and the process employed to arrive at these measures are steps in wrong direction. They violate four cardinal principles that institutions in a knowledge based society will have follow.

They are not based on assessment of effectiveness, they are in compatible with the freedom and diversity of institutions, they more thoroughly politicize the education process and they inject an insidious poison; that will harm the nation's long-term interest. These measures will not achieve social justice. I am as committed as anyone to two propositions. Every student must be enabled to realize his/her full potential regardless of financial or social circumstances "increasing the supply of good quality institutions, more robust scholarship and support programmes will go much further than numerically mandated quotas. The knowledge economy of the twenty first century will require participation of all sections of society. When we deprive any single child of any caste of relevant opportunities, we mutilate ourselves as a society and diminish our own possibilities… I can not help concluding that what your government is proposing poses grave dangers for India as a nation…

Experts are very angry with our govt. for everyone in the National Knowledge Commission is only Brahmins or upper caste. Atleast now let our dear Prime Minister know the true colour of caste thro' the NKCs. The experts said because our PM have studied in US he does not understand the venoms of the brahmins. The brahmins when it comes to reservations more so for OBC are very much against it. The politicians would have understood this when Mandal's recommendations were announced by our ex Prime Minister V.P. Singh. Anything that pertains to reservations for OBC be it job or in education the arya Brahmins are dead against it. The deepest venom in them is that they will talk in favour of other policy and project to uplift OBC other than reservations.

Experts said the views of Pratap Bhanu Mehta are perfect the govt. should have immediately defused the resulting agitations by terminating them form services, only like it was done in Karnataka. The govt. human approach made the medicos more and more selfish and arrogant. If the govt. would have taken stern actions, the agitation would not have been carried over a fort night claimed the experts.

Mehta claims that the govt. action is disappointing for the reasons claimed by them. What is more disappointing is that neither the govt. nor the so called Aryan policy makers have never in these 60 years raised any concern for our (OBC) people or any form of their upliftment says the experts; when everybody affected in the problems like; monsoon failure, weavers problem or closing of dyeing units and so on. Only when the govt. announced reservations for OBC they (arya brahmins) continuously write in the media as if they are concerned of the improvement of the OBC’s lot by injecting the insidious poison that will at large harm the nation’s development for if majority remain to be uneducated and poor; who happen to be the true workers for the nations real progress. The arya Brahmins who pen such articles with the media’s sole support hoard the nation wealth without any toil or hard work. It would have been better said by Mehta; says experts when we deprive any single brahmin child of relevant opportunities, we mutilate ourselves, the brahmin society and diminish the arya brahmin possibilities and not as falsely said by Mehta when we deprive any single child of any caste, of relevant opportunities, we mutilate ourselves as a society and diminish our own possibilities. Do they ever think how many millions of OBC/SC/ST children have been deprived of their relevant opportunities? When the very thought one child in their caste may be deprived makes them shun and protest OBC reservations. When the govt. makes some constitutional amendment to enable atleast a few from these millions of children to get educated the so called NKC member convener is so upset and want to render his resignation! What does this drama of Metha teach us? He cannot be neutral to one and all! Any policy helping the deprived majority is felt and described as danger to nation! What a great friend is our Mehta to OBC? What does he think of the intellectual capacity of our people? Let him know his place. What guts he has to make statements like that, reservations in central govt. institutions would … harm the nation's long term interest.

First our friend Mehta resignation is a blessing for us says the experts, for atleast now let our PM learn a lesson and put OBC/SC/ST in that post says those experts. Further we want to know the powers of NKC? Secondly we want to know is; it above constitution and PM? We want the NKCs to be dissolved before they do more harm to OBC!

Express news service The Karnataka government will dismiss from services those doctors who are in an unauthorized absence for a long spell.

Express news

With the debate on reservation of jobs in private sector summering commerce Minister Kamal Nath had said that having quotas in private jobs would not impact FDI flows into the country, which grew over 40 percent last fiscal to 7.5 billion dollars.

"There would be no impact of FDI flows" he told reporters.

Express news

Supporting the PM's suggestion for job reservations for the underprivileged in the private sector, the left parties on Wednesday demanded legislation to this effect while observing that industry would continue to reject proposal for such "affirmative action".

Mere appeal to the industry will not serve the purpose as they continue to harp on "merit and efficiency" the CPM, CPI and the forward block said adding that only enactment of the law would ensure reservation in the corporate sector.

CPM politburo member M.K. Pandh made it clear unless there is a law job reservations will not be implemented….

The New Indian Express

V.P. Singh could argue that fighting with his back to wall he had no time to consult any one. But in heaven's name, why make the same horrendous mistake a second time.

There can surely be a case for reservations even for the OBC's officially called Other Backward Classes but in reality Other Backward Castes. The trouble however is that this case has never been considered if the expression is permissible on merit, it has almost always been driven by crass political even personal calculation. Let the record speak for it self…

The issue of extending 27 percent reservations to institutions of higher learning such as the IITs and IIMs yet to be considered by the union cabinet has already been declared an irascible indeed immutable decision. Clearly the inner conflict within the congress party is at work and indeed visible in all its coarseness to the naked eye…

However, given the state of the Indian polity there is going to be no escape from the folly of 27 percent OBC reservations in medical colleges, IITs, IIMs and so on. These will be rammed down the country's throat, possibly through an ordinance some sweet words about increasing the number of seats in the institutions concerned to "protect everybody's interests". But only the mindless would believe that a change on such a large scale can be introduced in one go. Especially at a time when there are at least 20 percent vacancies in the faculties of the relevant institutions.

Kancha Ilaiah the invertebrate supporters of all reservations put his fingers on the heart of the matter when he said on a private TV channel that all children must go to neighbourhood schools with a higher enough standard of teaching and reservations would become irrelevant in 25 years. It is too tall a order for the country that has not been able to make primary education universal and compulsory in nearly sixty years nor is able to reduce the appalling dropout rate in primary and secondary schools, especially among dalit and adivasi students.

The article of Inder Malhotra feels sad for reservations of OBC. The following are the expert views on this article. (1) Why no article published in favour of reservation find its place in the media as frequently as articles against reservations! Is the media's stand only in the hands of the upper caste that too the arya Brahmins? Why still no Brahmin has come up with a suggestion of doing away with castes which is the root cause of all these problems? Secondly why do they deny entry of OBC to these govt. run institutions? Is it not an act of discrimination? Third, why not the arya Brahmins run institutions for them as they have surplus merit and economy? What is the difficulty? Our strong view is they totally self confidence and lack merit! What is their legal stand to say no to OBC reservations? Why majority of the population is denied the right to study in central govt. run institutions? It is high time govt. makes some laws to punish the few people who are against it. They can be treated as anti social elements. What is their ground, when they say nation's future is at stake? What is the nation achievement in these 60 years in terms of majority? Only the arya Brahmins have flourished both economically, politically and employment wise beyond any specific dimension! They have by their selfishness crippled the majority. They no more shed crocodile tears about nations prosperity or the merit of these institutions. They are bothered about nation not about the majority! What is the nation to them a place without OBC!

Experts said our friend says OBC does not mean other Backward Classes but Other Backward Castes; we ask what is wrong?, when you arya Brahmins stick so tight to caste! You shed the Varna in laws of manu this would naturally make reservation meaningless.

Express news

Charging that the AIIMS was being converted into a political hub union Health Minister Anbumani Ramadoss today warned of stern action against those "fomenting trouble". Responding to queries from media person on the AIIMS imbroglio Anbumani said the problem was a fall out of the anti reservations demonstrations that took place on institution premises. "They were a few faculty members who were instigating the striking students. Action will be taken against those making this medical hub a political hub". He said adding that there was no interference from him in the functioning of the institute. Referring to the posting of a dean, he said the AIIMS administration was not performing for the last two or three years and "the administrating posting therefore became necessary". I am the president of the institute and by an act of parliament and I have to do all I can in streamlining the administration he averred.

The minister however promised that he would talk to the staff on Monday on his return to Delhi and try to settle the issue.

Many found it intriguing how was it that the south particularly Tamil Nadu remained cool, when the protests against the move that provide 27 percent reservations for Other Backward Classes (OBC) in elite institutions have raged like a wild fire across north India?

Does not the plethora of fears raised over reservation - on sacrificing merit promoting inefficiency, likely fall in standards of the institutions and fuelling caste politics – have any relevance to the south too?

For Tamil Nadu which now adopts a 69 percent reservation policy, has been in all governments order, litigations, committees, debates protests, experimentation with creamy layer based on economic status and every thing that is debated upon now during the last 152 years ever since the then British government issued a standing order in 1854 (No. 128-2) urging collectors to divide the subordinate appointments in their districts among the principal castes. It was not the Britishers alone who felt the need for reservation in jobs. In the neighbouring kingdom of Mysore the king gave a similar directive in 1894 and then followed it up in 1921 by issuing an order after appointing a committee under Justice Leslie Miller in 1918 to study the issue. Then Diwan Vishweswaraya a brahmin resigned his job in protest against the order. Madras presidency however paved the way for a reservation policy through an order in 1921 clearly listing out job quota 44 percent for non Brahmins and 16 percent for Anglo Indians and Christians 16 percent for Muslims and 8 percent for Scheduled castes. That was during the regime of the South Indian Liberal federation, popularly known as the Justice party after the outfits' English daily news paper "Justice".

The G.O. remained only on paper till 1927, when Muthiah Mudaliar, a minister in the next government headed by P. Subarayan (grand father of Rangarajan Kumaramangalam) issued fresh orders which came to be known as the communal

G.O. to ensure job reservation in the registry department. Prior to the govt. had also issued an order in 1923 exhorting admission of SC in schools. The communal G.O. subsequently underwent changes – somewhere along the time the term "backward Hindu" was also incorporated and quotas were provided for in jobs and admission to colleges till 1950 when the GO was struck down by a full bench of the Madras High Court on the ground it was against the article 29(2) of the constitution. The verdict was given on two cases one filed by Shanbagam Duraisamy and other by C.R. Srinivasan. The grievance of Duraisamy was that she had applied for medical college admission and could not get a seat in view of the reservation policy and Srinivasan had contended that he was denied admission despite having the qualification in an engineering college because he happened to be a Brahmin.

Duraisamy's case was argued in court by Alladi Krishnasami Iyer a member of the constituent Assembly and former Attorney General and Srinivasan was represented by V.V. Seenivasa Iyengar, a retired judge. Though the state govt. went in appeal the Supreme Court upheld the high court verdict on Apr. 9, 1951. Later a lawyer P. Rangasamy has recorded that the during the hearing of the appeal in the Supreme Court it was learnt that Shanbagam Duraisamy had obtained her BA degree in 1934 and that she had not applied for any medical college in 1950.

The striking down of the communal GO led to wide spread agitation in Tamil Nadu, which impelled the Centre to bring in the first constitutional amendment adding a fourth clause to Article 15 (prohibition of discrimination on grounds of religion race, caste sex or place of birth). Clause 4 reads as follows "nothing in this article in clause (2) of Article 29 shall prevent the state from making any special provision for the advancement of any socially and educationally backward classes of citizens or for the scheduled caste and scheduled tribes.

Based on that, reservation resumed in Tamil Nadu, helping a few of backward class people forge a head in 1957, after the reorganization of states the reservation formula followed in Tamil Nadu was, SC 16 percent, BC 25 percent, OC 59 percent.

Later on DMK government raised the quota for BC to 31 percent.

After M.G. Ramachandra became Chief Minister in 1971 an annual income of Rs. 9000 was fixed as ceiling to avail of BC quota benefit – a sort of economic criteria and also a creamy layer concept. MGR withdrew the order on income ceiling after his party was routed in the 1980 Lok Sabha elections and also raised quota for BCs to 50 percent thus taking the total percentage to 68.

In 1988 when the DMK was in power it was raised to 69 percent with one percent added for STs. The BC quota to 50 percent was also split into two 30 percent for BC's and 20 percent for denotified tribes and most backward classes (MBC) a new group carved out of the BC's following the agitation organized by the Vanniyar Sangam led by Dr. S. Ramadoss.

Today Tamil Nadu provided 69 percent reservations in jobs and admission for all courses, including medical and engineering colleges. But a cursory look at the cut off marks for MBBS admission will tell whether merit takes any beating due to reservation and if at all how much. In 2005 the cut off in open competition was 294.83 marks while for BC it was 294.59. (The difference being mere 0.24 marks out of 300 or just 0.08 percent).

For MBC the cut off was 292.50 over and above that BC, MBC and SC student cornered 374 of the 433 seats in the open category.

Despite so many BC's becoming doctors over the years the states health care delivery system has won accolades from economists and social scientists including Nobel laureate Amartaya sen. The reservation policy is said to have enabled more people from the backward classes and scheduled castes to become doctors and go to their villages and small towns for practice.

The idea of reservation has sunk deep into the collective psyche of the people that they no more dread it though there could be individuals with grievances on that score. Unlike the medicos who went to strike as they were gripped by the fear of the unknown – not knowing what will be in store in the new quota regime – people in Tamil Nadu have accepted reservation since the state has seen remarkable upward mobility of the backward classes and scheduled classes and facilitated inter-mingling of different castes.

So in the prevailing social ambience it is hard to expect a small section to suddenly go up in arm against what is a mere extension of reservation to central institutions. Above all not only most of the caste groups in Tamil Nadu enjoy the BC or MBC facility, the socio political dynamics is such that reservation is not a dirty word in the state, where even the Brahmin association is demanding quota for its members.

The experts welcomed this article whole heartedly. They said this article is a strong feed back for one and all who think that reservation for the OBC will make nation to deteriorate from merit and quality. What is all the big fuss? Tamil Nadu is one of the leading states in medicines. HIV/AIDS patients from Andhra Pradesh and from northern states folk for treatment to Tamilnadu. The medical development of T.N. is unparalleled with any other states claims the experts.

So no more the topic of the merit criteria should be spelt out by the arrogant "meritorious" arya Brahmins warned the experts. Especially the medicos alone made a strike about reservations. This is beyond doubt for the experts that the strikes are induced by the BJP/RSS and certain Hindu Brahmin mutt and nothing more. That is why our esteemed Honourable Union Health minister Dr. Anbumani was targeted by the subordinate Dr. Venugopal. When Dr. Anbumani sat in the OP section in AIIMS, Delhi the expert said the ego of Dr. Venugopal was hurt. So he made a melodrama he would resign after taking into confidence that the resident doctors would also resign if he resigns. But Dr Anbumani did not heed to him.

So dear medicos; don't underestimate the capacity or intelligence of the OBC's. We know that you very well know they are intellectually superior to you that is; why you fear and shun OBC reservation in AIIMS, IITs and IIM said the experts. It is high time to bring to light that the so called NKC had failed to properly analyze the pros and cons of reservations for OBC for they saw it only in the arya Brahmin angle for their peril and only sorrow was that the number of seats enjoyed by them would be shared by OBC. This alone is inhibiting the minds of all the arya Brahmins say the experts. No rhyme or reason ever place any role claims the expert. If alone one takes the state T.N. which has the highest percentage of reservation say 69% is no way lagging behind any other state. Even today it can be marked as the best performed state in medical, engineering and education field in general. So it is high time we think of the good of the majority who was lagging behind mainly due to lack of opportunity and provide them reservations in central govt. institutions, which is our duty towards our brotherhood and a dire need to make India a developed nation and not a developing nation even after 60 years of independence claims experts.

Experts strongly feel it is always the legal field which is only dominated by the arya Brahmins which always gives verdict against reservations. Can any non Brahmin get justice from court, whatever be her / his merit? Courts play a role supportive of Brahmins and are always anti to reservations claims experts. So the court has always been against the non Brahmins, SC/ST. Experts said even recently the Madras High Court has given a judgment on the PIL filed by a retd. IAS officer Karuppan Dalit…. The very illegality about it was, the case came first under a bench which itself was unconstitutional. Since the PIL was whether reservations for SC/ST existed in teaching post in IITs, the honorable court dismissed the case by saying such reservations do not exist. So the very judgment in the opinion of the experts was against the very constitution. Thus in India says the expert the legal field always blindly favour only the arya Brahmins as majority of the law giver were only of the arya Brahmin community. Even if one of two judges were non Brahmins or SC/ST they once go to this profession mostly favour them! Thus according to experts when it is reservation the verdict was always against it.

That is why the experts said in the article the court not even properly verifying in the case of Shanbagam Doraisamy who had not even applied to the medical colleges put a case as if she was denied a medical seat this false case lead to the cancellation of the communal G.O. by the Hon’ H.C. of Madras. Such frauds were always practice by the Brahmins’ in courts.

The New Indian Express, Sachin Pilot the writer is a congress MP in Lok Sabha.

My understanding is that people don't really like reservations. It could be an empty table in restaurants, seats in a movie hall or in a train if these are marked reserved it doesn't leave the waiting lot very happy. And for me that is a natural reaction.

It is however important to see whom these are reserved for , let's say pregnant women or a berth in train is marked… then our initial disappointment or even mild frustration does change into a more compassionate between the cases in these example and those of reservations in educational institutions.

"There is need to extend our support system so that those side lined for centuries catch up. But that does not mean we have to necessarily compromise on merit"…

Experts say first the Lok Sabha M.P. needs to know the notion of merit. Does it really mean merit is scoring high marks in the school final exams and the entrance exam. Experts said a big no. Most of the arya Brahmins who score high marks in both have all comforts in their command but majority OBC/SC/ST do not even have 1/100th of these comforts; several children do job before and after school hours they have to walk many kilometers to reach school. Many do not have nutritious good food. They above all study in schools berth very poor infrastructure!

Now dear Sachin can we ever compare their marks and say the arya brahmin is meritorious and the OBC/SC/ST is not meritorious. Are we compromising by giving seats to the OBC/SC/ST by means of reservations merit? Never, infact these deprived children by all means are in all way more meritorious than the arya Brahmins. Secondly T.N. which has been always following a very high percentage of reservations in our opinion produces the best of doctors, engineers and so on.

They are those people who have come up due to reservations. So a reservation for OBC is the need of this hour. If they are left without any support; the majority will create a national problem. A revolution may break out. For all good posts are usurped by the so called arya Brahmins. They are the policy makers for us. How can any good be done to us? Might be you as a M.P. have power and means to get anything so you have no need to use reservation; think of the majority, they need it! You think not from the plane of an M.P. but from the plane of a majority living in rural areas with no education. Certainly only reservation can lift him for the bogus term merit used by the arya brahmins has made them permanently poor with labouring over 8 hours as per the laws of manu "for the riches of a sudra makes a Brahmin feel sad"!

Express news

The proposed 20 percent reservation for the OBC in Central Government institutions like Indian Institute of Technology (IIT) and the Indian Institute of Management (IIM) among others has gained more support in the state, with political parties and other organizations extending support and demanding it immediate implementation.

In a release PMK president G.K. Mani said that the party condemns the medico's agitation against reservation, adding that the stir was instigated. He also recollected the demonstration staged by PMK founder S Ramadoss recently in support of the proposal in Chennai….

The Joint Action Council (JAC) of college teachers Tamilnadu has also welcomed the move, brushing aside the claim that reservation would lead to dilution of quality. The organization also welcomed the Union Minister for Human Resource Development Minister, Arjun Singh's refusal to reconsider the proposal. The Tamil Nadu government Medical, Dental students, CRRIs and post Graduate Students Association in a release said reservations were not against merit and the existing 69 percent reservation in Tamil Nadu has promoted socio economic and health care delivery system of the state to a very great extent.

The experts accepted the claims of various organizations and once again pointed out the anti reservation groups to view the problem in a selfless and in an unbiased way. For the majority cannot always be denied every basic comfort for the sake of pomp and luxurious lazy living of the 3% arya Brahmins.

Express news

It is official quota for OBC from next year. The UPA left coordination committee has decided to fix reservation for the OBC in educational institutions at 27 percent by bringing in the legislation in the monsoon sessions and implement the decision from the academic session commencing in the June next year, 2007.

The committee resolved to implement reservation in letter and spirit and increase the seats in educational institutions under the purview of the central government. A student leader along with other students of the Allahabad Central University getting his face painted in favour of proposed reservation policy of the government in Allahabad on Tuesday.

Express news on the first page.

SC rethink on quota for PG medico’s. Vacates interim order on providing 10% quota for SC / ST.

The Supreme Court today vacated its interim order directing the centre to reserve 10 percent seats for scheduled castes and scheduled tribes in post graduate medical courses for the year 2006 – 2007.

A bench for Judges K.G. Balakrishnan, A R Lakmanan and D K. Jain accepted the submissions of additional Solicitor General Gopal Subramaniam that enforcing the 10 percent order for this year would create practical difficulties. He said admissions this year should be free of reservations in the 50 percent as this was not provided for in the prospectuses that has already been issued. The reservations will however continue for the 50 percent state quota.

The Bench has passed the interim order in April 24, 2006 on a writ petition filed by a batch of 14 doctors challenging the method of calculating seats reserved for SC/ST candidates in post graduate medical courses….

The experts blamed that always the courts be it a Supreme Court or the High court it was always in Brahmin monopoly. All the judgments were only favouring brahmins. The highest eyewash they do is to make the non brahmin or SC/ST judges to sit and make a resolution which is in favour of the arya brahmins and against OBC/SC/ST. Now also one SC judge, one Chettiar and one Jain were used by the cunning arya brahmins to vacate the very interim order of providing 10% quota for SC/ST in PG courses in medicine.

Experts still pointed out especially regarding reservations the OBC/SC/ST can never get any favourable judgment. It can be only against, experts claimed. Only a revolution in the judiciary can do some good, favouring reservation they said. Even Judges belonging to OBC, SC and ST are highly crippled of their rights by these over populated arya brahmin judges. They contended that judgments on reservations were always in keeping with the laws of manu. Until reservations for OBC/SC/ST are made in the judicial field and the laws of manu is totally banned it is impossible for us to get any justice from these courts claimed the experts.

Express news service

Chief Minister M. Karunanidhi on Monday urged Prime Minister Dr. Manmohan Singh to promulgate an ordinance immediately to ensure 27 percent reservations in IITs, IIMs and AIIMS for the other backward classes (OBC) Karunanidhi in his letter to the Prime Minister said.

“The very future of the down trodden communities is linked with this affirmative action. Twenty seven percent reservations may be implemented fully without any dilution or postponement.

He further said if urgent action was taken to promulgate an Ordinance to ensure the reservation for OBC living all over the country would be grateful to the union government.

The chief minister said reservation in educational institutions, controlled and financed by the centre had been one of the long pending demands of the people of Tamil Nadu, especially those belonging to OBC.

The level of expectations of the people of Tamil Nadu had been raised when the union government announced its proposal to implement 27 percent reservation for OBC in admissions to educational institutions like IITs and IIMs he added.

Thus the reservation for OBC/SC/ST in the central govt. run educational institutions all these years can be only construed as discrimination based only on caste.

Express news service

Joining the national protest against reservation, a group of students from various colleges in Chennai propose to organize an anti reservation protest in Chennai….

The centre should be more responsible keeping in mind the interest of all sections of the society the students said.

This statement made by them said experts would more fit for reservation for the majority who happen to be OBC are denied good education after their school final; the govt. should not have kept quiet and unconcerned for over 60 years after independence about the welfare of the majority. It should have implemented the reservations for OBC ever 50 decades before. If that would have happened says the experts, India would have been now declared as a developed nation. It is unfortunate the arya Brahmins who are so much self centred are making wrong statements about the centre. The experts further substantiated the joining of the IIT-ians' and software professionals is not really joining but they are the ones who are actually organizing it they added. Further the experts said it was the alumni of the IITs who are dominated by the arya brahmins has with them 10% of the nations economy and they are the fanatics who were very much against reservations and most of them are software professionals, headed by one Narayanamurthy INFOSYS; Bangalore. The experts claimed he had a strong backing of the Kanachi mutt, VHP and the RSS group. They have a lion share of India's economy. He is also the chairman of IIM Ahmedabad. He is a Brahmin to the core who is dead against reservations said experts. Infact the youth for equality is backed and part and parcel of the IIT alumni's and the RSS/VHP/BJP claimed experts.

Media had always been with them in projecting about any form of anti reservations. So before a protest after a protest and during the protest they will also simultaneously announce in the press.

Right Angle – Swapan Dasgupta. Let's be quite clear, the affected students and those likely to be most affected are in a woeful minority. After all how many can even hope to secure a place in the IITs, IIMs and the top medical colleges in India? The average student may take a pick of the lesser colleges, but even if a single student is forced to settle for lesser institutions because his positive grades are offset by a negative caste, it would be tantamount to a human rights abuse.

As is evident we are by definition talking about a real minusculity. Yet this minusculity is relevant. Those who have been agitating on the streets under the youth for equality banner are by far the best and the brightest of our youth, those who can hold their own anywhere. They have precious little need for either affirmative action or grace marks. They are India's undisputed creme de la crème. They are the one's who have transformed India from an over populated, third world back water into nerve centre of the knowledge economy. Yet for the last four weeks not one politician of consequence from any of the main stream parties has dared to be associated with this minority movement.

When it comes to the caste question and reservations an insidious political correctness over whelms India…

At the very heart of the protests against the new quotas is a frighteningly modern demand. The students are asserting their right to be treated first and foremost as Indians overriding all class caste and creed. That's a message we have been told constitutes the highest form of nationalism. At the same time that is a majoritarian argument no self respecting minority's will ever countenance.

Expert views were recorded on this article. These experts first asked these boosting youth. Are they enjoying any form of equality in terms of food, money or comfort, medical facilities or good education; asked the expert. Then what is the good thing have they ever seen in their lives except hard labour asked the experts. These people that is youth for equality as claimed by Swapan Dasgupta are "India's undisputed crème de la crème" what does he mean by this? After commanding all forms of comfort good food and good institute for school educations, very nice extra tuitions and super infra structured schools they get some 30-40 marks more than the downtrodden, poor, with no good school or any form of tuition, even no time to read, for after school hours they have to work for their stomach and family with no books, get about 30 to 40 marks less than crem ede la crème layer for a total of 1200, or 1100. Who is really bright? Who is meritorious? Who honestly need support? To ones surprise the later category is majority also. Still the nights spent by most of them are, sleepless, for their father in the majority of the cases of these students are drunkards and would pick up a quarrel in the nights added experts. Actually by the way of reservation if the nation would have developed these downtrodden, rural poor OBC certainly the nation would have been made progress in the field of technology, medicine and sciences in leaps and bounds asserts experts. It is the blunder the govt. did not implement reservations soon after independence. That is why they have become more and more downtrodden and now it needs more number of years to uplift them says experts. The youth for equality should now at least be happy the govt. is implementing reservations to save them, but on the contrary it is shocking to see they are green with jealous for the past four weeks, they are in the roads protesting a paltry percentage of reservation given to them in govt. run institutions like IIT's IIMs and AIIMS. The experts wonder how come they call themselves falsely as youth for equality, when they, themselves ought to have protested against the govt. for not giving them some percentage of seats as reservations. It is the boastful nature of them to shamelessly claim themselves to be Indias' undisputed crème de la crème when they are nothing of these sorts they are much lesser than our last of students, claim experts. The very isolation of main stream parties not associating itself is a strong sign to prove them that their protest are wrong and anti social. The arya Brahmins should learn to think in terms of majority and not only as 3% Hindu minority. Unless they change their attitude they have to be only striking in streets with the ever faithful media taking their photos and projecting in the best possible ways in the print and in TVs. Further they should have some sense to see the suffering of the downtrodden mass have been denied even the very basic need in their lives. They should have learnt a lesson from the silence of the main streamlined parties that their strike is wrong. Their claim they are the "creamy" is wrong. If the comforts and the facilities are given to OBC/SC/ST students their performance would be much more superior to them. What is their right to say OBC/SC/ST should not be given seat to study in these institutions which are run by the central govt.? If it is atleast their institution we can accept. Who are they to say no to the constitutional amendment which is favouring the downtrodden majority? All these amounts to two things apart from arrogance; (1) the media is in their monopoly more than this is that (2) the courts of law are at their disposal. The catalyst agent adding to these two factors is the Hindu religion. For the lawyer who argued for Dr Venugopal is the same lawyer who argued in favour of the Swami Jayendra Saraswathi of Kanchi mutt in an attempt to murder case! Can public and the OBC understand their strong network? They said their network is more finer than that of the spiders any ripple of any order would be monitored by them and supported by their religious, leader; the spider. So they would crush all things which trouble the network apart from getting sumptuous meals and etc. Whose property is IITs, IIMs and AIIMS asks these experts? Is it govt. or the arya Brahmins? Who is financing them is it the govt. or the arya Brahmins? Hundreds of crores have been wasted over these arya Brahmin all these years.

Only from 2007 June our people may be permitted to study lamented one of the experts. What a national waste to educate 3% of the selfish arya Brahmins with the national economy in hundreds of crores all these years? That too after getting education they leave India! Let the minority Hindus get 3% reservation as per their population! Let the govt. make one more amendment in all these govt. institutions the back log will be filled for OBC/SC/ST in all posts more so in teaching positions. If we can be more just and legal let us oust the extra people from arya brahmins so that they don't exceed more than 3% in all the teaching posts in IITs IIMs and AIIMS for this alone can guarantee of the OBC/SC/ST to get good education without discrimination and mental torture. Some experts said we would protest continuously till this is achieved. Let them resign their post on moral grounds if at all if they have some thing in their life as integrity. Let them not pollute the nation by publishing such caste fanatic selfish articles said these experts. Several of them complained that their article was not even acknowledged by the Brahmin media. Even in publishing they are so fanatic said experts. Even if we write in response to their selfish, improper articles they don't accept or have the heart to publish it said experts.

Sunday Spin by Bishwanath Ghosh, what's your caste? In the Sunday Express

Biswanath_g@yahoo.com

He speaks of how accidentally asked his Kerala girl friend of her caste which sparked a problem in his love life. As we experts visualize is that the question by Bishwanath as he says could not have been that superficial or on the surface! He wants to know the caste to build strong relations. For we have hundreds of instances were intercaste marriages have been a big failure. Either the bride or bridegroom burnt alive or killed or so in T.N and in India. Also it is very rare to see a Brahmin man even accidentally marry non Brahmin women. But it is matter of routine and surplusness a Brahmin women marrying non brahmin men mainly for the sake of his education, position and above all money. This is a very common feature in India!

And as IIM, Ahmedabad; board member N.R. Narayanamurthy asks can it be imposed from Delhi…

The experts put forth the following suggestions;

(1) The knowledge commission chief Sam Pitroda's idea of increasing more number of IITs is a foolish task for even a chief without any knowledge will not commission to such a stupid act (2) why should increase of seats be made, without any problem of disturbing the present serene, let us only give out from the existing seats; 27% for OBC who are nation majority. For over 4 to 5 decades they arya Brahmins have enjoyed now at least for 4 more decades from 2007 to 2047 let the OBC enjoy the central govt. institutions and their fruits their of. So the experts gave a very simple problem to the problem. Thirdly after all the teaching faculty of IITs and IIMs are so great they can teach OBC to a proper extent that no merit problem (which is however a false way of talking ill of OBC) comes, for if you deliver a good goods certainly the out come would be proper. Finally all their fear is the alumnis of IITs are the ones who are having a grip of 10% of Indian economy. They fear if OBC learn this trade how according to laws of manu can they tolerate economy in the hands of Sudra for it will pain a Brahmin [43 ]. Thus for sometime let things go as it is, but give reservation to OBC. Only thing we suggest is in the vacancies mentioned please appoint only OBC/SC and ST, so that OBC, SC/ST students do no suffer the pang of discrimination due to caste.

Deccan Herald

"Instead of uprooting casteism, reservations should not become a bane for society. Reservations as an aid to egalitarianism began in India much before VP Singh and Arjun Singh found the Mandal recipe useful for their own not entirely for unselfish reasons. The story began in princely Mysore around 1870 and this nonviolent method of protecting the vast majority from brahminical tyranny soon became popular in southern India…

The experts claim even of today the brahminical domination and tyranny is being suffered very silently by the OBC/SC/ST, for instance Dr. W.B. Vasantha of IIT (M) with her credentials is suffering silently, under the regime of Arjun Singh, and he has not taken any action, for he is an upper caste, and supports the Brahmins and the upper caste. If reservations are not announced by a govt. we may have a civil war some social scientists suggests and warn us. Thus they congratulate the rulers of the nation to take proper steps at appropriate time by saving the nation.

Times of India

Thirty five people, four of the media person's and some police men, were injured in a clash between police and pro-quota agitators here on Friday. The clash turned so ugly that even members of the press were not spared, beaten with sticks and lathis by the students.

Shouting slogans against those opposing the Union Human Resource Ministry's proposal to introduce reservations for OBC's, hundreds of pro quota medical students and junior doctors took out of procession from Gandhi Maidan. On reaching the busy Dak-Bungalow round about, they burnt effigies of anti-quota protestors. And without any provocation they turned violent and attacked journalists and camera men covering the agitation with sticks, lathis and other objects.

Photographers of the Press Trust of India, The times of India and Hindustan times, Camera men and Star New and Aaj

Tak TV Channels and some journalists were injured and are being treated in Patna medical College.

The pro-quota agitators also attacked policemen with stick forcing them to resort to a heavy baton-change, leading to injuries to about 20 protestors. In a bid to check the rampaging students the policemen began raining lathi blows on them forcing them to scurry for cover. In all around 35 people including students, media men and police personnel were injured in the clashes, a senior police official said.

The experts said that if the govt. does not take the appropriate step of giving reservations to the downtrodden OBC at the proper time certainly peace of the nation will be disturbed. Further lakhs of graduates from OBC/SC/ST remain unemployed. Unless govt. takes some steps to help the discriminated lot due to caste the nation will have to face dire consequences. The reservation infact announced by the govt. is much paltry. For it is unfortunate a substantial wealth of the nation is in the hands of the 3% arya Brahmins. One can see the luxurious life they lead. None of their protestors were furious because they do not need to strike to get anything more. What they enjoy from the nation is much more than their share. Such sort of inequilibrium cannot sustain for a longtime in any society says the social scientist. At one point of time it is sure to tilt. Before things become chaotic it is the govt. that should avert it by peaceful means. The OBC must be given reservation both in jobs and in education.

The Hindu

Questioning the very formation of the Nation Knowledge Commission, CPI (M) Polit Bureau member; Sitaram Yachury on Friday said the need of the hour was to ensure implementation of the 27 percent quota for OBC students in institutions of higher education.

"I don't know why the knowledge commission was even established at the first place and for whose knowledge it is "Mr. Yechury said at a seminar on "Debating education.

Hindustan Times

As the government finalises its plan on quota in elite institutions, a parliamentary panel on Monday recommended a graded fee structure with less for poor students than others. The aim to ensure fee should not act as an obstacle for poor aiming for higher education. Higher education needed to be supported by increased allocations by the state in coming years, a parliamentary standing committee on HRD said in a report…

The expert said almost all from upper castes will also come as poor to get permission to pay less for getting a pay certificate to that effect was not difficult in this nation.

Hindustan times

Unlike the attention seeking campaign by anti reservation activists faculty members in favour of the proposed reservations are organizing low profile demonstration at AIIMS every day.

Around 50 faculty members from AIIMS, Saldurjung, LHMC, MAMC and UCMS have come together under the banner of medico's forum, for equal opportunities and are demanding implementation of reservation in educational institutions. The forum has decided to organize daily lunch hour demonstrations at AIIMS and other medical colleges.

The pro reservation faculty members have also submitted a memorandum to the President blaming AIIMS, Director P. Venugopal for "aiding and abetting the agitation". Protests are not allowed inside hospital premises. This agitation is an exception to every rule and the director is to be blamed for it said Dr Anoop Saraya a faculty member; he said "Reservation for backward students in higher education destroying merit is nothing but bogus argument". The medical students should be equally sensitive to attending the hapless patients and they were duty bound to attend them round the clock.

This pro reservation activist campaign never found in media or in any of the dailies. Only Hindustan Times, New Delhi; published it. Thus the media, most of the leading dailies are trying their level best to project and create a feeling the whole nation is against reservation. So news about any campaign or protest pro reservation is never given any place in the media especially in dailies like the Hindu, The New Indian Express and the D.C.

Deccan Herald

Communication gap in time of saturation of communication! It appears that the anti-quota protesters and the quota advocates are talking through each other and not talking to each other.

And there is reason to believe that the media did not do its job of clarifying issues and that instead it followed the news in a blinkered fashion as it unfolded. If there was a protest march, the cameras and the reporters went there. The media did not bother to find out what the real problem was

“What India needs is equality of opportunity for all”

Meanwhile there is no doubt that everyone is contributing to the label. Union minister for Human Resource Developments, Arjun Singh told to Lok Sabha on Wednesday that he had not made any announcement regarding reservations for OBC in the higher educational institutions…

They should know that quota based reservations is not the right way of opening gates of opportunity to the SC, ST and the OBC. Reservations give the false impressions that the disadvantaged groups are not on par with the rest of the society.

What India needs is equality of opportunity for all. The government and the private sector must do everything to develop the intellectual energies of all sections of Indian society.

Reservation is not the way to go forward. Create a fair system of education where genius can blossom. Nothing more is needed. The issue is not reservation. The focus has to be on creating excellence. People then know how to make use of it.

The experts said media especially the leading dailies will immediately become paranoid once govt. talks or wants to implement reservation. Immediately statements like “what India needs is equality of opportunity for all” why all these days such articles did not find place. Did any Brahmin media ever said the downtrodden rural poor need at least better opportunity leave alone equality with a rich, city dwelling brahmin. No they jump in anger, frustration and above all insecurity and say, reservation is not the way to go forward, forgetting in good old times of British the brahmin begged for reservations.

Thus was reservation the only tool to go forward then! Or if it is brahmin’s reservations as they have asked on 24$^{th}$ Dec 2005 in their Brahmin Association held in TN then will it not imply it is the only way for creating excellence! What a double stand they have? They are really people with double tongue claims the experts. Does T.N. in comparison with other state ever lack in the criteria of merit or in merit, mind you they are the state with highest percentage of reservation in all state held posts and in education. In fact they are, one of the leading states in India. To be more particular most of the engineers working Japan in department of metallurgy are from Anna Univ. are working in Japan. Japan has more Anna University graduates and P.G. in this field than from any one of the IITs. So one cannot say reservations will ever dilute the standard or excellence or merit. Experts warned the writers that OBC's are much more meritorious than the arya brahmins; the only thing is that they had every where god fathers to lift them even if they lacked merit but real merited OBC has no god fathers to give them even their due credit added the experts. One of them said in those days Ramanujan the mathematician could see the viceroy and get help for his travel to U.K. etc.; can even today our OBC leaders or intellectuals, however meritorious they are, can see the P.M. or President and get any form of help, however meritorious they are asked some of the experts!

So the experts in a single tone expressed that reservation was the only means to go forward at least get, some sort of equality in this caste dominated society. For they all expressed only reservations can annihilate, the discrimination due to caste playing havoc in the development of the downtrodden people in education.

They all happily cited T.N. which has developed in almost all the areas which had the highest percentage of reservations. Thus for the nations development the central govt. too can raise the reservations upto 69% they claimed!

Times of India

Anti quota protests by medical students and doctors across the country continued on Wednesday even as union HRD minister Arjun Singh told parliament that, "the centre would take an appropriate decision without in anyway diluting it commitment to raise quotas for OBC".

Medical services were hit hard at all the four medical colleges in Kolkata as junior doctors stayed away from outdoor departments. The agitating medico's of Amritsar came in for a surprise when it found support from former cricketer and local BJP MP Navjot Singh Sidhu. Addressing medicos he said he supported them. “Through he was all for providing quality education to the downtrodden, merit should not be ignored at any cost he said…

Experts said the chief of the BJPs had no guts to stall the reservation for they were scared that they would loose the votes of OBC who formed the majority so they should have instructed the Sikh, BJP, MP to talk against the reservations. The experts flayed Sidhu for the statement “Merit should not be ignored at any cost”, for they are still unable to find where merit was compromised. Given identical surroundings the experts said certainly the rural OBC would fair better than the city arya brahmins; what was lacking in the rural OBC getting a little lesser marks in the school final was they had no teacher, no time and above all not even a single proper meal a day. They had to take up after school hours some job to save the family from starvation added the experts. Who is meritorious?

The Statesman

Joining the medicos in their protests representatives of the AIIMS Faculty Association began a day long fast without striking work. “We are supporting the students and doctors and we want the govt. to listen to them” said Dr K.K. Handa general secretary of the Faculty Association.

“We have burnt these notices as an expression of our solidarity said. Meanwhile, the govt. said several of the striking doctors in hospitals have resumed duties.

In Patna, clash between pro-quota agitators, policemen and media persons injuring 35 people when they took out a procession from Gandhi Maidan which became violent. They were shouting slogans against those opposing the union human resource ministry’s proposal to introduce reservations for OBC.

Infosys chairman Mr. N.R. Narayana Murthy, today said increasing the number of seats in central elite institutions was not the solution to the reservations issue. He said “we should concentrate more on primary education and provide students with more nutrition, books and other facilities”; reported from Bangalore.

The experts expressed their shock at the statement of the AIIMS faculty Dr K.K. Handa's statement, "we are supporting the students and doctors and we want the government to listen to them".

This show the peak of arrogance of Dr Handa was Dr Handa working under the govt. or was he running the private institution was their first question. They strongly suggested that govt. should take stern and appropriate action against him for he was breeding insubordination for he says govt. should listen to students and doctors.

Secondly the media was reporting the pro reservation strike in a very low profile.

Third they said the views given by Mr. N.R. Nararyana Murthy shows that he wields power due to his fat economy! The arya Brahmin says increase of seats in not the solution then we suggest in the present number of existing seats let the govt. order 27 percent reservation for OBC. The govt. should be stern in this and should not roll back.

Times of India

"Why do defenders of merit not find money power a factor against merit? Medicos forum for equal opportunities.

Pro quota forces have begun to assert themselves and oppose the anti reservationists several doctors have formed an organization called, Medicos Forum. For Equal Opportunities (MFFEO), while various dalit and backward caste organizations have called for a joint meeting on Friday to chalk out their course of actions. The dalit and OBC organizations include, Mahatma Phule Foundation, Akhil Bhartiya Gurjar Mahasabha, Prajapati Mahasabha, Yadav Mahasabha, Dehat Morcha, Kisan Morcha, All India OBC Railway Employees Federation and Indian Justice Party.

These organizations have also invited students from Delhi University, Jawarharlal Nehru University and the medical colleges.

The medicos forum comprising members from various medical colleges including some top doctors of AIIMS, Maulana Azad and Indian Council of Medical Research have in a signed statement condemned, the anti reservation agitation Terming it as an expression of "upper caste chauvinism' the forum said the anti-quota agitation would deepen the caste divide by exploiting the insecurity among the upper caste youth”.

The forum has also objected to a section of medical students and doctors claiming to represent ‘all doctors’ and ‘all medical students’ and ‘supporters of merit’. It questions why the anti quota people have not opposed reservation for NRI seats, management quotas, governing body quotas and paid seat. Why do defenders of merit not find money power a factor against merit” they asked. It called a lie claims that incompetent doctors got in due to caste based reservations.

“Even with quotas one has to get past the basic cut off percentage for selection. Secondly all passing out examinations are the same for all medical students once they have entered the medical college” the forum clarified.

Offering similar arguments dalit leader Udit Raj whose Indian Justice Party represents scheduled castes organizations said “upper caste doctors leading the anti reservation agitation were not really for merit or else they would have opposed the admission in medical colleges thro’ donation or NRI quota”.

Several factors mentioned by the pro reservation groups are very much welcomed and appreciated by the experts. The experts praised their missions and wanted them to make more protests all over the nation. The experts also wanted the govt. to analyze the points given by them. They whole heartedly thanked Times of India for recording this news without bias. They all agreed on the fact that money power was never taken as a factor which always diluted the merit. As most of the arya Brahmins had money power they never brought out the powerful enemy of merit; viz. money and power.

Times of India

Former Chief Justice of India S. Rajendra Babu has embarrassed the government by declining to head the National Commission for Backward Classes even after his appointment was notified in the gazette. But Justice Babu said he was not to blame for this unusual situation because the social justice ministry had failed to consult him before initiating the process of issuing the notification in February.

“I got the offer for the first time only a day or two before they issued the notification”. Babu said, “even then I told them clearly that having held the high test judicial post in the country I was not interested in taking up an assignment under any ministry”.

So why did the government still notify his appointment as chair person of the panel, a post that has been vacant since August 2005. Babu who was CJI for a month in 2004, denied his decision was related to the controversy over the proposed OBC reservations in central education institutions.

The New Indian Express

In New Delhi medical students belonging to the backward communities demonstrated to demand reservations in all sectors. Once of their demands was to fill in the complete back login in quota allotted to SC/ST in all sectors.

Express news

The group of ministers grappling with the reservation row appears likely to suggest to the P.M. a three pronged approach to defuse the snow balling; the plans to increase seats in existing institutions and set up new ones, and implement the quota in phases.

In Jawaharlal Nehru University tonight at least three persons were injured in a clash between opponents and supporters of the quota.

A pro reservation march in Patna by over 500 doctors turned violent as some of them clashed with the media persons forcing policemen to resort to lathe charges…

In Delhi, Congress leader Ashok Yadav denounced the strike by medical students to protest reservation for other backward class student and warned about a growing pro quota movement.

“Next Monday the nation will witness the strength of OBC. Students of 36 universities will arrive. It will be a mass movement. He told a meeting attended by Delhi University and Jawaharlal Nehru University students.

N.R. Narayana Murthy says, “Before talking about increasing the seats the government should talk to the directors of the institutions and find out if they have enough resources”.

The experts came very rudely on N.R. Narayana Murthy for the following reasons.

(1) they now wanted appointment new directors and chairmen from OBC for they felt almost all the directors and chairman of IITs, IIMs and AIIMS are only arya brahmins and only one or two may be non brahmins but they are the henchmen of them. They also contented with the fact that the majority of these directors and chairmen of these institutions were appointed by the BJP government and how can they ever be pro reservation, they will be brahmin fanatics who are basically against non Brahmins and so automatically will be against reservation. How can any meeting with these caste fanatics will ever give any form of solution? asks the experts. It is the cunningness of Narayanan to give such suggestions when he was not asked for it said experts. They all wondered why Narayanan was so much interested about reservations. Basically all his expressions show his stand which is much against reservations for OBC, for the experts strongly claimed his backing is not only BJP, RSS and VHP but above all, the Kanchi Sankaracharya mutt where even today Brahmins and non Brahmins do not dine together. They are made to eat in the verandah and Brahmins in the well furnished hall. If even in food which is very short lasting thing they discriminate Brahmins from non brahmins, will they ever accept for reservations which is not a days fair asks the experts. So our Narayana murthy if happens to be our PM in a minute will cancel reservations and say these institutions are only for Brahmins who are meritorious and how can the non brahmins study in them, totally forgetting the fact these institution are fortunately run by the government and not by them.

Hindustan Times

A tough talking Arjun Singh mounted pressure for quick decisions on the OBC quota on Sunday by putting the onus of deciding when to place the issue before the cabinet on Prime Minister Manmohan Singh.

At a press conference the Union HRD minister ruled out a relook into the decisions – to increase quota for backward sections to 50 percent in that has generated a lot on controversy.

“It is entirely of the cabinet to decide. As soon as the P.M. allows it to come before the cabinet a decision will be taken” he said when questioned if the bill would come up in the ongoing session of the parliament.

Singh dismissed as propaganda the media's description of the anti reservation agitation as "Mandal II. There is no mandal II, III or IV it is all propaganda to vitiate the whole. He also came down heavily on the National Knowledge Commission, headed by Sam Pitroda that recently voted 06-04 against the quota.

The major role in these reservations for OBC is played by the media first says the expert followed by the legal circle. The media is making a false propaganda as if reservation was condemned by the nation by only projecting and publishing articles against quota. None of the stir or protest for the implementation of the quota or which are pro-quota are given any coverage. The experts condemned the media for such discriminatory acts. They said the media was always in the hands of brahmins but it does not mean they can do such discrimination in publishing about reservation. It is high time they correct themselves. If we non Brahmins vow not to buy your dailies what would be your position asked experts? Why is the media playing such a fowl game they wondered! They are not even keeping the ethics of the media they added!

Hindustan Times

Emergency services at government hospitals particularly at the AIIMS and the LNJP were crippled on Saturday as junior and resident doctors observed a 24 hour strike to express solidarity with the medical students protesting against the OBC quota proposal. Till evening the strike had claimed the lives of two patients who needed immediate attention. But the doctors refused to budge from their stand as news came in from across the country similar protests.

In Delhi Simrankelkar a 30 year old thyroid patient bled to death while she was being shuffled between LHMC and LNJP. Her father in law said, we took her to LNJP first where the doctors told us to take her to LHMC. At LHMC they told us to take her back. She died on the way back!

In the second case a 35 year old man, who had met with an accident two days ago and was supposed to be operated upon on Saturday died at LNJP. There were several others including babies who have been left to suffer in the stand off between the government and the doctors…

The experts viewed these deaths as barbaric and the striking medicos and the doctors must be arrested and treated as terrorists under the terrorists act. The court has the conscience to pronounce that govt. should pay them in the striking period. The experts asks what is law or judgment, in India they say it is manu dharma so under all circumstances only the arya brahmins will be supported. Even death is not condoled or the courts have not now ordered them to resume work. Loss of life of patients is not a matter at all. The deaths due to their not attending their duty which comes under essential services are not condemned but the court has the guts to ask the govt. to pay the striking medicos and doctors even during their period of absence. What is the action the govt. is going to take, those reported to have died are very young just 30 years and 35 years, who is going to take responsibility for their death. The expert says these are the ones reported how many more have died because of these striking doctors. Who will be held responsible for it? The experts want the govt. to take appropriate action against them. Further the experts ask can any one go on strike and demand for salary from the employer.

Deccan Herald

The anti reservation stir intensified in the national capital and several parts of the country…

In Bangalore a scuffle broke out on Sunday between an IIT – Alumni led rally protesting against reservation.

The experts made it very clear that 10% of the national economy is in the hands of the IIT alumni, who are very strongly against reservations for OBC in the institutes of higher learning. Further the experts added that these people have strong fall back from the Hindu religious mutts. They are very strongly knit with the Narayana Murthy of Infosys and Hindutva organization who practice caste in full swing. The experts claim for most of the protests they are the funding agents. They don’t want to loose the seat which has made them fat in economy. The media also have a nice understanding or to be more precise a strong link with them to publish it with many photos claims the experts as most of the media is owned by them.

Deccan Herald

The middle class deserves what is gets, it needs to vote if it wants politicians to take note of them. Sushant Sareen;

Use of 'reservations' as a tool for increasing the political support base suffers from a triple paradox. One, even though reservations is a populist measure it doesn't bring new votes. Second the moment a person is empowered he will no longer beholden to the party or person who empowered him. Third the politics of competitive reservations is useful only until the policy is implemented. Once reservations are implemented they stop yielding any dividend to the advocates.

Now the experts feel that the writer is using a very wrong term viz. "reservations is a populist measure" for they feel that only reservations can improve the down trodden which is very well known to the Brahmins and other upper castes that is why they protest and the media has taken such pains to blow it up and show, not only that, the media is also very careful on the other side not to even record the pro reservation protests or stirs.

The second paradox clearly shows a Brahmin mentally. Till the needs are got they will serve the feet, once they have everything they will severe (cut) the head. So they view the majority OBC like them for after all according to abidhana chintamani OBC are only their half brothers! Thus the validity of paradox by sushant needs to be chanted for OBC. The very status of them shows their loyalty to the arya Brahmins and other upper castes which had made them as downtrodden and poor!

Once the reservations are implemented certainly the beneficiaries would be OBC. What is the yield the dividend or the advocates expect? The have done a great social justice for equality can be achieved, discrimination would be certainly reduced though not eradicated claims experts. What experts feel is what more is to be achieved by reservation than a good education for the majority! Sushant say "so before we start off on a new round of reservations let us pass a law that no one will be allowed to go abroad for medical treatment. Experts ask Sareen why such a rule was missing all these days and what is the need for the implementation? We first wish to say majority of the OBC that too living in rural areas cannot even get the medical aid from a health care leave alone a good doctor or from a costly clinic. If this is a case how is that rule going to affect the OBC who are to get their paltry right of 27% in higher education asks experts.

Press trust of India, The Statesman

Medical students today appeared sharply divided over the proposed reservations for OBC in premier technical institutions of the country as demonstrations stirs burning of effigies and road blocks, both in support and opposition to move were witnessed in north Bihar's Muzaffarpur and Darbhanga towns.

While slogan – shouting students of Sri Krishna Medical College and Hospital (SKMCH) in Muzaffarpur burnt the effigy of Union HRD minister, Arjun Singh their counterparts in Darbhanga Medical College and Hospital (DMCH) chanted slogans denouncing Mr. Singh's cabinet colleague Mr. Kapil Sibal for his opposition to move.

Raising anti-Singh slogans the medical students at Muzaffarpur burnt the effigy of HRD minister accusing him of trying to divide the society along caste lines.

The professors warned that if the govt. did not shelve its plan to introduce quota for OBC they would intensify their agitation. On the other hand nearly 300 medicos of DMCH took out a march from Karpoori chowk to their alma mater in support of proposed reservation for the OBC staged a sit in and blocked the road in front of the institution for four hours from 9 am till 1 pm.

They later submitted an 8 point charter of demands to the district magistrate which apart from reservations for the OBC in medical and technical institutions included lowering the cut off mark for eligibility to take CBSE medical entrance tests for SC, ST and OBC to 40% and exempting them from mandatory interview for admission to technical institutions.

The experts said media was avoiding pro reservation agitations or stir so that they want to impose in the minds of the public by saying reservation for OBC in institutions of excellence was wrong. Thus the media was playing a very ruining role for OBC and no dharma of the media was ever followed by them. The media clearly showed by it acts it was only arya Brahmin property. Except one or two dailies all other dailies were only protecting and propagating that reservation would ruin the progress of the nations. The experts ask; what did they do all these centuries to the nation? The only most visible thing is that the poor has become poorer with the life span of the rural people has drastically reduced as they do not get any proper medical aid and nutritious food. All well trained doctors go to work in big hospital or fly abroad leading none to work for the rural poor so if some rural poor are themselves educated certainly, they would be useful to their area at least claimed experts.

Times of India

Gurpreet Mahajan. The writer is professor in politics, JNU. There has been a spectacular increase in the number of institutions seeking minority status in the recent past. In all of 2005 the national commission for minority educational institutions received just 380 applications form institutions seeking minority status. This year with the government announcing its decision to reserve an additional 27 percent of seats for OBC in all centrally founded educational institutions, the number of applications received by the commission is already more than 200. At this moment it seems that schools and colleges that had not previously sought minority status are now seeking that special status to ward off government intervention.

If the proposed reservation for OBC come into effect we can expect the number of minority educational institutional to steadily increase. As the number of open seats in prestigious educational institutions declare with increased reservation less brilliant Hindus will be driven towards community colleges that will come up to fill their need… For sometime now we have focused on the real possibility of caste identities becoming fixed markers of our identity on account of caste – based reservations policy. It is time to recognize that steady increase in reservation quotas is likely to trigger a process of desecularisation along with the consolidation of majority community identity and this is going to have serious implications for the future of secular India.

The experts are not able to understand the professor does not mind about the caste and about the caste society due to the laws of manu. The still demeaning factor is scavenging is allotted to a particular caste, the burial of dead and the rituals to another caste and so on. No one is bothered about it but why are they bothered about caste based reservations when caste based employment (or profession) exists! The experts feel because of caste based profession that there is more threat to secular India than caste based reservations. What has the secular India so far given to these (OBC/SC/ST) caste people but have only made them still meek and deprived.

Deccan Herald

Union minister of state for commerce and industry Jairam Ramesh on Sunday blamed the media for fuelling the anti reservation agitation across the country and reiterated, the resolve of the UPA government to "fulfill constitutional obligations"

"That reservation will destroy talent and merit is a bogus argument and I find media has been fanning the summer of discontent in the country "Mr. Ramesh told news man here at the side lines of the news conference. Claiming that apprehensions about the so called adverse impact of reservations is more a myth than a reality he said.

Tamil Nadu is an ideal instance where 70 percent reservations is currently in vogue and the southern state has demonstrated that it could adapt to any hi-tech changes as far as other states can….

We can aim excellence by a process of learning and intensive teaching and I understand there is a need for this" he observed, 'but besides religion, caste is also taken seriously; the Indian community and it is nearly impossible to deny the practice of casteism in the country".

The experts welcomed the views of union minister of state for commerce and industry Jairam Ramesh. He has correctly stated that media fuels quota protests. Further the caste in India is cruel than religion in India. It is more caste clashes than religious clashes experts said.

Experts still wondered why the media was not working against the caste discrimination. In fact the double danger about media is, it was favouring both caste and religion but only pro reservation was its topic of hate. It has played a serious havoc these days in misguiding the public, experts added. Some rule must be made to punish the misguiding and misleading media, experts urged. Media no doubt was very steadily fueling the quota protests.

Times of India

To strengthen the anti quota agitation students of Indian Institute of Technology are requesting their alma mater to join the movement.

"We are a mixed bag of students from across the country. We are asking our school and college authorities to start an agitation for the cause of merit. We are writing to the IIT alumni to chip in whatever way they can" Said Kaushal Mittal a student of IIT Powai…, We can have reservations on the basis of economic backwardness but not on caste".

Why the OBC were defined six decades ago classification was made on the basis of economic backwardness. Today a lot has changed we should redefine what constitutes the OBC said another student.

The experts said the poor knowledge of IIT-ians about caste system. Castes (varnas) were codified from the times of manu (given in laws of manu). The cruelty faced by them under the rule of Brahmins can also be seen if they have read proper history laws of manu clearly shows how other castes were ill-treated, discriminated and we denied the right to get educated. The wealth if at all earned by a Sudra pained the heart of a Brahmin so the laws of manu said he (arya brahmin (can usurp) the wealth from the other castes. Also the IITians has no concept about castes that is why they say six decades ago based on economic backwardness OBC were defined! Backward castes are non Brahmins not defined but created by Brahma. Now we ask the super power IIT-ian to first relinquish castes then the very reservation would become meaningless!

The role of caste which has been very cruel has today forced for reservations. Why IIT-ians do not want to banish caste? If caste did not exist certainly quota for OBC would never arise said the experts. The very anti quota stir was infact organized, funded and lead by the alumnis. Experts said some of them who witnessed the anti reservation protest in Chennai said the participants of the silent protest were living a life of luxury.

They came in hired cabs or own cabs and they had all comforts at their command and not like our protests where we would buy water packets for Rs. 1/- added the experts. They distributed 500 rupee currencies to the protestors. These people (IIT-ian) do not know the caste humiliations suffered by the OBC. One expert asked does any one of these idiots know the meaning of “Sudras”. Sudra means prostitute son and he is a lower caste being born from the legs of Brahma.

(1) How do these engineers feel about the birth of Sudras?

(2) These people only duty is to serve others physically and live in poverty!

More than 50% of the OBC are Sudras. How do these IITians act to solve these degraded caste by birth asks the experts? “Reading a stipulated syllabi and scoring marks based on that securing, a seat in IIT is got and after getting the degree they fly to west for money. This ends the problems of you all IIT-ians” asks the experts have you ever pondered over the evil of caste by birth! The need of the hour is reservations for them so that they too can lead a life of dignity and equality. You all protest saying equality is last, who should protest the OBC or you asks the expert! If they start to claim equality and protest for it India would face only a civil war and this may lead to the blood bath in India, so the expert advice “do not anger the quiet people and continue to deny their rights”.

The New Indian Express

HRD minister Arjun Singh is a member of Other Backward Classes (OBC)…

That is according to caste certificates officially sealed and delivered to express for a fee ranging from Rs. 200 to Rs. 500. These certificates can now be used to get all benefits from job to colleges admissions available as per government policy.

The express traveled across the three districts Lucknow, Kanpur and Faizabad and obtained OBC certificates of well known public figures….

All these certificates can be obtained in a day from one place but the express team tested it for over a week obtaining one or two certificates on one day and tried the same the next day and the next, to see how easy it is.

From Lucknow itself the news paper obtained a pile of blank but stamped certificates of caste, motor vehicles registration you name it. On these blank form officials signatures can be forged. But the express got the certificates duly issued by a proper authority with serial numbers and actual signatures. Names of high caste prominent people were deliberately used to explain how vulnerable the system is and how blindly the certificates are signed…

Just contact the relevant tout give them the name, address, the caste and the cash. Certificates are home delivered with the three to five hours, no questions asked no one cares to check the documents no interviews, as padre process. Nothing….

Experts say these can be set right in a raid and by giving proper punishment to the concerned. So from this the Brahmin has no integrity so if he/she can get some benefit do not mind to cheat the govt. It only shows the two sides of the coin. But media giving such things and they are teaching to public to lead a life of immorality is highly questionable said experts.

The Indian Express

The issue of reservation raised question of inclusion and equality. We ask whether those earlier excluded from participation in society are being included in meaningful ways…. The effect of the polity has to be to create a system that allows for easier mobility and one which punishes those who unlawfully prevent others from surging ahead.

A Swedish Sociologists, Staffan Lingberg in mapping changes in Tamil Nadu; over 25 years period. He sees progress steady if slow. During the first stint of his study (1979-1980), the untouchables and lower castes were not allowed to sit in the presence of upper castes nor could they share water or sit in the same tea shops. Noting changes he and his colleagues conclude both lower castes and women were substantially better off. The reason varies from land reform… Other research tells us that if SC/ST student passes high school, his chances of getting into professional education is better than those of a non reserved category student IITs are coming close to fulfilling SC quotas and more students from this category are able to compete…

While the change in inclusion can be assessed by evaluating quantitative indicators of educational evolvement and performance what does one do to assess whether inclusion at social and psychological level has taken place? The real goal of affirmative action policies is to integrate the lower caste individual into the “general category”. It is not to empower him / her as a dalit or lower caste person. There is a stigma attached to being lower caste in Indian society the goal of policies of inclusion therefore has to be to get rid of this stigma and of differentiation and discrimination based on caste. The American civil rights act 1964 had a similar goal: citizens should not be discriminated on the basis of colour, race, religion, gender and other primordial identities. It was an anti discriminatory law.

In recent times, certain dalits may wish to valorize dalit culture and identity yet the majority would be happy if the stigma associated with low caste status and polluting occupations disappears. They do not wish to be known as a dalit engineer... The use of reservation as a tool to achieve equality often forgets this longer objective of obliterating stigmatized identity.

The structural changes ushered in by reservation policies need to dwell on how this can be achieved. Three possible routes are first, the empowering effects of education which makes people capable of competing with others on their own merit. Second the rapid introduction of Technologies do away with degrading occupations such as scavenging. Third may be most important and demolishing of the prejudice of the upper castes.

The experts said if reservations are to be stopped “Caste away the caste”.

The Times of India

8000 crores question high price for buying peace on quota issue. The only suggestion given by the experts is to keep the seats as it is and given 27% of reservation for OBC for the sake of equality, lack of knowledge.

Times of India

The national knowledge commission appointed by the Prime Minister in a 6-2 opinion had opposed reservation expressed “displeasure” that government did not consult them before going ahead on quotas for OBC in central institutions. If the six members felt that they had been ignored they should have had the self respect to resign immediately.

It is a common and respected convention world over that members of nay commission or body appointed by government who oppose state policy and publicly express their opposition must resign immediately. Clinging on the positions in government bodies and issuing statements deriding constitutional provisions are appalling and unprecedented. Two members have stepped down rather late in the day…

It is in this context one needs to understand the relations between knowledge and reservations. The majority of the commission has not provided any empirical evidence that reservations will harm the accumulation of knowledge.

The commission has showed complete lack of understanding of social concerns. Mere technical or business expertise does not constitute expertise to comment on social issues. The commission seems to be preoccupied with development of business rather than the business of productive activities that contribute to development. It is said that they have become opponents of provisions to ensure equal opportunity to weaker sections.

Their failure to understand Indian society and the immature way of expressing their own concerns have made the commission just another mouth piece for the anti reservation lobby, there by trivializing sides rather than advise government on basic issues relating to educational empowerment of the deprived classes. In any cases, they do not have a mandate or authority to comment on the time tested policy of reservations.

It is interesting that the commission expressed its opinion on reservations on the day when the results on this year's civil services exam were announced. In the open category an OBC candidate ranked their and there were six reserved candidates including an SC who were selected under the open category for IAS. It clearly demonstrates that SC, ST and OBC can compete with the best provided equality of opportunities is created.

It is public knowledge that fees in educational institutions have increased manifold over the past few years. This commercialization of knowledge needs first to be settled rather than harping on reservations. If access to knowledge is not available then; how can one ever think of others as equal? It is unfortunate that the six wise men have thought of making education accessible by regulation of fees and instead tom-tommed their position against reservations.

This article of Mr. D. Raja is actually a proper and an appropriate assessment of the so called National Knowledge Commission. They should atleast be now ashamed of their act and resign on moral grounds demanded the experts. The experts said they would make a strike or protest until they are removed from the post.

Hindustan Times

The write is a professor of Sociology, School of Social Sciences, Jawaharlal Nehru University. When the abbreviation OBC is spelt out it refers to other Backward Classes and not castes.

This is the meaning that the constitution of India gave to the phrase. There is a world of difference between the two and it is only by means of deft political jugglery that the impression has been conveyed that the constitution is interested in uplifting a section of the population that can be labeled under the totally spurious term called backward caste.

As the constitution was interested in a lot of backward classes it did not easily equate them with any cluster of castes nor did it privilege a solely a caste based criterion to determine them. Backward classes could mean the poor or village based artisans or the unemployed or those who are in remote regions or those who are educationally deficient or those who suffer from food deprivation. In all such cases the focus is on actual people and not communities. As those answering to such descriptions could come from any caste the founding figures of the constitution very purposively used the term OBC strictly in terms of classes and not castes.

Contrast this with the way the scheduled castes were classified by the constitution though the job of determining which castes were untouchables proved tricking in all such castes, only the caste criterion was the determining factor. Economic or educational concerns had no role to play when it came to deciding on the listing of the scheduled castes. So when the constitution clearly states classes and not castes as in the term OBC it is deliberate and not an oversight.

As there were too many competing factors to content against, the first backward class commission set up in 1955, under Kaka kalekar decided to fold up its proceedings after rounds of futile negotiations. Its members found no clear criteria for determining backwardness that would be free from counter claims. In spite of long deliberations the task proved hopeless and unachievable and the conclusion was more or less notified by the government of the day.

The mandal commission seemingly had no problems at all in designating the backwards. This is because it cleverly masked the term class with caste and put forward set criteria where economic and educational backwardness received very low points compared to social backwardness. As the definition for social backwardness was very fuzzy, it was possible to score 12 points on the four criteria with in this subset. There was no need to go further and consider educational or economics factors as only 11 points were needed to be eligible for reservation. Space does not allow a full elaboration of the various categories under which backwardness points can be awarded on the mandal formula but it needs to be noted that the economic criteria get the least emphasis. While on social backwardness it is possible to easily score 12 points the best one could do in terms of economic backwardness is only three. This is surely a travesty of what it means to be backwards.

The points system also contains other constitutional improprieties. When considering social backwardness, full three points are given if a caste practices child marriage. In this case you have a clear instance of actually being rewarded for breaking the law. What could be more unconstitutional? Further again with in the category of social backwardness a caste scores three points if other castes have a low opinion of its standing in society. As is widely known to specialists in this subject, no caste thinks well of any other caste. Even among Brahmins there is no unanimity. All too frequently charges of being imposters, pretenders and worse are freely traded among them.

This bias in the points system clearly reveals that the mandal commission did not want to draw attention to the fact that is intended beneficiaries were really not economically and educationally backward. This is to be exerted, as the castes forwarded for reservation by Mandal were never really discriminated in history as the scheduled caste were that they chose not to seek education and urban jobs as the so called forward castes still does not make them a once persecuted category or community.

It must also be remembered that those so called backward castes were often the rural dominant castes and they were perhaps on occasions were brutal in their persecution of the scheduled castes than the Brahmins were. The deprivations resulting from the odium of untouchability were the highest in the South India a region where the traditional dominant castes are in the main the designated backwards of today. Interestingly Mandal activists ignore this feature of rural castes interaction.

At the heat and dust of political controversy around the issue of reservations for OBC obscured from view that it was really all about class and not caste the mandalites were able to win valuable political points. They also spuriously occupied the moral high ground by claiming that they were adhering to the constitution, when infact they were calmly subverting it by making backward class synonymous with backward castes, the Mandal commission opened the gates for identity politics of the kind that would be very hard to reserve.

It is just as well that the Supreme Court has asked the government to clearly spell out what is means by backwardness and to whom the label should apply. Hopefully in the course of this exercise the unpardonable equation of backward classes with backward castes will be exposed and we can also breathe easily that the constitution has not been breached.

Experts some of them professors of politics and sociology who were also Ambedkarities critised that a JNU professor of sociology did not know the fact the very father of our constitution Dr Ambedkar himself has very clearly stated in the constituted assembly that the word classes means castes. That is why deprived classes became deprived castes and the author of the constitution has himself very clearly stated that the terms castes and classes are synonymous. How can D Gupta say that the constitution is breached by using the two terms backward classes and backward castes as synonymous which is clearly erroneous for Dr Ambedkar has clearly said both means the same. Here the expert once again state that legal untouchablility is being the powerful tool used by the Brahmins to win everything. For in their opinion law is nothing but laws of manu. That is why experts went on further to substantiate even today the OBC face educational and economic untouchability, for it is clearly said out in Manu [ ]. That is why in many of the centrally run educational institutions of higher learning that OBC are discriminated educationally harassed and tortured even if they are much more superior to manu brahmins.

Several of these OBC are paralyzed in their activities or development in these institutes as they are affected by the legal untouchability and educational untouchability.

Some of the experts felt that this sort of problems faced mainly by the OBC either makes them subdued or at times very much depressed or at times makes them go abroad. When this is the case of intellectual OBC working in these institutions one can imagine the plight of OBC/SC/ST students! That is why the dire need of the hour is 27 percent reservation for OBC students is these institutions. The govt. should also contemplate to give employment in the teaching and research posts for OBC, SC and ST then alone it will help the OBC/SC/ST students to learn without fear and discrimination. If they are the teacher with so much of hatred (seen from the anti reservation protests) they will ruin the OBC, SC and ST students by mentally torturing and harassing them and failing them. Till the educational untouchability coupled with legal untouchability is wiped out, the reservations in these institutions must be carried out faithfully. Thus from the times of manu majority of the brahmins never were economically backward. It would be very much surprising to know that majority of the nation wealth in the form of money, land and jewellery which can be roughly estimated as over 40% of the nation wealth are in the hands of the Brahmins claims some of experts. One of the experts said atleast one forth of the Madras city wealth is concentrated in the Kanchi mutt. More so to claim that economically well of brahmins are staunch supporters of Kanchi mutt and are anti to reservation for OBC. This group includes the alumnis of IITs.

For it is said 10% of the national economy is in the hands of the ex IIT ians and they are the group who are very much against reservation and some of the experts claim they were the ones to organize the anti reservation protests. In fact they help the media to blow up the protests beyond proportions.

They (brahmins) identify themselves as forward class of people and infact as claimed by gupta they have become the forward caste people or other classes or equivalently other castes. As said by gupta instances of caste clashes in rural areas! Expert asks which brahmin ever lives in a rural area. Maximum in every village their can be almost one family who work in the temple associated with the village and eat that temples income. Their close relatives would be big shots in the city. Secondly as observed by gupta, it is wrong OBC wish to stay in the rural areas, in fact the majority of the migrant labourer are only from OBC/SC/ST.

They with no other employment work, as daily wagers, employed as painters, mason, waiters in hotels or as road layers, sweepers and drivers. They are the worst affected sections of the people we (experts) ask Gupta; can he find any brahmin working as a scavenger and carrier of night soil in the head, but there are several backward class people who do this job. How many sweepers (including women) who earn less than Rs. 50 a day are among Brahmins (or Brahmin women)? But there are many from OBC women as well as men who work as sweepers and earn less than Rs. 50 a day without any other monitory profit. How many of the construction labourers and road layers are from brahmins? These OBC due to agricultural failure have taken up these professions. Further they have no education. Neither are their children educated so they have no other choice but to take up these jobs mostly as daily wagers to keep themselves away from starvation. Thus these OBC are not only backward castes in economy, employment but above all in education. Thus when the govt. says the backward classes or castes it implies these people and can never imply the brahmins. For the expert says one or two families which claim to be brahmins and who live in these rural villages eat every day nice food offered to god in temples with milk and ghee and get as pocket money which the people (public) offer for the gods. Thus can these brahmins be called backward they have good food, a good profession of archakas and above all their children invariably get good education. Do not miss the difference warns experts.

Further experts claim that in Tamil Nadu, mainly there is very less and no suicide of farmers the only reason attributed by the socio scientist for this is that Tamil Nadu is the only state which enjoys the maximum percentage of OBC and SC/ST reservation viz. 69% in India for many years. Thanks to the govt. for implementing such reservation and saving our agriculturalists claims experts.

Thus if the govt. has to save the farmers and get a way to see a developed nation the reservation for OBC is very essential. Further it is very important to note the nation for over 60 years even after independence has remained only a developing nation because the education and the economy has not reached the rural people above all OBC/SC/ST. The nation wealth is still only in the hands of the 3% Brahmins. Even when govt. tries to implement some form of reformation and means to improve the majority who belong to OBC the so called Brahmins using the media to full extent stall it. Thus experts have taken up to win the implementation of reservation to OBC in higher education in the institutions run by the central govt. even by blood shed. They are determined to get reservations.

Hindustan Times

Deveshkapur is professor at the University Of Texas, Austin, Sunil Khilnani is the author of “the idea of India”. No aspect of India’s developmental experience has been so marked by disparities between rhetorical ambitions and actual achievement as our educational system.

The importance of developing our human capital – not just for the intrinsic benefits it brings to our individual citizen, but also for our collective economic development – has been long and widely recognized. And yet we have been unable to realize this national imperative. While the reasons are multiple a very large share of responsibility must lie with those responsible for managing our educational policy.

Since independence the ministry of education has been one of the weakest ministries, for the simple reason that the opportunity for rents was limited compared to other ministries railways, industries and agriculture. Minister weakness has been reflected in the poor quality of its leaders. With the notable exception of Maulana Azad our education minister have been non entities or party hacks, lacking any vision or understanding of the vitally important charge they hold….

Take again the fact that the vast majority of India's college students over two thirds – go to arts commerce and science colleges places with poor facilities, abysmal teaching no accountability, which deliver a caricatual education or the fact that the vast majority of young adults (more than 90%) will never have access to tertiary education….

But the barely 3 percent of the higher education budget that goes for IITs, is a grim reminder of who we really care for. Indeed today we can't even get highly paid government teachers to even show up, to teach…

A moments reflection might lead one to think that treating the IITs and IIMs as navaratnas – insisting that they become self financing and this re allocating the resulting annual savings of nearly 700 crores to wards funding schools and colleges where ten thousands of depressed classes are today getting a pathetic education as well as offering more targeted scholarship and funding other access routes may actually serve social justice more effectively than simply indulging in a politics of symbols. The experts viewed this article very differently. They said we can immediately close down all the 7 IITs and all IIMs or sell it. Even the land must not be given freely to them they have to pay rent. Some of experts add when the govt. is trying to sell SBI, NLCs and what naught why should not the govt. sell the IITs and IIMs which are only blood suckers of the national economy with minimum benefit to the majority and the maximum benefit to the minority arya brahmins? The experts said the gradual close down of IITs or change them into state run technological universities. This can be easily achieved with in 5 years. Further the experts added those who get education also leave the country and serve other nations. They are also wondering why the nations hundreds and thousands of crores of rupees every year is spent for one class of people who do not also serve our nation and they have become fat in economy, they say 10% of the national economy is in the hands of the IIT alumnis and they are against reservation, in a way terrorists for they are the once who are the brain and money behind all these anti reservation protests said the experts. The only nice solution would be to close down all IITs and IIMs with in 5 years

The Telegraph

Whether or not it lowers the quality of students, Arjun Singh's quota plan is likely to hit the standards of teaching at the Indian Institutes of Management.

Many faculty members say they would consider leaving if the proposed 27 percent quota for other backward class students is applied at B-schools.

The best faculty members will leave and there will be an exodus said a professor at IIM Ahmedabad.

A young teacher one of 18 who joined in the past two years said he and other colleagues have already began exploring job options. They believe that once the quota is put in place academic standards will fall.

We may go abroad or join other institutes said another teacher, a Ph.D. from Texas University who has begum feeling frustrated over his decision to return to India. This young professor an IIT graduate who has taught for eight years at the London school of Economics said I am seriously exploring other avenues I am here by choice – I can leave if I feel that the purpose of coming to IIM Ahmedabad is not being served".

The Professor who did not wish to be named, said many of his colleagues share his feelings. "For most of us who have taught at foreign universities there is no dearth of opportunities".

A senior professor Ramesh Bhat confirmed than an exodus was a genuine possibility.

I joined IIM Ahmedabad because I saw growth opportunity here. This is a world class institute. But its character will change if nearly 50% students are selected on the basis of reservation, not merit which is the hall mark of any institution of excellence. A professor suggested that many people from the Other Backward Classes (OBC) are so well off that they don't need reservation.

“If the centre is so concerned about the backwardness of the OBC who does it not introduce reservation at Tata Institute of Fundamental Research or ISRO (Indian Space Research Organization), they cannot because quality will be compromised”.

A senior professor in the Ahmedabad campus Anil Gupta has a different view on reservations.

He argued that private schools are in effect reserved for the rich rather than the meritorious because of high fees they charge. So if reservation is wrong then this is equally unjustified.

I am all for affirmative action Prof. Gupta said but only first generation OBC (those not well off) should enjoy the benefit.

He suggested that all quota students at the IIMs should be provided free coaching by the government. They should also spend one extra year at the institutes.

All the experts welcomed that if they resign IIM – Ahmedablad, the govt. can contemplate in filling the arising vacancies by the OBC, SC and ST which would certainly be a blessing in disguise. They however said their view of merit will go down or any such bunkum is never applicable.

They also said as pointed out by Gupta that certain private schools are only for rich like wise he said these are certain schools which are only for caste religion which is also equally unjustified. Govt. institutes at least now should learn to respect the constitution and there by implement reservations for OBC/SC/ST without any grudge.

Once quota is put in place the experts said the brahmin population will not get admitted and not the academic standards will go down as they claimed. Again the experts wanted the faculties of IITs to resign so that they can make way for OBC/SC/ST at least now so that the filling of by OBC students and appointing OBC/SC/ST candidates in teaching posts go hand in hand. Still the expert questions will ever these fellows; who has enjoyed the benefits of faculty positions in IIT and IIM resign their posts? It is a big question mark they said!

The claims of Gupta that OBC are economically well off is a false statement said the experts for they said when OBC form the majority of the nations population, they do not own over 1% of nations wealth but is shocking to see the arya brahmins who form 3% of the nation populations have a command over 10% of the nation wealth or economy. Further they asked what is the percentage of arya Brahmins who live and work as plumber, scavengers, sweepers, and doing all types of menial works?

The Statesman

Quotas and reservations in India are merely populist measures to win votes for the next election. They would never churn up the backward classes to fore. Even if some people were put on the creamy surface they would go down unless they are intellectually prepared and enabled to stay on the top in a competitive environment.

The quota and reservation system is one of the biggest frauds being committed against historically deprived classes. This is the way to keep them down forever.

The best way to raise them high is through a system of challenge and response, by providing OBC access through merit and need based scholarships so that backward class kids can compete in the market place. Let's keep in mind that historically deprived people are not genetically deprived.

But neither corporate America nor any social group has ever advocated that diversity and social justice should be achieved through a system of quotas and reservations for historically deprived classes as it is being done in India.

The experts first claim that in no nation in the world like India there is caste system which is all the more based on by birth. The experts further said they are codified as laws in laws of manu. Let the govt. first ban all the castes and ban the laws of manu. Even to this the arya Brahmins would protest as they protest for reservations. Further when it is black versus white the discriminations is seen but in brahmin versus non brahmin one cannot find the discrimination in general appearance.

Further as said by Batra, OBC are not genetically deprived; "Brahmins still claim merit, will go if reservations are given, what is the equations put by these cunning arya brahmins" asks the experts. When only the talks about reservations for OBC the brahmins act as if they are bothered about the society. Are any of the arya brahmins ever bothered about the society? What do they mean when they say society? The brahmin society! For when a OBC becomes his competitor, who has been deprived of everything in this unjust rule of brahmins when some paltry quota is implemented! This pains each and every Brahmins and deep concern over the nation (brahmins alone form their nation, for nation does not include OBC or SC or ST for to them nation is not as important as the people especially deprived class). If genetically they (OBC), are not deprived then how will reservation have any impact on merits, these brahmins always are contradictory. The expert claims such statements are surplus in the laws of manu.

Dear Batra all your claims and comparisons cannot fit in the Indian contest for in the Indian contest a non brahmin is brainless cannot possess money, cannot be rulers if they don't take the ministerial assistance from the brahmin. OBC/SC/ST have no merit and above all according to laws of manu, born to serve the Brahmins meekly. But fortunately they are not the sons of the soil! They are not deprived but their opportunities denied by high caste to be more precise their opportunities are snatched!

The Telegraph

Pockets of unrest against quotas began mushrooming in several states today sucking in some of the brightest among the young and raising the possibility of a throw back to the divisive Mandal agitation day… In the age of 24 hour television that relentlessly replays striking images, a missing phenomenon during the Mandal agitation/protests is that even small groups could create an impression that a mass movement has begun to take shape.

Human resource development minister Arjun Singh has attacked the Knowledge Commission Chairman, San Pitroda for opposing reservations in education.

The chairman of the knowledge commission says he was unaware of the constitutions amendment (on reservation). If he is that ignorant how can he be the chairman of the knowledge commission? Arjun asked a Television channel.

The experts said the media more so the television channel is playing a very destructive role in blowing up as if the mass movement has began to take shape against reservation. This was very much condemned by Dr Ramadoss. In Tamil Nadu certain

T.V. channels more so only backed by the arya Brahmins was only showing anti reservation protest several times a day none of the pro reservation protest were ever telecasted by them. Thus the arya Brahmins are trying to create an impression that there is lot of people against reservation. They cannot get new photos repeated only a few photos plays a role but the repetition is very frequent trying to make a mass impact against reservation said the experts. Print media in majority of the cases showed only anti reservation trend. Any stir or protest made, pro reservations were schematically hidden by the media. The media is for getting the very fact those who are for reservation happen to be the majority in comparison with the paltry percentage who are against reservations.

The New Indian Express

… According to recent international rankings, only the IITs are noted at all and even these are not any where near the top of the charts…

The tiny quality sector in Indian higher education is severally undermined…

Serving the needs of mass access and social mobility for disadvantaged groups is important but it is not the only goal of higher education…. Research universities everywhere have some common characteristics.

Meritocratic values, “A meritocratic university seeks …”

India is truly at a turning point. If the nation is to fulfill its economic and technological potential in the 21st century, it must have an elite and internationally competitive higher education sector at a top of large and differentiated higher education system…

The author statements that India international competitiveness would be spoilt or lessened or lowered cannot be accepted at any cost argued the experts.

Further when even 60 years after independence they have not achieved any thing solid says the experts and all institutes of education have been dominated only by them. In spite of this what merit they have exhibited asks the experts. So by giving reservations the situation is not going to be changed and his high claims are absolutely biased rooted on prejudice says the experts. When majority of the nation is denied good education how can any one expect the development of the nation that is why the nation still remain under developed!

The Times of India

Doctors who withdrew their anti quota stir following a supreme court order banning all kinds of strikes and dissent on the reservation issue on Monday expressed a grievance before the apex court that the govt. did not apply the ban on pro-quota agitators we with drew the strike following the May 31 order of the apex court. But exactly four days later, pro quota supporters held rally at India gate and the govt. did not take any action against them on the basis of the SC order counsel for the pro-quota doctors.

Can there be two sets of norms – one for the anti quota and the other for the pro quota on the basis of a single order of the SC, counsel asked and the govt. was making it appear as if it was tacitly supporting the pro-quota elements.

This very article shows how the anti reservation protesters are biased over the protest of pro-quota people. They go to SC for clarification. “Further govt. move is supported by pro quota people they have not paralyzed for 19 days, the out patient services and hence the usual functioning of the govt. hospital” added the expert. They also pointed out, these are the places were law is not in the hands of the courts but very much with the anti reservation protesters. If reverse was taken up or challenged by us the court would have clearly said they are with the govt. and it is not against constitution! Thus the laws of manu is still ruling us it is a pity added the vexed experts.

The Statesman

The government should use the services of teachers and students who have not benefited from the quota system for its important projects writes; Chittlabrata polet.

There has recently been quite a storm over the HRD ministry’s decision to increase quotas in govt. jobs for SC/ST and OBC candidates. Such quotas however should not prevail even after 60 years of independence. Jobs should be open only to talented people. There should be no compromise over quality. Inefficiency should not take precedence over excellence. Institutions would collapse if under qualified persons are employed by virtue of their quota privilege … Quotas promote caste and segregation and contribute to a policy of national disintegration….

Experts feel that the author who is professor of History, Jadavpur University claims that marks cannot be the criteria in case of some people who mainly do not fall under any of the quotas. His article supports them strongly and recommends them for teaching jobs.

Expert claims that; those who get lesser marks and do not get the stipulated marks needs consideration. What is the injustice they see on giving those deprived castes reservations asks experts. Further these deprived classes or castes are the ones who are not only deprived of basic good education but also they are the one’s who are deprived of employment, economy, and basic food. They lead a life of misery and poverty and a rural non Brahmin working over 8 hours a day and cannot have a square meal.

If not under quota privilege but under the caste of Brahmins even if 50% are under qualified the experts says the institutions would not collapse. Only if quota privilege is under even if the employed are above qualified still according the arya arrogant brahmins the institute will collapse.

The Statesman

The reservation policy has a long history and a natural cause but a politicization of caste and caterisation of politics lent largestly to a large number of educated urban Indians becoming vociferous against the prevailing reservation policy of the union govt.

The protest drew more doubt from educational urban Indians which would explain why some have preferred to call it a movement of (rural) Bharat against (urban) India…

The anger and frustration among a large number of educated deserving students directed against the pious reservation policy for Dalits has not prompted friendly relations between the advantaged and disadvantaged. This should not have happened merely because of a nations’ pledge to provide solace to its backward communities…

The author Samitkar who is a teacher of Sociology at Presidency College, Kolkota, view is not appropriate says experts and give the following reasons.

(1) Caste always exists not only due to reservations. For even in apartments if there are say 25 owners 2 OBC or SC/STs and 23 brahmins. The two OBC or SC or ST will be socially segregated from all activities of the apartments. They would be treated as untouchables, this occurs even today in the heart of the cities in T.N.

(2) On the contrary if 23 OBC/SC/ST live in an apartment with two brahmins. They will not freely mingle with OBC/SC/ST. These two houses will be together though there is no reservation policy in the apartment where they live.

So according to the experts caste code in India is so strong be it in a city or in a rural area that free unreserved intermingling is very rare between brahmins and non brahmins. This does not occur due to reservations. So atleast if govt. gives reservations for OBC/ST/ST especially in education that will atleast help the deprived class from the sufferings of the social stigma.

Not many are aware that 12 years before V.P. Singh introduced reservations for OBC in 1990; the issue had cropped up in a virulent form in Bihar in 1978…

A full fledged meeting of the committee was convened and the proposals put before it. The formula was unanimously accepted and Bihar government in 1978 issued a notification. In this manner reservation for OBC was introduced with cooperation of Janata leaders belonging to all castes. The lesson was that a fair compromise always finds takers. If V.P. Singh has handled the issue more tactfully when announcing his policy of reservation for OBC in 1990 it would not have caused so much turmoil in an emotionally charged atmosphere different sections of the society are not able to appreciate the views of others. If however questions are debated calmly people can understand and appreciate views advocated by the other side.

In the past weeks we have seen a raging debate on the proposal for reservation for OBC in institutes of higher education. This was yet another example of unnecessary hysteria blocking minds from looking at issues in a balanced and objective way…

The experts view the article in the following way.

(1) The word reservation itself makes the upper caste especially Brahmins loose their balance of mind.

(2) Almost all articles written by them would state either equality is lost or merit is at stake or nation will face a degeneration or merit and quality should not be compromised. All articles view OBC reservation will make merit less in the educational institution. Experts added that they have become fed up with this. They said look at T.N. which has the highest percentage of reservation no one can ever say the merit has been compromised to the least degree. They said equality of OBC with other upper castes has been lessened to a very great extent. Further among students very rarely we see the caste discrimination.

Deccan Chronicle

Members of the Madurai Medical College students association staged a demonstration demanding that the reservation for OBC should be implemented in higher education and requested their counter parts in New Delhi to give up their protest.

Dr Gerald, a member of the association said the support for the reservation for OBC in higher education was growing day by day.

How can they say that only the so-called upper class students are qualified to become outstanding doctors".

"It is only in the six states where the quota system is in effect, the public health system is good. Tamil Nadu which has 69 percent reservation has an outstanding public health system. So how can the medical students in Delhi say that reservation in higher education will ruin the system, he said. Infant mortality is least in these states he added.

Similarly in the IIT Madras among the 450 faculties only two are dalits and 50 are from backward classes.

"If doctors in Delhi continue their protests we will show our protests more vigorously when college begins in June, said a protesting doctor.

The experts said by giving OBC reservation in 10 years time certainly India would become a developed country when the arya brahmins had made India remain a developing country for 60 years, as OBC are the sons of the soil they would make the nation a developed one in a decade!

Deccan Chronicle

Human resources development minister Arjun Singh made a first announcement about reservation in India's higher educational institutions. It was a side comment and would have hardly elicited a major response atleast not at that stage, had it not been for the Mandal – phobic media response that wasted no time in stoking burnt out fires. A small group of medical students (please look at the newspapers from day one to trace the movement) came out atmost reluctantly to protest against reservation. The television channels swooped down at them – Eureka! Here they are! and each statement, each movement was reported over and over again, 24 hours of a day. News papers even published the photograph of one student as their discovered icon and carried her interview with her to establish the point.

Cameras bring out the best and the worst in us ordinary mortals. One can see that when the political leader walks by print journalists with a "no comment" but stops suddenly when he sees the television cameras and waxes eloquent on all possible issues. In fact he often forgets that what he had just told print reporters off the record, was now being aired by him live on television; but what the heck, how else would he become a celebrity? So it was no surprise when young students started coming out with the running media commentary being a major incentive. Effigies were burnt, students photographed facing water cannons (that by the way is part of the protest as all who have walked on the streets on issues know) and a general atmosphere drummed up in the country; unscrupulous politicians, vote banks and justice. Merit they say must count. It must but does money determine merit? But those who are able to study in the best of schools and have a better chance than the other in clearing entrance examination! Or by those whose parents have the money to pay the capitation fees and get their wards into medical and engineering colleges that have opened all over the country. Are the doctors and the engineers coming out of such colleges; meritorious, or have they got admission because of the money they have? Is it not true that the student passing out from a post school with the best laboratories and teachers have a decided advantage over the students passing out from a poor run down building that passes for a college in Banda. So if the first gets through the examinations for higher education, is he really more meritorious than the other?

But no one cries when one student pays the money and gets into a medical college even while the more meritorious ones line up outside simple because their parents cannot afford the capitation fees. Instead such colleges living off hard cash are eventually recognized by governments and their degrees produce doctors who failed to clear the entrance examinations for the competitive institutions, but had the money to get admission into these institutions. But when it comes to reservations without money, but on the basis of social equity, there is a hue and cry as if all the opportunities for a handful of students – when compared to the millions who do not even have the luxury to dream – are being taken away.

When India became independent, the foundation fathers decided in their wisdom to introduce reservations for SC and ST. It was felt necessary at the time to give the dalits a level playing field, and for the country to atone for the centuries of discrimination and oppression. It was decided that unless the government stepped in to give the victims of gross untouchability a hand forward, they would find it very difficult to climb up the social ladder and join the rest of India in its march. Reservations for the most visibly backward section of the society were thus introduced. Several years down the line it was realized that the backwards in North India were not being able to keep up as well. Victims of their social status as well as OBC had been largely ignored in North India because they were essentially sandwiched between the visibly oppressed Dalits and the upper castes. The resurgent Congress did not really take their plight on board, and while Muslims, Dalits, Brahmins and Thakurs were attended to as influential vote banks. The OBC divided with in were left out of the government sights. I do not want to go into the social political history since independence…

Sometimes the social indications are more powerful and in the case of the reservations policy for both Dalits and the OBC quotas have brought in a level of empowerment, a desire to be educated, a willingness to move out and seek jobs, a greater awareness of their rights….

This government has allowed itself to be dragged into an unnecessary political controversy when the need of the hour is a clear statement that the forgotten India needs a helping hand, that caste based oppression continues to exist that merit is relative in a country where millions go without schools, teachers, laboratories, and that the government role is to strengthen existing policies removing discrepancies, and shaking out the cobwebs for a more responsive and dynamic reservation policy.

As for us in the big cities and the medical students neglecting their patients, get real there is a big world out there that does not revolve around the "us" but belongs to "them". Merit can have meaning only in a reasonably level playing field, for otherwise it is just a hollow cry signifying nothing.

The experts feel happy over some of the factors recorded by Seema Mustafa. Why are they talking about merit which happens to be a hollow word according to them? They have not said what really merit means? Merit is itself relative for it is judged on the beholder! Secondly the experts claim because of reservation alone the farmers suicide in T.N. is very less or nothing in comparison with the farmers suicide in the north where reservations are protested, happen to be very high.

Hindustan Times

As the debate on OBC reservation spreads. I've unearthed a few interesting facts which raise pertinent questions…

Beyond the issue are reservations necessary is the question; do they work? A study done by two IIT directors (and quoted by S.S. Gill, Secretary to the Mandal Commissions) shows that only 50 percent of vacancies at IITs for SC and ST are filled up of the remaining half 25 percent drop out before completing the course. Much of the same is true of Delhi. According to the 2001 – 2002 report of the Parliamentary Committee on the welfare of Scheduled Castes and Tribes between 1995 and 2000 barely half the SC quota and merely one – third of the ST quota of undergraduate seats was filled up.

Medical Colleges fare no better. Last year at Kanpur's Ganesh Shankar Medical College as many as 50 percent of SC candidates for the final year MBBS exam failed. A further 17 percent got only through grace marks.

For some the conclusion seems to be quotas don't work. They either are not filled up or many candidates don't finish the course. But you could also argue the opposite. If even 50 percent of the quota is filled that's proof of efficacy because without quotas the percentage of SC/ST students in higher education would be a lot less. Next what do reservations mean for those they don't cover? Seen through their eyes they mean students less qualified get preference over those more qualified solely because they are OBC. That's penalizing the upper castes for being upper castes. Thus reservation favours caste over merit. This is why students have risen in protest…

The experts gave the following comments on the Karan Thapar Sunday sentiments.

(1) According to experts the views given by the two IIT directors cannot be fully accepted for in the first place majority of the IITs do not wish to fill all the seats, so they don't send the very selection card and in other times give a cut off which will cover only 50% of the reserved quota. This is very intelligently done by the "merits" of India.

(2) Granted they join the institution like IITs, IIMs or even in medical colleges, which are mainly dominated by Brahmins and upper castes the SC/ST castes students are discriminated tortured and harassed by the faculties and ultimately failed by them. For the experts say in many of these Brahmin dominated institutions, the attendance register itself carries behind the name the term SC or ST; so it is very easy for them to be discriminated. In some classes the brahmin teacher harass them by asking them questions and when they do not give the correct answer they are humiliated before the class and so on claimed experts. One of the basic questions to our commarade Karan Thapar is; are the teachers in these institutions so useless that they cannot make a student pass? What is the great merit in these teachers who die more mediocre than teachers in other institutions. Experts also added the situation is very different in T.N.

(3) Experts asks when college teachers and university teachers could make SC/ST students fair so well in T.N. how is it when the “merit” (brahmin prof.) teaches and they fail? Does it imply meritorious will be taught only by merit?
(4) The experts said invariably the brahmin and upper caste teachers not only ill treat SC/ST students but see to it they don’t get first class. Unless these students are taught by their own community it is very difficult to help these students. Thus parallely govt. should bring reservations in the teaching posts in these institutions.
(5) If this is this case who will give such high pay for them.

Finally let them first give reservations for OBC for atleast 10 years and see the development of the nation said experts. Unless and until they are developed it is impossible for India to become a developed nation claims experts. Experts demanded the govt. to fill at least the existing teaching post vacancies by OBC/SC/ST candidates.

The Telegraph

The author is former director general, National Council for Applied Economic research.

Education must not become as politicized as it has.

… None of his politicking has any relationship to the genuine problems of getting a good education by the poor, the socially backward and the minorities. Educational standards must not be devalued to give them the pretence of an education. College, higher and professional education demands the underpinning of a good high school education. Government must perforce improve the quality of school education for all…. The object must be to improve the capability of the weakest in society not just give empty opportunity.

The experts view that none of the schools which produce cent percent passes and some ranks induce any form of intelligence but only by memory they get these marks. This has nothing to do with merit. What one needs is a regular coaching and an appropriate surrounding to study claims experts. So if OBC/SC/ST also are given these they can easily outshine the upper caste says experts. Without any comfort they score good marks so even granted we give them little more attention their performance would be unparalleled with the upper castes. So it is not the school which is to be improved we should give them opportunity after their schooling and while at school also. Even 60 years after independence if we go on talking about improving the schools it clearly makes us ask the question what have we done in these 60 years? So it is high time we give them opportunity in institutions of higher learning and see developments. Conjecturing is wrong in the part of upper castes and these elite groups who go on writing articles in the print media as if it is in their hands; about merit and other hollow things. These only project the selfishness and concern for them and absolutely no concern about the majority.

Hindustan Times

Jamia Milia Islamia is the first university in the country to increase the number of seats for students after the HRD ministry decided to introduce reservations for OBC. This July we will increase 699 seats across the board in the undergraduates and professional courses said Prof. Mushiral Hasan, Vice Chancellor.

The experts not only welcomed Jamia Milia Islamia's act of increasing the seats but also wanted the govt. to implement reservation for Muslims.

Now not in phases allies to govt.

Hindustan Times

At its meeting with the left, the UPA made a strong case for quota based on the GOM's recommendations.

Justification

- Quota will help to create a more equitable society enable greater upward and social mobility.
- OBC have not benefited from the expansion in education opportunities in the last six decades.
- Constitutional amendment that came into effect on January 20 this year allows state to make special provisions for socially and educationally backward classes SC and ST.

Way forward

Widening quota must be done in a manner that's just and efficient while not diluting quality. The catch

- If the number of general seats is to be preserved, the number of seats will go up by 54 percent.
- Capacity expansion will take at least 3-4 years.
- Getting more teachers is not easy possible way out reemploy teachers till the age of 65.
- The UPA allies have begun to mount pressure on the Manmohan Singh government to give immediate to the 27 percent reservation of seats in central institutions and private unaided centers.

The Hindu

The central government more to introduce reservations for other backward classes (OBC) in elite institutions of higher and professional education – popularly known as Mandal II seems to be heading towards a stalemate. In this article we propose a possible solution that might take us beyond the debilitating stand off between merit and social justice.

This is clearly an ambitious and optimistic agenda, especially because Mandal II proves that some mistakes are destined to be repeated. Once again the govt. appears set to do the right thing in the wrong way without the prior preparation. Careful study and opinion priming that such an important move obviously demands. It is even more shocking that students from our very best institutions are willing to re-enact the horribly inappropriate forms of protest from the original anti-mandal agitation of 1990-91. As symbolic acts, street – sweeping or shoe polishing or group inequalities in higher education. ST, SC, Muslims and OBC lag behind.

Sikhs, Christians, Upper caste Hindus and other far ahead.

Table 1

| Table shows percentage of graduates in population aged 20 years and above | | |
|---|---|---|
| Caste/Communities | Rural India | Urban India |
| ST | 1.1 | 10.9 |
| SC | 1.2 | 4.7 |
| Muslim | 1.3 | 6.1 |
| Hindu – OBC | 2.1 | 8.6 |
| Sikh | 2.8 | 25 |
| Christian | 4.7 | 23.7 |
| Hindu – Upper Caste | 5.3 | 25.3 |
| Other religion | 5.4 | 31.5 |
| All India Average | 2.6 | 15.5 |

Source – Computed from NSSO 55th Round Survey – 1999 – 2000.

Shoe shining send the callous and arrogant message that some people – castes are indeed fit only for menial jobs, while others are "naturally" suited to respectable professions such as engineering and medicine. However, the media do seem to have learnt something from their dishonourable role in Mandal I. By and large, both the print and electronic media have not been in cendiary in their coverage, and some have even presented alternative views. Nevertheless far too much remains unchanged across 16 years.

Perhaps the most crucial constant is absence of a favourable climate of opinion. Outside the robust contestations of politics, our public life continues to be disproportionately dominated by the upper castes. It is therefore unsurprising, but still a matter of concern that the dominant view denies the very validity of affirmative action. Indeed the antipathy towards reservations may have grown in recent years. The main quotas and the like as benefits being handed out to particular caste groups. This leads logically to the conclusion that power – hungry politicians and vote bank politics are the root causes of this problem. But to think thus is put the cart before the horse.

A rational and dispassionate analysis of this issue must begin with the one crucial fact that is undisputed by either side - the overwhelming dominance of upper castes in higher and especially professional education. Although undisputed this fact is not easy to establish, especially in the case of our elite institutions, which have always been adamant about refusing to reveal information on the caste composition of their students and faculty. But the more general information that is available through the National Sample Survey Organization clearly shows the caste – patterning the relevant data are shown in Tables 1 and 2.

Table 2

| Over – and under Represented groups ST, SC, Muslims and urban OBC severely under represented. Sikhs, Christians Hindu upper castes and others over represented. | | |
|---|---|---|
| Table shows group share of all graduates and a percentage group of share of 20+ population. Numbers below under indicate under representation, above 100 indicate over representation. | | |
| **Caste and Community groups** | **Rural India** | **Urban India** |
| ST<br>SC<br>Muslim | 43<br>47<br>52 | 71<br>30<br>39 |
| Hindu OBC | 82 | 56 |
| Sikh | 107 | 169 |
| Christian | 200 | 154 |
| Hindu UC | 205 | 164 |
| Other religions | 200 | 200 |
| Source computed from NSSO 55th Round Survey 1999 – 2000. | | |

Table I shows the percentage of graduates in the population aged 20 years or above in different castes and communities in rural and urban India. Only a little more than 1 percent of Scheduled Tribes, Scheduled Castes and Muslims are graduates in rural India, while the figure for Hindu upper castes is four to five times higher at over 5 percent. The real inequalities are in the urban India, where the SC inparticular but also Muslims, OBC and ST are way behind the forward communities and castes with a quarter or more of their population being graduates. Another way of looking at it is that ST, SC, Muslims and OBC are always below the national average while the other communities and especially Hindu upper castes are well above this average in both rural and urban India. Table 2 shows the share of different castes and communities in the national pool of graduates as compared to their share of the total population aged 20 years or more. In other words the table tells us which groups have a higher than proportionate (or lower than proportionate) share of graduates. Once again with the exception of rural Hindu OBC and urban ST the same groups are severally under represented while the Hindu upper castes, other religions (Jains, Parsis, Buddhists etc) and Christians are significantly over represented among graduates. Thus the Hindu upper castes share of graduates is twice their share in the population aged 20 or above in rural India and one – and – a half their share in the population aged 20 or above in urban India. Compare this for example to urban SC and Muslims whose share of graduates only 30 percent and 39 percent respectively of their share in the 20 and above population.

It should be emphasized that these data refer to all graduates from all kinds of institutions countrywide – if we were to look at the elite professional institutions the relative dominance of the upper castes and forward communities is likely to be much stronger, although such institutions refuse to publish the data that could prove or disprove such claims.

Although it is implicitly conceded by both sides, upper caste dominance is explained in opposite ways. The upper castes claim their preponderance is due to solely to their superior merit, and that there is nothing to be done about this situation since merit should indeed be the sole criterion in determining access to higher education. In fact they may go further to assert that any attempt to change the status quo can only result in "the murder of merit". Those who are for affirmative action argue that the traditional route of caste dominance namely an upper caste dominance namely an upper caste monopoly over higher education - still remains effective despite the apparent abolition of caste. From the perspective, the status quo is an unjust one requiring state intervention on behalf of disadvantaged sections are unable to force entry under the current rules of the game. More extreme views of this kind may go on to assert that merit is merely an upper caste conjuring trick designed to keep out the lower castes.

What is wrong with this picture? Nothing except that it is only part of much larger frame. For if we understand merit as sheer innate ability, it is difficult to explain why should it seem to be an upper caste monopoly? Whatever people may believe privately, it is now beyond doubt that arguments for the genetic or natural inferiority of social groups are unacceptable. If so how is that roughly speaking, one quarter of our population supplies three quarters of our elite professionals? The explanation has to lie in the social mechanisms through which innate ability is translated into certifiable skill and en cashable competence. This point to be intended or unintended systemic exclusions in the educational system, and to inequalities in the background resources and education presupposes.

It is their confidence in having monopolized the educational system and its prerequisites that sustains the upper caste demand to consider only merit and not caste. If educational opportunities were truly equalized the upper castes share in professional education would be roughly in proportion to their population share, that is between one forth and one third. This would not only be roughly one third of their present strength in higher education, it would also be much less than the 50 percent share they are assured of even after implementation of OBC reservations!

If the upper caste view needs an unexamined notion of merit that ignores the social mechanisms that bring it into existence, the lower caste or pro- reservation view appears to require that merit be implied of all its content. While this is indeed true of some militant positions, the peculiar circumstances of Indian higher education also allow alternative interpretations. In a situation marked by absurd levels of hyper selectivity 300000 aspirants competing for 4000 IIT seats for example merit gets reduced to rank in an examination. As educationists know only too well, the examination is a blunt instrument. It is good only for making broad distinctions in levels of ability. It cannot tell us whether a person scoring 85 percent would definitely make a better engineer or doctor that somebody scoring 80 percent or 75 percent or even 70 percent.

In short it is only a combination of social compulsion and pure myth that sustains the crazy world of cut off points and second decimal place differences that dominate the admission season. Such feticides notions of merit have nothing to do with any genuine differences in ability. The caste composition of higher education could well be changed without any sacrifice of merit simply by instituting a lottery among all candidates of broadly similar levels of ability – say the top 15 or 25 percent of a large applicant pool.

But the inequalities of our educational system are so deeply entrenched that case inequalities might persist despite some change. We would then be back where we started – with the apparent dichotomy between merit and social justice in higher education. How do we transcend this dilemma? Is there a way forward where both merit and social justice can be given their due?

Experts view that the term merit is not explained by the arya Brahmins and the upper castes. They have always claimed themselves to be meritorious! Let us look back in these 60 years the merit has not made any form of significant changes in the living conditions of the majority. So it is high time it is banished and those people are given at least some opportunities. They have not even come up to fight in full strength. Before they gather momentum and stage their demonstration we give them what is their right says experts. Inequality cannot rule for a long time for by the natural force of equilibrium it is sure to loose it stability. The experts feel from now on words let us burry deep the mythical word merit and proceed on to build a developed India with the progress of the majority!

The experts claimed they can give hundreds of instances in which OBC/SC/ST merits were buried for the only reason they did not belong to the upper castes. This sort of discrimination is very common in jobs and higher education. The meritorious among the OBC/SC/ST are very much tortured than an average and are driven mad by the upper castes. The meritorious OBC/SC/ST are very much targeted by the upper caste that is one of the reasons why the real merit talent and capabilities in them are never given any recognition or representation by the state or nation says experts. They said they can show how meritorious doctors, lawyers and teachers were humiliated harassed and tortured by the upper castes. If it is not an exaggeration OBC high capabilities and talents were their (Brahmin) enemies which hindered any form of worldly progress!

The only reason for the upper caste arrogance “merit claim” is the surrounding, economy and above all the comfort which they can command to achieve after all marks; these three things a rural OBC/SC/ST cannot even dream of as round the clock he struggles for the very living which does not even mean a square meal a day. So the question of nutrition or comfort is only a myth in his life. That is why the upper caste is able to speak of the hollow merit possessed by them says experts without defining the very notion of “merit”!

When genetically it is proved that OBC/SC/ST lack noting and when mathematically it is clearly brought out by the authors well tabulated dalas (thanks to them) why should any one have any form of hesitation to give 27 percent reservations for the OBCs asks experts?

Times of India

In a move that might fuel another round of anti quota protests the ministry of social justice and empowerment has announced that it is favourably considering in raising reservations for SC from existing 15 percent to 16.23 percent following the increase in the population.

Social justice and empowerment minister Meira Kumar admitted that her ministry had issued directions to NGOs aided by it to follow reservation rules for SC, ST, OBC and physically challenged said, she was satisfied with the response from NGOs on this count but had given more time to them to employ people from these sections.

The Telegraph

Ashok Guha. The author is the professor of economics at the School of international studies JNU.

In the recent polemics reservation has been projected primarily as an issue of caste conflict. It was well defined profiles of beneficiaries and of losers and the losses of one are the gains of the other. The govt. claims that its negative effects are confined and an optimal consensual solution achieved by increasing university admissions so as to maintain the number of open seats. All this sound and fury has drowned out debate on three key questions. What explains the persistent backwardness of the scheduled castes and scheduled tribes more than half a century of SC/ST quotas?

How will reservation affect economic growth? And what will be their long – term impact on the social and economic status of the groups they seek to protect?

A starting point of an exploration of these issues is a piece of data from a major university. For several years now the fractions of SC, ST and other backward classes applicants for admissions to this university has averaged around 90 percent, 6 percent and 17 percent respectively as against the existing and proposed quotas of 15 percent, 7.5 percent and 27 percent. A immediate implication is an ST applicant will have at least a 50-percent higher probability of success than one without quota benefits while SC and OBC applicants the chances of success are more than double those of non quota candidates, regardless of merits and entirely because of caste.

This may of course be precisely the objective of the quota regime …. It waives eligibility criteria for application for quota candidates and substantially reduces their passing requirements in its admission test. Once a quota candidate is admitted it gives him preferences in scholarships and hostel seats. Finally it awards him a degree even if his performance is well below that required of a non quota student…

In educational institutions, the sizeable presence of such a group compels lowering of standards of instruction; this retards the academic progress of non quota students as well. Since caste quotas are imposed in faculty recruitment too, the capacity of the institution to teach its students – including those not admitted through quotas is impaired. A vicious circle of progressive deterioration in quality is set in motion.

The dilution is quality of admissions, instruction and evaluation of credibility and value of degree…. With such a decline in quality, India's comparative advantage in human capital – intensive activities is certain to disappear. Employers are bound to turn increasingly to other countries where the educational system is not similarly constrained. Not only will they recruit fewer graduates from Indian Institution, they will also curtail their investments…

This is the likely outcome of the vast expansion of quotas in educational institutions alone… There is little doubt that reservations will extinguish the growth surge of the Indian economy. In the long run the groups is supposedly protected will not be able to escape the consequences of the resulting stagnation. The nation at large needs to ponder over this before it lets itself be steered into a bottomless pit by selfish and short sighted politicians.

The experts were against the IITs and IIMs students going abroad. Today over 10% of the nations wealth is in their hands and most of them get good education from these institutions for which govt. spend thousands of crores to it and finally as clearly said by Guha serve the other nations. If the quota can cut down such brain drain it is a welcome thing so that they work for the nation from which they have become educated. It is very good if they recruit fewer graduates from Indian institutions.

Experts said the last resort on people who speak or write in print media falsely project about low quality of teachers or students would be, we would claim damages for their derogatory remarks on our talents. Experts warned them of dire consequences. Further only these arya Brahmins must ponder over their baseless fear which is truly short sighted and highly selfish.

Also the degrees of IITs in post graduates and doctoral levels needs major scrutiny; for they (experts) say at least over 40% of these thesis or dissertations submitted by the students are not original. Most of them copied. That is why in IIT Madras the Ph.D thesis submitted by students, are kept under lock and key in the library. For instance if any staff/student wishes to refer a Ph.D. thesis submitted by to IIT Madras they have to get a willingness letter from the guide of the student under whom he / she has submitted the Ph.D. thesis. This sort of closedness is not in vogue in any other universities. Further this shows if the guide fears that thesis contains some copied maternal naturally he / she will not sign the willingness letter. This sort of practice is very unique to IIT Madras for it fears it would be clearly caught for plagiarism. The experts challenged if IIT Madras has no inhibition or fear let it first keep all thesis in open, second let it use the latest software which can detect plagiarism before and just at the time of synopsis defense said experts. The experts feel the main reason of the claim of Guha that students do not come to study is mainly due to the fact that the feed back they receive from the other OBC/SC/ST students the horrifying or problematic environment for them; so they shun to study, and secondly the faculties who handled them. Experts requested the govt. to fill the faculty positions in these institutions by OBC/SC/ST for their own survival success and give them a week crash course to these faculties to encourage the OBC/SC/ST student and not to join hands with upper caste faculties and torment the OBC/SC/ST. The experts said most of the OBC/SC/ST faculties in the institutes of higher learning think they are appointed on merit and they look down upon the OBC/SC/ST. This is also one of the reasons why the percentage is not properly available about by the OBC/SC/ST.

The New India Express

Look at the basic principles of social policy… The second reason forged caste identifies can't be dismissed as another usual Indian vice is that, these pieces of easily obtainable paper make the point sadly; few have done in the quota debate : the assumption of total congruence between deprivation and caste has been questioned authoritatively but to absolutely no effect as far as political action has been concerned. As this newspaper has always argued and as opinion polls have conclusively shown a greater majority of the better off wants targeted action for the socially disadvantaged – people are just not sure that caste is the right marker. As caste identity becomes increasingly purchasable those doubts will grow.

The expert view that this article is also to project only reservation is wrong. Now they have in this article just projected the myth all will get caste certificates which according to their polls opinion is easily purchasable. How ever in this article the usual slogan of merit will be destroyed is not used!

The New Indian Express

The medical students of Amritsar Government medical college protesting against the proposed reservation policy in premier central institutes on Saturday, they held a poster "Kick this Arjun Singh Monster out" was the wordings written on it.

Services in state run hospitals across the country were hit on Saturday as resident doctors struck work to protest the proposal for reservations for OBCs in elite educational institutions…

The poster says the experts only shows how the arya brahmins have the arrogance to call say “kick this monster out”. They all asked only one simple question. Can the non brahmins make a protest of his form against the brahmin Union Minister? Granted they ask such a question what would be the plight of the non brahmins who have asked to kick the monster. Why is media highlighting this!

The experts said this is the strong evidence to show to the public, how our legal system is caste dominated and only the laws of manu i.e., the Aryan law rules India.

The New Indian Express

The Supreme Court on Wednesday took strong exception to the governments decision not to release the salary of the doctors who took part in an anti quota stir saying it should act like a model employer.

The government says something at one point in time and … When government came on May 29 you were on your knees. After the strike is over you say ‘no work no pay’ a bench of Justice Anjit Pasayat and L.S. Panta observed.

The court said that at the time of passing the order the government had said no action whatsoever would be taken against the doctors and that it simply wanted them to resume work as patients were suffering.

The court said but the assurance given by the government it would have passed an order on this point also.

The court wondered if the concept of ‘no work no pay” was not known to the government till May 28 and observed that it was not expected from those government. “Ask your government to be a model employer” said the bench.

The court ordered the doctors absence from work shall not be taken note of for the purpose of their completion of training and internship and appearance in post-graduate examination.

The court fixed July 17 to hear the doctors’ plea after Solicitor G.E. Vahanvati sought time to seek instruction from the government on the question of releasing the salary for the strike period (May 14 to June 3).

Complaining that the government was taking punitive action against doctors who took part in anti-quota stir, resident doctors associations of AIIMS and Maulana Azad Medical College had moved the court on June 30.

They had also sought the leave of internship and residency during the strike period from May 14 to June 3, be condoned so that they were not rendered ineligible to appear in the coming PG entrance examinations, PTI.

This article says experts, clearly shows that the complete legal department i.e., the court and law are in the hands of the arya Brahmins more so to say in the hands of BJP/RSS. That is laws of manu are practiced by SC and HCs claims the experts. They (the doctors) who abstained from work during May 14 to June 3 by their very conscience have no right to get pay why should court give any law? When one does not work how can one expect pay? This is illegal the experts said for when they have made a strike that too against govt. if they have any shame they should forgo the pay even if govt. gives them. Experts said arya Brahmins are too much inclined only after money and they may loose their self respect but not money! Suppose non Brahmins (medicos) have made a stir against the Brahmin govt. and the brahmin Union Health Minister who is their institutes president can they ever expect the court give verdict in favour of the medicos. The medicos would have faced only dismissal or suspension said. So he called the non Brahmins to join together in unity and fight for social equality and justice. The experts demanded SC and HCs should not practice laws of manu, they wanted reservation in all legal posts as per their strength so that they can get justice not favour like the brahmin who get it from SC and HCs in all their cases. Experts said any case which is brahmin verses non brahmin very rarely a non brahmin gets justice.

Sunday Express

Salary during strike reflecting the differences between health minister Anbumani Ramadoss and AIIMS Director, Dr P. Venugopal on order issued by the eminent doctor was described as “null and void” by an official junior to him.

In escalation of the feud Medical Superintendent Dr D.K. Sharma on Saturday rejected an order issued by the Registrar V.P. Gupta which had said that the Director has ordered release of salary to the anti- reservationist resident doctors during their strike period last month.

It is clarified that the orders issued by the director AIIMS only pertains to the completion of the scheduled academic term of the residents.

"If any such order regarding sanction / release of salary to the resident doctors for the strike period has been issued, it is null and void since no administrative / financial approval was granted for the release of salary". Sharma said in a communication; Gupta had said that Venugopal had also ordered adjustment of the strike period of 17 days between May 15 and May 31 last when the resident doctors struck work.

The experts said the very paying of salary to the doctor during the 17 day strike period was a violation of natural justice, for this would set a very bad example to the employee for by citing this situation they will carry on a strike or suspension of work but claim for salary. How is this justifiable? Experts said if OBC director had made the strike certainly if the OBC doctors were involved in it against a brahmin in govt. the law would have been favouring the govt. The experts asked in surprise; what a law? What a brahmin domination?

The New Indian Express

Oversight committee Chairman M. Veerappa Moily who was asked to prepare a road map for implementation the 27 percent reservation for OBC in elite educational institutions on Saturday said with the right infrastructure the increased quota could be put in place without lowering standards.

"By including more people, I don't think the standard will be lowered provided we give equal support in infrastructure and faculty" he told NDTV news channel.

As the common entrance test was administered at the same standard for general and backward caste categories, the question of merit was misplaced, he said.

"I introduced the system of the Common Entrance Test (CET) for the first time in the country. Consequently every thing transformed said the former Karnataka Chief Minister, Moily said that the standard was never lowered in the CET. The difference between general category and backward class students may be 1.5 to 2 percent and between SC and backward class students it may be between 2 to 3 percent. He also said that there can be 50 IITs in the country.

Asked how the money would come for such an expansion Moily said funds should be available. He had discussed the issue with Prime Minister Manmohan Singh Planning Commission Deputy Chairman Montek Singh Ahluwalia and the Finance Secretary who has assured that the money would be available.

The New Indian Express

The functioning of AIIMS which was the nerve centre of anti reservation agitation and whose director P. Venugopal is at loggerheads with Union Health Minister will come under scrutiny of a high level expert committee.

The four member committee headed by M.S. Valiathan, which has been asked to give its report within three months, would examine whether the AIIMS has achieved the purpose and objectives for which it has been asked to make recommendations of efficient utilization of manpower resources to attract the best talent, retention of faculty, provision of better opportunities for utilizing the talent available and optimization of scientific, technical and non-technical man power.

The New Indian Express

The issue of doctor's strike, against OBC reservation and against the subsequent removal of AIIMS director Dr P. Venugopal goes to the Supreme Court. The apex court has listed a public interest litigation (PIL) filed by a social health activist organization, Peoples for Better Treatment (PBT). It will be heard on July 24.

The PIL seeks imposition of "Punishment" on the doctors who went on strike, pointing out that under the oath of Hippocrates, no medical practitioner could turn away an ailing man.

The petition filed by PBT president Dr Kunal Saha, an Indian born American and a leading AIDS specialist in the University of Ohio, US also sought punishment in the case of death of two patients who were denied emergency medical treatment of the AIIMS during the strike period last month.

The secretaries Medical Council of India (MCI), Union of Health Ministry and the Director of All India Institute of Medical Sciences (AIIMS) have been arrayed as respondents.

The petition sought the intervention of the apex court to bring an end to the needless loss of lives and sufferings of patients due to the abrupt strike of the resident doctors at different hospitals in Delhi and other cities in India including a shut down of many "emergency" departments resulting in denial of life saving medical therapy for critically ill patients.

The PIL seeks cancellation of the "license to practice" of doctors involved in strike.

This is the first time the apex court would be taking up a case!

Experts said that the honourable Supreme Court is so sympathetic towards the striking doctors and wants to pay for the 17 day strike period. The Supreme Court superior law does not bother about the patient death or the closure of the intensive care unit or the emergency unit of those hospitals.

A well renowned doctor of the Ohio University, USA on the other hand files a PIL to cancel the license of the striking doctors. Law itself is different in USA and in India law is very different especially when it concerns brahmins. The supreme duty of any doctor is to attend to patient, their negligence in their duty due to unwanted strike has resulted in the death of two patients. One has no idea how many deaths due to medical negligence in the strike period have gone unreported. When so much advertisement is given to the anti- quota protest the media has no heart even to publish the death of the patients and warn the striking doctors. But the media has the guts to warn the govt. and criticize the Union Health Minister for working in the OP dept. in AIIMS when the doctors were on strike.

The media's attitude is absolutely shameless. The experts condemned the legal stand in supporting the striking medicos and the Director of AIIMS for his careless behaviour since he was morally responsible for the death of the patients. They demanded that Dr Venugopal should resign because he was the cause of the death of innocent patients. They said that he should be banned from practicing medicine because he had violated the Hippocratic Oath. They also wanted compensation to be given to the families of persons who lost their lives due to the doctor's strike. The experts demanded all the medicos who went on strike should be likewise punished for they are the root cause and responsible for the death of patients and the immeasurable inconvenience they are causing to the innocent patients.

The New Indian Express

Anbumani Ramadoss brings to his portfolio as Union Health Minister a significant qualification. He is that rare incumbent of the office to actually have a degree in medicine. It would be normal to expect therefore that he would be sensitive to the surgeon's greatest concern, procedure. At AIIMS however he had made that attribute conspicuous by its wholesome absence. His very public meddling at the institute threatens much more than the tenure of its director, P. Venugopal has proceeded on leave to protest the ministers intervention to place in key faculty positions candidate other than those recommended by the institute's search committee. The committee incidentally was headed by the Union Health Secretary. This very unseemly episode highlights the problem at the country's premier institutions our medical research and the further deterioration that could be ahead.

It is an apt moment to contrast the objectives with which the All India Institute of Medical Sciences set up in 1956 and what it is currently mired in. It was to offer world class medical education, serve as a referral point for tertiary care and undertake cutting edge research. It has in the decades since retained its status as the premier medical college in India comparable with the best in the world. All the same it has failed to record the kind of growth that status demands. It is shocking that today just about 50 seats are offered each year at the MBBS level. To be a worthy determinant in upgrading the quality of doctors produced by a country with India's size and challenges, that is just too paltry. Those challenges in terms of the huge numbers of patients that come to it for routine care, have meant that AIIMS must make for the sorry state of government hospitals. In research despite the much publicized fight of faculty, much has been done. But the sub-continents unique health profile it has been argued, needs enhanced focus on indigenous work.

If AIIMS is to play the role in the country's health sector that was envisaged by parliament and that the institute for all the talk of resource crunch has shown it is capable of it must be freed from political interference. Ramadoss is right to be worried about reviving public opinion in favour of AIIMS. But by meddling and there by over turning administrative procedures he is becoming part of AIIMS's problem. At this advanced state of crisis, the onus is on him to rescue AIIMS from politicization. That can only be done by respecting the autonomy of the institute and its director.

This article says the experts in another instance or proof of the Brahmin arrogance. The union health minister Dr Anbumani must respect the autonomy of the institute and its director. Who is running AIIMS? Who should question? This clearly indicates the union health minister, who is a non-brahmin Sudra should respect a subordinate brahmin who still director of AIIMS Dr. P. Venugopal. Thus the experts feel that medicos agitation at the AIIMS and other medical institutions are well contemplated and performed, because of caste identity it has been getting the over whelming support of the media.

Express news

Another patient who was terminally ill died due to lack of treatment at AIIMS on Thursday evening as resident doctors at the premier medical institute are on strike in protest against director P. Venugopal's dismissal.

Ramjeet (45) a resident of Bihar has been under treatment at AIIMS for over a month and on Thursday the striking doctors did not attend on him Mahavir said.

Two patients have died allegedly due to lack of treatment at AIIMS since Wednesday.

The experts lamented that neither the politician nor the media is bothered about it. The value of life has become so cheap when compared to the medico agitation anti-reservation.

More than the agenda against reservation the medicos strike was very much deep rooted in the autonomy and power of the AIIMS director which became shaky when the President of AIIMS Dr Anbumani Ramadoss made some appointments and changes in the administration.

**Today it is AIIMS**

P.V. Indiresan – former director of IIT, Chennai. The governing council of the All India Institute of Medical Sciences has condemned the dismissal of its director, Dr. P. Venugopal. Such is the contempt that politicians have for the professional class that it has no fear of treating its most reputed members in anyway they like. Politicians have no reason to worry about the professional class, because the Indian professional class has generally been servile. In the council itself members of the professional class would have voted for the removal of Venugopal, mortgaged their independence merely to please the minister. It would not have occurred to them that the same humiliation imposed on Venugopal can be inflicted on them too.

There was a time when people accepted the principle of the Divine Right of Kings. They accepted that by the accident of birth alone. God had granted king the right to rule and rule absolutely.

Professional institutions need autonomy because they have to pursue long term goals whereas politicians operate with short term priorities. There are few politicians who can think beyond the next election, the next bye election even. On the other hand professional institutions take long years to mature. They need to be insulated from the temporary pressures of the political class and of the dictates to the ephemeral bureaucrats too who, like butterflies, flit from one department to another.

Autonomy does not sanction liberty without responsibility. Autonomous institutions too are bound by discipline, the discipline imposed by peers. In the case of hospitals, the Indian Medical Association plays that role. In the present case if the minister had his complaints, he should have asked the IMA to advice him, and accepted its advice. That is no reflection on the authority of the minister. Truly speaking the authority of a minister depends on how great his institutions are not how puny they are. Such autonomy is good for politicians too. Their prosperity depends on the quality of professional institutions, India would not be what it is today but for what the IITs and the IIMs have accomplished. No hospital which has been perpetually interfered with by politicians will have the capability to treat the same politicians when they suffer a heart attack or brain damage. Professional excellence does not thrive under a political shadow. There is not a single outstanding professional institution in the world that has operated successfully under political control.

Having said this, this situation would not have arisen if top professional of the country had stood their ground. The country does have super scientists who commanded exceptional influence. It is unlikely that they will take a public stand. They will bemoan the state of affairs in private. But in the presence of authority they will behave like Dronacharya. In India today we are having the situation the guardians of the constitution are systematically destroying it. The political culture today is getting close to the decaying days of the Moghuls when local satraps became a law into themselves. An old saying goes: Vinasha kale Viparita buddhish (when doom looms intelligence runs contrary).

The experts first claimed, P. Indirasan as the Brahmin fanatic who was always very much against the non Brahmins. They portray him as a Brahmin supporter and medias' favourite. He always writes rubbish they said. He is above parliament and below the laws of manu for he has been says some experts. He is the root cause of bramanizing the IIT Madras.

Does he (P.V Indaresan) ever know the basic; that AIMS is only under the govt control and not under the brahmin control asks experts, respect any of the association it is strange says experts he wants the negotiation throw' the IMA might also feel a brahmin should not be dealt badly or punished by sudra. Experts say the conflict is as per the laws of manu for Dr Ramadoss is a Sudra (he may be the union minister and president of AIIMS who is bothered about it) so he intervening in the administration of Brahmin (Dr Venugopal the director of AIIMS) cannot be accepted? It is the crux of the laws of manu. That is why all the brahmins, brahmin media and above all the BJP and its sister parties strongly objected it. Finally the supreme law SC also supported only the brahmin say the experts. Thus it is high time we non Brahmins learn a lesson said experts. Thus as per the laws of manu a Sudra however highly placed he may be he should not penalize a Brahmin even if it is legally justifiable!

The Indian Express

The centre today suffered an embarrassment when Delhi high court stayed the controversial move by Health Minister Anbumani Ramadoss to sack AIIMS director P. Venugopal.

Hours after justice Anil Kumar granted the stay, AIIMS resident doctors and faculty members who went to strike on Wednesday to protest the move to dismiss Venugopal suspended the agitation late to night and resumed work immediately.

Responding to the High court order Ramadoss said “all legal course of action” including an appeal in the Supreme court will be taken against the stay on the removal of Venugopal while a buoyed AIIMS director announced he would attend office from tomorrow and vowed to fight on the autonomy of the institute.

Granting a stay on Venugopal’s petition challenging the recommendation by AIIMS apex body headed by Anbumani to dismiss him, Justice Anil Kumar directed the centre and the AIIMS to file their replies within two weeks on the matter.

He also issued notice to the election commission on Venugopal petition seeking Anbumani’s disqualification from Parliament on the ground that the minister by being the president of the institute was guilty of holding on “office of profit”.

Justice Kumar ordered the stay in a packed court room after three hours of arguments by additional solicitor general Gopal Subramaniam, who appeared for the centre and senior counsel Arun Jailtly who represented Venugopal.

The ASG disclosed to the court at the fag end of his arguments that the government had earlier asked Venugopal to gracefully resign on his own and only when he refused to do so was a decision taken to terminate his services.

Apparently peeved over the manner in which the institutes affair are being run, justice Kumar during the hearing said, “it appears either the President (minister) or the director is running the whole show”.

The court made the observation after the ASG defended the govt. move to appoint Dr Deka as Dean of the institute without the approval of the search committee despite the incumbent being relatively junior to several other contenders.

What ever decisions the search committee takes, can it be ignored? There must be a reason to ignore the same", Justice Kumar queried, when Subramaniam admitted that Deka's appointment did not have the approval of either the director or the search committee which had prepared a list of probables in which Deka figured only as the 17th preference.

Health minister Anbumani Ramadoss today clarified that it was the collective decision of the institute body of AIIMS to improve administration of the premier institute which has become a political hub instead of medical hub.

Experts felt that the situation faced by the Health Minister, was unfortunate. They all showed their mental agony over the subordinate Dr Venugopal display of arrogance which is due to his caste.

Several experts described the situation as a clash between the sudhra employer and a Brahmin subordinate employee. They all said that laws of manu played middle man role was always favouring the brahmin. They have only put forth the following questions for the reader to ponder over for they call this as the time of crisis.

1. Can ever in the history of India a Brahmin union health minister be used in the court by a Sudra director of AIIMS?
2. Will the court ever accept the W.P or WIP filed by a subordinate OBC against the Brahmin employer?
3. Will court ever reverse the order in favour of the OBC subordinate verses the brahmin superior?

Experts say because the court is just laws of manu they have the guts to stay the termination made by the Sudra union minister against his subordinate director Dr Venugopal. Is the credential of Venugopal important or his code of conduct important as an official!

Can the court say if some one has high credentials he can behave in any questionable way?

More so can a sudhra with high profile behave like Dr Venugopal?

The experts have put all these material before the reader so that they can contemplate what is law? Is it relative or absolute or purely depending on caste?

The New Indian Express

The BJP on Friday demanded the resignation of Union Health Minister Anbumani Ramadoss and urged the Prime minister to reinstate Dr Venugopal as the director of AIIMS with honour and dignity condemning Ramadoss attitude towards Venugopal. BJP spokes man Ravi Shankar Prasad said that it was unfair to humiliate doctors and professionals in the country. He said that Ramadoss must resign owning up responsibility for the complete mess in the Health ministry.

In a letter to the Prime Minister BJP senior leaders including former Prime Minister Atal Behari Vajpayee and leader of opposition in the Lok Sabha L.K. Advani, expressed serious concern over the arbitrary manner in which Venugopal was dismissed. They said the Prime Minister must clarify his stand immediately since his name and the authority of the office of the Prime Minister were involved in the matter. Pointing out that the dismissal of one of its outstanding cardiac surgeons has shocked the country, the BJP leaders said that the decision not only humiliated Venugopal but also the professional community in the country. They said it was an irony that during the golden jubilee year of AIIMS, the autonomy of the institute was being totally shredded "what could have been a year of celebration, renewed dedication and true autonomy has now been turned opposite" they said. The BJP leaders recalled that though there have been occasions earlier when directors or senior professional differed and disagreed it did not result in any punitive or vindictive action. They also reminded that once during the NDA government tenure a very distinguished scientist was suspended from service and as soon as the matter was brought to the notice of the then Prime Minister he was reinstated with honour.

Express news

Terms the appointment of Union Health Minister Anbumani Ramadoss as the president of the AIIMS "illegal" Panlit Foundation for Rural Transformation (PIFORT) in a petition to Prime Minister, Manmohan Singh has charged the government with eroding the autonomy of the premier institute.

The petition demanded the removal of the minister as president of the institute and called for the appointment of a competent scientist or visionary to the post.

The petition further argued that combining the two positions of minister and head of the institute leads to day to, day interference in the affairs of the institute their by completely subverting the norms of good governance and eminence.

Is the country shocked by the dismissal of Venugopal? Experts said we were also very happy and infact welcome it for he (Venugopal) cannot compensate the lost lives of the patients. What is his moral responsibility in the death of the patients and the suffering which the institute has caused to the innocent patients due to the17 long day strike which was utterly unnecessary and illegal! So experts warned the BJP not to use terms like the dismissal of Venugopal has shocked the country. Only the BJP and the related brahmins must have been shocked by it. This clearly proved Venugopal is politicalizing the issue as basically he belongs to that party.

Further experts said (some of them doctors) that his dismissal was in no way affecting all the doctors of the country for what the govt. did was only legal and proper!

Further experts made it very clear that even the anti quota or reservation protest had its ulterior motives in the feud between the Sudra master and the brahmin servant they said.

The New Indian Express

Any social process if it has to achieve its intended consequences requires the benefit of deliberation and due diligence. The impatience with detail displayed by the proponents of reservation politics to push through 27 percent OBC quotas in one quick stroke in the forth coming monsoon session of parliament is disquieting. It should be disquieting even for those who are otherwise open to the idea of a social welfare intervention of this nature, because the impatience on display points to political rather than social agendas at work. Since political rhetoric does not feel the need to submit itself to the impediments of measured discourse or the lines drawn by the constitution there is now even talk of reservation for the judiciary. Such competitive mouthing of the reservation mantra could end up giving affirmative action a bad name.

In the fury of words the original intent of OBC quotas (to ensure access to higher education among those denied such access due to social asymmetries) is almost forgotten. The expansion of the existing capacity in our institutions of higher learning – a move print medias has firmly supported cannot be achieved by the wave of arm. It is to actually benefit students its demands painstaking even laborious effort and planning. The observation made by the authorities of the Indian Institute of Management Ahmadabad to the Veerappa Moily committee is worthy to close attention. It look the institution 20 years to get the land, two years to plan and five years to implement its most recent expansion plan, they pointed out. Similar hurdles are sure to surface in other spheres as well.

We would therefore counsel the promoters of OBC quotas to be more considerate in their approach. S. Ramadoss is not known for his nuanced view on the issue and CPI's D Raja can say what he wishes about reservations as he has been won't to do. But the other raja- V.P. Singh should atleast as the discernment of a former Prime Minister; direct the discourse in a more rational direction. The monsoon session of parliament promises to be full of sound and fury over the issue, but let it at least signify something.

This article only shows their (Brahmins) mental agony says experts. Unlike us they have the media to vent their feelings claims the experts. Their very mischief of calling V.P. Singh the other raja is sufficient to depict their internal hatred and aversions towards reservations. Experts feel that Brahmins can never feel for social equality or social justice when it comes to non brahmins as they are very selfish and their minds is still only ruled by the laws of manu.

Express news service

Swami Dayananda Saraswathi is a renowned teacher of Vedanta and has founded a number of organizations throughout the world for spreading the Vedic message.

The spiritual teacher gave a series of lectures on the topic “Vedanta and psychotherapy” in the city last week.

In an interview to express he shared his views on reservation issue that hogged headlines recently and other social evils.

Excerpts from the interview:

Q What is your comment on reservation issue which is being talked about?

A There should be a complete study on the problems of the students where are the students in which areas they are coming up and how many students are going to really benefit from the move etc. As the reservation is a pre-election promise and a populist proposal they have not really gone into the matter in detail and studied it properly. There should be an indepth study on the issue and based on that a very detailed report could be prepared to enlighten the people on the quota implementation.

In our country the student population is struggling to get seats. Even those who have secured high percentage of marks are denied seats for professional courses. We have to see that the student population is happy. At the same time we have to do something to help the poor and economically backward students. They should come up in life. But this should be done in an enlightened way not just for getting votes.

Q: On the existing system of reservation …. ?

A: The politicians say that they want to abolish caste. But in reality they are perpetuating it at one time and expressing the desire to abolish it at another time. By adopting this double stand they are taking the country backward and not forward.

Q: Despite the growth of spiritual movements and growing religious fervour among the people social evils like untouchability and two tumbler systems continue to rock the rural masses…

A: We just allow the people to drop all these notions and come to a state that treats people equally. If the caste system is not perpetuated it will go automatically. As long as you have caste based policy and politics there will be no answer to this. We have to give them up in order to bring about a change in the minds of people. And gradually it is happening and it will go.

Q: But religious leaders have a specific and special role in achieving this task…

A: They are doing it to the extent possible they can. The only culture where oneness of all is existing is talked about is our culture. We do not look upon anything as secular. Money knowledge success or any other accomplishment; everything is associated with spirituality here. We do not have a culture that divides. Ours is a culture that unites. We should make people realize the essential aspect of our culture i.e. spiritual content and get rid of all differences.

Q: What you thing as panacea for casteism?

A: Caste based politics and caste based vote bank are the real problems. Actually we do not need the caste system. It will go automatically. Over a period people have changed and should not change certain things over night. Casteism had served its purpose and it can go now. I think it is going. As long as these people do not perpetuate casteism it will go. It will not last. Any system is not absolute. A system serves its purpose for which it came into being for a length of time and when the people go out of it goes away; it withers away.

Q: There is a notion that Hinduism propagates fatalism what is your comment on that?

A: It is wrong propaganda.

Hinduism makes you responsible for what you are. Hindus are the ones who have made even gods responsible for what they do. I take responsibility for what I am doing. Many of my past actions have led to my present state and for that I take responsibility and that is not fatalism.

In fact in the name of fate Hinduism has a shock absorber which others do not have. When my efforts fail I do not blame any force for the failure and go ahead with the confident that I can write my future. I accept what I have and go ahead. Indeed it is a positive attitude.

The experts highly condemned the notion of karma in Hinduism which was the root cause of all inequalities faced in India. Today the caste in India is based on birth and birth is based on karma. It is not as said by swami caste based politics and caste based vote bank are the real problems. If the so called fictitious vicious circle is true where is caste for a brahmin may be born as a sudra, a sudra a panchama and a panchama a brahmin and so on. Thus the fatalism in India is that people when they are sick even do not go to doctor, instead pray to gods, visit temples or fortune tellers. This is very common mainly among the lower social strata of people. They blame their poverty, lack of education to fate. Until the propagation of fatalism is stopped experts feel it is not possible to improve the society. Experts contented a brahmin is very cunning and does not believe in fate. Who is perpetuating caste systems asks experts. It is mainly the brahmin mutts that practice caste and infact grow the caste. Even today in Kanchi mutt Brahmins and non brahmins cannot dine together! Have any of the religious leaders ever questioned this asked experts including Swami Dayananda! Experts further said all brahmin religious leaders keep up only the brahmin caste. They may talk universality, equality and what naught but only practice casteism. Can a black complexioned baby from a slum ever get the blessings of achariyas?

She / he with a black complexioned baby in their arms can never even get a prasadam! What is casteless India? Experts asserted that the brahmin owned or headed religious mutts are the breeding place of casteism! In Kanchi mutt only brahmins are allowed to stay. There is not even in a single non brahmin who is housed in the premises of the Kanchi mutt said experts. When the religious mutts are so mean to propagate caste why should others falsely say politicians practice caste, which is an utter falsehood with the exception of BJP/RSS/VHP.

These parties BJP/RSS/VHP can be mainly termed as brahmin parties, for the representation of non brahmins in these parties in very meager. Even if they are present it is impossible for non brahmins to be a decision maker! It is still surprising to see the Kanchi mutt chief has not opened his mouth over reservation for he is against reservation and some of the protests are carried out only with his backing!

Express news

The centre should immediately implement 27 percent reservation for OBC in higher education institutions through an ordinance and all political parties should support the move, Paatlali Makkal Katchi (PMK) founded president S.Ramadoss said on Monday.

Ramadoss on a mission to unite the forces of social justice in favour of the move told reporters that he had already met Prime Minister Manmohan Singh congress president, Sonia Gandhi and Union Human Resources Development minister Arjun Singh and had impressed upon them the need for giving effect to 27 percent reservation to OBC students in higher educational institutions.

Ramadoss said he also had meetings with ID (U) president Shared Yadav, RJD Chief, Lalu Prasad, LJP Leader Ram Vilas Paswan and spoke to Bihar and Uttar Pradesh Chief Minister Nitish Kumar and Mulayam Singh Yadav, respectively, besides former Prime Minister V.P. Singh and sought their cooperation. I am happy to report that the UPA – left parties supported my views the government will be implementing the policy of 27 percent reservation from the next academic session and an enabling legislation would be brought in the ensuring monsoon session of parliament, he said.

"With a view to mobilize support for reservation I visited Mumbai, Hyderabad, Bangalore and several other parts", he said.

Ramadoss said according to article 340 the centre constituted the first Backward Classes commission in 1953 which submitted its report in 1953 recommending 70 percent of seats in all technical and professional institutions to the OBC and job reservation. The second backward classes commission known popularly as Mandal commission, which submitted its report in 1980, recommended 27 percent reservation to the OBC in higher educational institutions.

Ramadoss said if Kaka kalekar commissions recommendation had been accepted, OBC would have received the benefit of reservation in the last 51 years and if Mandal Commissions recommendations had been accepted the benefit could have come in the last 26 years, PTI.

Express news service

PMK founder Dr. S. Ramadoss on Monday urged the DMK government to enact a law enabling reservation in private sector in the forth coming budget session.

"The proposed law should be a model for other states" he told journalists here.

Ramadoss further said reservation in private sectors should encompass all those companies that had obtained bank loans or other concessions like power tariffs and subsidies from governments. Companies operating on their own resources could be exempted he suggested.

The New Indian Express

Union Human Resource minister Arjun Singh on Friday said the issue of the rights of Dalit Muslims and Dalit Christians needs serious consideration.

Religion should not be the ground either for giving or depriving anyone of reservation, said Singh when asked about the demand of Dalit Hindu converts to Islam and Christianity to be given this benefit.

I personally feel that no section of society should have a feeling that it was being deprived of anything on the ground of religion” he said.

Pointing out that several Muslim castes had been included as beneficiaries under the Mandal formula, Singh said if they were not being given this benefit they should fight for it.

Replying to a question he expressed reservation over the formation of a front, comprising Muslim organizations in Uttar Pradesh.

“Creation of such a front was a matter of concern for secular parties’ he said.

He however felt that it should be examined what led to the new political development in the state.

Indian Express

Medical students on Sunday turned down an appeal from Prime Minister Manmohan Singh to end their stir even as the Union Health Ministry gave resident doctors at the All India Institute of Medical Sciences (AIIMS) 24 hours to vacate hostel rooms if they are continuing with the strike.

After a meeting with the Union Health Secretary P.K. Hota, in the evening representatives of students groups said that the strike would continue as the appeal did not take into account their demands. The meeting was attended along with medical superintendents of the hospitals and AIIMS authorities.

Dr Kumar Harsh spokesperson for students group youth for enlightenment said that the striking medicos had three demands – roll back of the proposed OBC quota, review of the existing reservation policy by a non-political committee and a concrete statement by the Prime Minister on the issue, “The appeal does not address these criteria” he said….

Senior residents at AIIMS were put on a 24 hour notice to vacate their hostel rooms. But, residents said after their meeting with the Health Secretary that they will not back down ….

This is a fool proof of Brahmin arrogance and all this shows, Dr. Venugopal the director of AIIMS is strongly backing the medicos/ students that is why they say even no to our PM. Thus we see the irresponsible behaviour of the medicos. They are least bothered about the patient death. Neither media nor the govt. is ever projecting the death of the patients. The experts ask the following questions.

(1) Who will pay any compensation to these patients death as the death of the patient have caused only due to the denial of medical attention as the doctors / medicos were on strike foregoing their very duty!

(2) Why the media did not display the death of patients as it has been displaying the strike so many time in all media both print and TV? Why was this news not placed before BBC? Did they display the death of patients due to denial of medical attention because of protests?

(3) It is a pity that when the politicians, media and the medicos have left the patients in lurch the public atleast should have condemned the irresponsible behaviour of the striking medicos and doctors. Who will help these relatives or kith or kin of these patients who are in the prime of their age? Further how many patients really died due to the negligence of these striking doctors? Did any Brahmin patient die due to this strike asked experts?

Experts said atleast the govt. should take the moral responsibility and suspend the license of these doctors atleast for a period of 6 years. Then only such type of hazarderous mistakes will not repeat. Further the govt. must pay some compensation to the nearest kith and kin of these patients who have died due to negligence and this strike.

The New Indian Express

How can you support the extension of reservations to the OBC while moping up the resentment among the upper castes? The two national political parties are being incredibly disingenuous in playing the quota issue for maximum electoral gain… When a small group of health care workers striking duty in blatant abdication of their primary responsibility to patients can so easily appear to hold the government in breach of due diligence, someone has to be held accountable… Add the economically weak among the upper castes. It is so evidently an appeal to various sections of the social coalition it has successfully wooed in the nineties in Utter Pradesh which goes to polls next year…

Excellence and equity are twin objectives India compromises onto its own peril in this phase of globalization. But can excellence be maintained by expansion of castes (a key need) without adequate investment in educational infrastructure? And is equity really best achieved through reservation on the basis of an old caste census….

The editorial shows their hatred towards reservations at large says experts who reviewed this article. What ever be the problem it is unethical on the part of the few upper caste medicos to strike which has largely affected the innocent suffering patients who are poor and who need medical attention at large immediately. This is an unforgivable offence by the striking medicos claimed experts, they continue their protests as they get the support of Dr. Venugopal, the AIIMS director.

Secondly it is an eye wash to say economically weak among the upper castes, for the upper castes get always the monetary support and all political support from their finely knitted network. The upper caste especially Brahmins are never economically weak. The percentage of economically weak in reality may be less than 0.01% in their population. For the major portion of the nations economy is only in the hands of the brahmins claims experts. Further experts warned the arya brahmins and the pro media for their loose talks about excellence would be lost. It is high time they check their tongue otherwise we would forced to check them said a few experts. Why the media is not coming up to question the brahmins? The discrimination mainly based on castes which has been in existent for over so many thousand of years. Why increase the seats ask some experts? With the same number of seats the 27 percent reservation must be given, so that the govt. need not spent any money on improving the infrastructure. They all said only reservations alone can establish equality in this society as it is Varna that, continue to rule us!

The New Indian Express

In an acrimonious controversy about the 27 percent reservation for OBC in centrally owned institutions a repeated justification for the move pleaded by HRD Minister Arjun Singh and later by others including the PM is that parliament as recently as Dec. 2005 gave its mandate for this move by passing the 93rd constitution Amendment Act.

In the 93rd Amendment was really designed and enacted by Parliament in Dec. 2005 to enable the govt. to reserve a quota for OBC in centrally administered institutions there would be much to be said in justification for the govt. move on OBC quotas. But this is not the case. What ever be the merits of the move it is vital for the lay people to know the 93rd constitution amendment Act passed in Dec. 2005, cannot be its justification. Only a legal understanding of the problem of reservations can reveal that the government plea is disingenuous. Why was the 93rd Amendment to the constitution made in 2005. It was made not because the central government did not have the power to make reservations for OBC in educational institutions till 2005. If the centre intended to do that since 1951, it had the power as a specific amendment was made in Article 15 (4), enabling the central and state governments to make, "any special provision for the advancement of any socially and educationally backward classes of citizens or for the Schedule Caste and Tribes "The socially and educationally backward classes including OBC". With this specific power the central govt. made reservations only for SC/ST and not for OBC. This was despite the fact that several states like Karnataka and Tamilnadu had made reservations for not only SC/ST but also for OBC in institutions of higher learning and not only in state owned and state aided institutions but also in those that were privately owned. Even after the Mandal judgment of the supreme court in 1992, which gave go ahead for 27 percent reservations of OBC in public service, for over 14 years in the central govt. did not think of OBC reservations until Arjun Singh's sudden announcement in April.

Why did successive central government not make such reservations all these years and why did the present government think of them only in April? One can only speculate on the factors that caused the govt. to make this move. But it is significant to note this reservation with the help of Article 15 (4) do not even require a law to be made by parliament. They could have been instituted by mere executive order.

Now coming to the 93rd constitution Amendment passed by Parliament in 2005. It was enacted by Parliament for a totally different purpose to overcome the Aug.12, 2005 decision of the Supreme Court in the Inamdar case. In that case the apex court had struck down the existing reservations made by state govt. in private unaided medical and engineering institutions, because it held that such institutions had a fundamental right to occupation guaranteed by Article 19 (1) (g) of the constitution and the states did not have the power to impose such reservations on them.

Prior to the Inamdar case several states had for over 15 to 20 years made reservations in private, unaided professional institutions for SC/ST and OBC. The Inamdar case for the first time declared these reservations illegal. Parliament therefore rightly reacted and made the 93rd Constitution Amendment by a new Article 15 (5) in 2005, giving specific power to make reservations for the advancement of any socially and educationally backward classes of citizens for SC /ST and declared that nothing in the fundamental right of a private educational institution to prevent making such reservations. Thus the 93rd Amendment had a very limited objective – to restore and status quo ante to the Inamdar's case and enable the states to make reservations even in private unaided educational institutions.

Given the background the 93rd Constitution Amendment passed by parliament cannot be a justification for the central government current move to reserve 27 percent seats for OBC in centrally owned institutions. It is only fair that the govt. should not make use of the so called parliamentary mandate of 2005 for it and candidly state the real justification for this step when for so many years it had not felt the need to take it!

Experts strongly felt the court of law be it the SC or HCs never give any judgment against brahmins and their rule. Experts feel they are forced to see the reservations as thought by the brahmins are against them, the reasons are they are highly selfish and never want any form of equality enjoyed by the downtrodden for very clearly the laws of manu states that the sudras that is the majority of the OBC are born only to serve the brahmins meekly without having any form of luxury, they should only live in poverty for if they possess money or riches or lands according to the laws of manu it makes the brahmin feel sad and there is some law which gives the right to the king or brahmin to plunder or forcefully take the wealth from him [ ]. So by chance if in India they impose reservations for OBC they will not only get educated consequence of the good education they will be economically better off which cannot be tolerated by the brahmins.

So both the legal as well as the media are one with the brahmins so they will protest it in all possible forms and ways claims expert. Experts want to make it very clear that the constitution and that specially the parliament is even above the courts for the simple logic parliaments consists of the peoples representations where as the courts the brahmins ideology. They recalled with woe how each time the Brahmin force stalled the reservations for OBC. The brahmins they (experts) feel can never think of equality especially when it comes to sudras that too in education. They blow up everything beyond proportion as they have media and law in their hands said experts. They clearly known the law will never help the sudras it is the laws of manu which has very badly discriminated the sudras, all the more the laws of manu clearly says a sudra should not learn even if he tries to secretly learn he should be punished by either cutting of his tongue or by pouring molten lead into his mouth and ears; [43] if that is the case how will these brahmins accept the reservation for OBC.

The New Indian Express

Amid opposition from private sector to reserving jobs for SC and ST the govt. today said it will not impose quotas on India inc. but only wanted more “affirmative action” from the industry for the socially underprivileged. “The government is not talking of any enforced action for the industry that would be run down your throats". Minister of state of commerce and industry Ashwini Kumar said at a CII conference for the especially abled people here. "You (industry) are already into affirmative action through your corporate responsibility initiatives. It does not need any enforced action, you are already doing it", he said referring to the on going differences of opinion between the government and the industry on reservation and affirmative action which he termed "spurious" PTI.

Experts are very much bothered about the govt. soft dealing with these brahmins, while they sit on the head of them. Why not make a strong statement they have been living for centuries discriminated and dominated by you brahmins it is high time they get their due!

Experts very clearly wanted to inform the protestors that their vote was never the deciding authority of win or loose. This was not the case with the OBC if the OBC did decide that a party should not rule it will be clearly shown by their votes. So the experts warned the politicians to do justice to the OBC atleast after 60 years of independence. Experts were very upset that the social equality cannot be attained even after 60 years. The great difference between the haves and the have nots in India is based on caste which is very unfortunate they said. The sons of the soil in India, majority of whom are OBC/SC/ST happens to be the have nots of the nation. Only caste in India makes the difference between the haves and have nots.

The New Indian Express

Dalit panthers of India (DPI) on Monday staged a demonstration in the city urging the centre to implement 27 percent reservation for OBC in higher education without any further delay.

The agitators including large number of women led by party general Thol Thiurmavalavan demanded the centre not to succumb to the pressures of the anti reservation groups. But to proceed further with its decision to implement the 27 percent reservation in higher educational institutions like IITs and IIMs.

Experts claimed if the OBC are joining hands with SC/ST certainly the nation will not be in a position to even tackle the situation. In Tamil Nadu the largest party for the dalits have already; joined hands with the OBC and have shown their demonstrations. So it is a clear indication that the govt. has to face a bigger problem if it succumbs to pressure by the paltry 3% population!

The New Indian Express

The centre on Monday contented before the SC that is interim order for 10 percent reservation for scheduled castes (SC) and Scheduled Tribes (ST) in medical colleges would open the flood gates of litigation as it was not envisaged at the time of holding the entrance examination.

Let us not carve out quota after the examination and result. The provision for 10 percent reservation is not there in the brochure and hand book.

The New Indian Express

Two senior minister in the union cabinet on Monday came out in support of reservation of seats in elite educational institutions, while the congress adopted the Rahul Gandhi line of middle path on the issue … Chidambaram who is a member of the group of ministers on reservation of jobs in private sector, also rejected the argument that reservations would dilute the quality of either student intake or the quality of students passing out. Union Agriculture minister Sharad Pawar speaking at a function in Pune, also favoured reservation, saying it was just a misconception that merit was the domain of a privileged few and asserted that people representatives must make laws that provide reservation for the marginalized sections.

“It is good if standard educational institutions made attempts to improve quality, however to keep away students belonging to neglected strata in the name of merit was unacceptable” Pawar said.

Deccan Chronicle

The confederation of Indian industry on Thursday strongly objected to mandatory reservation of jobs in private sector for the socially underprivileged, while BRP leader and grandson of late BR Ambedkar Mr. Prakash Ambedkar asked to govt. to make public the non- performing asserts of industrial champions opposing the quota. Mr. Ambedkar, a former Lok Sabha member also demanded that the government to bring a legislation making reservation mandatory in the parliament session next month. But the CII view of such a law is brought in, the MNCs might say why they should set up shop here and might move out. If the legislation is only for domestic industry even then it will be detrimental. The govt. Mr. Ambedkar said should spend 20 percent of the total budget on the economic development of the SC and ST every year, CII president R – Seshasayee, said “Mandatory reservation in any form is not conducive to competitiveness in the industry. It is not acceptable. The industry need to take positive action to empower the backward classes to join the main stream and help them develop abilities to compete and empower them with education and employment skills.

“To ensure that SC, ST and OBC are able to complete and acquire job skills and as our response to the governments call for reservation in the private sector we are setting up a task force on affirmative action under the chairmanship of J.J. Irani, Mr. Seshasayee said, who is also the managing director of Ashok Leyland in his inaugural press conference.

Asked if the government brings in a legislation making reservation mandatory, he said, “I do think the government will be responsive to a dialogue. We are committed to a dialogue. We are committed to positive action. I do not think legislation is a solution. If there is a law you have to abide by it.

“The industry should take active steps to develop the technological skills of SC and ST and support their education, he said adding that the private sector should address issues related to employment and trade apprenticeship to ensure employability to the backward classes.

There should be no discrimination while competing for jobs” he said.

The experts feel that Mr. Prakash Ambedkar’s views on reservation is very welcoming for he has proved beyond means that he is closely following the foot steps of his grand father Dr Ambedkar. When such great people are for reservation experts asked why we should ever react to the 3% arya Brahmins and a fraction of the striking medicos against reservation! They once again upheld the views and suggestions of Mr. Prakash Ambedkar. Experts asked industries in which reservation is not given to the OBC then if the OBC SC/ST population does buy the goods supplied by them what would be the industries plight? Secondly we (OBC) SC/ST will start such industries so that only the industries run by the Brahmins will have to face a slow death said experts.

The New Indian express

The government today softened its stand towards student s indulging in anti quota agitation by saying it will protect the interests of every one even as it is ruled out any dilution on OBC quota in higher educational institutions.

Replying to a calling attention motion in the Lok Sabha, Arjun Singh today indicated that the govt. was in no hurry to bring the bill. The government will have collective view taking into consideration all aspects. No one will have any objection if efforts are made to address the cause of tension. That alone will be the wise thing he remarked.

There have been suggestions to increase the number of seats in the general category.

Softening his stand towards the striking students opposing the legislation he said 'the students are our children, we should also listen to them. Arjun brushed aside reports that government was speaking in different voices on the reservation issue.

"This is false. Forget that there is any division in the cabinet or there is any possibility of it, he said, even as Kapil Sibal entered the Lok Sabha and occupied the chair behind him. "Kapil Sibal had something to say on the subject' shouted the opposition, MPs taking a dig at the claim made by HRD minister.

However Sibal chose to keep mum on the issue. Arjun said the main reason for the delay in implementing the constitutional amendment providing OBC quota was due to the assembly election in five states which had brought the model code of conduct into force.

According to Arjun Singh "The central govt is aware of all the views expressed and shall take an appropriate decision without in anyway diluting the commitment arising out of the constitutional amendment in Article 15 (5).

The issue of a central legislation to give effect to provision of Article 15 (5) covering central institutions was under government consideration. He said a decision would be taken in the 'interests of society as a whole and in the interest of all"

"We have given a commitment to the nation. The constitutional amendment had been passed by parliament. There is no question of going back on it" he said. If some tensions were simmering in society on the issue it will be sorted out he said.

Attempts are made to twist the whole issue. Even the election commission was misled by such canards.

The DMK and PMK, MPs wanted action to be taken against the students opposing the reservation. Devendra Yadav JD (U) wanted the govt. to set up a time frame for the reservation of OBC in educational institutions, while accusing the students for paralyzing the system. Santosh Gangwar (BJP) was at pains to clear the misconception that his party was against reservation. "My party never opposed reservations", he remarked.

A caption of a aged pro-reservation activist trying to break the police cordon at the All India Institute of Medical Sciences during a rally in New Delhi on Wednesday, PTI.

The experts wanted to present the following points before the public /reader.

In case of pro-reservation protests held all over India we see the protesters are not only students. But in case of anti reservation they belonging to the only organization viz. youth for equality backed by the upper caste and brahmins. Of course they get the media's complete support. The same news is only carried in all print media but also repeated several times. But the pro reservation protests which involve old and young no particular group or organization, no specific caste or creed is never given any importance either in the print media or in TVs. Experts said the public / politicians and all the brahmins who are backed by the Hindutva force must know that by denying the majority, the equality, the nation is sure to face a drastic problem. So the expert said it is high time some justice is done to them. The paltry reservation of 27 percent is a drop in the ocean. But even if this is denied after so many years and after the 93rd amendment of the constitution the nation will face a very big problem which cannot be put down easily!

The New Indian Express

The government had passed orders reserving 50 percent of total seats for SC, ST and OBC in private professional educational institutions other than minority educational institutions.

While issuing order regulating the admission to private professional educational institutions (provision of reservation admission of students and fixation of fees) on May 25 the govt stipulated that a private aided or unaided professional educational institutions should reserve seats for candidates belonging to SC, ST and socially and educationally backward classes other than the economically weaker sections of society in the sanctioned intake to such extent as may be notified by the govt from time to time….

The experts were very happy Pondicherry govt. is the first one to implement 50 percent reservations for SC ST and OBC.

The New Indian Express

In a written reply in the house, minister of Human Resource Development Arjun Singh said that the proposed legislation would be in pursuance of the constitution (93rd amendment act 2005).

It enables the legislation to make law for advancement of weaker sections in matters of admission to all educational except minority institutions.

This law had already come into force on Jan 20 last, he pointed out. A draft of legislation in pursuance of the above amendment is under preparation and is expected to be ready for introduction during the current session of parliament" he said.

Meanwhile the oversight committee set up by the govt. for preparing a road map for introducing quota for the OBC has its interim report ready to be submitted to the Prime Minister in the next few days…

The panel formed after country protests against the move announced by the HRD minister has been asked to suggest ways to accommodate the OBC without harming the interests of general category student, UNI.

Experts have shown their protests to the govt. for always stating that the politicians are trying to accommodate the OBC without harming the interest of general category? To be more precise why not government feeling for the majority of the deprived classes! Experts feel that money alone wield power!
The New Indian Express

Union govt evades direct reply. Avoiding a direct reply to the supreme courts' query for details on the basis on which the number of OBC were calculated to justify the proposed reservation policy, the centre on Tuesday sought indicated that his question had been gone into the mandal case also after which a list of socially and educationally backward classes and OBC was drawn up.

SC query on quota. In compliance with the SC's direction that the govt. of India shall specify the bases applying the relevant and requisite socio economic criteria to exclude socially advanced persons / sections. (creamy layer) from other OBC, the government set up an expert committee consisting of Justice Ram Nandan Prasad (Retd.) as chairman. And as per its findings the govt decided to exempt creamy layer from govt. employments, the centre said in its affidavit in response to the petition before the Supreme Court, challenging the controversial reservation move.

The petitioner, Ashok Kumar Thakur had sought quashing of the quota move and termed the Mandal commissions finding about there are being 52 percent OBC population in the country as fictitious. He claimed that the National Sample Survey Organization (NSSO) and National Family Health Survey (NFHS) statistics put the figure of OBCS at 36 percent.

In response, the counter affidavit said no final decision had been taken on the policy. Any decision on this could be implemented only by a law enacted by parliament, it said adding, writ petition seeking to challenge a proposal even before it is fully formulated is therefore premature.

Experts feel that the Brahmins seek the support of the brahmin law even before the implementation of the policy. The functioning of the courts and its relation to the right of brahmins according to experts are inseparable. Whatever form of justice they want will be granted by all the courts. Unless the courts are first cleansed and laws of manu being practiced in these courts banished it will be impossible for the OBC to get any form of justice especially when the case is against them or involving them. Experts sympathized the OBC and were very sad to bring to light the bitter truth that they would never in their life get justice. The only solution to get justice is to have OBC courts or courts for OBC in which there is not even a single brahmin judge who gives justice to the OBC or passes verdict on the cases.

The New Indian Express

In the midst of his battle with AIIMS Director P. Venugopal, Union Health Minister Anbumani Ramadoss now faces a petition seeking his disqualification as MP.

The disqualification petition against Ramadoss on the ground that he was holding an office of profit as president of AIIMS was referred by President APJ Abdul Kalam to the election commission a few days ago sources said.

The petition filed by a retired bureaucrat was under the process of the commission sources said. Venugopal who is at loggerhead with Ramadoss had earlier sought the Minister disqualification in his petition before the Delhi High Court accusing him of holding an office of profit, PTI.

The New Indian Express

Against the backdrop of question being raised about its autonomy, a high level expert committee is confident of making the All Indian Institute of Medical Sciences (AIIMS) “as immune as possible from authorities”.

We hope we can come out with some recommendations which will make AIIMS as immune as possible from authorities “M.S. Valiathan who has been named the head of the four member expert committee on AIIMS affairs said from Manipal.

The committee is beginning its task against the backdrop of a bitter row between Health Minister Anbumani Ramadoss and AIIMS Director P. Venugopal whose sacking was recommended by the institutes governing council and stayed by the Delhi High Court.

Observing that the committee had been given a “broader assignment” through “comprehensive’ terms of reference encompassing all aspects the former director of Kerala based, Sri Chitra Tirunal Institute of Medical Sciences said the panel’s recommendations would “redefine” the functioning of the premier institution.

“Even the Prime Minister and the Health Minister have favoured autonomy of the institute. It is only how you interpret it “Valiathan said.

“We will meet in the next four five days in Delhi…

Admitting that the issue of retaining the best talent with the institute is a “difficult job” for the panel, Valiathan said when AIIMS was set up it was only such body in the country. ‘Today the picture is different. Corporate hospitals have better scope, much better scope, much better facilities. Many go there and even politicians go there so you cannot stop from going” he said PTI.

Experts still claimed that they are yet to find the meaning of “autonomous” when the institute is getting aid from the govt. They said autonomous colleges exist but they don’t get the total aid from the govt. These are central govt. institutions and are purely run by the tax payer money. Majority of the tax payers are non brahmins because brahmins hoard the money in mutt which has escaped from all taxes. So experts said these institutions have to abide by the govt. and in such situations court cannot unnecessarily intervene and pass a judgment against the govt that too when it is a policy made in a parliament with full majority. Certainly the court which is a brahmin monopoly (3% people) cannot rule the decision of the 97% of the population especially, the parliament asserted experts.

The Indian Express

Most major changes have roots in some crisis. It required a huge financial crisis in mid 1991 to deregulate imports and exports and dismantle the license raj there by ushering in an era of economic reforms. In 2006 we again have a crisis of sorts. There was a nation wide agitation mainly against but also for the numerical based quotas in central higher education institutions. The crisis almost went out of control it required the Supreme Courts intervention.

An over sight committee has been appointed to draw up a plan for implementation of the quotas. Three out of the five groups tasked to provide inputs have submitted their reports. Huge investments are recommended. The groups on IITs and IIMs have suggested phases in implementation. Some genuine concerns have been raised. The exercise is restricted to the central higher education institutions. These institutions enroll less than two percent of the student population. The bulk of the higher education is with the state governments. The basic private higher education is also large and expanding. If would be sad if these exercises is restricted to merely increasing seats in a handful of institutions overlooking the weakness in the system as a whole.

The systemic problems of higher education in India have been analysed in the recent study (ICRIER WP No. 180). Some of these problems such as unwieldy affiliating system, inflexible academic structure, uneven capacity across various subjects, eroding autonomy of academia, and low level of public funding are well known …

Higher education in India comprised nearly 18000 institutions. The majority of them are affiliated colleges that enroll 90 percent students at undergraduate level and 66 percent at the post graduate level. India has the highest number of higher education institutions in the world – almost four times that in US and entire Europe and more than seven times the number in China. Many of the Indian Institution is non viable understaffed and ill-equipped, two thirds do not even satisfy the minimum norm of UGC. All this makes the system highly fragmented, scattered and difficult to manage. There is a strong case of consolidation.

The distribution of capacity across subject areas and at different levels is uneven. For instance facilities for post graduate education in medicine are grossly inadequate. While there is heavy demand for some courses for many others there are no takers.

Public funding for higher education is small and unfairly distributed. Nearly one third institutions do not get any government funds at all. Of the remaining about half get some funding form central govt. with only a handful of central institution that cater to less than two percent of the students getting 85 percent of central funds, bulk of the higher education system depend on the state governments, most of them facing financial crunch. Pumping government funds in relatively better funded central institutions to accommodate additional quotas will further skew the public funding of higher education.

Faced with financial crisis all institutions other than central institutions and state institutions in Bihar and UP have raised their tuition fees. Higher education in India in now increasingly expensive, beyond the reach of the poor… The standards of academic research are low and declining.

A comprehensive review of higher education is required. In the present crisis we have a chance to bring about a knowledge revolution through comprehensive reforms in higher education. Let this opportunity not be wasted.

None of these institutions other than the one run by the central government get an aid of 700 – 1000 crores for a year says the experts. Further atleast because of such institution some of the OBC get a little education claims the experts.

The New Indian Express

Social justice and secularism the soul of our nation, had to be safeguarded even by sacrificing our lives, said Union Minister Ram Vilas Paswan.

Chairing a meeting on social justice to celebrate DMK president M. Karunanidhi's birthday organized by forum of educationists in the city on Wednesday, Paswan urged the people to fight for social justice on their own.

"For how long can we expect our leaders like Karunanidhi and other to fight for us?" he asked, urged people to think between living as slaves and as a life to self respect.

Even in the present anti reservation struggle Paswan said the UPA government had not "encroached" upon "others" seats. The centre had only decided to add 27 percent more seats to enable the same number of seats available at present for students under open category, he remarked wondering why the forces against reservation were still intolerant.

On merit being the argument he said to southern states which has more offered quotas than northern states are much more developed. Raising some fundamental issues like reservation for Dalits in minority communities he said reservation could never be compromised.

"More importantly when we are fighting casteist elements we should not forget those forces working to divide the nation on the basis of religion "he pointed out. In a nation with minorities as majority in borderstates he had secularism had to be safeguarded. He urged people to fight for social justice and secularism if necessary even ‘with blood’.

PMK founder Dr. S. Ramadoss in his loud and clear voice for reservation said there was no place for those with “compromising views’ when it came to social justice.

Demanding reservation in judiciary which was constitution he said article 312 read with 236 paved the way for the creation of a judicial commission at the national level that could well lead to a proper representation for SC /ST and OBC in the judiciary.

Quoting Vekatachalaih Commission report he said there were only 10 judges from scheduled castes out of 600 odd judges in high courts of the nation. When cases of social justice are heard in these courts how can we expect justice to be delivered? he asked. Further Ramadoss demanded reservation in media also.

On implementation of reservation in private sectors, he said a group of ministers had been formed but nothing had been done for two years. “Chief Minister Karunanidhi can do it. He will do it I will be the tool’ Ramadoss said.

Dravidar Kazhagam (DK) president K. Veeramani also highlighted the lack of adequate representation for backward classes and scheduled castes in the judiciary. More into praising Karunanidhi, the DK president said reservation is a “birth right” which is “non-negotiable”.

CPI state secretary D. Pandian said the backward class was still a “prisoner of thought’ even after centuries students from backward class has always excelled, provided given opportunities he said, citing numerous examples.

Reservation in private institution was inevitable as everything was being privatized said W.R. Varadarajan member of CPM central committee, demanding people to unite and fight for social justice ‘if we fail to do it now future generations will not forgive us’ he argued.

Terming “reservation” is nothing but “affirmative action” Union Finance Minister P. Chidambaram ever a diplomat and confining himself with in the “cabinet dharma” said no one was “hesitant” or “confused”.

It was the government duty to implement 27 percent reservation for OBC he said, emphasising that the centre was firm in its commitment. Reservation was the “only door” towards social justice in Indian society he said. But the strength of the opposition would also have to be evaluated for a final victory, he added.

Accepting that answers have to be found to difficult questions related to reservation in private sector as villages has more BC than others. Chidambaram said solutions could come only step by step.

Experts seconded the statement of K. Veeramani who said reservation is a ‘birth right”. Infact the expert said since caste is based on birth it is natural that reservation is OBC birth right.

Experts said they are yet to understand government needs to make any form of compromise with general category students when the OBC students have been denied their right and degraded and humiliated even after 60 years of independence. Experts argued govt. should be courageous enough to put an end to the anti reservation protests and its close associate the court of law!

Express news service

The centre will soon convince the students who are striking against the contentious issue of reservation in premier institutes, EM Sudarsana Nachiappan MP, Chairman of the Parliamentary Committee of Law and Justice said here on Tuesday.

Speaking to a mass demonstration organized by the all India confederation of OBC employees, welfare Associations, Nachiappan said that the backward classes which remained underprivileged and downtrodden for centuries, deserved support from the government.

It was centers duty to fulfill the aspirations of all sections of society. “But instead of understanding this and settling the issue in an amicable manner, a section of the student committee is getting impatient and fighting it out on the streets; he deplored.

The govt. is spending huge amounts of money on students in IITs and IIMs from the tax payer money. But most of them look up to the west and apply for a citizenship there”.

The father of Indian constitution B.R. Ambedkar had clearly said that those belonging to socially and economically backward classes needed support from the authority.

Contrary to media reports, union minister of Human Resource Development Arjun Singh who introduced the Bill had full support of Congress Chief Sonia Gandhi and other ministers Nachiappan said.

Due to the model code of conduct put in place during the Assembly elections we could not pursue the matter we will soon convene a meeting with the students and convince the striking ones" he said.

The government would also take up the issue to backlog in filling vacancies in jobs meant for reserved candidates.

The Indian Express

As government proposes to enhance reservation in institutes of higher learning a number of seats reserved for Scheduled Castes and Scheduled Tribes students in post graduate courses in the different institutes of technology.

The prescribed reservation for SC and ST is 15 percent and 7.5 percent in post graduate courses offered by Indian Institutes of Technology.

But only 11.91 percent and 3.95 percent of the seats could be filled in the SC and ST category respectively in 2005 – 06. Minister of state for Human Resource Development MAA Fatmi told the Rajya Sabha today in a written reply to a question.

In 2004 – 05 while 10.62 percent seats were filled by SC students, 3.96 percent seats were filled by those in the ST category, the figure standing in 9.88 percent and 4.2 percent respectively in the previous year. As much as 11.29 percent of seats were filled in SC category and 4.18 percent in ST quota in 2002 – 03 compared to 10.1 percent and 3.03 percent respectively in 2001-2002.

Fatmi said the inability of fill the reserved seats in Indian Institutes of Technology was attributable to various social and economic factors. Joint entrance exam, there has been a sizable jump in the number of candidates clearing the Joint entrance

Examination (JEE) to the IITs in the first attempt after the implementation of the new exam format.

In reply to a question minister of state of Human resource Development MAA Fatmi said while the number of candidates who cleared the exam in the first attempt was 28.49 percent in 2005 it jumped to 43.50 percent in 2006.

The total number of women candidates almost doubled from 29,291 in 2005 to 58, 997 in 2006.

In the new JEE format the application fees for women candidates is half of that from the male candidates in the general category. Increase in seats. The All Indian council for Technical Education (AICTE) has sanctioned an increase of over 42,000 seats in technical institutes this year.

The AICTE has approved an increase of 42, 277 seats in 529 technical institutes with a maximum number of seats – 31, 111 being in engineering and technology, minister of state of Human Resource Development MAA Fatmi said.

As many as 4, 785 seats will be added for MBA and PGDBM courses, while there will an increase of 641 seats for pharmacy and 740 for MCA.

Experts said in IIT Madras the brochure for PG courses does not contain any information regarding reservations for SC/ST. Even the advertisement for M.Tech, M.S. and Ph.D. admissions published in the dailies does not contain any information about reservation of SC/ST only the cost of the application form is less. Apart from that the same committee interviews them. They also write the same entrance test said the experts. Infact some of the experts said that if SC/ST the committee only takes extra caution not to select them. They said the liaison officer invariably favours the institution and never the student said the experts. One of expert is the father of an SC student who had done Ph.D, in IITM. The experts said SC/ST students are uniformly harassed by IIT (M) in one way or other. The fall in percentage in the year 2004 – 05 said the expert may be due to the impact of the BJP govt. said the experts.

The New Indian Express

While ruling out job quotas in top management levels in the private sector, the CII Assocham committee has assured that it will expand its employee base via a variety of measures.

The report on affirmative action has come up with proposals including mentoring 100 SC/ST entrepreneurs a year, providing for 50,000 scholarships in reputed institutes of higher leaning and setting up special coaching centers to improve the quality of schools in select district.

The Prime Minister has maintained that he is against imposition of quotas through law.

Sources said that CII and Assocham will set up a committee to chalk out a detailed code affirmative action for companies to follow. For instance the balance sheets of Indian companies will hence forth disclose another vital statistic. The number of SC/ST employees hired during the year.

In a closed door meeting with Prime Minister Manmohan Singh on Friday, members of the committee headed by JJ Irani presented the report seen as crucial to save India; inc for imposition of job quotas.

So far the Prime Minister maintained that he is against imposition of quotas through law. The committee has dropped its demand for a separate judicial commission and a separate law to make discrimination in the workplace based on reservation is punishable.

Experts were tensed to learn from this news that the company can discriminate in employing based on caste but the down trodden majority cannot have any claim over the jobs. Further the experts were vexed to note the stand of PM. They all said PM was only for the minority rich upper castes and brahmins and never for the majority downtrodden non brahmins which includes OBC, SC and ST. The experts viewed any judicial commission as a total waste of time and money for they give the most unlawful findings!

Express news service

The two major industrial bodies the CII and Assocham on Friday made it clear that the recommendations of their task force on affirmation action was not an attempt action to pre-empt any move on the govt. to bring legislation on reservation in private sector.

Spelling out the “cornerstones” of the recommendations at a joint press conference, CII, Chairman R. Seshasayee said that the private sector industry was against any legislation that would compromise the sanctity of its non – negotiable freedom of choice in employment.

As private sector accounted for only two percent of the total workforce, enforcing reservation would not solve the problem. However acknowledging their responsibility; the private sector industry wanted to be a part of the solution. The entire initiatives would be voluntary and not forced on the private sector by legislation as legislative compulsion would only be counter productive.

The industry also stated that whatever was done would further its competitiveness and not erode it.

Seshasayee said that the private industry had been taking some affirmative steps with regard to economic backwardness. “But we discovered that there was some disconnection’. While we looked from the economic angle, politics was looking at social backwardness. So we decided at both economic and social backwardness” he said.

“It is not as issue between the government and the industry. Both have the same objective to bring up SC/ST people to the higher echelons of the industry by imparting good training, education and equal opportunity to complete and come up on merit, he added.

CII task force chief said TATA son’s director Jamshed J. Irani who did not want to comment on the response of the Prime Minister to whom the report was submitted however said that the prospects of legislation on private sector reservation had now receded. Detailing the recommendations, he said that the industry would strengthen its HR systems to enhance opportunity for SC/ST candidates with equal qualifications and capability is terms of new recruitments. They would also be encouraged to bring more number of SC/ST to higher echelons.

Irani and Assocham task force chief K.I. Chugh said that the two organizations would encourage their member companies to start quality educational institutions in municipal areas with the help of NGOs.

Stating that it might not be possible for the smaller industry units to take up the affirmative action, they said that the major companies would also help entrepreneurial development of SC/ST through their supply chain management.

The larger companies would “mentor and create” atleast one entrepreneur each from the SC/ST category a year and cost and quality being equal and companies would give preference to enterprises and promoters, partners and proprietors from SC/ST category.

To help improve the employability of SC/ST candidates, coaching programmes would be conducted for 10000 students in 10 universities in the first year which would be ramped up to 50000 in 50 universities by 2009. The two bodies would come out with a come of conduct for their member companies which would come into effect from October this year.

The CII would set up a council to coordinate these efforts by its member companies and appoint an ombudsman to look into the observance of code. A similar set up would also be created in the Assocham.

Both Seshasayee and Assocham president Anil Agarwal said that the two bodies would apply the creamy layer principle while implementing the affirmative action and also made it clear that they would go by whatever definition the Supreme Court gave for the “creamy layer”

Agarwal said that effort now being taken for inclusive growth would ensure that there was a tangible difference in less than 10 years.

What CII is for and against on Reservation in Pvt. Sector proposed.

- Concrete steps to give better opportunities to SC/ST in employment.
- Companies will be encouraged to provide more executive positions to SC/ST.
- Create 100 entrepreneurs from SC/ST through supply chain management.
- Establish coaching centers in 10 universities for 10000 students in the first year and 50 universities and 50000 students by 2010.
- Sponsor training programmes at CII’s centers of excellence for SC/ST candidates.
- Provide a large number of scholarships for higher studies.

- Provide support to SC/ST students to prepare the entrance exams to professional and technical courses. Ten centers for 5000 students in the first year.
- Partner with NGOs to improve primary level education in government and municipal schools.
- Adopt a code of conduct for companies.
- Appoint ombudsman and council to oversee voluntary adoption of the initiatives.
- Creamy layer as defined by the Supreme Court will be left out.

Opposed

- Oppose legislation on reservation in private sector.
- No compromise on non negotiable freedom to choice in employment
- No compromise on competitiveness of the Indian industry.

Experts only made the note that the two major industrial bodies never spoke about OBC. Further they pointed out how many Brahmins are economically backward? They were still vexed to see that all the five Anil Agarwal, KL Chugh,. Jamshed J. Irani, R. Seshayee and Sunil Bharti mittal were all Brahmins or upper castes and none of them were OBC/ SC/ST but they were making policies for reservation in private sector for SC/ST. The pity is the majority of the people OBC had no representation and above all their problems about their position was not even opened. They said such neglect of the majority will one day or other find itself in the form of an out burst of a revolution. They only place before the reader why in the first place no representation is given, secondly why when reservation for OBC in educational institutions is discussed why even a mention about reservation for OBC in private sector was not dealt with.

Express news

The intensive media coverage of Mandal II may have succeeded in making the national emotive issue, but it does not seem to have made much difference to the established, pro – reservation national mood on this subject. Public opinion evidence of the last decades had always shown that the Indian public offers solid support to all measures of affirmative action including reservation.

It is no different after Mandal II, if one goes by the latest “The Indian Express CNN IBN poll conducted by market research agency AC Nielsen. The popular verdict in this poll is 57 percent to 37 percent in favour of the governments decision to extend quotas to higher education for the OBCs (the remaining 6 percent did not have an opinion).

Given the over sampling of urban educated and well off respondents in this survey it can be guessed that if all sections of the population were to be fairly polled the verdict will be 70 percent to 30 percent for the quota. This is exactly what it was six months ago when a national representative sample were quizzed on reservations in the CNN – IBN – HT state of the nation survey.

When asked to choose whether reservation in higher educational institutions will lead to equalizing of opportunities or a loss of quality and merit, 63 percent opted for the former only 34 percent for the later. It seems that in this instance the politicians and political parties read the public mood right.

But the political class seems to have under estimated the maneuvering space available to it for fine tuning the system of reservation.

While the respondents in this poll had little sympathy for proposal for doing away with all quotas and reservations, they showed a clear desire that the system should be improved to target the needy.

Given an option a two thirds majority favoured reservations on an economic criterion rather than use purely caste based criterion. The idea of giving reservation to the poor found favour with a majority of the SC/ST and OBC respondents also. This implicit critique of the government could have become a major issue if the media had generated more light than heat.

The poll reveals a narrowing funnel of information seventy nine percent had heart about the controversy, 61 percent had heard about the government decision, 42 percent knew it is going to benefit only the OBCs and 36 percent knew the decision was to apply to higher educational institutions. Taking everything into account only 14 percent of the respondents really knew what is controversy was all about.

The political calculations of swinging the “vote banks” also do not seem to have materialized so far. Only 12 percent of the respondents said their voting preference is likely to be influenced by this decision.

The ruling UPA appears to have an upper hand in the calculations of votes. It is more likely than the NDA to retain its old votes and to snatch voters from its opponents. But these are very small gains and the sample of this survey is too small to arrive at any definite conclusions in this regard.

The poll was conducted between June 3 and 6 among 776 adult persons in and around five major metropolitan centers (Delhi, Kolkata, Mumbai, Chennai, Nagpur). Half of the samples were selected from rural locations at least 30- km away from city centers. Quota sampling was used to ensure that the social profile of the sample reflected the country’s social profile, 25 percent SC/ST, 40 percent OBC besides 11 percent Muslims.

A survey of this size and kind in subject to a 5 percent standard error and therefore can only be broadly indicative of the larger population. Therefore the poll can be seen as the first broad indicator of the national mood after the Governments’ decision to implement the OBC quota in higher education.

The experts said it is high time the politicians of India take some care about the majority of the down trodden. They do not live in India according to the experts they only exist. Where is time for them to think of anything they have to battle for one meal a day. They work from morning to night yet live in absolute poverty and are denied every comfort of life!

Express news

While serious efforts are under way to implement the quota regime in central Universities, a parliamentary committee has said that there is a serious shortage of faculty at the current level of demands.

As per the information available in the 16 central universities there were 1,988 vacancies as on March last year with the Banaras Hindu University and Delhi University having as many as 687 and 396 vacancies respectively the Parliamentary standing committee on HRD said in its report.

"Situation is also far from satisfactory in JNU, Visva Bharti and Mizoram University", it is said noting that out of total vacant 1988 posts 1056 posts belong to lecturer category.

Expressing surprise over such a large number of posts of lecturers lying vacant which was the entry point to a university the committee apprehended that situation might be worse in state Universities,

In order to attract and retain the qualified and highly motivated teachers the committee felt that there was a need to supplement their salaries with an attractive package of perquisites and support for academic activities coupled with appropriate recognition with outstanding achievement.

The experts gave the following suggestions;

With the UGC scrapping the UGC net and asking qualified people with M.Phil and Ph.D to be directly appointed to the posts of lecturer it would be now possible for the government to fill all the vacancies.

When so many qualified people are unemployed it is very shocking to know that thousands of posts are vacant. The government can adopt UPSC to help in selection of a list of candidates and post them, then and their in universities.

Whenever vacancies arise! This would maintain a standard and quality in selection. The experts said even in famous institutions like IIT (M) the selection are not based on merit but on caste. Only upper caste that too arya brahmins are the maximum who get appointed some of them are from industries who do not know the abc of teaching or the subject or above all how to handle the fresh students said one of the experts whose son is studying in IIT(M). He says IITs sustain to get name and fame mainly because of the students who can memorize anything and everything for majority of the facilities are not good teachers. IIT faculties are always busy finding one or other means to go abroad mainly short term to earn more. The experts say this craze is not only with the faculties but also with each and every student. Most of the students who come to do Ph.D. want to visit other countries at least once. They mainly work so that they can visit abroad said one of the experts who is a student studying in IIT (M) for his doctorate degree.

The New Indian Express

A depressing calm is discernible in government circles on enlarging the scope of reservation in educational institutions. Having won the consensus on quotas for OBC there appears to be a sense of fait accompli in government policy. An Indian Express CNN-IBN opinion poll showed this weekend the public mood is over whelmingly in favour of fine tuning the reservation policy beyond the automatic installation of quotas along familiar axes. To recap the findings of the survey published in the Sunday Express a strong majority of respondents in a representative sample expressed support for reservation in higher education.

However an even larger majority – 67 percent recommended that economic criteria, irrespective of caste, must be included in identifying beneficiaries.

It would amount to rebuffing a historic opportunity if the government hastened its current programme to install quotas in higher education purely on the recommendations of the Mandal commission report. Affirmative action is the key to making India a more inclusive society. It is also crucial in realizing its knowledge potential. Contrary to the false choice often offered between excellence and equity, affirmative action in fact enlarges the pool of excellence by enabling people to overcome social prejudice and economic deprivation. The question is this: Is automatic identification of caste groupings alone adequate? Deprivation in India - as in every society – operates along various axes. Caste is one. But to believe that it alone should be a marker in programmes of affirmative action may be a poor reading of India’s Complex socio Economic realities.

It bears repetition that the sample for the opinion poll included the urban educated in greater proportion to their percentage of the total population. It is from the elite that opposition to quotas has come. Even so the survey is a plea for affirmative actions but in the programmes with more optimal targeting of beneficiaries. It would be a pity if the government ignored this sentiment and stayed with the easy course-this is simple adherence to the Mandal report.

The experts in first place flouted the very statement, “Equity isn’t the enemy of merit but mandal is. They said equity when it comes to OBC they feel merit is lost. Who is poor in India, everyone knows the poor are only from OBC, SC and ST. No brahmin is really poor, even if he happens to be poor he can get every form of support from the other brahmins. If brahmins are actually in the economically backward strata! the experts asks how many brahmins are doing the job of a scavenger or a sweeper or a vettiyan (working in the grave yard). There is not even a single brahmin employed in these posts! So the expert asks how can we say an upper caste or a brahmin is economically backward.

The experts firmly believed by chance the government adds the word economically backward immediately most of the brahmins will get a pay certificate to show themselves as economically backward and no OBC or SC and ST will be in a position to use the right of quota. Further as the brahmins can see any one will use that influence to get the seat under the economically backward reservation. Thus in India the very meaning of reservation will loose it meaning claimed the experts!

The article says caste is one to believe that it alone should be a marker in programmes of affirmative action may be a poor reading of India’s complex socio economic realities.

The expert said as far is India is concerned caste is everything. It someone is a brahmin he/ she gets more than their due; merit or real worth in the case of brahmins recedes back and only caste wins every thing be it a post or an admission to desired colleges or a prize or anything and everything the only reason for it is they have the major economy of India in their hands and above all they are the policy makers of India. For instance take the case of S.C. or H.C.s or Media, NKC or CII or the directors of IITs, IIMs, AIIMS… they are Brahmins. In certain places there is no representation for the OBC or SC. If this is the case the expert asks, how can we ever get justice? If OBC or SC or ST comes immediately they make a very big scene and say merit is lost and make street demonstrations and media which is their monopoly blows it out of proportion! To no policy they will nod, all plans should give them cent percent benefit other wise merit would be lost and the nation will be leading to crisis this is their usual statement said the experts. Unless we become policy for ourselves it is impossible for any betterment of the majority (OBC, SC and ST). What is even passed in parliament will be made null and void by the Court.

It is a pity in India the minority Hindu brahmins are ruling the majority!

The New Indian Express

The on going controversy over OBC reservations raises some important questions is reservation on the basis of the decades old OBC list really in favour of socially and educationally backward classes?

Is reservation at the post graduate level justified? Is reservation in favour of OBC in private unaided educational institutions established by non minorities violative of equality? Does the 93rd amendment make it mandatory to provide such reservation? Is the unanimity among political parties on the issue the equivalent of the constitutional of the proposal?

In my humble opinion after an indepth study of the provisions of the constitution and various Supreme Court judgments, the answers to all the five questions have to be in negative. Equality to all citizens is the founding faith of our constitution as enshrined in its preamble. It also constitutes one of the elements of the basic structure of the constitutions which cannot be amended.

However given the existence of inequality among the citizens an enabling provisions was inserted in Clause (4) of Article 15, which conferred power on the state to make special provisions for the advancement of socially and educationally backward classes (not castes) of citizens or for Scheduled castes and tribes.

Accordingly provision for reservation was made in professional colleges over 50 years. Now that a substantial number of persons who belonged to OBC have made tremendous progress thanks to these and other measures, they should not be considered as belonging to the backward classes any longer. In this context it is appropriate to quote what the Supreme Court had said in Peria Karuppan as early as in 1971. “Government should not proceed on the basis that once a class is considered as a backward class it should continue to be backward class for all times... Reservation of seats should not be allowed to become a vested interest". Once this elimination is achieved the percentage of reservation would also have to be proportionately reduced. Without doing so extending the benefit to the same list now and to the same extent is unconstitutional.

Regarding the second question in Dr. Preeti Srivastava 1999, a constitution bench of the Supreme Court held that at the level of super-speciality courses, no reservation is permissible since it is contrary to national interest.

As regards reservation at the post graduate level, the question was left open. However the reasons given for holding that there shall be no reservation for super-speciality courses.

As regards reservation at the post graduate level, the question was left open. However, the reasons given for holding that there shall be no reservation for super-speciality courses apply with equal force to post-graduate courses. The Supreme Court has held in Indra Sawhney that once a person gets selected to a cadre by direct recruitment on the basis of reservation, he belongs to one class along with others directly recruited in open competition and a second level reservation there after for promotion is not valid. This principle holds good for reservation in post graduate courses also. After a candidate secures admission to a degree course by reservation, all the students so admitted belong to one class. Therefore selection to higher courses should be according to merit obtained in the degree course.

As far as the last question is concerned clause (5), which was inserted into Article 15 by the 93$^{rd}$ amendment enables the state to provide quotas for backward classes in private unaided professional colleges established by non-minorities is ex-facie discriminatory in violation of law laid down by 11-judge bench decision in TMA Pai. Secondly it is also discriminatory towards backward classes as they are deprived by securing reservation in minority colleges.

Now for the question, is it correct to say that the 93$^{rd}$ amendment makes introducing OBC quotas in higher professional courses mandatory? It certainly is not. The constitution bench of the Supreme Court as early as in 1963 in M.R.Balaji held that Article 15(4) is only an enabling provision and that it does not impose an obligation to provide for reservation. The same holds good for the 93rd amendment.

Finally the question, is the unanimity among all political parties in this issue proof of the constitutionality of this step? Again it is not. Such unanimity only indicates that no political party wants to oppose it for electoral considerations.

Having said this, let me also state that any reasonable reservation made with due regard to the provisions of the constitution and Supreme Court judgments would be acceptable and in national interest.

The experts very strongly oppose the intervention of court for they said in all matters provisions made by the constitution was always made negative by the Courts judgments. They place before the public the only question. Is people's representation made in the form of constitution powerful over the courts or the courts (the laws of manu) can over rule and ride the constitution? Why no upper caste or a brahmin is raising the fundamental and the basic question when equal opportunity is denied to majority by a one digit minority and discriminated even at the entry point; what type of faith can these majority have in the founding faith of our Constitution as enshrined in its Preamble. We have lost complete faith in the Courts for law for even as of today it is only the laws of manu that is practiced when it comes to any case a brahmin verses a non-brahmin it is only the brahmin who win even if the non-brahmin has all justice in his side? What change should be made to get justice asks the experts? We never get justice from these courts says the experts. Innumerable examples and citations were given by them.

They (OBC) have not enjoyed any form of reservations in the central institutions for 60 years so the experts demanded immediate implementation of reservations in IITs and IIMs for which the OBC taxes are used in large amounts to the very functioning of these institutions. Even if they (OBC) have enjoyed reservations in professional colleges over 50 years, they are only run by state and not by centre, so the experts demand now reservation in central educational institutions.

The New Indian express.

The Supreme Court on Thursday again refused to stay the implementation of the government order providing 27.5 percent reservation to OBC candidates in admission to educational institutions.

Additional solicitor General Gopal Subramaniam informed the vacation bench comprising Arjit Pasayat and C.K.Thakker that the work in all hospitals has been resumed following withdrawal of strike by the doctors last evening in compliance with the court order directing striking medicos to call off their strike forthwith....

So when the topic is reservation comes unlike the OBCs all brahmins would be awake and media and court would act in support of them said the experts. Even before an order is implemented the brahmins show their hyper awareness and go to court for stay added the experts.

The New Indian Express

Agitating students on Thursday met President APJ Abdul Kalam and submitted a memorandum, demanding roll back of the proposal to provide 27 percent reservation to OBC in Central Government educational institutions.

“We are against caste based reservation. Higher education must be merit based. Only those students who meet the merit criteria should be in these institutions” said the students which called on Abdul Kalam at the Rashtrapathi Bhavan.

“Caste-based reservation is just a ploy of politicians to divide us. It is only vote bank politics at our expense” he added...UNI.

The experts were very much upset over the statement and our President meeting them. They recalled with owe how the President who is from T.N. a Periyar land in the first place gave these students an appointment. In the second place they strongly flouted the statement of these students they are against caste based reservation when all educational institutions of higher learning are ruled by the caste with brahmin domination and brahmin atrocities on the meritorious non brahmins. Further they attacked the very statement that caste based reservation is just a ploy to divide us. We OBC are always divided from the brahmins for the atrocities they do on us said the expert so it is not the government or politicians dividing, but it is the very caste which has already divided us they said. Further the experts said our President who gives the students appointment never carters to take any of the burning issues of the OBCs or agriculturist and has never given us any appointment however series it may be. All our problems are turned down by the Chief Secretary who is a brahmin they said. Because it was a brahmin issue and the agitators were brahmins the Chief Secretary has fixed an appointment for them the experts said! The only valid suggestion and question is will Brahmins leave behind caste?

The New Indian Express

The reality of this government's brave new world of higher education, as reported by this news paper should give pause to even those who have argued that 27 percent additional reservation is a blow for social justice. Take the Tata Institute for Fundamental Research, a seriously a world class institution that attracts global scientific talent and which figures in Arjun Singhs 100 plus list. Will the TIFR now have to desperately try and fill up its quota, restructure its cost structure and deal with UGC bureaucrats now vested with greater discretionary powers? Most probably in that case what will happen to its real agenda, the search for and nurturing of those with exceptional aptitude in some of the most abstruse branches of science. Even Arjun Singh should hesitate to answer that things will be exactly as before.

And the cabinet that Arjun Singh belongs to must scrutinize his proposed implementation plan of the 27 percent quota with the knowledge that on them lies the responsibility of not encouraging institutional sabotage in the name of egalitarianism. Singh has slipped in dozens more institutions, including those that receive no government aid under the excuse that their status that the UGC confers makes them ready for direct government interference. This is a bad faith of the worst kind, especially since it was never officially indicated that the extra quota was going to apply to institutions that are all bad private.

For these and other institutions Singh's dictum that admissions must be fair and transparent could very possibly mean a deliberating loss of academic freedom. The power arrogated to education bureaucrats to check for compliance and impose penalties for deviations on the other hand is reminiscent of the license permit raj.

Clearly Singh has decided that the 27 percent quota policy is going to be his answer to liberal policy making. The space for liberal policy making has been won after a long political fight that saw intelligent leaders of both the national parties realize that quality and efficiency in most fields cannot simply be mandated by fact. Should a politician whose personal ambition far exceeds his political support base be allowed to so seriously challenge that?

The experts feel what has merit and quality which the upper caste has professed and practised all these years have achieved, except they have made the poor poorer and the rich richer nothing more!

Unconstitutional way of functioning is the academic freedom argues the experts. The experts view that giving equality for OBC/SC/ST is loss of standard and merit and global scientific talent is destroyed. In these over five decades what has TIFR achieved? How is its functioning useful to humanity? Except the Iyengars domination in TIFR what have they concretely achieved? All these pure maths stuff is not going to serve any one claim the experts! The biggest Ayyangars regime is TIFR the experts said. He has over 90% of the Iyengars lists who have been in TIFR for these years, the experts added. What have they achieved in terms of society and nation? Western universities have scraped the maths dept which do not have practical use! Some of the western universities have merged maths department with Computer Science department claims the experts.

The courts cannot over rule or middle the policy made by the two houses for in India the courts are occupied only by the upper caste more so only by the brahmins. The experts further feel the laws given by these courts are nothing but laws of manu! In such a terrible situation it has become useless to take anything to courts for it is equal to their rule and power in the hands of 3%!

**Beware of the three Singhs**

*By Kumar Ketkar – Kumar, Kethat@expressindia.com*

Not many people remember Charan Singh the man who appointed B.P. Mandal to head a commission in the late ‘70s to collect data on the other Backward Classes (castes) belong to the Congress Party till mid-1960s.

He left the congress to form Bharatiya Kranti Dal because he felt the agriculturists did not get justice and a fair deal in the Nehruvian systems.

The Jat leader of Uttar Pradesh’s farmers was actually reflecting the aspirations of the “Middle peasantry” Translated in Caste terms, it meant that the so called OBC (Sometimes known as the middle caste and also “forward castes”) were left out and impoverished by the Congress capitalism. Almost all Lohaiaite socialists and some communists at once empathized with Charan Singh. While the capitalists used to criticize Nehru for being socialist the socialists and comrades would often condemn him for being capitalist.

Charan Singh was not advocating the cause of the agricultural labour, surviving far below the poverty line, nor was he championing the demands of the landlords. He was articulating the political and economic ambitions of the OBC, who were not part of the congress system. The congress system was a rather vague manifestation of aspirations of all classes and castes. Till the mid-sixties vagueness helped the party to appeal to all sections of the electorate. But two decades after independence there were murmurs of resentment in all the rural north. The so called Hindi belt had gone through social movements and the production relations remained steel framed by the tight caste structure.

The congress culture had been able to co-opt the Scheduled Castes and tribes by offering them affirmative action in the form of reservations. Indian Society took these reservations in its stride. There was a feeling of guilt among the upper castes which actually were at the helm of all the parties of left and right. The conservative opposition to those reservations did not cause protest like the one we are witnessing today.

The OBC may have been poor but were not untouchables. Also a large number of them owned land and cattle. Many OBC

were skilled craftsmen. The OBC occupied a vast socio-economic space in rural India. They were hostile to the former untouchables. Indeed caste relations between the OBC and the BC were far more vicious and violent than between the upper castes and the dalits. Most dalits were landless labourers and lived at the mercy of the OBC land owner. There were a few upper caste land lords, but by and large the agricultural means of production were controlled and owned by the OBCs.

Mrs. Gandhi knew of the OBC rebellion was aimed at destroying the congress culture. She shelved the Mandal report submitted to her government in 1980.

The upper castes had learnt quite early almost by intuition that the future belonged to the knowledge industry. They had begun to invest in education. Therefore the primary contradiction in the rural society was not Brahmins versus the BCs but the BC versus the OBC. Indeed hardline casteism prevailed more among the OBC and even the BC as they observed the hierarchy very strictly. Whether it was drawing water from the village well, going temple or marrying a lower of higher caste partner with in BC or OBC ranks, it generated hatred murder or mass violence Landlordism was formally abolished in 1951 but it had taken other forms. The restructured landlord and his henchman, the money lender and their allies in the bureaucracy represented the actual ruling class. The well-of OBC represented that ruling class.

This class began to become more and more assertive and aggressive after acquiring larger economic clout in the 60's. Charan Singh became its spokesman.

This OBC aggressiveness was predominantly a northern phenomenon. Maharastra and the southern states had their versions of middle class political movements but the violence and viciousness was not so intense there.

The congress political culture received its first jolt in the 67 elections as it lost power in the entire north. The newly formed Samyukta Vidhayak Dals were nothing more than the castes OBC fronts. Indira Gandhi had come to power just about a year ago. It was a blow to her and the rivals within the party decided to challenge her leadership. The party was split in 1969 and she had to once again take a kind of "hold all" approach to appeal to all sections of people. The slogan “garibi hatao” reflected the old congress style. The OBC fronts and their left allies could not hold on to power as they could not mobilise allies. Also the OBC base was not large enough. Their opportunity came in 1977, when the Janata Party came to power. The so-called JP movement had actually mobilised this class and given it political gains after coming to power, the Charan Singh faction confronted Prime Minister Moraji Desai and the compromise was the Mandal Commission.

But the collapse of Janta Dal government and Indhra Gandhi’s return stalled the OBC advance she knew the OBC rebellion was aimed at destroying the Congress culture. She shelved the Mandal report, submitted to her government in 1980. But the mobilisation continued and was challenging the Congress in every village in the north. But for Indra Gandhi’s assassination in 1984 and remobilization of Indian’s transcending caste and religion, the congress would not have come to power with such stirring majority under Rajiv Gandhi.

Rajiv was politically naive and did not initially understand the craftiness of his deputy V.P Singh.

It was not merely Bofors and Ayodhya that defeated Rajiv Congress. It was defeated by the OBC armies. How many people recall today that Mulayam Singh Yadav and Lalu Prasad Yadav gained clout after 1989? Both are products of the JP movement.

V. P Singh was going to complete the cycle started by Charan Singh. When he announced unilaterally the implementation of the recommendations of the Mandal Commission, he was consciously destroying the mass base of the congress. He knew that the congress could be kept out of power by appealing to the new OBC power and destroying the delicate balance that party had been able to maintain. He may have lost power but he ruined the rainbow coalition that the Congress had represented from the days of the freedom movement. Ironically Arjun Singh used to hate V.P Singh then not for political but for personal reasons. But today the same Arjun has picked up the same bows the V.P had collected to attack the congress.

If the congress system succumbs and surrenders to the mandal conspiracy today it will have committed political suicide. It seems that the party leadership is paralysed and does not know how to deal with this Mandal bolt. Arjun Singh knows what he is doing just as his predecessors V.P Singh and Charan Singh, knew what they were doing. All these Singh have only one aim, to wipe out for ever the congress coalition, which represented the whole of India.

If the Congress Coalition is weakened or wiped out the day would not be too far when the Indian state too would wither a way like the Soviet State did. The experts wishes to straighten the wrong record given by Kumar Ketkar, that Congress system was a strong manifestation of aspirations of the brahmins and rich society and not of the aspirations of all classes and castes. Secondly Ketkar is not aware of the fact that OBC and BC are not different classes and it is wrong to state that caste relations between the OBC and the BC were for more vicious and violent than between the upper caste and dalits.

The experts feel it is unfortunate Ketkar does not know the basis that OBC (other backward classes) and BC (Backward classes) and one and the same. With his wrong notion of the very caste if he is going to speak of reservation it is going to be nothing but a very big blunder. The experts feel as the media is his and he is in the media he should not make historical blunder and the New Indian Express should not publish such an erroneous article.

The New Indian Express

**Some Reservations it could extend congress and BJPs banishment**

How can you support the extension of reservations to the OBC while moping up the resentment among the upper castes? The two national political parties are being in credibly disingenuous in playing the quota issue for maximum electoral gain. There is a striking casualness in the way something is being offered for everyone without even the pretence of informed ground work. This is why the government appears embarrassingly exposed each time its emissary offers a settlement to striking medical students. Where is the precise plan of action for increasing the total number of seats in higher education by the promised 54 percent they claim? It is a legitimate query. The UPA government has committed itself to that increase by the next academic year. When a small group of healthcare workers striking duty is blatant abduction of their primary responsibility to patients can so easily appear to hold the government in breach of due diligence someone has to be held accountable.

But if this is the state of government strategy, the opposition too is being rather carefree in identifying beneficiaries that should be targeted. The BJP on the first day of its National Executive meet in Delhi welcomed the extensions of quotas to OBC in higher education with a strong of caveats. Remove the creamy layer among the OBC from entitlement to affirmative action. Add the economically weak among the upper castes. It is so evidently an appeal to various sections of the social coalition it had successfully wooed in the nineties in Uttar Pradesh – which goes to polls next year. One could argue that this was the game begun by the Congress for the same states but he BJP too jettisons its responsibility as the primary party in opposition by failing to keep focus on the more pertinent questions of infrastructure and delivery mechanisms.

Excellence and equality are twin objectives in India. Compromise on to its own peril in the phase of globalization, but can excellence be maintained by expansion of seats (a key need) without adequate investment in educational infrastructure. And is equity really best achieved through reservations on the basis of an old caste census? It is difficult to escape the suspicion that the Congress and the BJP are walking out of these questions with an eye on UP. Experience yields a sobering lesson. Mandal is actually an instrument that has banished congress from the political matrix of UP. And for the BJP gains in the Hindi heartland have come on the more inclusive bijli, sandak pani palnk. So what accounts for this political amnesia?

The experts feel that the Brahmin arrogance is at its height when they ask, can excellence be maintained by expansion of seats? Our only answer to them is first provide education achievements of the so called directors and deans of IITs and IIMs and that of the OBC working in the same institution? One is made to wonder what has made them directors and deans except the caste or the holy thread they wear. Why are they not putting a data of this and ask. Is not the excellence buried because of these meritless appointments? What can be their answer? Why has no media published this? The only answer is caste! The defect of the mind is caste! They are worried only they will not get any opportunity if OBC are given reservation. For they fear OBC are much more intelligent than them so they want to keep them below inspite of their qualifications and superiority. No medicine can be applied when the defect is in the mind of the Brahmins!

**Ethnic quota to be introduced in UK firms**

Companies that bid for lucrative government contracts in the UK are set to be rejected if they do not employ enough black and Asian workers under a pilot scheme.

A government committee has drawn up plans to question competing firms about their attributes towards race before awarding contracts a media report said on Monday.

The firms will be asked to provide figure showing the number of blacks and Asians employed. These will be compared to the proportion of people from ethnic minorities living near the company's offices and will be a factor when deciding on the winning bid "The Times" said here.

It follows the release of statistics showing that people from ethnic minorities are twice as likely to be unemployed as the white majority. The plans were approved last month by the ethnic minority employment task force which brings together seven government departments and comprises as many ministers.

According to Iqbal Wahhab, a task force member and Chairman of the ethnic minority advisory group a government blacked think tank the plans were moving a head quickly.

These new procurement policies are required to assist employers in making more enlightened recruitment decisions" Wahhab said "It may be unpopular in certain quarters but the fact remains that we should not have been in this kind of position in the first place".UK has planned to give reservation in firms!

The New Indian Express

**OBC quota to be implemented in phases: Panel.**

The Veerappa Moily Committee has suggested that reservations for other backward classes (OBC) in institutions of higher learning and technical institutions be implemented in a phased manner from next academic session.

The oversight Committee in its interim report presented to Prime Minister Manmohan Singh. On Tuesday, is understood to have given indications suggesting that the implementation be carried out in a phased manner in view of certain issue and constraints expressed by educational institutions through the five sub-groups.

However Moily's recommendation is at variance with the statements of Union Human Resource Development Minister Arjun Singh, who has been insistent that the quota should be implemented in one go.

The experts feel when the government is not interested or does not want to take the very action formulated by it they put panels, committees to perform the delaying tactics in a very organized way.

The New Indian Express

**PMK flays quota for OBCs in phased manner**

Expressing concern over the introduction of quota for OBC in all Central educational institutions by the Congress led UPA Government in a phased manner, PMK founder leader Dr.S.Ramadoss said that it was the constitution and the tenets of social justice.

Speaking to journalists in Thailapuram near here on Tuesday, he said the Central Government should introduce a legislation effecting 27 percent reservation to OBC in Central educational institutions, deemed universities and institutions aided and managed by the Central Government in the monsoon session and implement it in the current academic year.

The introduction of reservation in a phased manner would raise apprehensions among the people over the commitment of the UPA government to social justice he said adding that the move had disappointed his party.

Instead of implementing reservation in a phased manner, the government should do it in a single stroke he urged.

Welcoming the non-inclusion of the creamy layer in the Quota Bill Ramadoss said it was the only positive aspect of the Bill.

However, his party would oppose if the Bill was to be tabled in the Parliament on the last day of the Monsoon session and referred to the standing Committee Ramadoss said.

The PMK leader said despite the pressure exerted by Chief Minister M.Karunanidhi, the UPA government had staggered the reservation system for OBCs.

If the DMK leader decides to oppose the move as he did in the NLC disinvestment issue the PMK would extend support to him he added.

Barring Finance Minister P. Chidambaram, all other ministers from the state had been rendering support to immediate introduction of reservation to OBC Ramadoss quipped.

He further said that reservation must be introduced in appointment of judges from lower courts to the Supreme Court and in print and visual media as well.

He charged that the anti-quota stir in the capital was instigated by a section of the media. To oppose the party had planned to stage a demonstration in Chennai on Friday, Ramadoss added.

Ramadoss was very well aware of important factor that most of the anti-quota stir was in most cases instigated by a section of media. Media is playing a spoilt role to propagate anti reservation.

The New Indian Express

**Anti-quota stir resumes**

Hundreds of students and medicos today took to the streets here and clashed with police as they resumed their agitation against the UPA government’s decision to provide reservation in government aided elite educational institutions.

A day after the Union Cabinet decided to introduce a Bill in Parliament providing 27 percent reservation, students from several colleges and medical schools in national capital joined under the banner of “youth for equality” in flaying the move and demanding its immediate reversal in scenes reminiscent of the earlier medical agitation students converged in large numbers at Jantar Mantar, some 500 meters away from Parliament, chanting anti-Arjun Singh and anti-government slogans, waving the tricolour and singing the patriotic songs.

Attacking the government for taking the decision, “under political compulsion”, the medicos representatives said they would go to any extreme, including going on a strike, this time bound to force the government to recorder its decision.

AIIMS resident doctors have decided to go on mass causal leave on Thursday to protest the bill introduction.

“We will continue the agitation with the same intensity. Earlier we had called off the stir after the Supreme Court intervened. The apex court should now tell the government as well to reconsider the decision as the matter is subjudice”. Anil Sharma spokes person of AIIMS Resident Doctors Association said.

He said the medicos were also meeting legal experts and a indicated that they would move the Supreme Court against the decision to provide reservation for SC/STs and other backward castes in elite Central Educational Institutions including IITs and IIMs, from next academic year. The experts said their protest should have been first condemned by the media for the court has already banned all types of protest be it against or for reservations. Their disrespecting the very order of court is not mentioned by the media, claimed the experts. The media playing the role of supporting the anti reservation protestors are giving it multi dimensions and colours added the experts. Even before the Bill is tabled for its final approval the arrogant upper caste are making big protests which is given greater prominence by the media.

The New Indian Express

A delegation of union Minsters and MPs from Tamil Nadu on Thursday called on the Prime Minister Manmohan Singh and demanded immediate implementation of the proposed 27 percent reservation for other Backward Classes (OBC) in central government run higher educational institutions including IITs and IIMs without any dilution.

In a memorandum to Singh the DMK led DPA MPs urged the Prime Minister not to dilute the proposal in any manner and sought its implementation without any delay...

The MPs also submitted a copy of Chief Minister and DMK Chief M. Karunanidhi's letter to the Prime Minister last week and of a resolution adopted in the Tamil Nadu Assembly on May 31 in this regard, UNI.

The New Indian Express

Amid hiccups over the consensus on OBC quota, Oversight Committee Chairman, Veerappa has ruled out any "bulldozing" the proposal for 27 percent reservation in elite and central educational institutions.

He also asserted that protecting the excellence of educational institutions would be the "main theme".

There is no question of bulldozing the implementations of 54 percent. We are clear about it... I will examine institution – wise and decide whether they can go in one go or in phases", he told Karan Thapar's "Devils advocate" programme on CNN-IBN.

The Committee set up to prepare a road map for the 27 percent OBC reservation would take a holistic approach, he said adding that he has "an open mind" about introducing the 54 percent increased seats in a phased manner rather than at one go.

"It's an open-ended decision. I am keeping my mind open, this side or that ... Some institutions may come forward and say that with one kind of proposals made by the over sight committee we can go at one go. That is likely to happen. And some institutions will say it is impossible even if all the facilities are lying there, Molly said. Asked specifically whether the committee would allow individual institutions to decide on the implementation, he said "very true very true".

Molly admitted that all the five sub-groups which have recently submitted their interim reports to the committee, have said it is not possible to implement at one go. It has to be phased.

The former Karnataka Chief Minister said that he would not attempt to persuade the subgroups to change their mind. Instead he would hold extensive meetings with all the institutions as well as the experts from out side to understand then reasons before the Oversight Committee takes a final decision, would be based on "a holistic approach".

He said that he has convened a meeting of governing council Chairman and apex council chairman on August 17 and invited experts from all over the country on August 28 and August 29 to “ask them what should be the process”.

The committee was asked to submit its final report to the government by August 31.

Dismissing the dissent by Planning Commission member Bhalchandra Mungekar and UGC Chairman Sukhdev Thorat that question of staggering the implementation was beyond the remit of the committee, he emphasized that protecting the excellence of educational institutions was “the main theme” PTI.

The experts highly objected to the statement protecting the excellence of the educational institutions was “the main theme” or vedas.

They said if we the sudras learn or repeat slokas in those days by Manu we were punished, in that way we were denied opportunities of education and learning. Now if we are given education the educational institutions would loose their excellence. This is in the days of Modern Manu. Their Upper caste) minds have not been changed. The crux of the information is sudra should remain uneducated. What have they achieved in these 60 years by their excellence asks the experts?

The New India Express

Prime Minister Manmohan Singh is said to have asked HRD Minister Arjun Singh to exclude “unaided educational institutions” from the proposed 27 percent reservation or other backward classes (OBC) in higher education and incorporate in the Bill a provision that allows these quotas to be staggered.

Consequently, the HRD Ministry is redrafting the Bill to be introduced in the current session of parliament. its earlier draft Bill had under the definition of “Central educational institutions” included unaided institutions deemed to be universities such as Tata Institute of Social Sciences BITS Manipal Academy of higher Education and Symbiosis Pune.

The New draft Bill will only extend quotas to Central Government funded institutions. The earlier draft had proposed to fix the quota at 27 percent whereas the new draft is expected to say” up to a maximum of 27 percent” allowing for a smaller quantums in initial years, consistent with the Moily Committee's recommendations There two decisions are bound to create ripples in the UPA alliance. PMK leader and Health Minister Anbumani Ramadoss father S. Ramadoss is reaching Delhi on Wednesday to meet the congress president Sonia Gandhi and PM. Said a senior PMK leader. Along with the DMK, we can't compromise on extending reservation to private institutions... This is the core issue of our alliance in T.N.

The New Indian Express

There could be no clubbing together of Backward Classes (BC) and Extremely Backward Classes (EBC) for the purpose of providing reservation, the Supreme Court said on Tuesday.

Setting aside a Jharkhand High Court order upholding a state government notification of Oct. 10, 2002 ordering amalgamation of the two classes, a bench Justice AR Lakshmanan and Justice L.S. Panta said it violated article 14 which mandated that only equals should be treated equally. It is well settled that to treat unequals as equals also violates Article 14 of the Constitution, "It ordered.

The bench said that inclusion or non inclusion of a class in the list eligible for reservation was the job of the Backward class Commission as was held by the Apex Court in the Mandal Commission Case.

In the present case, the bench noted that the reservation for Extremely Backward classes had been granted on the recommendation of a Commission constituted for this purpose.

"The recommendation is based after a detailed survey. On the other hand, when the amalgamation of the categories took place, there were no material or empirical data to indicate that the circumstances had been changed other than a mere bald statement to the effect, it observed.

The SC also asked the state to constitute an expert committee to enquire into the recommendations/ complaints of over inclusion and under inclusion of the various groups.

Jharkhand which emerged from the erstwhile Bihar had adopted a Bihar law providing 73 percent reservations for appointments in government services.

The New Indian Express

The quota debate is often presented as a conflict between equity versus excellence. The chorus of the anti-quota protests, Getting in more OBC/ SC/ST, you dilute the standards of an institution since they are “undeserving”. The facts however are not black and white. Using the Right to Information Act The express obtained data for what are considered the two most competitive examinations in the country the AIIMs and the IIT entrance examination for 2006 and 2005.

An estimated 75,000 students vie for the 50 seats in AIIMS while 3.5 lakh for the 5000 seats in the seven IITs, the Institute of Technology at BHU and the India School Mines, Dhanbad. The results show that for both the AIIMs and IIT entrance exams the gap between the toppers in the general and reserved lists in large between 20 and 30 percentage points. But when it comes to the last successful candidate in both the lists for both the institutions, the gap narrows to almost half of that.

AIIMS entrance results for both 2005 and 2006 show that SC/ST students would not have got in at all had there not been quotas for them. Because of the topper in SC list for both the years scored lower than the last successful candidate in the general list but just 8 percent in 2005 and four percent this year. This is significantly not the case in the IIT entrance exam in both 2005 and 2006. Many SC/ST students would have got in even in there were no quotas – several have scored higher than the general students. Given that the AIIMs pie is tiny- only 50 seats each year of which 34 are for general and only seven and four for SC and STs respectively as compared to almost 5000 seats in the seven IITs this also lends credence to the view. That the larger the pie the higher the chances that the principles of equality and excellence won’t come into conflict. This also makes a strong case for an increase in the number of seats in higher educational Institutions says K. Srinath Reddy, Head of Cardiology at AIIMs”. I am not surprised that the marks secured by the last qualified candidates are high across the reserved lists as well.

IIT JEE

2005 2006

| | Seats | Top | Last | Seats | Top | Last |
|---|---|---|---|---|---|---|
| Gen | 4633 | 82 | 45.60 | 5527 | 82.25 | 27.90 |
| SC | 405 | 63.30 | 28.00 | 700 | 52.90 | 16.67 |
| ST | 76 | 53.90 | 27.80 | 156 | 48.19 | 16.67 |

People forget that in addition to these entrance exam marks there are qualifying marks so all candidates admitted have 60% in their school leaving exams and at least 50 percent in case of SC/ST. Failed or below par SC/ST candidates or General Candidates are not chosen.

Similarly in IITs All those who have made in general SC or ST have scored at least 55 percent in their 10 + 2 qualifying exam.

AIIMS's

| | 2005 | | | 2006 | | |
|---|---|---|---|---|---|---|
| | Seats | Top | Last | Seats | Top | Last |
| Gen | 81.0 | 66.17 | 34 | 76.33 | 68.00 | 34 |
| SC | 58.33 | 54 | 7 | 64.33 | 57.66 | 7 |
| ST | 54.00 | 51 | 4 | 59.16 | 56.83 | 4 |

The experts view these marks are also not a very fair representative of their performance for even the answer papers for the SC/ST are of different colour from the general category, then how can one expect any form of justice in the evaluation of these answer papers which is invariably done by the upper caste people. Can you not see this from how much they are against reservation? Will they ever be fair in evaluating the answer paper asks the experts? With so much of hindrance the marks which they get is a great achievement claimed the experts! If the evaluation of these papers are done by OBC/SC/ST teachers and answer papers are not of different colours certainly their performance would be different. Because of the discrimination only they loose marks claims many experts!

The New Indian Express

The Bills to give reservation to Other Backward Classes (OBC) in higher educational institutions are expected to come up for approval before the Union Cabinet on Monday.

This comes after Human Resource Development Minister Arjun Singh met Prime Minister Manmohan Singh on Friday, where he argued for fulfilling the UPA Governments commitment made in May last.

Sources said the meeting was scheduled for the evening and the only item on the agenda was the, consideration of the OBC quota Bills. They said the Cabinet note and the proposed Bill were likely to be circulated, on Monday itself for the Cabinet's consideration.

While the Bill on aided institutions deals with IITs, IIMs, AIIMS, IISc and Central universities like BITS (Pilani) and Manipal Institute of Higher Education,

HRD Ministry, sources said is ready with both bills as Arjun has been keen to implement the quota from 2007-08 academic year. The HRD Minister wants the Bills passed in the on going session of Parliament which ends on August 25.

There are two view points on the fate of the Bills which are going to come up for Cabinet approval. First it is expected that the Bills would only be introduced and sent to Parliament Standing Committees or a Group of Ministers (GOM) for further examinations.

The New Indian Express

Amid controversies and pulls and pressures from its own allies DMK and RJD, the congress led UPA government was forced to introduce the OBC bill providing 27 percent reservation for OBC in central aided educational institutions including IIMs and IITs.

The Central Educational Institutions Reservation in Admission; Bill 2006, which will be made operational from the academic year 2007 and staggered in its implementation makes no mention of the creamy layer, which has turned into yet another bone of contention between the Congress and its allies.

While Congress Union Ministers Kapil Sibal and H.R. Bharadwaj strongly opposed the inclusion of creamy layer in the Bill at the Union Cabinet meeting their DMK and RJD counter parts forcefully argued for its inclusion.

Introducing the Bill, Union HRD Minister Arjun Sing said, Millions of Backward Classes all over the country have been looking forward to get an opportunity for equitable access to the institutions of higher learning maintained by the state. To this end a small but irrevocable step is being taken today. It is now for the wisdom of Parliament to provide shape and substance to their dreams, which they have nurtured for decades".

The bill provides for a mandatory increase of seats in Central educational institutions, which would be attained over a maximum period of three years, beginning with the academic sessions in 2007.

The bill states that the Central educational institutions will increase the number of seats so that the number of seats excluding those reserved for the persons belonging to the SC, ST and OBC is not less than the number of such seats available for the academic session immediately preceding the date of the Act coming into force.

Arjun Singh said the Bill was brought in order to benefit millions of students belonging to socially and economically weaker sections of the society.

The New Indian Express

Anti-quota protests flare up. Demonstrations held near SC and India Gate.

Police use water cannons and teargas shells.

Demonstrators detained and some sustained injuries.

7 lathe charged in Kolkata Doctors, students may strike work from Friday if the OBC reservation Bill is tabled in parliament.

The experts are wondering at the blackmail used by the doctors and students. They said unless the government takes stern steps to stop such type of protests by the anti reservation cadres this would lead to other types of complication disturbing the very peace of the nation.

The New Indian Express

The intake of students in IITs across the country would increase by 54 percent once the centre gives. The green signal for implementation of the proposed 27 percent reservation in elite educational institutions, according to the director of IIT Bombay Prof. Ashok Misra.

However he made it clear that the number of students under general category would be maintained at the current level. "Keeping the general category intact we will increase the number of intake of students" Misra told the media here on Thursday... Most of the experts said IIT's teaching faculties do not do any work. They go abroad, show some input by conducting some workshops which benefits no one and so on.

The number of Teaching hours by the IIT faculties is very less even when compared to deemed and central universities. Further all the research scholars working for Ph.D in IITs and HTRA (Half Time teaching research Assistants) so they (facilities) also exploit the teaching labour of the students. Thus the experts suggest with the same infrastructure more number of students can be taken and the faculties can work in two shifts and two shifts for student classes can be taken!

Now the government should ask the Working hours of facilities in all IITs and take necessary steps to put the teaching staff to work more hours or the stipulated hours for the high pay they receive.

The New Indian Express

Even as Prime Minister Manmohan Singh is engaged in evolving consensus on the quota for OBC the DMK has demanded that the reservations should not be "diluted in any manner".

"I strongly insist that the UPA government should respect and respond to the peoples mandate by immediately extending reservation to the OBCs at all levels in all institutions under the government of India without any dilution", DMK. Chief and Tamil Nadu Chief Minister M. Karunanidhi said in a letter to the Prime Minister.

The Chief of the DMK, which is a key ally of the central government ruled that reservations were denied to the OBC since 1950 in educational institutions, state run units and the central administration.

"Thus the social justice to the OBC was denied and delayed for more than four decades. Delayed social justice is denied human rights for the OBC" he said.

Despite Kaka Kalelkar Commission and the Mandal Commission emphatically recommending reservations for OBC in education and employment in consonance with Article 340 of the constitution, “We are still making in irrational debate over implementing the reservations for the OBC” Karunanidhi said.

Under there circumstances, Karunanidhi said the Centre should not give in to the “unjust” demands of the anti reservation lobby who form five to 10 percent of the country’s population.

“They suggest to implement the reservation for OBC in a phased manner to apply creamy layer concept and to increase the seats for general category. I am of the firm opinion that, if these demands are not accepted it would amount to preferentially; empower the empowered at the cost of oppressing the oppressed.

It is strange to note that the creamy layer concept is not applied for the unreserved category till date” he said.

Karunanidhi said the reservation for OBC in education and employment is the “surest route to empower the hitherto socially oppressed sections of the society”.

The common minimum programme of the government has accepted the 27 percent reservations for OBC” without caveats”, he said adding consequently the people especially the OBC and Dalits solidly “stood behind us and voted us to power at the Union Level”.

He pointed out the Tamil Nadu Assembly has unanimously passed a resolution for immediately implementation the 27 percent reservation for OBC, PTI.

The New Indian Express

The lush green campus of the Indian Institute of Technology Madras (IIT-M) reverberated with the cries of “We want Merit” on Tuesday as a cross-section of students took out a protest march against the proposed caste based reservation in admission. The Union Minister for Human Resource Development, Arjun Singh had recently announced that the government would introduce 27 percent reservation for other Backward Classes candidates in admissions to Central Institutions.

The proposal as and when it is implemented would be applicable to the Indian Institutes of Technology and the Indian

Institutes of Management among the other higher educational institutions.

A cross section of the IIT Madras Students which had expressed their dissent over the proposal, walked a distance of about 3km from the Jamuna Hostel to the main entrance raising slogans against the reservations. They sported black arm bands and carried placards and banners denouncing Arjun Singh's proposal. Nearly 300 students took part in the march.

The participants took a stand that introducing communal reservation would dilute the academic standards at the IITs quality was very hall mark of the technical institution.

A note explaining the stand of the protesting students said, we are not against the concept of reservation but against this gross dilution of merit. Reservations should not promote elitist tendencies amongst the backward classes".

The students stressed that the union HRD ministry should use its resource in strengthening the primary education and said that such proposals of reservation would only lead to social inequality. They also criticized the proposal contending that in the absence of a well-regulated and unbiased system, the more will end up undermining the concept of democracy itself.

The students felt that the unity among students itself was under threat as it will lead to division among students and undermine the multicultural environment present within the IITs. Some of them even wanted to see the "Politically motivated" proposal taken back at the initial stages itself.

The experts were very harsh with the IIT (M) students community for calling the OBC reservation as "politically motivated" move. They felt these students leading a selfish, protected comfortable life do not think of the deprived majority that is why soon after graduation go to western countries either for further study or for earning. They have no mind to even think that they have exploited nation best education and wish to serve another nation without any concern for their own country. They are a set of educated rouges said the experts who are talking about merit or other false things when they are unconcerned of the nation or its improvement.

Further the experts said their view about their statements we are not against reservation, but against dilution of merit" is like shedding crocodile tears. If they have even a little concern over the nation will they go to other countries for money after exploiting our nation they asked! Now we OBC are concerned about the waste the nation is making by educating these unworthy selfish students! Secondly given the comfort and opportunities to the OBC they are certain to outshine these students added the expert which is evident from the marks they get in the school final!

The experts viewed this protest march as made by the IIT (M) facilities, that is why such a protest could be organized with in the campus and the media was also allowed to take a complete coverage of the protests! They further said the director of IIT (M) is morally responsible for it and he should come forward to give his stand and explanation for such a protest when they are functioning under the very MHRD with Shri Arjun Singh as the Central Union Minister Even as of now the experts said their was no absolute unity amount upper caste students and the OBC/SC/ST.

An article by Ravinder Rao an associate Professor, IIT Delhi which has appeared in Editorial page of the New Indian Express dated under the title.

Revenge of the creamy layer: Selection purely through lower cut offs is detrimental for the student and the institution. If bright students are the ones to be admitted to IITs and AIIMS or similar institutions they have to be actively identified at an early stage and given a good pre college education...

The experts ridiculed the views of the Prof. of IIT Delhi for he said his first feeling or statement that the students selected on the basis of entrance exams marks cannot be called as the bright and deserving students. For according to the experts they are the students who are economically very well off so that they can afford to go for the best of coaching classes for the entrance exams having good coaching they get good score by shear hard work.

So according to the experts the question of dilution of merit or any such sort is only a fiction (story) made by the arrogant upper castes!

The hard truth once reservation is given they (upper castes) are frightened, for they have to face a tough competition in employment and further studies. The experts gave their concern and fear how these teachers would behave with the OBC / SC/ ST students in class and in the evaluation of their answer papers. So they demanded the govt. to see to it that OBC/SC/ST faculties are recruited in IITs and AIIMS with out this the OBC / SC/ ST students are sure to suffer discrimination.

The New Indian Express

A law has been proposed to require central institutions to reserve 27 percent higher education seats for socially and educationally backward with in three years while adding seats to protect the prospects of unreserved candidates.

Besides 22.5 percent such seats already reserved for scheduled castes and scheduled tribes, in the central educational institutions (reservation in admission) Bill 2006, introduced in Parliament reserves for other Backward Classes 27 percent such seats.

But the Bill says that for reasons of "financial, physical or academic limitation or standards of education" the government may permit an institution to increase the intake "over a maximum period of three years". The Bill requires such permission by the government to be, "by notification in the official gazette".

The reservation is such cases "shall be limited" so that seats available to the OBC for each academic session are commensurate with the increase in the permitted strength for each year, it says.

The Bill introduced by Human Resource Development Minister, Arjun has triggered students agitation in the capital and else where, but no political opposition.

The New Indian Express

A Bill providing for 27 percent reservation for OBC in unaided higher educational institutions too will not have to wait long.

This was indicated by Human Resource Development Minister while talking to news persons here. "That Bill too will be brought in Parliament in due course of time", said Singh.

The legislations for introducing OBC quota in government aided institutions of higher learning like the IITs and IIMs was introduced on the last day of the monsoon session of parliament which concluded on August, 25, 2005.

The Bill has been referred to the department related standing committee which will debate its finer points. It will be introduced in its final form in the winter session of parliament - UNI.

The New Indian Express

The quota bill neither faced any opposition from the UPA allies, nor has any problem anticipated according to Union Human Resource Development Minister.

Talking on the sidelines of the launch of an educational initiative at the Indian Institute of Technology Madras (IITM) here on Sunday, Singh said the matter was now before the parliamentary standing committee and the cabinet. “Cabinet has the last word on the issue" said the minister.

When asked about the arguments from various quarters for and against the contentions bill, Singh said, "I did not notice any opposition yet and I do not anticipate any problems in the future

The New Indian Express

The supreme court on Friday directed the union health Ministry to keep 10 percent seats in each subject in abeyance for the reserved category students in admission to postgraduate medical science courses in the All India Institute of Medical Sciences (AIIMS).

The bench comprising Justice, K.G. Balakrishnan and Justice D.K. Jain issued the directions when the matter was mentioned before the Bench contending that the first course long for admissions to PG medical sciences courses was to start on Friday and unless and until interim relief was granted by this court, no seat will be left for the reserved category. The matter has been listed for hearing on Monday.

Earlier the centre had reversed its decision to cancel the all-India entrance test after allegations of wide spread copying in the test when the Supreme Court showed its disinclination against such a slop. UNI.

The New Indian Express

Union Human Resource Development Minister on Friday said implementing 27 percent reservations for other backward classes in Central Educational Institutions at one go was not possible.

Singh who made a surprise visit to Chennai two days ahead of his scheduled official engagements here on Sunday told Journalists at the air port that, the decision to implement the reservation in a staggered manner was taken after consensus in a Union Cabinet meeting attended by representatives of all allies including those from south India. It is not possible to introduce at one go".

According to him the UPA government cannot be said to have failed in implementing the quota at one go. "It is not a failure. It (the staggered implementation) is a practical solution for a difficult problem" he reasoned…

The New Indian Express

Congress legislature party leader D. Sudarssanam on Sunday expressed dissatisfaction over the reported remarks of PMK founder Ramadoss "creating" an impression that the Congress and its leaders were not in favour of reservation for OBC in higher education.

"It is unfortunate and regrettable that Ramadoss had made such remarks… No one can undermine the commitment of leaders like party president Sonia Gandhi, Prime Minister, Manmohan Singh and in Charge of party affairs in Tamil Nadu, Veerappa Moily towards the cause of reservation" Sudarssanam said in a statement issued here.

Claiming that the Congress stood first in working for the cause of social justice and social emancipation of minorities and backward communities Sudurssanam said, the Congress could say with pride that it was the one which had brought the first constitutional amendment when there was a threat to communal reservation in the country.

It was the Congress which constituted the Mandal Commission to uphold social justice, he added. "Also the Congress was instrumental in adopting a resolution in the State Assembly supporting the Centres decision on reservation" Sudarssanam said. "It was the Congress which constituted the Mandal commission"

“Like minded people act unitedly by realizing the problems of each other the rule of both the UPA at the center and the DPA in the state cando wonders" he added. However Sudarssanam said the party was proud to place on record the remarkable role of Chief Minister M. Karananidhi, Ramadoss the Union Ministers, MPs and MLA s in the state on the issue.

## Chapter Four

# ANALYSIS, SUGGESTIONS AND CONCLUSIONS BASED ON DISCUSSIONS, QUESTIONNAIRE, INTERVIEWS AND MATHEMATICAL MODELS

This chapter has five sections. Section one gives a view and analysis by a group of educationalists about the role of media in OBC reservations in institutions run by the central government. The present functioning of the media as described by the experts' forms section two. Section three gives the analysis by socio scientists on media and OBC reservations. In section four we give the origin of reservation in India and the related analysis by the experts. The final section gives suggestions and comments and views of political leaders.

In this chapter we analyse the role of media on (against) reservation to OBC in the educational institutions run by the central govt. of India.

### 4.1 A View and Analysis by Group of Educationalists about Role of Media in OBC Reservations

By the term media we mean both print media and T.V. Print media includes dailies, weeklies, monthlies as well as books. We give a complete analysis of how the news is projected which includes like the page in which it appears, the size of the head lines the size and the importance of the photograph and so on.

About things like who writes about it and his relative position in society and his social concerns selfless or selfish and so on is also analysed.

1. Student protests in anti reservation as covered by media (print and T.V).
2. Pro-reservation protests of students as carried out by media.
3. Coverage given by the media about other anti reservation protests.
4. The other pro reservation protests like human chain, seminars, conferences, etc-the importance given by media.
5. Essays or research articles against reservations for OBC … given in print media.
6. The pro-reservation articles or essays for OBC carried out in the print media.
7. The display of the difference of opinion between a Sudra authority and a Brahmin subordinate (A union cabinet minister versus a director of a central govt. run institution for which the union cabinet minister is the president) by the media (studied in chapter two using fuzzy models).
8. The upper castes versus rural poor, comfort, pomp and luxury versus sufferings, poverty and discrimination – coverage by media.
9. Sufferings of uneducated youth from OBC/SC/ST versus achievements of Brahmin youth projected by media!
10. Problems faced by OBC/SC/ST to get even basic education or needs displayed in media.
11. Reservation in industry!
12. Employment OBC/SC/ST in private sectors.
13. Reservation for minority in education and employment.
14. School education – media’s concern.
15. Constitution versus Law – media’s projection.
16. Role of Hindutva political force as covered by media.

Further it has become pertinent to mention that the views given are from the experts and authors have only put them and analysed it using fuzzy models. The views and news do not pertain to inflict anyone only to uplift India’s majority from the problems they face so that social equality reaches the last man.

The study is more of social concern over the majority who lead a deprived life in spite of scientific advancement. Till date living conditions of the last man in these six decades that too after independence from the British has only drastically deteriorated where as the luxury and pomp to which the upper caste especially the Brahmins are living have disproportionately developed beyond words. This sort of disparity will certainly disturb the equilibrium beyond repair! Thus it is high time politicians or public take up the issue and does laws or amendments so that some sort of balance is struck. We now proceed on to describe the role of media in blowing up the anti reservation protests carried out by students.

### 4.1.1 Students Protests in Anti Reservation as Covered by Media

1. The news given in the print media is briefly given in chapter three of this book. The views given by the (experts) are also included in chapter three, then and there. The attributes or concepts taken for analysis are only from the experts. Some of the articles have been disproportionately blown up by the media about the anti reservation protests are given in chapter two and three for the reader to know the true struggle! It has become pertinent here to mention that only one or two captions involving protests mainly done by female medicos found their place many times both in T.V. and print media. It was recorded by 98% of the experts it was like stills used in movies to attract the public. Is this justifiable asked the experts? Further they said every single anti reservation protest carried out by the students was given more importance that is basically wrong even along the very lines of media dharma. Public was very desperate for media was overdoing things which was against social justice. It is appropriate that the majority gets a chance for good and standard education in at least institutions run by the central government. The very act of disproportionately blowing up anti reservation protests is a fool proof to know how selfish the Brahmins are; when it comes to education for sudras (reservation for majority) they mainly make such protests because it is against the laws of Manu to give education to the sudras. That is why they want to deny the Sudra education in institutions run by the central government; is only dominated by them. The socio scientist said how can we get education in the institutions where they get education, that is why even the govt is favouring them for the govt. is run by them as they are the mouth piece of the govt. and as they occupy all high positions and more so are the policy makers for the majority. This is evident from the acts of NKC. Why is not even a single non-Brahmin in NKC? Who are they to make decision about reservation for OBC? Who are the Brahmin and upper caste medical students to protest against reservation and media to blow it up? The anti reservation protests were so much projected and given undue importance in comparison with the socially relevant protests carried out by SFI and others by media, both print and T.V. For anti reservation protests was anti social but others were socially relevant and so on. The media should answer this, as it is their moral responsibility.

2. The student ant reservation protest according to the expert was not spontaneous. Certainly it is proved beyond doubt that some one was making them to do the protests. The media has failed to investigate the brain behind these protests? Their failure to do so only proves that this protest is well supported by them. The journals/news papers have purposefully failed to analyse or investigate the people behind these protests.
3. The T.V media conducted a poll about anti-reservation thro' SMS in the NDTV and the results were declared as 70% against reservation. Is it a proper representation of data for who see English T.V and send SMSes. Once again what is the role of the majority? The bias of media; in reporting, was heavily condemned by the majority of the experts.
4. Almost all the anti reservation protests were carried out by the photos in the first page and the media had tried to show a big crowd in comparison with the pro reservation protests and above all other important social issues like suicide by farmers, death of patients and so on.

5. In all the anti reservation protest only women were projected on T.V and print media. Is it to get the mass appeal? Only a particular photo was carried out many a times.
6. Only the anti reservation protests were given that importance for nearly over 40 days and no importance was given even about the death of the patient due to these protests.
7. Why no action or investigation was carried when in the student protests only a gutka seller was forcefully for money made to immolate? Still date why no legal or police action for no falsely portraying the gutka seller as a student has taken place!
8. As reservation was main threat to Brahmins all media's irrespective of their party joined together to blow up the protests, which shows that media is run by Brahmins so only the Brahmins have organised it and given importance for the anti reservation protests.
9. Usually the media avoids any form of students only because it is an anti reservation protests it was given undue importance.
10. The media failed to carry out proper analysis (1) Why the protest was only by medicos and (2) Why it was carried out in the extreme northern states?
11. As anti reservation protests were only caste based and media is in the hands of the upper castes the anti reservation protest by the student was blown up beyond proportion claimed the experts.
12. Why these students who made anti reservation protest first make anti caste protests if they were genuine about the issue, for once there was no caste the question of reservation becomes null and void. So as rightly said by Periyar they want to claim themselves to be of high caste by birth and so automatically make protests for steps taken to achieve any form of social equality.
13. Experts claim that the anti reservation protests were well brained by the high Brahmin officials in New Delhi that is why these protests were given so much importance by the media.

14. The experts claim that the guts with which they stayed out of classes and work (duty) for over forty days, cannot be carried out unless they get the support of the head of the institution. This is evident from the fact the court too said that full pay must be given and their period of absence should not be taken into account! What a beautiful imposition of laws of Manu said these experts. Even now only the laws of Manu is in vogue. Both the orders of the court as well as the protests of the students are not only well covered, but repeated many times by both print and electronic media.
15. The media once again tried to use the good old policy of the Britishers viz. divide and rule policy. The student protests were well schemed by some, Hindutva brains who tried their hand in bringing a misunderstanding between the SC/ST category students and OBC students; for both would come under reserved category and they (media and the protesting students) did not know OBC were eligible for reservation. They tried to befriend the SC/ST students, were called by names like KATA students and so on. Even soon after the agitation some of the SC/ST students were forcefully asked to shift from the neighbouring rooms of upper caste students for the only reason they were dalits. Even they refused to dine together. This is the power of caste. The media has covered only a section of it.
16. Media both print and electronic never carried out the pro reservation protests of students properly.
17. The experts blatantly said that the anti reservation protests which sparked in many places in the city had a strong support from the Brahmin administrators. In actuality students were forced to protest. In IIT Madras a flash protest against reservation was carried out by students, with a complete support from the authorities. This is very much evident from the fact that no action against them was taken (for the court very recently put an order that no anti reservation or pro reservation protests should be made in the campus) in spite of it they made such anti reservation protests.

18. The organization behind all the anti reservation protest was the youth for equality. The experts claim that the youth for equality is an organization built and backed by the Hindutva force. They further claim the director of AIIMS, Dr. Venugopal has an upper hand over it. Are we fools to believe students or medicos protest on their own? No, they are only bothered about their future and their life, their security only on the basis of clear understanding and under taking they have taken to protest with the protection of the court which says they will be paid even in their strike period for which they did not go for work. What is law? “No work no pay” or No work full pay”. How can law be made special? When anything is special can that be called equality? If anything is unequal or partial can it law? This is the Manu law which is ruling India even today one law for brahmins and another law for non-Brahmins.
19. Finally the media which should have been just and work for national development, social equality and above all teach students to make one and all live in a better society without poverty is doing work in the reverse gear!
20. When the anti reservation protests by students, the minority had the guts to project and write in print media that reservation would stall the nation progress and speak of merits! What is the status of the majority?
21. It is claimed that since the media did not have BC/OBC/MBC/SC/ST journalists and reporters in chief executive positions, they have covered so many times the same protests and photos, there by giving so much importance to the students carrying the anti reservation protest.
22. The T.V. media took interviews about the reservation and this was very much biased for all most all the reservation interviews were only given by anti reservation high class people so the media was only protecting anti reservation aspects. Further as the interview was in English and many students from rural area who were really for reservation could not participate and air their views. A picture that reservation was not welcomed by one and all, was falsely projected by the T.V media to the world.

23. In the student anti-reservation protests the media showed very cheap portrait of politicians and the central ministers in power in an effective way. They again and again published these cartoons.
24. The Indian media that was always taking the politicians views in all matters failed to highlight their views about reservation. Very surprisingly the politicians were united and were all for reservations so they went all the way to project only the few anti reservation protests made by the students!
25. Youth for equality advised the authorities for the improvement of school education when the reservation for OBC were talked about. Now over these years have they ever talked about the improvement of school education?
26. In Indian history the politicians so called elite or educated ones and above all the NKC and the policy makers have always been against any social change that too when it pertains to OBC/SC/ST. In that situation when the politicians have unitedly come to give reservations for OBC in the govt. run institutions how come the media had the courage to speak against it. That is why for thousands of years India has never seen any change towards the real development of the majority. This is one of the main reason why the students were made to protest and it was covered by the media says the experts.
27. In the institution of AIIMS they cannot protest in campus as per the law, but with the help of the authorities, for they could put a shamiana and air cooler and carry out anti reservation protests. Do they respect law or law gives way to them?

### 4.1.2 Pro Reservation Protests of Students as Reported by the Media

We in this section analyse how the media has carried out the pro reservation protest.

The main discussions and arguments given by the experts were OBC/SC/ST students, socio scientists, doctors, politicians,

lawyers, public and some OBC/SC/ST union leaders are as follows are as follows:

1. The students said neither the T.V media nor the print media ever did any justice in reporting the pro reservation protests by them.
2. The anti reservation protests on the other hand was given the most prominence by both T.V and print media.
3. The main false reporting by the media was that pro reservation agitation was supported by the govt. but in truth no one ever supported them.
4. The real number of participants reported by the media and the disciplinary action taken by the police as claimed by the media was always false. The photograph never carried our total mass or number.
5. T.V. never telecasted our agitations or protests.
6. Our protests were never carried out under the shamiana with air cooler and our arrival was never by car or hired cab. We were left to the mercy of sun or rain. We cannot even afford for an auto fair and we travelled by bus or cycle to the place of protest. Our banners were our shouting with no mike.
7. As we have no backing from any one or party our protests were built on poverty and not on prosperity.
8. We were not even interviewed for our views as in case of anti reservation protestors.
9. The media tried to split us form SC/ST students. This only resulted in our (OBC) intimacy and most of the SC/ST students participated whole-heartedly in pro-reservation protests.
10. Only a few print media run by D.K and PMK covered our protests.
11. We did not receive encouragement even from our parents as the Brahmin students, even the authorities did not support us, they only found fault with us; said the students “Then how will we get the support of the media asked the desperate OBC/SC/ST students. We don’t have AC rooms or temples to discuss how to carry out our protests we did carry out our meeting under the sun in the corporation play grounds.”

12. The experts said if the media worked really for the uplift of the socially deprived and for the real nations improvement it should have carried out both the protest and should have shown to the world real differences between economic and social conditions of the anti reservation protestors and the pro reservation protestors and hence the real need for reservation. But on the other hand the media only showed the anti reservation protests and tried to project that reservation would harm the nation but in truth reservation alone can develop the nation and above all the majority. This biased way of projecting views has only shown the fact that the media being owned by the upper caste is against reservation and nothing more so only the importance of pro reservation was not-taken up as they viewed it against their own development.
13. Some of the pro reservation protests by law college students did not even appear in the print media but the legal views and consequences to protect the students who made anti reservation protest was repeatedly given importance in the print media claimed the experts.
14. Some of the experts said that medical students pro reservation protests in south was never given any media coverage. Some dailies reported it in few lines; some did not even make a mention.
15. Almost all pro reservation protests photos were very small in size or did not find the photos in the inclusion of news items. This was the biased action of the media claimed majority of the pro reservation protestors.
16. Even weeklies, bimonthlies and monthlies failed to include the pro reservation protests said the experts.
17. Some of the pro reservation students claimed as they were not glamorous the T.V. media did not take their protests into account. For who want non-glamorous rural poor, our very appearance irritates the upper caste how will they not protest our reservations though it is not their any right protest for they are not givers said some of the frustrated youth.
18. Some of the OBC/SC/ST youth said after all we are a burden to our parents as we are unemployed or

underemployed even after our graduation or post graduation or our professional degree (like B.E or B.Tech or Diploma) so they were willing to sacrifice even their lives for a good social cause. They asked how many brahmin educated youth are unemployed? Why Govt. is not taking a survey and implement reservation and put behind the bars the anti reservation protestors they asked.

The anti reservation protestors must be dealt as anti nationals or terrorists for their protest amounts to acting like traitors of the nation as they are hindering the development of the majority, they claimed. They blamed the Hindutva backed parties for they were the brain behind the anti reservation protests.

This is very evident from the print media and Sidu a upper caste cricketer, a B.J.P politician openly said he was against reservation for OBC; i.e. why we again and again say "equality for youth" is floated by the fanatic Hindutva forces.

### 4.1.3 Coverage Given by the Media about other Anti Reservation Protest

1. Most of the experts said that the coverage given by the print media and T.V regarding the anti reservation protests by other people was more in comparison with the pro reservation protests.
2. Actually it is a pity claimed the experts that not many conference or seminar or symposia regarding the anti reservation protests were organised. Youth for equality always made only hunger strike and nothing more. It was also only the medicos who carried out such protests in reality these hunger strikes were supported by one or two times by IT officials, software engineers and others who were Hindutva backed or labelled.
3. It came out as an interesting impact when the bar girls/dancers based in Bombay extended their support and said they would give a performance and collect funds for the striking medicos. The experts wondered the inter relation. This was very well projected in the media.

Prominence was given to this news in a leading daily by publishing it in the first page of it.

4. The experts wondered why the anti reservation protests was carried out mainly by youth for equality which always projected the photos of only two or three persons. Why these protests were not carried out by engineers and other student like lawyers?.
5. A few Brahmins from all educational institutions formed the human chains under the organisation of youth for equality, which was very well covered by the print media and T.V, claimed the experts. Apart from this there was not even a single political leader who ever had the courage to talk against or organise a single meeting against reservation in comparison to the many dozens of such conferences addressed by several non brahmin political leaders favouring reservation. Even several central cabinet ministers participated favouring reservation. It was a pity that both the print media and T.V never even reported it and in fact these central cabinet ministers even made fun of the media to write at least some lines about these protests and telecast it in T.V. It is but pertinent to mention both the print and T.V journalist covered the issues by recording the speeches and by taking photographs but only they have a strong inhibition to publish it claimed the experts. They record it mainly to know our views and pass on these views and material to their Brahmin boss said the experts. The experts further claimed that the reporters were non Brahmins and many showed enthusiasm to publish it but only their Brahmin boss where against publishing it. The Brahmin editors scanned the news and all pro reservation news or photos were just thrown out, be it a central cabinet minister or a veteran politician. They were basically against publishing of any pro reservation activities but gave prominence to all anti reservation activities. This was the secret policy followed by them. This is a common deal among all Brahmin dominated media claims the experts, but the reporters were forced to cover the complete programme mainly for them to know what was happening in support of reservations.

The experts said the fear driven Brahmins had no courage to come out in public and organize protests against reservation, for they knew fully well their act of protesting against reservation was very wrong so only they gave their brains to youth for equality, who made all protests and they could get only other than a gutka seller to be burnt against reservation. Till date no proper enquiry or investigation has taken place about how a gutka seller came for self-immolation against reservation when in the first place he was not a educated. The identity of the gutka seller as well as the persons who forced him to do so has not taken such an extra care to project anti reservation protests. The experts wondered why they have not given any importance or done any investigations over this. They claimed the media utterly failed in its duty by not publishing any information about how the gutka seller has appeared in the staged drama of anti reservation. The only conclusion drawn by the experts is that on one side the media is in their hands and secondly the politicians do not want to make a change towards real nations development or progress or is it fear of losing support from the policy makers who are Brahmins? The experts further claimed that the fear driven Brahmins to make sure that pro reservation protests are annulled and it does not take place once again took refuge from the laws of Manu (i.e. the supreme court?) that no protests about reservation should be done either against it or in favour of it. In spite of such orders, the experts said several symposia, conferences and seminars were held in support of reservations presided by central cabinet ministers, veteran politicians, powerful socio scientists in fact organized by students of backward and SC/ST castes. But the Brahmins (and upper castes) had no guts to organize any such ones but once a while to show their existence held some human chains and sported themselves with black ribbons in their left hand sat in silent protests with a collection of a dozen and odd students.

It is still a pity to record that these protests were well covered by the media but the pro protests symposia, seminar or conference was not given any coverage in the leading dailies or weeklies or monthlies or by T.V.

### 4.1.4: The Pro Reservation Protests … Given by Media

The experts claimed with vexation and depression that pro reservation protests even organized one day seminar by doctors/by intellectual group was given the least coverage. In that case they said how can any coverage with anti reservation protests comparison be ever made.

The following list of over dozens of seminars/ conference or public hall meeting or closed meeting presided by great shots like V.P Singh, Dr. Veeramani, D.K leader, Dr. Ramadoss, Dr. Deivanayagam, Dr. Ravidardan and other social movements were conducted. They claimed uniformly that no coverage was made in most of these and even if coverage was made they said it in a line or two in an obscure place in the dailies. Further they compared it with the anti reservation protests which found their place in the first page of the leading dailies with big colour photographs. This is sufficient for the reader to understand the media bias in India.

Some of the experts said even some of the dailies run by the OBC/SC/ST never gave repeated propaganda as made by the leading dailies dominated by Brahmins which were in fact the mouth piece of Brahmins and not of the nation but however all these information was only projected as the nations views which was the criminality unquestioned rights which they wield over the ignorant masses who form the majority of the nation.

The media of any nation which ought to be balanced in giving information in an unbiased way has utterly failed to give any proper information about the facts pertaining to pro reservation; conducted by way of conferences / seminar / symposia by intellectuals or socio scientists or great people of the nation. The only reason the experts attribute for it is not only the theme pro reservation but also they are not Brahmins so they see to it that no propaganda about them is made in their Brahmin dominated media.

The experts said media which ought to have carried out the analysis of pro reservation conferences / symposia and seminars or say even protests which always started and ended by speech given by VIPs was not even reported. An unbiased media should have done this simple work; for this analysis alone can educate one and all about the need of reservation as portrayed by OBC/SC/ST leaders. It is not voiced by them for their own well fare for they can lead a life of comfort and pomp without it, but a majority need it for the main need of the hour for the nation is to become developed. The media in their (experts) opinion has not analysed even a single speech given in any of these protests. They say all these only make one understand that the role and the very functioning of the media was to promote and establish the rights and the comforts of the so called Aryans and nothing more. It would be improper if media is termed as the mouth piece of the Aryans supremacy and it never gives the true news of the OBC/SC/ST say even of the socially deprived. The experts said even if a single conference would have been held by them against reservation how many times and in how many different angle and by how many different category of people the analysis would have taken place? So before we make pro reservation protests it is high time we make protests against the biased media or it is proper we make it as a policy not to buy their dailies or weeklies or monthlies so that they would be left with no other option to close it down. Before all these the govt should now at least not give any of the govt advertisement to their print and T.V media as they are not functioning to promote the nations interest but the Aryans interest. Unless such policies and plans are made in years to come India would not only face some revolution but also the difference between the haves and the have not would be so widened that the nation would not only loose peace but face uncontrollable terror for the majority of the have-nots are only from the OBC/SC/ST.

This is very much evident from the fact that the crimes, thefts, murders for money are becoming more and more!

The experts meekly claim that the very fact media shuns the news and views shows that the news and views would put the Aryans who have been enjoying the resources of the nation in an unlimited and disproportionate way that is unethical. That is why they do not want to speak about the reservation but in the most cunning way want to speak and give propaganda against reservation.

### 4.1.5 Essays or Articles Against Reservation in Print Media

Experts said the print media only in this period when the topic of reservation for OBC in the educational institutions run by central govt was to be made after passing the bill in both the houses the print media has published several articles against reservation. Most of the articles were penned by Brahmins. In fact several of the non Brahmins said these print media had made special requests from these persons to write articles against reservation. The main points carried out by almost all these persons, are against reservations.

The merit in these institutions would be at stake.

Reservations given to the OBC would stall the nations progress (what progress has been made even after 60 years of independence asks these experts)

By giving reservations to the majority who are poor and from rural areas the nation is acting in a discriminative way! (what about the equality the majority are denied asks the experts who is discriminated by whom? Let the nation as well as the Brahmin fanatics first settle the issue who is deprived by whom?) Can media if it was unbiased publish such scrap? Because the media is itself is totally in their hands and as the policy makers of the nation they have the heart to write articles in this form and make false propaganda about reservations said the experts.

Most of the articles suggested that the govt. before giving reservations in these institutions develop the standards of the schools and give every one a good standard school education so that merit would not suffer, (Here the experts put forth only. The following question why after 60 long years of independence only now they are thinking of developing the standard of all schools so that good education is given to one and all. What has made them suddenly realize this big truth that too after 60 years? What are their true ulterior motives asked these experts?)

Thus the only way to deny and at least delay the process of reservation for OBC with an immediate effect was to speak of improving the standard of schools, unfortunately the govt. has also become fooled and acted to this effect! When OBC intellectuals wrote a tight fitting reply those articles were not even acknowledged. None of them were published! This is the instance of media monopoly says the experts.

Another song they sing against reservation is that by reservation their (Brahmin) dreams are killed. The experts say for them only their bloody dreams are killed but in case of OBC/SC/ST, they are killed, for they are deprived of their right and the very opportunity to be educated in an institution run by the central govt. These institutions being funded by the govt. are run by the tax payers so it is just, all of them get an equal opportunity to get education from it. But OBC/SC/ST cannot enter the iron door of these institutions. Is it just? In the false name of merit they (OBC/SC/ST) are denied good education in these institution and they form not only a majority over 95% of the population but what more is they are the natives of India.

One of the articles that has appeared in the New Indian Express says that if reservation is given to OBC they would get education form Pakistan, Arabia, Kuwait, who is against it; asks the experts. It is better they quit our soil once for all for they have done enough harm! Further in that article they said already they are taking such steps. All the experts said in a single tone that if the Brahmins leave India; India would become a developed nation by five years.

The experts said the media has failed to project the real rural India for it always projects India as a nation full of prosperity, where as majority are deprived of every form of comfort, good food, clean water, education, etc., but the media never cares to project or publish the truth of poor sufferings. Unless the media is unbiased the experts claimed that any nation would never become a developed nation. That is why India only remains as a developing nation for over 60 years. In India media is not only biased it is totally concerned only about the 3% Brahmins welfare and indifferent about the 97% people basic needs. The reality is in these 60 years the rural poor has only become poorer but has lost their health and hence their longevity. They due to the evil effects of modernization face health hazards. Even the drinking water they get is polluted by chemical. In these conditions how can India boast of its development? Still the last man carries the night soil on his head for his living. Will this not be sufficient to know the real status of the majority? Why is the media unconcerned about all this? The experts said till a Brahmin carried the night soil on his head he would never know any pain of the poor. They said we should make them carry the night soil on their head, as equality; for they say equality is lost if reservation is given and even if one of them is denied an opportunity in their caste it is not equality. So if the nation is to apply equality in everything it is high time they should be forced to carry night soil on their head.

Thus the print media carried out such stupid arguments against reservation with no relevance to any equality and their writings only showed their haste, fear and insecurity for they know their identity which is in a waning state would be lost if reservation is given to OBC said the experts, that is why none of their reply to these articles was published!

The other attributes are already analysed in chapters two and three of this book.

## 4.2 The Present Functioning of the Media as Described by the Experts

The description given by an heterogeneous set of experts, many well educated have given the following short comings of the media.

1. Media by all standards have lost its media ethics
2. No psychologically trained or balanced person runs the media
3. Media is used for complete profit and for the development and propaganda of only a particular community
4. No distribution of power in any media (Daily news papers, T.V, etc) among all castes and religion.
5. Projection of Hindutva by inculcating in the young minds superstition and vulgarities.

6. Irresponsible way of display of stills advertisements etc; which will pollute and disturb the minds of the youth in their adulthood.
7. Media unconcerned about the welfare of the nation.
8. Anti reservation stand is exhibited at each stage.
9. No respect for people representatives only talks of court and law which is basically the barbaric laws of Manu.
10. No concern shown about rural poor or about their welfare.
11. They are so selfish and fear of losing their power due to reservation in this fear the media has lost its reasoning.
12. Media has failed to project the true statistics, the percentage of upper caste faculties in IITs, IIMs and AIIMS
13. Media unconcerned of the poor rural OBC/SC/ST youth and their problems.
14. Media has failed to project the corruptions of Venugopal but has disproportionably projected as if Dr. Anbumani was holding office of benefit that is purely caste based.
15. Media even now projects only their (Brahmins) problems or about a Hollywood or Bollywood in their day to day life and is least concerned about the social issues of poor or needy (for money is their main criteria not improvement or development or upliftment of the majority).
16. Why is not media giving the list of industrialists and the caste they belong to in India, who own half of India's wealth?
17. Media's lack of social concern.
18. No evidence is given to prove that merit will suffer but continue to project merit will be lost due to reservation.
19. Lack of business expertise to think how their propaganda does not take into account the social issues related to reservations.
20. Media does not give a complete projection to the scandalous life of the swamiji but only projects the govt. actions as if the govt has erred by taking action against swamiji, which amounts to anti national activities.
21. When govt did not accept the NKC views that issue was blown up out of proportion by the media and they made a big fuss and resigned from it. This act is again an act of anti nation as it was clearly, caste based support.

22. The media is more preoccupied with development of business rather than the business of productive activity that contribute to the development of the nation as a whole.
23. The media is very much concerned over a particular communities welfare and is practically blind about the even the basic needs of the majority but at every step talk of reservation as a means which will stunt the nations progress.
24. The media has not talked about the commercialisation of knowledge that has to be first settled but on the other hand harping on reservations. If access to good education is not available to majority then how can they talk of equality? This point is totally missed by the media to give any form of exposure.
25. The media has failed the project the weakness of the NKC, which has taken stand against reservations forgetting about the majority, but on the other hand thought of making education accessible by regulation of fees which is a major bundler by the so called elite members of NKC.
26. Has the media ever projected about the child labor which is at its peak or only among the non Brahmins with any care? Rarely the news about it that too when forced by some association or when the daily has no news to tell they would put some scanty news about child labour. Even the T.V media has drastically failed to bring to notice about this.
27. The media is unconcerned about the adults, that is, especially youth when it projects murders done by students. It is a pity to note that a murder of students which was carried out in A.P was given such coverage in the leading dailies that with in a month a similar type of murder was carried out in the daylight in the city college campus in the staircase by a student. First who is morally responsible for this? the media, the daily and the journalist who carried out this news are to be punished. So the media has behaved not in keeping with the interest of the nation.
28. Did the media ever provide to the public about the percentage of highly placed Brahmins and upper castes in the lucrative jobs viz; like professor, associate professor, directors, V.C., chairman etc; both in the industrial and educational fields?

29. Has the media ever worried about the unimaginable difference between the rural poor and rich? Should not this be the first concern of the media? Should not experts like economists, socio scientists write in the dailies about this social issue?
30. What are the steps taken by the media to propagate the school dropouts who are mainly from the OBC/SC/ST?
31. Why did not the media(both TV and dailies write articles to educate the school dropout's parents about the importance of education).
32. Why mandal name is dragged? How many negative articles about V.P Singh have appeared who has been the victim of Brahman's arrogance? He became a martyr by relinquishing the position?
33. Has the media ever professed the public to leave the notion of caste/Varna inculcated in the laws Manu to make the nation for advancing towards development.

    On the other hand all forms of superstitions and other irrational practices are imposed through the media very regularly in an insidious way poisoning the society!
34. Why the pro reservation protests held through out the nation was never reported by the media them as they have reported in an exaggerated scale the anti reservation protests? Should not the media act justly for both if it has some concern for the real development of the nation majority of the people.
35. The media is not acting in the justifiable way when the same photos of striking medical students against reservation was many times put, in T.V channels and also by the B.B.C (for world wide publicity) on the contrary pro reservation protests were never properly covered or never covered by them.
36. No empirical evidence that reservation will harm accumulation of knowledge proved by any commission scientifically in this Technical age.
37. Media's immature way of expressing the concern of the upper caste has only made one understand their anti reservation lobby.
38. Media by its activities has shown least social concern and utter selfishness of the Brahmins.

39. The media has proved beyond doubt its failure to understand the Indian Society for the brahmins have proved to be outsiders and have not become one with the natives.
40. Media has severed as an opponent to provisions to ensure equal opportunity to weaker sections.
41. Media has failed to understand that it does not have a mandate or authority to comment on the time tested policy of reservation, but on the other hand it has indulged in false propaganda.
42. The media has so contemplated with the NKC to express its opinion on reservations the day when the results for civil services exam were announced. This shows how biased the media is functioning.
43. Media from the day of its existence failed to expose the poverty of people for the people who die of starvation are only from backward castes/ST/SC.
44. Media has miserably failed to give any space or propaganda about the unemployed/bonded labours for they are mainly and only form backward classes and SC/ST.
45. The media has now also failed to provide any form of coverage to educationally deficient for this class is also from the OBC/SC/ST but has multifold times magnified and given false propaganda against reservations.
46. Media is trying to confuse the majority of the people between the two terms backward castes and backward classes, when the very father of our constitution Dr. Ambedkar himself has very clearly stated in the constituted assembly that the word classes means castes. The media is now projecting as if there was confusion on the terms ‘castes’ and “classes” in the very constitution. The media is unable to get any support against reservation so lives in trivialities and trying to create organized confusion.
47. Why is the not media projecting the legal untouchability and educational untouchability?
48. Why is the media not giving any proper and true exposure to the life of scavenger and give a proper propaganda in these 3000 year? Not even a single Brahmin has served as a scavenger and carried the night soil on head!

49. Why has the media failed to express about the majority for whom the life challenge is getting a meal a day not even a square meal? How can they raise high through a system of challenge? Why is the media, which is to be highly concerned and sensitive to the issues so blind of their problems? The only answer is, caste based selfishness as the media is mainly dominated by upper castes!
50. Why has media not taken any step to ban caste by birth? and start to profess there is no caste or ban castes. If caste is not present naturally the very word reservation based on caste becomes a meaningless term!
51. Why has media never openly criticised the laws of Manu? What does this imply? But on the other hand the revolutionary Periyar is portrayed as most of the media in a unfriendly way!
52. Dr. Ramadoss has condemned the T.V media as well as the print media for trying to project a wrong message to public by blowing up out of proportion the anti reservation protests and for hiding all protests of pro reservation! Media has failed to act justly but acts based on caste!
53. All these acts of media has annulled the possibility of any healthy discussion on reservation.

## Observations given by various categories of people on media

The media in India is dominated by the upper caste Brahmins and Banyas. The worst part of the functioning of the media is that they function in such an atrocious manner that only news in favour of them are published and given coverage beyond need and normality. This abnormal and noxious functioning of the media has resulted in the risk of the non Brahmins and minorities for their abilities and their normal need is totally suppressed. This was very glaringly seen in the reservations policy. The media only projected the anti reservation protests that too carried out by a very small set (two or three) of women medicos. The photo was repeated several times in the news and by leading news papers where as the pro reservation rally or protest or conference or seminar carried out did not even figure out in two or three lines in leading dailies and the T.V media. Of recent times after Periyar this has revoked and evoked a lot of controversies and faithlessness in media in the minds of non Brahmins. For only the rich upper caste alone can give the representation of his feelings through media. It has become impossible for non Brahmins to express their thoughts through media however great or powerful they may be. Under these problematic conditions we feel it essential to discuss the flaws of media. During our discussions some of the views expressed by majority of the four categories are included.

1. Almost all the non Brahmin students be it doctors or lawyers or engineers or other professionals very sadly said none of their pro reservation protests was ever given in the media be it T.V or otherwise, even the media run by the non Brahmins which just was very neutral to both pro and anti reservation protests with a few exception.
2. The Brahmin students were very much satisfied for their anti reservation protests were not only given due importance but the same protests were repeatedly published in dailies and also aired several times even in foreign T.V. Just they gave lot of publicity to project the anti reservation protests so as to create a blown up or exaggerated view. Infact some of the Brahmin students proudly commented "we are the darling of the media so even if we sit for few minutes with placards we will be given a lot of coverage". This was not just a sweeping statement, by a student but a very true fact about the functioning of the media. The media tried to create a wave that reservation was ruining the meritorious and infact it could not influence the majority.
3. Some of the non Brahmin educationalists said that several of their articles written to the leading dailies and other magazines which was pro reservation was never published infact not even acknowledged in several cases which were in response to the articles against reservations written by some fanatic Brahmins some of whom were industrialists, doctors and educationalists. This sort of media's bias or preferential stand is questionable against any norm of the media dharma; several of them who made these statements were holding high posts, some of them retired IAS officers and some cabinet ministers.

4. A very big seminar was held recently which was basically a pro reservation seminar in fact attended by over 500 delegates majority of them doctors presided by a central minister. Even in that seminar this minister said the media will not cover this truly not even a line was recorded in many of the leading dailies. So one can understand how this media functions only as a mouthpiece of Brahmins and Banyas. So even a cabinet minister has no say over them (media). Thus the media is acting in an autocratic manner giving false and biased information about reservations claimed many people.
5. Political leader like Dr.Ramadoss has very clearly said the media is trying to create problem that is hindering the peace of the nation by their false reporting. They are trying to spread false things, about reservation to OBC in education in the institutions run by the central government.
6. Even the D.K president Dr.Veeramani said in his statement the very news published by the media was to instigate unrest in the nation. These acts were very systematically carried out by the BJP and pro BJP parties like RSS and VHP. As majority in these parties are Brahmins and as they cannot say no reservations to OBC, for fear of loosing their vote bank which are only OBC, as the Brahmins form not more than 10%, they put up a dual form only externally which is pro reservation which is a total fraud and anti reservation which is projected in media which is their true colour.

## 4.3 Analysis by Socio Scientists on Media and OBC Reservations

What is the basic or the fundamental right the students have to say not to reservations for OBC who form the majority of the nation population? Their very act is unlawful. It is still undemocratic on the part of the media, which was encouraging the forward caste students to launch an ill-conceived agitation. Who are they to say no to reservation for OBC? The students have not taken the right to feed one and all and save people from starvation in their hands or as their duty. The students have not sought for the “right” to educate one and all” but suddenly have taken rights to say no to reservation for the depressed and oppressed majority when the government is imposing reservation. The government should have first put an end to their selfish act and secondly the media should have ignored their act of protests; but on the contrary the media which was dominated by forward castes, with extra ordinary indulgence showed the anti-reservation agitation of doctors and medical students as well as fully blacked out the pro reservation protests made by the doctors and staff under the banner of Medicos Forum for Equal opportunities in the same AIIMS campus.

The media has projected the wrong rights of the students protesting reservations for OBC as the right to protest; which is very dangerous to their own development. The media now at least should accept its selfish way of projecting and protecting wrong information for selfish ends? Media by these acts has not only ruined the nation development but its image in public.

The students by their faulty agitation, deprived the poor patients right to get treated; for over a month. Without any shame took the protection of the court for their salary without working. Can anyone with any little self-respect accept the pay when they have not earned it by their work? Even if government pays them, is it not their moral responsibility not to accept it, but they have heinous gone to court’s protection to get their pay which they don’t deserve and the law is so unlawful to order the government to pay the salary in the strike period of over 40 days. It looks very much for a balanced mind that all things that has taken against reservation including media’s support and legal support projects everything is occurring in a topsy-turvy way! All this show glaringly one thing, absence of OBCC/SC/ST in the news room and absence of OBC/SC/ST judges and lawyers who can do independent judgments and above all the absence of OBC/SC/ST policy maker for the real progress of the nation or the majority.

All these in the long run in going to impoverish the media, law and above all the nation.

The media which has blacked out all pro reservation protests is an unpardonable one for they by this act it has denied the public of any means of open discussions of the problem of poor Indians (rural poor) who form the majority of the OBC/SC/ST, which consequently will continue to keep India only in the list of developing nations even after 60 years of independence. This act of media has stopped the people to debate or discuss over the real need of reservation for OBC.

The media should have heavily condemned the insensitive and casteist forms of protests like - the “symbolic” sweeping of streets; the shining of shoes, the singing of songs of warning OBC and others to remember their place; all of them only portrays their arrogance. The media both print and T.V put them without any comments about it. Here also the media has absolutely failed in its duty. What these “meritorious” students were likely to achieve if they had such arrogance and contempt towards their own fellowman, who are more than half of the population of India. Is this ‘heroic’ or ‘shameful’? media till date has projected this as heroic which is an heinous act of the media!

It is an act of open shame for the anti reservation protestors for not even mentioning about seats given or purchased by the rich in these institutions. Will not this affect the merit? If they purchase the seat for money next they will purchase al marks also for money. So merit is money for these Brahmin students. Only the pro reservation protestors have asked this question!

Shamefully both the media and the anti reservation protestors have hidden this fact or ignoring it for mostly those people who get such seats are not only rich but they also belong to upper castes; one or two may be from other castes. Thus the court, media and the anti reservation protestors are blind to this question. What about NRI; seat in IIT, for several years survey show that these are the students who fail miserably in their exams and their foundation in basics is questionably poor; but magically they complete the course even at times with a first class. Why no one has so far questioned this! But if OBC/SC/ST is given reservation they trumpet merit is lost, discrimination is made and nation progress is hindered. What is nations progress first let them understand that the majority progress alone is the nations true progress. The nation has failed to become truly progressed for the majority in these 60 years, that is, after independence has become only poorer in all sense so how can the nation boost of its development. In truth what has really happened is that the rich more specifically the upper caste and the Brahmins has become not only rich but also well educated getting education mainly from the premier institutions of India or from abroad. Now when the government announces a small percentage of reservation for the majority of the people they are making protests which are mainly the after moth of the hysteria induced by the media coverage and some indulgent authorities who support and nurture such type of anti reservation strike.

How is it the High Court now remaining so silent about the protests carried out in AIIMS campus against reservation, when the very terms of a High Court order is, no protest or demonstration is permitted with in the AIIMS campus. Yet nobody demurred when the anti reservation students occupied the lawns put up shamiana and air coolers and received the "Solidarity" of traders, event managers and IT employees (whose employers usually ban their own staff from ever striking work). Not only this; the strike was not only carried out in shade with coolers but the protestors were lying down in comfort playing cards or reading novels in contrast with the pro reservation protestors who found it very difficult to buy the pa water packet for Rs1/- under the scorching sun. No media found it interesting to cover it as after all the protestors belonged to OBC/SC/ST or Muslims and not the high caste Brahmins, but the media has always forgotten that it is giving absolutely no space for the majority of citizens of the nation which is bound to end sometime. If the media was doing some justice it should have compared it and contrasted the agitations carried out by the pro reservation protestors and the anti reservation protestors, and their display of social and economic conditions. No article or analysis in this direction was ever carried out by the media. Here also the media has miserably failed. Only one more fact which we wish to place before the readers which was very much pressed by our experts is that the anti reservation protestors were (1) given right to protest under the shade in the campus itself with air coolers. (2). If at all these anti reservation protestors were to do protests in streets they came in cabs with water canes and in many cases 500 rupee notes were given to the protestors as beta for their food (3). The pro reservation protestors made protests under the sun found it difficult even to buy water for Rupee 1. They either walked to the place of protest or maximum reached the spot by cycles. No media ever cared to cover their protests. In all cases the protest was very well organized as it had nice introductory talk and also concluding speech given by some VIP unlike the anti reservation protests which was always absent of by any such speeches by VIPs.

Media has failed miserably to cover these contrasts. Even if they have covered the pro reservation protests our intelligent people would compare them and come to some conclusions about the importance reservation for OBC in central govt. run educational institutions. But the media has miserably failed to bring out even this simple message to the public.

The media has very conveniently avoided any further discussions about the heinous acts of the anti reservation symbolic protests like sweeping the streets shinning shoes and singing songs of warning but has only proudly published them in the first page with big photos. The articles sent questioning their casteist and insensitive form of protests as a news items was never published. Is this media ethics? Can media in the first place publish it without condemning such acts? When articles were sent to them condemning such acts, is it not the duty of the media to publish it. Why did the media not publish the articles sent in protest of these? Nevertheless the experts have asked us to record the following:

These symbolic protests adopted by them shows the casteist element dominated in the minds and hearts of these young protestors, which has forced them to be un rationalistic. The have absolutely lost their reason and rationalism. They have acted like forest animals insensitive to the feelings of the coexisting sons of the soil who are in fact doing all these heinous jobs to make the nation clean. It is high time they realize it or they will be forced to realize said the experts.

Further they said as they have come as beggars to our nation, they do not even know that it is still indecent to beg than sweep the road or shine the shoes. If we as a counter publish big posters with begging bowl in their hands and their middle fingers tied with dharba grass with a tuft adorning their head and write vulgar songs what will they do? We even in our penury, we do not display such emotions said the experts. They also said our poverty had made the beggars more insensitive and it is high time the politicians under stand it and do not woo them for compromise for who are they to give anything to us. By begging these thousands of years they have got everything from us said the experts and now we have become beggars and they are sitting on our head and protesting against reservation getting strong support for this from media and court claimed the experts.

Let them not forget the fact history repeats itself so once again they will be forced to take the begging bowel in order to make their living in India warned the experts. What the experts could not comprehend was the failure of the media to analyse their fanatic and low level protests which was basically casteist and this has wounded the feelings of the majority. Will reservation ever exist if caste alone was wiped out of the nation? Why no media persons are writing in this direction to educate people? They do not want caste to go they want caste and above all they want to claim they are superior to others based on caste which is based on birth. Let the politicians fearlessly proclaim we will have reservation in every thing as long as the nation has caste that too based on birth.

The Media should have criticized the way anti reservation protestors especially the medicos were protesting in the AIIMS campus lying on lawns under the shamiana and air coolers reading novels, playing cards or sleeping. The expert asks whose money is spent for all this, our tax! The honourable court ordered that as a special case (and cause) they would be paid for their protests during the strike period. Can we ever get justice even in genuine cases, even in any of the prima facie cases boiled the experts? So with no work, for lying in shade under shamiana with coolers they were fully paid for over 40 days! What a law? Why no media ever discussed this? Do media and thus law or judgement go hand in hand in caste together? Here also the media has utterly failed to bring in the truth about the agitations done by the medicos against reservation and getting complete judicial support. If this was ever done by the non Brahmins, that is the OBC/SC/ST the govt and the Brahmins would have not only denied everything using law but they would dismiss the students and would be forced to discontinue their studies or would be behind the bars for disobeying the High Court order that no agitation or protests can be done in the AIIMS premises claimed the experts. Further the experts said that they would have written about the national waste, our disobedience of professional ethics our irresponsibility of not taking moral responsibility of the patient death in the period of agitation and so on and so forth. They would have portrayed us (OBC/SC/ST) as if we have been in that place as anti nationals, terrorists and as criminals. But why is the media supporting them for these acts? The experts claim only caste and only caste of media in India is Brahmins! Their failure to be just shows their true colour and caste.

Another fact that media was casteist is clear from the fact that the effigies of the congress leader Sonia, our PM Manmohan Singh and our HRD central minister Arjun Singh were burnt. Why these students and the protestors were not put behind the bar? Why law did not act on them? Why media never bothered to speak about their wrong and corrupt actions? The only answer the experts give is that the protestors are Brahmins and the victims; Sonia, PM and Arjun Singh are non Brahmins. This is the difference which has led to all police tolerance, political tolerance, legal tolerance and media's support said the experts.

The experts said because of the hasty and imbalanced actions of the media which was encouraging forward caste students to launch an ill conceived agitation forcefully, that has spoilt all the possibility of a rational debate on what the best way of building an inclusive education system really is! Now how is the media going to redo this knot? Or atleast will the media stop criticizing the govt. and the public about the fact that govt. did not give way to open discussions and debates over this matter; when the media is the root cause for the closure of all possibility of any rational debate over reservation. Further the media has taken a very harsh stand by criticising the pro-reservation scholars. Several of their views was not even published by these biased media but on the contrary it took very many hard steps and procedures to publish all anti reservation articles in fact the experts claim that many non Brahmin intellectuals were also requested to write against reservation and their article were published in the print media. After all these non Brahmins want to be in the news, so they easily tame their views to set in the (Brahmins) grew! The expert further claim that not only this pro reservation issues but also several of the caste prejudice in higher education is never given any importance or to be more precise never published. Even in the dining hall of the paramacharya of the Sankara Mutt at Kanchipuram, a non Brahmin however high caste he may be cannot dine with the Brahmins. Disgracefully even today the Brahmins are served food in the room where as the non Brahmin devotees are served food in the open corridor which is adjoined as the pavement of the dining room. So even while the non Brahmins are eating constantly the servers and others walk for they do not serve the food, mind you on tables! So when the Hindu religious mutt is holding the caste difference as Brahmin and non Brahmin, it is not something strange when the Brahmin students do not sit with the non Brahmin students in the dining hall in their hostels. How will ever the media write this for the instruction and the practise is very clear from their religious mutt! Has the media in these several years of existence ever written against this dining process in this particular mutt? (Baring some news papers like Viduthalai etc). When will the media ever accept any form of social reformation? Here it has become pertinent to mention that Periyar closed down some schools in his time for giving food in an discriminatory way. This special schools closed down by Periyar was the Cheranma Devi gurukula run by Vavesu Iyer. This has never been recalled by the media even today the medical students do not like to dine with dalit students. At that time Periyar was living to close down the schools now there is no Periyar that is why the Kanchi Mutt is functioning peacefully claimed the experts. Thus the media that should build in equality in society and lead the nation towards social reformation is not doing so. It favours such heinous acts of discriminations, that is why for over a month-long all anti-reservation agitation by the doctors and medical students of AIIMS and other colleges was showed and given more propaganda with extraordinary indulgence claims the experts. Except for the medias support, such agitations against reservation would not have been carried out, only their encouragement made them agitate in an insensitive and casteist way. Thus it is not an exaggeration if the experts claim media to be the major culprit for the anti reservation protests which has taken the lives of innocent patients. Has the media ever talked seriously about the death of these patients or has it even published a single article condemning the acts of the doctors who are the root cause of the innocent patient death, who were in their prime of life. Has the media ever written an article to make the doctors and medical students realize their irresponsible way of life, leading to the death of the patients so that they realize their fault. Infact their failure to do their duty has resulted in the death. Why has the court of law not inter feared with this? Does it mean non Brahmins life is less important than a Brahmins luxury? When Dr. Anbumani interfered and he went to serve the patients the media shamelessly wrote that Dr. Anbumani wants to show his face in the media and become popular with the medical coat on and so on. Why did media in the first place publish this article? Is this act of the print media justifiable?

Media in case of anti reservation protests always exaggerated the number of protestors and reported falsely the number of protestors. In the case of pro reservation agitations it never reported it properly. This is the questionable way in which media has been conducting itself claims the experts. They (experts) gave hundreds of such instances. Several of the experts said that the media was never pointing out the faults of the anti reservation protestors. The highest mistake made by them was they abstained from duty over 40 days resulting in the death of a patient even this was not brought out properly by the media. The media had the mind to encourage the placards written by the medicos said kick out the monster and the monster was drawn as Arjun Singh. It is a pity neither the govt.

nor the public had taken any action against these students. On the other hand the President, P.M. and others have met them and pleaded them to stop the agitation and return to their duty? Why this soft approach?

In the history of student agitations the govt. has never taken such an approach. They could afford to meet the president-Kalam and have the guts to demand the roll back of quota proposal. On Saturday May 6, 2006, The New Indian Express first page reported in bold that SC rethink on quota for PG medicos. "Vacates interim order on providing 10% quota for SC/ST".

The experts showed this instance that the courts and media go hand in hand unless, the media acts with independence it cannot do any act with integrity claims the experts. However the experts were wounded for both the court and media being dominated by the Brahmins they were like the two sides of a coin. If this is the case how will even the media do any form of justice in reporting the need for reservation and pro protests for reservations! They all uniformly claimed that unless the media is taken by the govt. and the court had reservations for OBC as per their percentage of population and the Brahmin domination is totally annulled in courts and their indirect ruling is stopped; the nation can never in its life time become developed and SC/ST/OBC would become as endangered species living in India!

The experts also in a single tone raised the very question raised by our Honourable Minister Arjun Singh; "Why have parliament if we cannot honour our commitment on quota? "The experts said it is very illegal when the issue of reservation for OBC in central run institutions was passed in both the houses of the parliament with over whelming majority how come the courts have even the right to analyse about it or ask questions like.

"OBC quota can divide nation tell us how you decided?" What is the basis of the norms for fixing OBC quotas? What is the rational for this? If the proposal is implemented what are the modalities of such implementation and the basis for the modalities?

The PIL filed by Supreme Court senior advocate Ashoka Kumar Thakur and Shiv Khera, sought a direction quashing the amendment to the Constitution Article 15(5) paving the way for enactment of a law providing for a further reservation of 27 percent apart from the constitutionally sanctioned SC/ST reservation of 22.5 percent. Further in a brief order the judges said the policy if adopted will divide this great country on caste basis" and has severe social and political ramifications". The judges said their "need not" be a stay at this junction…

The questions raised by the experts are as follows:

1. In the first place who are the S.C. and H.C judges and who appoints them?
2. When both the parliament houses which function as the mouth piece of people the court which is appointed by the people representative cannot act against the people wishes.
3. When any reform made by the two houses of the parliament for the well fare of the majority it cannot be a wrong move leading to any form of social or political ramifications.
4. The petitioners are according to the experts well backed by the Brahmin and Hindutva fanatics to file the PIL. If the courts were mouthpiece of the nation certainly even the admission of the PIL would have been not carried out. This instance is only a foolproof for the politicians and the majority to learn that the courts` are only mouthpiece of the Brahmins and never a tool to uplift the national interest or the nations majorities welfare or do justice claimed the experts.

If the Hindutva Political force had the guts it should have made their protests in the houses of parliament and not in the court. This clearly shows us as the vote banks are OBC and if they make any form of protest they would loose their stand so only they made this nasty move for they are always certain the courts are theirs and the law which it proclaim is the laws of Manu. So that they can over ride the people wishes and the amendments passed in both the houses with over whelming majority. The experts condemned the sneaky approach of the Brahmins using both media and the courts that were in their monopoly.

The experts charged the media for its malfunctioning in every move regarding reservations for OBC. The media and the courts can now join together and say that caste cannot be used and the court can pass a law to ban all castes in the nation. This is the most reformative step for the nations progress. For if no caste, exists we need no reservation and no discrimination. For the experts claim in India, the caste is not so much in the external appearance only it is in the mind of the people more so in the mind set of the upper caste and the Brahmins. This is hundred percent in the heads of the Brahmin religious leaders and the temples which function under them said the experts. That is why for the Brahmins caste is very important and they will do any thing to protect caste.

Further when the MPs in the upper and the lower house of the parliament did not open any thing about reservation but favoured reservation, but the BJP cricketer Sidhu and a few others talk against reservation outside the parliament. This is the double stand they always hold one can recall in the year 1948(5-4 -48) in Delhi a Brahmin Anatha Shayanam Ayyangar spoke in the parliament one should find a unifying fact other than religion; towards this goal. He said one should build non religious government. The next day he in Madras before the High Court-Chief Justice, Rajamannar, he give a upanyasam about the coronation of Sri Ram; this upanyasam was published in the Hindu dated 5-4-1948. This instance is placed before the readers mainly to show not only the double stand and the double tongue of the Brahmin but their view or their behaviour in the Parliament; distinctly is different in real life. If one does not think it as an exaggeration the views are sometimes diametrically opposite. They have not changed in their basic character. They pretend to be some personality in the parliament and do and behave right opposite in public as well as in private.

But they have changed their media strategy this time some of the experts claimed. The experts said before the media when some important meeting or protests by the OBC/SC/ST are held they will never cover it but by chance during these protests some mishappening occurs like burning of the bus of some damage to the public property then the media will flash it so well, not about that protest but about that mishap in such a way that the public image about the persons who organized these protests are brought down to dust. This was the medias usual strategy. But in the pro reservation protest the media has changed its strategy the experts said. If some mishap occurred during the pro reservation protest they did not blow it up; for two reasons (i) The public should not think that the pro reservation protestors are so powerful and so it is better we do not protest against reservation (ii). If they have to cover the mishap they have to atleast mention in a line or two such protests for pro reservation has taken place; now a days they avoid both very cleaver, said the experts. Thus the experts said not even 10% of the pro reservation protests found their place in the media be it print or T.V. thus all the pro reservation protests were black out by the media. This act of the media has only stalled the nation from progress. the views of the people who stood for reservation is put in black and white before the public one can further develop a debate over it and a healthy understanding between the upper castes and the non Brahmins could have developed; but by their total blackout they have only made the misunderstanding between these two sections of people sharper!. This is also the vital place where the media has failed to do its duty claimed the experts.

Another important factor discussed in the media especially print media against reservation was that by giving 27 percent reservation to the OBC in the central govt run institutions of higher education that the 27 percent quota would cause institutional sabotage in the name of egalitarianism.

Media also is making a propaganda that this 27 percent additional reservation is a blow for social justice. What do these media imply by social justice, the rich who are mostly from the Brahmins and upper castes should continue to live in comfort and the poor majority from the OBC/SC/ST must continue to live in poverty and even work over 8 hours a day still cannot afford a square meal. Is this the social justice the media dictates to us? When govt say just a paltry percentage of 27% for OBC it is saying it is blow to social justice or is it a blow to the rich upper caste and Brahmins asks the experts?

Has the print or the T.V media ever made open the sufferings of the rural poor unemployed graduate. The experts said for instance the first generation graduates from the rural India remained unemployed in a very large percentage. No body ever projected the mental and psychological problems faced by them.

Those first generation learners after graduating with lot of economic problems, when remain unemployed, they face of lot of humiliation in their villages not only by their kith and kin but by the entire village who speak of graduation as a unless one, but only an handicap for the one who possesses a degree. These graduates either flee the village due to humiliation or migrate to other far off states to work as coolies or drivers or sales man and ultimately land up many a times as HIV/AIDS affected persons.

Has media ever given any article about it?

The media has the guts to claim that the agitating youth under the banner of youth for equality are by far the best and brightest of our youth who can hold their own anywhere. They are India's undisputed crème de la crème.

The media is falsely propagating that for the last four weeks not one politician from any of the main stream parties has dared to be associated with this minority movement. The experts answered their claims in several essays and pages but we here put only the crux of their reply. First if these students can hold their own anywhere why do they want to get education from IITs, IIMs and AIIMS? They can as well sacrifice their seats right now to a OBC/SC/ST and choose to get educated from a normal (average) college claims the experts. This act if has been carried out by any one of them it would certainly show their broad mindedness as well as their sensitivity to the oppressed and depressed majority and above all the very best. They but on the other hand are like the Dhronacharya asking the student who learnt all science of archery by himself his right thumb. Rightly from those days till date the Brahmin remains a Brahmin Dhronacharya taking away from the OBC/SC/ST what they have earned and making them powerless, helpless and hopeless. The cruel tact in them has been passed on from Dhronacharya to the last brahmin who is in the youth for equality in such a powerful way, they want to deny the OBC/SC/ST a good education that too from the govt run institutions.

Now we are not as meek and foolish as Ekalaivas to give our right thumbs. We will fight for our equality claims the experts. You (Brahmins) may not have changed but we have realized our faults and we have changed to make the society a better place for the majority said the experts. So even in their story a sudra, Ekalaiva could reach an excellence without a teachers help yet the jealously teacher Dhronacharya asked for his right thumb to make him invalid in that very art through out his life. This is the villainy of the Brahmins claims the experts. The Dhronacharyas are forgetting that todays Ekalaivas are not that budging as the good old Ekalaivas only the Dhronacharyas remain the same or still more cruel than before.

They (media) have further claimed that no mainstream potilician has ever turned for their supports. The experts say that the youth for equality banner is itself the BJP and top Brahmins brain child as they give as much of protection as the brain gets from the skull. How the print media above all lie so blatantly asks the experts. Because our President Kalam on 23rd May 2006 asked the doctors and students to end strike, they to appease the striking students constituted the Moily panel, to end to anti reservation agitation made by a few doctors and students Valiyathan panel was constituted, NKC came to the protection of the anti agitators by saying quota would pose grave dangers for India as a nation, PM met the striking medicos inspite of their arrogant behaviour. Youth for equality is nothing but a banner created by Dr. Venugopal together with the BJP, RSS and Sankara Mutt supporters said the experts. The media's nasty stand of being very openly lying and partial is very condemable claimed the experts. BJP to keep its vote bank tact, made several announcements in this period like, they support the reservation policy except to avoid creamy layer and so on. This was also shamelessly published in the media. So the claim as "youth for equality" as RSS and BJP banner with so few faces most of them 'women' appearing in the media in such postures and reading novels, sleeping, etc. was reported several times to buy the sympathy of the viewers (Public) in contrast with the pro reservation protestors majority of them being rural poor boys did not attract the media said the experts. The media in this period has not only behaved in an biased way more so in an antinational way by protecting, promoting and propagating antisocial activities blamed the experts.

The experts were very much pained to see the over reacting Narayana Murthy of Infosys, who said the very idea of increasing the seats was not at all acceptable, he said on the sidelines that too in a release of a newsletter of ace athlete P.T. Usha. Thus if he is a Brahmin and if he can voice against reservation the media will welcome and write the issue by creating occasions. What is the relevance for Narayana Murthy to talk of reservation in a release of a newsletter function asked the expert? Further they said he had the guts to ask the govt. to first consult the directors of these educational institutions. If the directors need faculties they should ask the government for the same; why should Narayana Murthy of Infosys, mouthpiece of them? Narayana Murthy who was very well known as an IT expert, an Infosys chairman has suddenly become an expert to talk about reservation for the OBC (when he is a staunch Brahmin); as he has media to express and the brahmin media now has place only to write about anti reservation factors that too falsely and foolishly claimed the experts.

The experts said it is high time a G.O. be made that if they talk ill of OBC or about reservation they should be penalized, for many a times they talk and portray OBC as useless by saying merit would be lost; India as a nation is at danger and so on and so forth. It is high time they are sensitive to the feelings of the OBC/SC/ST. They cannot be more casteist and racists said the experts.

Another person who objects reservation is the CII president R.Seshasayee said the experts, Is he the CII president or a OBC reservation expert? Who is he to object reservation? It is well known from Periyar's time that all Brahmin fellows are against reservation! Now this only show they have not changed but time has changed! If they want to maintain something as respect they should keep quiet said the experts.

They have no business to talk about the OBC reservations. One can very easily see the difference between a OBC person and a Brahmin talking about OBC reservations said the experts for instance, Uma Bharathi the Bharatiya Jan Shakthi (BJS) president and former chief minister of Madhya Pradesh who is a OBC called for a national debate on quota issue but we see the Brahmins in power only say merit is at stake, the nation would be ruined and so on and so forth. The media will write pages about anti reservation, which is unscientific and fictitious said the experts.

The media being in their monopoly shamelessly publish these articles said the experts, this is again very clear from Shiv Khera who was all these days never seen on media to discuss about any social issues suddenly is seen in the print media issuing a statement that he opposes caste-based quota. He further advocates that the only solution is to give reservation on economic basis he further says that people do not like casteism. But the Brahmins and the Brahmins students who agitated are so casteist that they made a ‘symbolic’ protest by sweeping the streets, shining the shoes and singing songs against OBC which was casteist. The expert asks, when the medicos made such casteist protest why did not Shiv Khera oppose it? What does this mean? Let him say no caste should be practised then without their saying reservation will become meaningless. It is very clear only reservation alone cannot be caste based even in the central government run educational institution; why the media did not publish any of the articles written by the non Brahmins against the casteists symbolic protests made by the medicos and medical students.

Thus it goes without saying a Brahmin is a Brahmin when it comes to OBC reservations. Only the non Brahmins are fools to think of them other wise said the experts. If the Brahmins are so concerned let them give away their caste and caste identity said the experts and let the media make it public and make a propaganda against caste and caste practices.

The non Brahmin central cabinet finance minister Chidambaram says we won’t review quota. He clearly says as a member who is in the group of ministers looking into reservations for OBC in higher education said the government had no intention to review the caste based reservation policy. “As I understand there is no ground to review whether there should be reservations or not…

Again the Brahmin Malhotra says with his Brahmin arrogance that why the same horrendous mistake a second time? The only question is how can one ever say V.P. Singh's move was a horrendous mistake? The reader is requested to contemplate over the venom in these statements!

The media has very indecently blown up the Anbumani's and Venugopals issue in a very nasty way. The only view one can see is how a brahmin subordinate wants to show off his power using some political forces like BJP, RSS and so on to put down a non Brahmin sudra authority. One knows very well the Union Health Minister Anbumani is also the President of AIIMS by the very fact he is the union health minister and Dr. Venugopal the director of AIIMS. Clearly Dr. Venugopal holds an office subordinate to Dr. Anbumani. The BJP demanded the resignation of Dr.Anbumani the basis of the issue being that Dr. Venugopal was removed from the post of director of AIIMS on the allegations he had misused the power as well as money.

The minute the termination of Dr. Venugopal was made the media was highly upset, one can see big head lines on these issues. P.V. Indiresan voices for the brahmin Vengugopal in the main editorial in which he says; our parliamentarians have vehemently espoused the principle that " Parliament is Supreme". There is little to choose between this idea and the Divine Rights of Kings. The experts said the IQ of the former director of IIT Chennai is very low for the Parliament is nothing but the peoples mouth piece elected by the people but the king is hereditary or he is appointed by the brahmins, it is a pity even this difference is not known to a man who held such posts. What has IITs and IIMs have so far accomplished in the international level? Once again government seeks explanation from AIIMS director was written in head lines; suddenly a plea to disqualify Anbumani as MP was made to ABJ Abdul Kalam our president. Ultimately the docile Anbumani formed a panel to probe into the AIIMS functioning for in the last three to four years, the quality and functioning of the AIIMS was not satisfactory. The experts put forth only one matter to the public about the media. Did any letter appear in the major editorial page of the New Indian Express or Hindu favouring the actions of Anbumani. If this would have been published then certainly one can say that the print media is working in a judicially balanced way for it puts the feeling of the Brahmin director of AIIMS Dr. Venugopal as well as the problems of the Sudra president of AIIMS Dr. Anbumani. But to once surprise all the experts said almost all the articles writing by the leading print media was only supporting Dr. Venugopal beyond limits and at the same time decrying the Union health minister Dr.Anbumani who is also the president of AIIMS. This instance is a fool proof to show how biased the media is functioning.

One can further say how a Brahmin subordinate is lifted and how a non Brahmin authority is badly put down. Every news in this direction projected Venugopal as a great cardiac specialist of World standards with pure and blemishless character on the contrary, Dr.Anbumani a sudra is projected as a corrupt and a dull person craving for power. Why is media acting in this biased way asks the experts? Is caste ruling the media? Will they ever say no to caste system? In fact they are the ones who grow casteism for without which they cannot live!

A clear cut act of media projecting 27 percent reservation for OBC as a crime done to the nation is brought out! Did the people who run the media ever think that every problem has atleast two sides? Is nation development majority’s prosperity or a minority Hindus prosperity that too who are not even the sons of the soil? Is nation development measured from the status of a majority or from the fact how many IT industries are in India? Is the nation development remains only in Brahmins getting education as doctors from AIIMS and not giving any chance for the majority of the non Brahmins to study in AIIMS?

Is the nation prosperity lies in IIT and IIM educated graduates leaving India and earning in lakhs by serving other countries? Is it just for them to leave India for money or serve India in poverty? Is this nation development? If any one talks in favour of reservation they blow only one trumpet that merit is at stake! Will they come out and define what they mean by merit? Denial of opportunity for the oppressed rural majority is merit according to them! What have they achieved by their merit except their purses are full their bellies filled all the time and the last man fights to have a square meal a day inspite of toiling over 8 hours in a day in sun or rain? Why no media is coming forward to air their problem? It is because they are not Brahmins. Economic ill balance can lead to social inequality. For if a man is uneducated inspite of being economically sound he is not socially well treated! So to achieve social equality one needs good education in a good institution. If alone our rural OBC/SC/ST people get education in these institutions certainly the nation can witness some social reforms! In the first place media has failed to produce even a single article in these directions. Secondly the media has never done its duty by reporting both sides. Ever since its functioning the media has been one sided alone and it never felt there was other side of reporting. These very many long years of bias reporting has resulted in the media becoming an unthinkable and non analytical tool which was only promoting one small section of people, viz, the Brahmins.

All the more when the issue was reservation for OBC, the foolish media is showing it vulgar blazed colour in all angles. That is why every Srinivasa Iyengar and Natarajan Iyer's article on anti reservation is put. Why do they write? What does the nation gain? Because they write such things, will the nation become developed over night? It is but pertinent to mention at one point of time or other that such bias reporting is going to have its bad effects! They by such acts are going to loose their own identity and power in due course of time. People especially younger generations are not only intelligent but also very critical to understand every thing. Atleast a day will come when all the non Brahmins will decide not to buy the news papers; dailies or weeklies or monthly and they would be forced to close it down. It may also happen that the forth coming governments may initiate to run a print media of their own in such a case these Brahmin media's may be forced to close down.

Now it is still important to mention that S.C had the fear, reservation will split the nation on the basis of castes; but the experts wants to record the drama they have played for the paltry 27 percent reservation for OBC that too in the educational institutions run by the central government alone, will be the seed sown by the Brahmins for the spilt among people on the basis of caste. We saw caste in every reporting in all possible angles said the experts. What do they loose by adopting the policy of live and let live? Why they want every thing good and give nothing to the non Brahmins? Unless they change their mental make up they are going to create problems for themselves! They are so intolerant when it comes to reservation! They are restless, hasty, unthinking! Why now the yogas, meditation, etc; do not help them to over come this? Why they want to treat the majority and more specifically the sons of the soil as inferiors to them? They do not want OBC to get education in par with them! What does this reveal? How will these Brahmins who get educated as doctors treat the non Brahmins? Will they really treat the illness in the non Brahmin patients or just for money give treatment? When there is no sense of equality in their mind, the experts strongly doubt about their honesty in their profession! So it is high time we try to isolate ourselves from these selfish Brahmins not in hatred but to save ourselves and protect ourselves from them! They have the guts to call a SC-doctor and describe in a movie the ignorance displayed by him. Their arrogance knows no limit. Let us remain united and fight for a social cause said the experts. Let the situations forbid us from taking any sort of medical treatment from these Brahmin doctors said the experts, for they in their hatred can kill, us this is clearly seen in their protests for they have no conscience about the patients lives; and the patients who died due to their protests; only protest against reservation was so important.

If they had been bound by any professional ethics they would have kept their protest against reservation in the off hours, not in the working hours. If they ever had some thing as mind for rights will they ever accept the pay even if the government gives them during the strike period when they have not worked for it? These two instances are sufficient to judge what a Brahmin is that too in such a small age in the time of student hood and practising period. How can they do any act to help the nation or for the development of the nation? What is uppermost in them is self and according to them development of themselves is the development of the nation and any thing favouring the majority which they construe as denial to them and will go on to say that will stop the nation from its development and they would voice from NKC or as a ex director or as a CII president or any Brahmin that the reservation for the OBC is denial of merit and above all because of it nation is going to face drastic and irreparable consequence; nothing more. Above all to write all these Brahmins are readily available and the media is in their monopoly. Only to write anything favouring a non Brahmin, the media will never do for they are so united even if a non Brahmin is at his height of glory and achievements the media will not write a single line about him. In fact even an average in their caste would be projected and pampered in all medias. They are very particular in their policy of spreading Brahminism and never give any voice to up lift a non Brahmin. Thus the experts said nothing is wrong even if we mention the media in India to be a Brahmin dominated media or as a Brahmin media when they have the arrogance to call a dalit doctor and disgrace him. Thus the media till date has not come out of the clutches of the caste that is why truth or real facts about many of the things occurring in the nation is never debated.

We in this book have restrained the media used in movies about reservation but only keeping record a simple thing. That a movie entitled “gentleman” was scripted and shown; the main theme was two school mates who were friends one Brahmin and another non Brahmin aspire for medical seat and the non Brahmin gets the medical seat on basis of reservation (it is shown) so the Brahmin boy vexed by jealous commits suicide. In the same movie it is shown how a Brahmin girl studies in the quota (reservation) of SC and faces some problems. In the first place the experts ask how many brahmins in reality have committed suicide? There are only instances that brothers kill each other in their community for property. This is in vogue form the time of Mahabaratha and Ramayana. Now we have seen this repeating very recently in the Mahajans family, the elder brother was shot down by the younger brother; that too they are from the very fanatic orthodox party, viz; BJP. Thus the movies against reservation screened by them have only given doubts in the minds of the non brahmins how far this is true. The experts further add that young non Brahmin students see this movie, they being very sensitive in that age, unlike the Brahmin boys will start to think against reservation and so on.

Thus this movie itself is a poison tree to poison the minds of the young non Brahmin children against OBC reservation added the pained experts! Why they do not take movies in which because of reservation a non Brahmin boy or girl, (first generation learner) up lifts the whole village or so, for such real instances are also existing and true. Thus a Brahmin does something means it always has deep ulterior motives to poison the non Brahmin youngsters mind from any form of social reformation but at the same time serve their ends by making them their slaves or work for them. Many a times they have said if every one is studying who will do the menial job of cleaning the road, burying the dead, washing our clothes and so on. The experts advice them only one question we have done this job for many thousands of years let our Brahmin youth for some time do our job and let us study so that now we will see whether the nation becomes a developed one! Thus they have also used movies but we in this book do not discuss about this very special mass media due to several constraints but say this topic is beyond the scope of this book for many technical problems.

Now we mention how, the media is biased by speaking about the quality of education in school at the eleventh hour when so many years have passed in which they never had any mind to even make a mention of it!

In the span of two or three months soon after the reservation bill for the OBC was passed in both the houses with the over whelming majority several articles in the print media has appeared about the improvement of the schools and their standard. The so called CII tried to set rolling of the ball by making the government to partner with NGOs to improve primary the level education in schools run by the govt. and municipality.

Provide support to SC/ST students to prepare for entrance exams on professional and technical course. The same CII-committee now as if concerned of the nation says that government should improve quality of schools atleast in select districts. IIT had asked the union HRD ministry to strengthen the primary education so CII and IIT go hand in hand. An editorial in the New Indian express dated 26th May 2006 says that Government schools have insufficient black boards and toilets and teachers don’t turn up to teach.

## 4.4 Origin of Reservation in India and Analysis by the Experts

Before we talk about reservations based on caste in the field of education we must know the origin of reservations. Was reservation an Aryan ideology or an Indian (non Aryan / non brahmin) ideology. In the pre Aryan period Indians did not know anything about reservations in anything, any one who wished to develop in any field had the full freedom to opt for the same and work for it. But only after the entrance of the Aryans who entered India bought the notion of reservation in India through the very Khyber pass. The concept of reservation is found in Laws of Manu. Infact the origin of reservation can be traced to the period of Laws of Manu; only in laws of Manu they not only reserved education and religion for them (for Aryans) but also very severely punished those who tried to learn or practice religion [43]. So they made such reservations based on caste by birth. The term “Caste by Birth” itself is a myth. For caste was assigned according to the birth of people from the god “Brahman”. People born of Brahman’s mouth were called brahmins i.e., the Aryans the twice born. People born of the arm of Brahman were Kshytriyas; people born of the thigh of the Brahman were the Vaishyas and Sudras were born of the leg of Brahman. People born below the naval of Brahman are termed to be impure and are untouchables according to the Laws of Manu. So by the laws of Manu majority of the population became impure. Say over 87% of Indians are impure and only the 3% Aryans who were alien to this very land became twice born. They had reserved teaching, studying, that is education and religion for them.

Their reservation policy severely dealt with those (non Brahmins) who learnt or practiced religion or learning. The work of the ruling of the nation was reserved for the Kshytriyas. The Vaishyas were supposed to do trade. The menial job serving the Brahmins with humility and sincerity was assigned to the sudras. Now they (nomads) were father of the reservations in India. When we ask for reservation in education no wonder they become hysteric and make so much protest using the media and court which is in their monopoly. All the menial jobs were reserved to the fourth caste. They declared caste by birth and we say reservation based on caste. The term “Varnashra dharma” was used instead of reservation. So they (twice born) till today have not done any menial job of scavenging or shoe polishing or burying the dead or removing the dead animals. They have also protected and followed the laws of manu very truly in reservation. So when we ask for reservations for the non Brahmins; it naturally becomes illegal and undemocratic and they agitate and burn a non Brahmin gutka seller; lying, he is a Brahmin medico later on whose identity was revealed only by the police investigation. This is their honesty. All this drama of medico, self immolization was staged and medias projected lavishly but lastly to find that the person was not a medico; but only a purchased gutka seller who was paid heavily for this acting.

Now when they promulgated reservation, we sons of the soil very meekly accepted it from the time of Manu. However the caste of SC/ST was not mentioned, to have born from Brahma. Now even after 5000 years that too after 6 decades of independence when we say remove the "Varnashra dharma" reservations and permit 97% to become educated in standard institutions that too run by the central govt. which is monopolized by them (brahmins) for these many years they are saying no and protesting for it and shamelessly our political leaders and people in power meet them for compromise! What is the need for it of what is there to compromise with them? Is not the nations development the development of the majority? More than the Muslims and the westerners the 3% Aryans who had come not as invaders and colonizers but as beggars has swindled India’s heritage and wealth. Now the major economy is in their hands. They by force are the power holders and policy makers of India.

Next we see how they begged the British for reservation of posts in Electrical Engineering Departments. Here also they only brought into the soil of India the notion of reservations.

Thus the reservations for Brahmins started in (Vedic times) times of Laws of Manu continued even in the British rule. Only because of these reservations these Brahmins are nearly 90% in all highly paid posts both in education and administration.

Only now the non Brahmins have thought about reservations and that too only for the past few decades it has found a place in constitution and that too the rule or order is only on paper and has never been in practice till date. Though the reservations for SC/ST are said to be in practice they have not implemented it in practice.

Thus the main fear in the Brahmin mind is that if reservation is given to OBC like the Brahmins in the British time, they too would become powerful for they alone are fully aware of the reservations for they have developed only because of reservations from the time of the laws of Manu. That is why they carry protest so strongly. Here also we see how prejudiced and discriminative they are, in their attitude that is even after 5000 years of exploiting and enjoying all comforts from the Indian soil. Some of the socio scientist feel that when they entered India they in order to equip with the Indian civilization that existed in those days made very strong rules as the laws of Manu to cope up with our civilization, that is why they choose to be our teachers; and unfortunately made all others poor and have usurped the soil. This they could do tirelessly in a period of 5000 years. They were very cunning to assign to themselves the ministerial post. They feared the strong community for they were not only physically strong but they must have been capable administrators and should have been the master of all arts, science and technology with sound economic background so only they were thrown down as sudras and dalits, hence laid down very strong rules to put them down, which was very cunningly done and executed by these brahmins by writing laws of Manu. That is how the third and the forth caste must have been born says the experts who are a socio scientist. According to them the more the rebel the people were they were more discriminated by lower caste.

Next we wish to state the other expert opinion on this matter. Some of them strongly urged that the concept of caste based on birth be banished and law be made so that no identification of caste will be made in name and no sign of caste like holy thread tuft and marks on the head be put then we will not ask for reservation; they claimed. We will have no caste, no religion and all are equal then why do we need any reservation said the experts. The experts were of very strong view that they (the Brahmins) had no voice to say any thing against or about reservations. For some of them recalled that the Brahmin Associations of Tamil Nadu held a conference on $24^{th}$ -$25^{th}$ Dec 2005 and in that they made a very clear resolution to request the govt. to give them, 15% of reservations in all posts and education. So when they are only a little less than 3% of the population they have the guts to ask for 15% reservations and when we are over 50% of the population, when only 27% of reservation was made as a constitution they are against it and they make all protests and our political leaders are so sympathetic they go out of the way and console and cajole them to come back to work. What a drama! When our farmers committed mass suicide none of these political leaders even had a heart to visit them, says these experts.

Some of the experts came very angrily down on the suggestions made by certain people in the leading dailies that the standard of the school must be improved or some other rubbish to that effect.

To this they, all reacted saying with their social background and the poor facilities, studying in such poor schools they (OBC and SC /ST students) not only clear the exams and get over 75% some of them over 85% and even some in 90% for they are writing the same papers the Brahmins write and still they are denied entry into institutions like IITs, AIIMS and IIMs.

If reservations are given certainly they will outshine Brahmins who boost of their merit and quality they said. For these Brahmins and upper caste by getting education in the best of city schools obtain only such marks with so much facility and encouragement. By shear hard work and no other merits they pass JEE or any entrance exam. What is the greatness or any merit or achievement in them, asks these experts. So they feels that not only 27% but 50% reservations must be given to OBC and 27% is enough for the open quota they feel.

Any one with little intelligence will certainly accept who is meritorious? A rural boy with no good/nutritious food working off hours of school for family, walking to school and studying in a ill furnished school say scores 40%, it is more than a city bread boy with full of encouragement with good and nutritious food, caring parents, goes in a car to school who studies in the best of schools with access to all types of tuitions gets over 90%, he is by all means inferior to this rural boy. Thus reservation is the only means to save them and help them of good education, especially when these institutions are solely run with the govt. aid; which is the tax of the majority.

Some other experts argued that reservations must be made upto 50% taking into account the following facts.

All these days the IITs, IIMs and other central institutions of national importance did not even heed to take any steps or even talk about improving the standards of all schools in the rural areas or the ones run by the govt. Have they ever thought over these factors? They had no mind even to adopt under NSS scheme to help these schools. Now to avoid the real issue they are so much bothered to improve the standards of these schools. These suggestions do not come out of concern or social awareness, they come out only to sabotage the reservations for OBC. If so concerned it would have been extended not when govt. is making some rules to help/uplift them. What is the meaning of their going on anti quota protests and strike?

Their concern to improve school education amounts to "pinching the baby and rocking the cradle "claims many of the experts. What have they done all these 60 years about the improvement of the schools infrastructure or its standard? It may be recalled with disgust that the then chief C.M. Rajaji (a Brahmin) of T.N. when saw the increase in the literacy rate after closing down of the liquor shops, immediately closed down several schools. It was the intelligent- Periyar who not only observed it but has recorded it.

Some of the experts expressed the following facts:

They asked how many of the brahmin children are rag pickers? Can there be one among a lakh or one among 10 lakhs? How about OBC? Statistics shows over 60% of them (rag pickers) are from OBC. Also they added how many brahmin children spend their lives on the platform? How many of the children from the brahmin community are child labourers? They ask why reservation for OBC? What is then right or responsibility towards the nation's upliftment? They are bothered only of their comforts and powers! They are the true traitors of the nation for they are agitating against the government policy, which tries to bring social equality. If the OBC have protested against any social equality for brahmins what would have been the plight of them? Their very agitation is not only illegal but against all norms of humanity, feels the experts. Their show of protest and the media's loyalty in carrying the news go hand in hand! What ever we write to these news papers or dailies or weeklies are never reported only anti reservation factors are shown out of proportion!

## 4.5 Suggestions, Comments and Views

The issue that merit would be lost if OBC students are given reservation was proved to be false from the model given by the experts. For it was always in the off state. Thus there is no relation between the merit of the institution and reservation for the OBC in institutions of importance. It is an absolute absurdity to have thought of linking both. Further there was the on state of them only when the expert was only Brahmin students. If the arrogant opinion of the brahmins to have linked reservation with merit of the institute. Many OBC outshine the brahmins in most of the fields but because media does not exist, to boost the OBC enough recognition or even in many cases sheer recognition is not even given to them.

The reservation given to minorities alone can save the nation from terrorism is the conclusion derived from our study. For in our nation except Hindu minorities i.e., Brahmins alone who are not the sons of the soil and all other minorities following other religions are all only sons of our soil so it would not be proper on the part of the govt. to say other minorities cannot have reservation. The main question before us is that how come the Hindu minorities who form such a small percentage of foreigners in India happen to be occupying nearly 70% of all the seats of higher learning in all the institutions run by the central govt of India. One is yet to find any proper answer except the fact that in India only this Hindu minorities are the policy makers of the higher education and they dominate all the teaching posts in institutions run by govt like, IITs, IIMs and AIIMS. So they have their own undetectable means to absorb their own caste people said the experts.

Secondly all the powerful departments like law, media and advisors to govt are not only dominated by them but ruled by them so even any amendment passed by both the houses in full flown majority are stayed by the Brahmin courts. Unless we have non Brahmin courts we would never get justice and the nation would never be free from terrorism (Gujjars recent protest).

The only question put by our honourable C.M. Kalaignar M. Karunanidhi is that the fate of 100 crores persons are decided just by two, that too they are not even the direct representatives of the people. He has very boldly said that if this is done the nation will be at a risk, for in India alone we have several divisions of caste and those caste at least at one time or other must be given their due for the caste name they carry are imposed by the Brahmins (laws of manu) for which they are responsible and now denying them the right for higher education that too in the institutions of higher learning by the very government is atrocious and would result only in a turmoil in the nation; he has warned. Thus the speech of T.N.C.M. “Karunanidhi” said the world must be made to realise that a situation where the horoscope and future of 100 crores people can be determined by two or three persons is a big threat to democracy. We are not seeking to march ahead but merely asking for space to keep going on our path “even that space is being denied to us” he remarked.

Whenever the Sudra community raised its head to establish its rights there were browbeaten into submission attacks coming out of the blue. The reservation policy had always faced opposition and hurdles right from the days of the justice party; the chief minister said.

The courts attitude has agonised the people for even the single drop of blood that oozed when the body is pierced reflects the pride of communal representation and reservation which are from Periyar's policies, Karunanidhi said. Karunanidhi recalled that the seeds of reservations were first sown in Tamil Nadu and protected by Periyar, C.N. Annadurai and K.Kamararj. We are seeking reservations to ensure that our future generations crossed hurdles to progress in education and the country was strengthened.

With the Tamil Nadu assembly reverberating with a battle cry against the Supreme Court for its refusal to vacate the stay on 27 percent reservations for OBC in higher educational institutions, Tamil Nadu Chief Minister M Karunanidhi today urged the centre to take immediate steps to prevent a nation wide volcanic protest.

"A time will come when the marginalised, oppressed and those pushed into the fourth and fifth castes" in the country will erupt like volcano. Without giving room for that the central government should immediately take steps to find a way out. There is no need to pass a resolution in this house pressing for this demand. The fiery speeches delivered here is enough, the Chief Minister said in the Assembly.

Thus the C.M of Tamil Nadu has very clearly shown the Tamil identity in reservation. Another matter of importance is that we see the denial of reservation which is nearly pending for over 57 years would only lead to some sort of blood shed. For the Brahmins cannot stand as an hindrance to govt policies trying to wipe out the social discrimination based on caste and caste alone.

They project marks as merit. Can anyone in the world agree to the fact because one person scored cent percent in all the exams he is a meritorious person. Merit is different, mark is different, mark are obtained by memorising all the subjects. They would not be even able to answer any implied answers from their books or even answers for questions which are not scheduled to appear. Thus the Brahmin student like memorizing the meaningless terms of Vedas will memorize these subjects. How can one call this a merit? Merit is a very term like Maya for if we accept their argument what have they done in these 60 years of their domination. They have not made any significant contribution to science or technology or medicine or literature.

Why are they boasting on some non-existing things like merit which only is based on shear memory test. Let them stop this! Already the editor of the Dalit voice has clearly wrote in an article in which he has said, “merit is my foot”. In Tamil Nadu our C.M. Dr. Kalaignar Karunanidhi had abolished the common entrance test to admissions to professional courses. This was challenged in the Madras H.C. The Madras High Court on Friday upheld the validity of a Tamil Nadu law abolishing Common Entrance Test (CET).

In a division bench comprising Justices P.K Mishra and Sampath Kumar said that the standards in professional colleges would not be lowered by abolition of CET. “It cannot be said that the Governments conclusion that forms the basis for normalisation process is so arbitrary as requiring interference by the court on the ground of being volatile of Article 14 of the constitution”, the bench observed while dismissing writ petitions from S.Aswin Kumar a student and four others challenging the Tamil Nadu Admission In Professional Educational Institutions Act Enacted by state legislative in December.

“While upholding the validity of the impugned statue we clarify that for admission to architectural courses aptitude test as prescribed in the regulations (pertaining to architectural council laws) has not been dispensed with”, Justice Misra said.

The judge pointed out that both Medical Council of India (MCI) and the All India Council for Technical Education had also argued that the petitions had to be dismissed in views of the TN act. “It should not be forgotten that the MCI regulations also authorise selection on the basis of qualifying marks”.

Dismissing absolute equality as a myth, the judge said there was an inherent possibility of ticking some answers more by guess as in the KBC (Kaun Banega Crorepathi),TV programme rather than any conscious selection of the right answers.

There would be no objection if a reasonable amount of equality could be achieved by any other method he added.

The mindset / psychology of the students has to be taken note of before taking further action in a matter, as the students are under the impression that there will be no CET this time and have planned to enjoy the vacation during May and June. Justice

JAK Sampath Kumar a member of the two judge bench that upheld the state enactment abolishing the CET.

Adding his views to that of the other judge in the bench Justice P.K Misra Sampath Kumar said, "Desire of the crores of people (to abolish the CET) cannot be denied by quashing the act but their desire should be fulfilled". Stating that the Act is a social welfare legislation to meet social justice, he said. The enactment is to prevent harassment and hardship to the socially, economically backward and weaker sections of students from both rural and urban areas getting admission into professional course. If that be so can it be stated to be bias, arbitrary and volatile of article 14 of the constitution. Certainly not, the claim of the petitioners is devoid of any merits".

Interestingly, Sampath Kumar also listed out five sections of people each of which would gain advantage and disadvantage respectively because of the hitherto existence of the CET explaining it by putting a tabular column the judge said; CET would only be advantageous to aristocrat schools, coaching centres and children of elite people and highly qualified parents as they can devote full time to studying with comfort.

And it would be disadvantageous to children of unqualified and illiterate parents, students who cannot afford to go to coaching centres and those who cannot devote full time in studying, parents who cannot take of their children's studies because of their work to eke out a living and students studying under greenwood tree with mosquito bites" he said. April, 27-04-07, Indian Express.

When just a state entrance exam CET is itself so much dependent on the social and economic conditions of parents, one can imagine how it would be for preparing JEE! What is the cost of the coaching centres? What is the caste that can afford for them with command? What can a rural Sudra or SC/ST poor do? Will they get or have money to buy any one of the competitive books for JEE leave alone the dream of coaching centres. Can they afford to buy paper and pen to work some problems?

When the govt after pondering over all these practical difficulties faced by them (sudras) to avoid further social discrimination announces a paltry percentage of seats to be reserved for them these arrogant brahmins had the guts to take this issue to S.C and S.C which always remained loyal to these brahmins by following the dictums of the laws of Manu has stayed the reservation. It is not an exaggeration if one says the democracy is dead in India; for we see the term equality has no meaning for the 3% are enjoying more than the 90% of the seats of higher learning. Where is social justice? Where is social inequality? Only the Brahmins made in the name of caste the social discrimination and social equality, now to establish at least to some degree the social inequality the only tool the government can use is reservation for the socially deprived and none others for that alone can save the nation from turmoil.

The suggestions the experts gave to the govt was if “Total reservation for backward classes to society should not exceed the limit of 50 percent as held by the Supreme Court and as contemplated by the frames of the constitution says the report …

1. In the first place when the constitution was written there was not even a single OBC representative; so as said by the CM of Tamil Nadu a through review and restriction of the constitution must be made as majority of the nations populations are OBC and naturally they need to be properly represented for it impinges even as said by the constitution “equality right of other sections (OBC) of society”.
2. The only way to solve all these problems is representation as per their strength as said by Periyar.
3. Unless these inequalities are properly settled at an appropriate time without much delay (as even now it is 60 years have passed) India is bound to face a revolution which cannot be easily put down, losses may be heavy or unimaginable taking into the simple fact they (sudra) form the majority of India’s population (Gujjar’s recent protest). They would be clearly supported by SC / ST. It is high time the arrogant Brahmins understand the grave situation and act reasonable. This wave is seen in T.N and already functioning in Rajasthan, A.P and Assam. The naxals problem which the nation faces is nothing, for the can no more accept the wide gap between the upper castes who are the haves and the lower castes (Sudra, SC/ST) who happen to be the have not for so many centuries.

Reservation should cover all areas like agriculture and non agriculture lands, capital, employment input and product markets including private housing and educational institutions Thorat said.

The scheduled castes and tribes suffer without access to the factors of production, both in agriculture and industry. The government should develop a policy of having reservation and preference in supply of certain inputs of production by government private parties he said. Thorat suggested setting up of a common land pool where lands could be acquired by the government and distributed to landless among the backward classes so that conflicts could be avoided.

Talking about reservation for OBC in various areas he said even if concept of equality by birth was legalised, it could not be implemented. Reservation, he contended was the ultimate solution for the OBC to invigorate their abilities. “The upper caste people should sacrifice their rights for those belonging to the lower castes through reservation. This is because their forefathers had enjoyed all the benefits and deprived the lower caste people of their rights. This is only solution for the country to improve”, he said.

The experts very warmly welcomed the views of the University Grants Chairman Thorat. They all said if the suggestion of the UGC chief was followed; certainly India would become a developed nation in two decades. Further the land and money owned by the Hindu religious mutts should be taken by the government and distributed to the OBC/SC/ST. The Nation is full of burglary, murders, etc mainly on the basis to have basic needs for the divide between the haves and have-nots in these 56 years has only increased very steeply. The upper castes have everything and the lower castes have nothing; how long can the same situation pull along asked the experts! They said as per Thorats’ order the social discrimination which was practised over many thousands of years should be at one point or other settled. They boast of modernity, they take up to western culture live in pomp and luxury and hate to think of OBC/SC/ST, we agree but what right they have to protest when govt takes any positive steps. The experts said at this point they are crossing their limit, which is first a danger to them. Now the govt is also surprised at the hysteric acts of the upper castes said the experts. The situation and the problem at hand is very clear said to experts. The day of settlement for equality is not a one distant said the experts. We are respecting our politicians who are our representatives said the experts that is why we have not so far taken the law in our hands waiting calmly for the final verdict and the act of the govt. Then we would act accordingly said the experts.

Further experts said our patients should not be mistaken for ignorance. We are capable than the upper caste who have denied us for so many thousands of years our just share in education, land, economy, industry and above all the basic equality. They all said we are awaiting our representative action only!

It has been observed with deep despair by several of the experts that when the SC put a stay the truth about it is that we are not able to find out from the judgement given by the two bench judges of the supreme court on 29th March 2007. What ever be it, the immediate reaction on the part of the upper caste authorities working in these govt institutions of higher learning was that they immediately want to release the selection list or to be more specific they were vexed to wait for their view was if there was further delay in the judgement they (authorities of these institutions) said it would speak on the students and in the completion of the syllabus for their students. The only question what the experts want to put in this book is that the OBC has been delayed nearly 60 years their right to education in these institutions but for the upper castes a few days or months of waiting is so much a problem! Now at least the reader and one and all will understand the selfish arrogant behaviour of upper castes. Are they ever going to think in terms of national development (that is development of the majority) or think and act as if only their development is the national development. Why no one is ever bothered about the majority of the students but only the few students who have joined or going to join for their courses in these institutions of higher learning? Their only argument is you all (OBC reservations so OBC students) have waited for over 57 years why not wait for one more year so that they can think about it in the coming academic year 2008. So when it is the upper caste student problem they are sensitive but when it concerns the problem of the majority of the poor sons of the soil the OBC they have the heart to say they can wait. In fact both the houses passed the amendment in Dec 2005 itself. Why is no one questioning about over a year delay? What happened to its implementation in the academic year 2006? A few months delay is so much projected by the print media in the first page with bold letters. Did these print media ever write any article or criticise the govt for being slow. Only one thing the experts ask as mentioned by the leaders suppose the OBC launch a nation wide protest paralysing the very functioning of these institutions (as these medicos did it over 40 days and successfully got salary with no work with the support of the judiciary for that is the power of the upper caste which rules and dominates the media and judiciary), what will then they do for the delay in imparting education to these students in these institutions. So first the experts requests them to think as responsible citizens and authorities of these institutions before they speak! The politicians also cannot cheat the majority by playing caste politics or vote politics. It is high time they too become serious of the issue. Already the nation is facing several problems if the OBC join to protest; the nations peace would be lost and the functioning will be a big question mark warned the experts! Because people trust their political leaders it is possible for them to run the nation in peace, once even in some corner faith is lost the nation would only face a turmoil which may not be easy for them to tackle they added.

Several of the experts viewed the government (centre) was trying to act too smart by blaming the court. In truth they were not taking reservations for OBC very seriously. Might be they are trying to play vote politics! The experts said the politicians are forgetting the biggest fact that people (public) also have very many educated who can feel and know their true mind! This time if OBC are fooled they can be sure of loosing their govt and the politicians cannot act dual characters. Still it is important to mention here that in U.P starvation deaths has taken place. In the last 2 years 174 have died due to starvation.

Majority of them are SC and Muslims. It was reported in Puthiya Jananayagam on May 2007. The following questions are put forth by experts.

1. How many Brahmins die of starvation in India?
2. How many upper castes die of starvation in India?
3. Why no print or electronic media dominated by the Brahmins never publish this news?
4. When the Brahmin medicos were striking against reservation for OBC our honourable PM, and President consoled them; several other VIP showed their solidarity. But when the starvation death in UP has taken place why none of the VIP's like to visit them or console them or announce immediate aid and save them from further death? What does this show? The only answer is that they are after all agriculture labourers from dalits and Muslims? Why should the VIPs be bothered about them? They are not shining Brahmins, which the medias like to project. They are not anything to the nation! They were given identity cards which guaranteed them that they would be given some job but they were not given any job by the govt. This village Belva in U.P is in Varnasi dist. They die of starvation because they do not have any job to do. Due to unemployment and no agriculture labour the people die. When there is no deficit in the yield and we have sufficient food grain why such death ask the arrogant people in power.

These poor people were categorised as living above the poverty line so they were not even given the ration. These poor people have pledged every thing. One Lakshmana from the Belva village said he had nothing to pledge except the marriage sari. So to save his daughter Sema he was going to pledge the sari but Sema died before it he said. She was just nine years old. In the months from May to September alone 5 children have died of starvation in the Belva village. These people are fighting against starvation. They have no means to earn. When people spend lakhs of money on trivialities, so many die of starvation. The only reason that can be attributed to this is caste based starvation. Even in the starvation death caste alone plays the role! Is this human? The divide between the poor and the rich has become so large that people die of starvation on one side and on the other side people visit temples and put in lakhs in hundies for God!

Session court judge Ajit Mishra resigned. Alister Pereira killed 7 persons on the spot by driving his car over them who were sleeping in the pavement. The judge Ajit Mishra had not taken any action and the accused. Alister Pereira was left free by Ajit's judgement. When this case went to appeal in the Bombay High Court. Ajit resigned. This is the way the upper caste judge Brahmin Ajit has given the judgement. He has given only 6 months imprisonment. Now the beauty is that the person Alister Pereira (21) who killed seven of the sleeping labourers on the pavement in the early morning hours is the son of a rich upper caste businessman of Bombay. Further he was in the drunken state, and the car that he was riding had no registration number. The seven persons who died were migrant labourers from Andhra Pradesh. From this the reader can easily understand that; judiciary in India, is for the rich and upper caste. For poor and from lower castes judiciary will be against them. This is also evident from the judgement of 29-03-07 on the stay for OBC reservation. No one is bothered about the millions of lower castes, two upper caste judges give a stay on reservation. What nice democracy is practised in India based on caste? Caste is everything and caste is based on birth. What will the poor sons of the soil do for the political leaders also are with them as they have economy! Now having seen two shocking news, we wish to state that Brahmin media has remained silent in these issues for it does not concern them. It has become very much pertinent to mention here that the media which so well on all the forty days projected the anti reservation protests had naturally failed to record the poor and lower caste people problems. What is the ethics followed by these upper caste media can now be well comprehended by the reader.

Finally only one thing is clear if in the just undergraduate courses like B.Tech. and MBBS etc., merit means marks! Can any one with little brain and rationalism ever accept merit means mark? What does the school final mark or even for that matter the marks in any competitive test imply? Does it imply merit?

It only implies that student has been coached with private tuitions coaching centres and books and becomes well trained like a donkey to do the work of answering the questions and scoring marks.

It means nothing. Given time comfort and above all undivided attention of the parents that too their parents are well placed in society so can command any thing to get marks for their children. They coach him and they do to get the first or the stipulated mark. But on the other hand the children from the OBC/SC/ST have no time even to attend the regularly the school; how can they ever dream of going to coaching centres or private tuitions.

When they cannot afford to get their textbooks or even second hand or third hand torn text books for clearing their school final.

They are from mostly rural areas and as the rich and upper caste they cannot command time for they have usually work assigned by their parents after school hours. If they are in the revision holidays these children for most of these days will spend their time on other work without touching their books. This is not all the suffering faced by millions of OBC/SC/ST students in the night after labour their father in most cases a drunkard will be in the drunken state and quarrel with his mother.

Even this he has to witness; so even the paltry hours of sleep, which they get, will be denied on majority of the days.

Thus these children with their sheer intelligence and merit score good marks. These upper caste children would be given good food nice coffee or tea and snacks to help them to read with energy and remain awake. Now who is meritorious the OBC/SC/ST children or the upper caste children? Who should be encouraged? Will not these children need a paltry percentage of 27% reservation in higher education in these institutions that too run by their own govt.

How long will this social discrimination pull on in the name of merit? How long is the Brahmin media going to publish anti reservation articles? How long is the judiciary going to give verdict against reservations? How long are the ruling parties going to put up this drama?

Will this nation have peace? Will not the nation break out of the usual gear, then what will the Brahmins do? Why should the majority suffer in their own land without good higher education when their govt is spending in crores every year to run these institutions?

What is the use of these institutions for OBC/SC/ST? It is high time they think seriously about this for any out break of revolution will only have evil impact on the nations peace! Let the govt boldly say if reservations cannot be implemented in these institutions it will be closed down! Then what will the supreme court media and the upper caste do? The only reason they (govt) can give is that these institutions do not in any way benefit the poor sons of the soil. What is the answer they will give?

As we authors do not have right to hide any views or edit their views we have no option except to put them with the maximum toned up language. So it has become our ardent duty to write their views, for any form of our editing would not only offend them they would fight with us and say we too are much biased.

So we authors have tried our level best not to edit, only in several places we have polished the language and in many places we have tried our level best to express views in a sophisticated language.

In fact at times the argument resulted in violence authors have always maintained a calm composed, balanced, democratic, approachable equanimity of temper in all the discussions. We also at the end of the discussion would list out the points given by them and got their approval.

The authors are forced to mention all this because the experts are from all walks of life poor, rich, educated uneducated any one who wished to discuss with us. Even when we discussed with the upper castes and Brahmins several time we feared that some form of problem would come up so we tried to have discussions separately in many occasions. We also wish to state that several of the rural uneducated openly asked, can ever our children enter the premises of the schools run by Brahmins and upper castes? We enter only as labourers, masons, scavengers or as coolies. This statement has moved us a lot.

# FURTHER READING

Daily Newspapers referred (April 2006-March 2007)

1. The Hindu
2. The New Indian Express
3. The Times of India
4. Deccan Chronicle
5. Deccan Herald
6. The Statesman
7. The Times

# INDEX

## F

## M

## O

## P

## S

## T

## About the Authors

**Dr.W.B.Vasantha Kandasamy** is an Associate Professor in the Department of Mathematics, Indian Institute of Technology Madras, Chennai. In the past decade she has guided 12 Ph.D. scholars in the different fields of non-associative algebras, algebraic coding theory, transportation theory, fuzzy groups, and applications of fuzzy theory of the problems faced in chemical industries and cement industries. She has to her credit 646 research papers. She has guided over 68 M.Sc. and M.Tech. projects. She has worked in collaboration projects with the Indian Space Research Organization and with the Tamil Nadu State AIDS Control Society. This is her 43rd book. On India's 60th Independence Day, Dr.Vasantha was conferred the Kalpana Chawla Award for Courage and Daring Enterprise by the State Government of Tamil Nadu in recognition of her sustained fight for social justice in the Indian Institute of Technology (IIT) Madras and for her contribution to mathematics. (The award, instituted in the memory of Indian-American astronaut Kalpana Chawla who died aboard Space Shuttle Columbia). The award carried a cash prize of five lakh rupees (the highest prize-money for any Indian award) and a gold medal.
She can be contacted at vasanthakandasamy@gmail.com
You can visit her on the web at: http://mat.iitm.ac.in/~wbv

---

**Dr. Florentin Smarandache** is a Professor of Mathematics and Chair of Math & Sciences Department at the University of New Mexico in USA. He published over 75 books and 150 articles and notes in mathematics, physics, philosophy, psychology, rebus, literature. In mathematics his research is in number theory, non-Euclidean geometry, synthetic geometry, algebraic structures, statistics, neutrosophic logic and set (generalizations of fuzzy logic and set respectively), neutrosophic probability (generalization of classical and imprecise probability). Also, small contributions to nuclear and particle physics, information fusion, neutrosophy (a generalization of dialectics), law of sensations and stimuli, etc. He can be contacted at smarand@unm.edu.

---

**Dr. K. Kandasamy** is a staunch Periyarist. He worked as a guest professor in the Tamil Department of the University of Madras till recently. He holds postgraduate M.A. degrees in Tamil Literature, Political Science, Saiva Siddhanta, Defence Studies, Education and History. He can be contacted at dr.k.kandasamy@gmail.com